\documentclass[12pt,a4paper, twoside,english]{book}

\usepackage[utf8]{inputenc}

\usepackage[italian,french,english]{babel}

\usepackage{indentfirst}

\usepackage[signatures,swapnames]{frontespizio}

\usepackage{amsmath,amsfonts,amssymb,amsthm,bbm,mathtools}

\usepackage[dvipsnames]{xcolor}
\usepackage{hyperref} 

\usepackage[active]{srcltx} 
\usepackage{color,soul} 
\usepackage[OT1]{fontenc}
\usepackage{authblk}

\usepackage[bottom]{footmisc}
\usepackage{tabularx}
\usepackage{makecell}

\usepackage{enumitem}

\usepackage{graphicx} 
\usepackage{subcaption}
\graphicspath{ {./images/} }

\usepackage{geometry}
\usepackage{setspace}

\newtheorem{theorem}{Theorem}[chapter]
\newtheorem*{theorem*}{Theorem}
\newtheorem{corollary}[theorem]{Corollary}
\newtheorem{lemma}[theorem]{Lemma}
\newtheorem{proposition}[theorem]{Proposition}

\newtheorem*{hyp*}{Hypothesis}

\theoremstyle{definition}
\newtheorem{definition}[theorem]{Definition}

\def\N{\mathbb{N}}
\def\Z{\mathbb{Z}}

\def\R{\mathbb{R}}
\def\C{\mathbb{C}}

\let\e=\varepsilon

\let\.=\cdot
\let\0=\emptyset

\newcommand{\norm}[1]{\| #1 \|}

\usepackage{fancyhdr}
\title{Evolution equations with applications \\
	to population dynamics}
\author{Elisa Affili}

\numberwithin{equation}{chapter}

\begin{document}
	

\frontmatter

\newgeometry{top=25mm, bottom=10mm, left=20mm, right=25mm}
\begin{titlepage}

\noindent
\begin{minipage}[t]{0.47\linewidth}
	\includegraphics[scale=0.33]{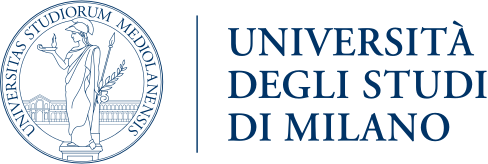}
\end{minipage}
\hfill
\begin{minipage}[t]{0.47\linewidth}\raggedleft
	\includegraphics[scale=0.17]{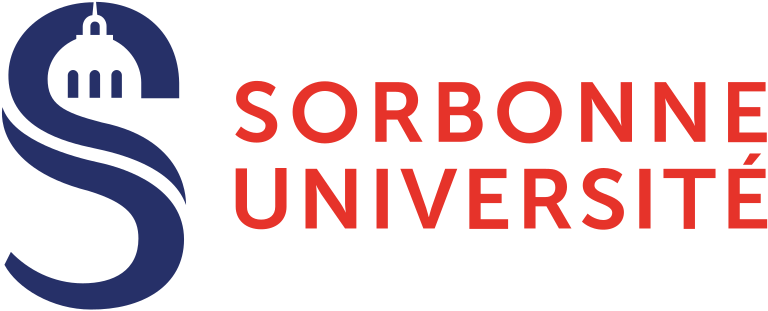}
\end{minipage}

\vspace{2em}
\noindent 
\begin{minipage}[t]{0.6\linewidth}
	{\large \textbf{Università degli Studi di Milano} } 
	
	\vspace{0.5em}
	\noindent
	\textsc{Corso di Dottorato in \\
	Scienze Matematiche \\
	Ciclo XXXIII \\
	MAT-05: Analisi Matematica}
\end{minipage}
\hfill
\begin{minipage}[t]{0.4\linewidth}\raggedleft
	{\large \textbf{Sorbonne Univeristé} }\\
	\vspace{0.5em}
	\textsc{
 	École doctorale de \\ Sciences Mathématiques \\ de Paris Centre \\
 	Specialité : Mathématiques}
 	
\end{minipage}

\vspace{1em}

\par
\noindent 
\begin{minipage}[t]{0.47\linewidth}
	\textsc{Dipartimento di Matematica }
\end{minipage}
\hfill
\begin{minipage}[t]{0.47\linewidth}\raggedleft
		\textsc{Centre d'analyse et de \\ mathématique sociales}
\end{minipage}

\vspace{1em}

\hrule

\vspace{8em}

\begin{center}
	{	\onehalfspacing
		  \huge \textbf{
	Evolution equations 
	with applications 
	to population dynamics} 
	\singlespacing
	}
\end{center}

\vspace{-2em}

\begin{center}
	\textsc{{\large 
		Doctoral thesis in conjoint program} }
\end{center}

\vspace{13em}

\noindent
\begin{minipage}[b]{0.47\linewidth}	
	Dottoranda \\
	\textbf{Elisa Affili} \\
	Matr. R12038
\end{minipage}

\vspace{2em}

\noindent 
\begin{minipage}[b]{0.47\linewidth}
	Relatore \\
	\textbf{Prof. Enrico Valdinoci}
\end{minipage}
\hfill
\begin{minipage}[b]{0.47\linewidth}\raggedleft
	Directeur \\
	\textbf{Prof. Luca Rossi}
\end{minipage}

\vspace{2em}

\noindent 
\begin{minipage}[position=b]{0.47\linewidth}
	Coordinatore del Corso di Dottorato \\
	\textbf{Prof. Vieri Mastropietro}
\end{minipage}
\hfill
\begin{minipage}[b]{0.47\linewidth}\raggedleft
	Directeur de l'école doctorale \\
	\textbf{Prof. Elisha Falbel}
\end{minipage}

\vspace{2em}

\hrule

\vspace{1em}

\begin{minipage}[b]{1\linewidth}
\begin{center}
	\textsc{Anno accademico 2019-2020}
\end{center}
\end{minipage}

\end{titlepage}

\restoregeometry

	\selectlanguage{english}
	\chapter*{Abstract}

The main topic of this thesis is the analysis of evolution equations reflecting issues in ecology and population dynamics. In mathematical modelling, the impact of environmental elements and the interaction between species is read into the role of heterogeneity in equations and interactions in coupled systems. In this direction, we investigate
three separate problems, each corresponding to a chapter of this thesis. 

The first problem addresses the evolution of a single population living in an environment with a fast diffusion line. 
From a mathematical point of view, this corresponds to a system of two coupled reaction-diffusion equations working on domains of different dimensions, which is called as in \cite{brr} a ``road-field model''. We introduce a periodic dependence of the reaction term in the direction of the fast diffusion line; in the ecological interpretation, this corresponds to the presence of more and less favourable zones for the growth of the population. 
Necessary and sufficient conditions for persistence or extinction of the population and the effects of the presence of the road are analysed through the study of a suitable generalised principal eigenvalue, originally defined in \cite{romain}. By comparison with the literature about reaction-diffusion equations in periodic media, we show that the presence of the road has no impact on the survival chances of the population, despite the deleterious effect that is expected from fragmentation.

The second investigation regards a model describing
the competition between two populations in a situation of asymmetrically aggressive interactions -- one is {\em the attacker} and the other {\em the defender}.
We derive a system of ODEs from basic principles, obtaining a modified Lotka-Volterra model relying on structural parameters as the fitness of the population and the frequency and effectiveness of the attacks. 
The evolution progresses through two possible scenarios, where only one population survives.
Then, the interpretation of one of the parameters as the aggressiveness of the attacker population naturally raises questions of controllability. With the aid of geometrical arguments we characterise the set of initial conditions leading to the victory of the attacker through a suitable (possibly time-dependant) strategy. 
Indeed, we prove that bang-bang strategies are sufficient and sometimes necessary over constant controls.
Finally, we treat a time minimization question.

The third and last part of this thesis analyses the time decay of some evolution equations with classical and fractional time derivatives. Carrying on an analysis started in \cite{SD.EV.VV}, we deal with evolution equations with a possibly mixed Caputo and classical time derivative.
By using energy methods, we prove quantitative estimates of polynomial or exponential type; the different behaviour depends heavily on the choice of the time derivative.
The decay results apply to a large class of diffusion operators, comprehending local, nonlocal, real, complex, and even nonlinear ones, of which we provide concrete examples.

	\selectlanguage{italian}
	\chapter*{Riassunto}

Il principale argomento di questa tesi è l'analisi delle equazioni dell'evoluzione che riflettono questioni di ecologia e di dinamica della popolazione. Nell'ambito della modellizzazione matematica, l'impatto degli elementi ambientali e delle interazioni tra le specie viene studiato mediante il ruolo dell'eterogeneità nelle equazioni e nelle interazioni nei sistemi accoppiati. In questa direzione, indaghiamo
tre problemi distinti corrispondenti a tre capitoli di questa tesi. 

Il primo problema riguarda l'evoluzione di una singola popolazione che vive in un ambiente con una linea di diffusione rapida. 
Dal punto di vista matematico, lo studio riguarda un sistema di due equazioni di reazione-diffusione accoppiate, che lavorano su domini di dimensioni diverse, chiamato come in \cite{brr} un modello ``campo-strada''. Introduciamo una dipendenza periodica in direzione della linea di diffusione per il termine di reazione, che nell'interpretazione ecologica corrisponde alla presenza di zone più e meno favorevoli alla crescita della popolazione. 
Le condizioni necessarie e sufficienti per la persistenza o l'estinzione della popolazione e gli effetti della presenza della strada sono analizzati attraverso lo studio di un adeguato autovalore principale generalizzato, recentemente definito in \cite{romain}. Tramite il confronto con la letteratura in mezzi periodici, si mostra che la presenza della strada non ha alcun impatto sulle possibilità di sopravvivenza della popolazione, nonostante l'effetto deleterio che ci si aspetta dalla frammentazione.

La seconda indagine riguarda un modello che descrive la competizione tra due popolazioni in una situazione di aggressione asimmetrica, in cui una popolazione aggredisce una seconda.
Deriviamo un sistema di ODE da alcune assunzioni fondamentali, ottenendo un modello Lotka-Volterra modificato che si basa su parametri strutturali come la fitness della popolazione e la frequenza e l'efficacia degli attacchi. 
L'analisi della dinamica mostra due possibili scenari, in cui una sola delle due popolazioni sopravvive.
Dopodiché, l'interpretazione di uno dei parametri come l'aggressività della prima popolazione solleva in modo naturale un problema di controllabilità. Tramite argomentazioni geometriche caratterizziamo l'insieme delle condizioni iniziali permettendo, con un'adeguata strategia eventualmente variabile nel tempo, la vittoria della popolazione che attacca. Infatti, dimostriamo che le funzioni di tipo bang-bang sono sufficienti a raggiungere l'obiettivo e talvolta sono necessarie rispetto a funzioni costanti.
Infine, trattiamo una questione di minimizzazione nel tempo.

La terza e ultima parte analizza il decadimento nel tempo in equazioni di evoluzione con una possibile derivata temporale frazionaria. Proseguendo un'analisi iniziata in \cite{SD.EV.VV}, trattiamo equazioni d'evoluzione con una combinazione di derivata temporale di Caputo e classica.
Utilizzando metodi d'energia, dimostriamo stime quantitative di tipo polinomiale o esponenziale; il diverso comportamento dipende principalmente dalla scelta della derivata temporale.
I risultati di decadimento si applicano ad una vasta classe di operatori di diffusione, comprendendone alcuni locali, non locali, reali, complessi e anche non lineari, di cui forniamo esempi concreti.

	\selectlanguage{french}
	\chapter*{Résumé}

Le sujet principal de cette thèse est l'analyse des équations d'évolution reflétant les questions d'écologie et de dynamique des populations.  En modélisation, la compréhension de l'impact des éléments environnementaux et de l'interaction entre les espèces dépend de la compréhension du rôle de l'hétérogénéité dans les équations et les interactions dans les systèmes couplés. Dans cette direction, nous étudions
trois problèmes indépendents correspondant à trois chapitres de cette thèse. 

Le premier problème concerne l'évolution d'une seule population vivant dans un environnement avec une ligne de diffusion rapide. 
L'analyse porte sur un système de deux équations de réaction-diffusion couplées, travaillant sur des domaines de dimensions différentes, qui est appelé comme dans \cite{brr} un modèle  ``champ-route''. Nous introduisons une dépendance périodique dans la direction de la ligne de diffusion pour le terme de réaction, qui, dans l'interprétation écologique, correspond à la présence de zones plus ou moins favorables à la croissance de la population. 
Les conditions nécessaires et suffisantes pour la persistance ou l'extinction de la population et les effets de la présence de la route sont analysés par l'étude de la valeur propre principale généralisée appropriée, définie pour la première fois dans \cite{romain}. Par comparaison avec des études similaires dans des environnements périodiques, nous prouvons que la présence de la route n'a aucun impact sur les chances de persistence de la population, malgré l'effet délétère attendu lié à la fragmentation.

La deuxième étude porte sur un modèle décrivant l'interaction compétitive et agressive entre deux populations. Nous dérivons un système d'EDO à partir de principes de base, en obtenant un modèle Lotka-Volterra modifié reposant sur des paramètres structurels comme la fertilité de la population et la fréquence et l'efficacité des attaques. 
L'analyse de la dynamique donne deux scénarios possibles, où une seule population survit.
Ensuite, l'interprétation d'un des paramètres comme étant l'agressivité de la première population soulève tout naturellement des questions de contrôlabilité. Grâce à des arguments géométriques, nous caractérisons l'ensemble des conditions initiales permettant la victoire de la première population avec une stratégie appropriée éventuellement dépendante du temps. En effet, nous prouvons que les stratégies de bang-bang sont suffisantes et parfois nécessaires face à des contrôles constants.
Enfin, nous traitons une question de minimisation du temps.

La troisième et dernière partie de la thèse analyse la décroissance dans le temps pour des solutions d'une classe d'équations d'évolution avec dérivées temporelles fractionnaires et classiques. Poursuivant une analyse commencée dans \cite{SD.EV.VV}, nous traitons des équations d'évolution avec une combinaison linéaire des dérivées temporelles Caputo et classiques.
En utilisant des méthodes d'énérgie, nous prouvons des estimations quantitatives de type polynomial ou exponentiel ; le comportement différent dépend fortement du choix de la dérivée temporelle.
Les résultats de la décroissance s'appliquent à une large classe d'opérateurs de diffusion, comprenant des opérateurs locaux, non locaux, réels, complexes et même non linéaires, dont nous fournissons des exemples concrets.

	\chapter*{Ringraziamenti}

Per primi vorrei ringraziare i miei relatori, Luca Rossi ed Enrico Valdinoci, senza i quali questo lavoro non sarebbe stato possibile.
Luca mi segue già da molti anni e le esperienze positive con lui sono state determinanti sulla mia scelta di intraprendere il dottorato. Ha continuato a seguirmi con attenzione e meticolosità durante questi tre anni.
Enrico fin dal primo giorno ha avuto fiducia in me e mi ha suggerito le migliori possibilità per lo sviluppo della mia carriera, incoraggiandomi, finanziandomi e aiutandomi quando il compito mi risultava troppo difficile.
A entrabi devo i miei ringraziamenti più sinceri.

Un ringraziamento speciale va anche a Serena Dipierro, che è stata molto presente nel mio dottorato, come collaboratrice e quasi come un terzo relatore, si è impegnata a coinvolgermi e promuovermi nella comunità matematica. 
Ringrazio moltissimo anche Henri Berestycki, per la grande ispirazione che mi ha dato, per i suoi preziosi consigli e per aver aiutato me e i miei relatori a formare l'accordo di cotutela. 

I would like to thank the two anonymous referees for their precious time and their nice comments on the report. I am also happy to thank the members of the defence commission, Luis Almeida, Sepideh Mirrhaimi, Fabiana Leoni and again Henri Berestycki, for accepting the task and devoting their valuable time to me.

Ringrazio anche per l'accoglienza il Dipartimento di Matematica dell'Università degli Studi di Milano, in particolare nelle persone di Vieri Mastropietro, coordinatore del corso di dottorato, e di Stefania Leonardi. Sono molto riconoscente anche a Daniela Lipari per avermi aitato a stringere l'accordo di cotutela. Un grande ringraziamento va anche ai miei colleghi dottorandi, per avermi condiviso con me gioie e dolori del dottorato e per aver reso molto più divertente il tempo a Milano.
Dedico un pensiero particolare agli altri dottorandi (ormai dottori!) e postdoc di Enrico e Serena, che mi hanno accompagnato in numerose conferenze in giro per il mondo e nei due mesi in Australia, facendomi sentire sempre a casa: Pietro, Luca, Claudia, Giorgio, Matteo, Matteo, Julien.

Je suis reconnaissante également à Sorbonne Université et au laboratoire CAMS pour m'avoir accueilli, déjà depuis mon stage de M2.
Un grand merci à Sandrine Nadal, Nathalie Brusseaux, Jean-François Venuti, Corentin Lacombe et Patricia Zizzo pour leur travail administratif, et à tous les personnes du labo  pour la belle ambiance, les discussions passionnantes, et pour m'avoir appris beaucoup sur la culture et la langue française.  
Ici aussi un grand groupe de doctorants et postdoc m'ont aidé avec leurs conseils et leur amitié, pendant ma thèse et le stage: merci à Romain, Samuel, Charles, Alessandro, Benedetta,  Federico, Julien, François, Noemi, Imke, Jérémie, Milim, José, Elisa.

Vorrei ringraziare tutti gli amici anche al di fuori dell'università che mi hanno sostenuto in questi anni di viaggi, conferenze e traslochi sfrenati tra Padova, Milano, Parigi e Stoccolma, e tutti gli amici che c'erano da molto prima che il mio dottorato iniziasse. 
Sono divisa tra la felicità di aver conosciuto così tante persone meravigliose e il rammarico di non aver passato abbastanza tempo con ciascuno.

Per la mia famiglia, il mondo dell'università e della matematica sono sempre stati estranei, ma questo non ha impedito loro di sostenermi e di avere fiducia in me. Grazie a mio fratello Andrea per avermi insegnato le sottrazioni e aver sopportato le mie domande sugli strani simboli che apparivano nei suoi libri di matematica. Grazie anche a Emanuela, Nicolò e Matteo per aver arricchito la nostra famiglia con tanta gioia e affetto. Grazie a mamma e papà per avermi dato tutto, senza mai chiedere niente. 

Ad Andrea potrei dedicare intere pagine di ringraziamenti, ma credo sia meglio farglieli di persona.

	\selectlanguage{english}
	
	\tableofcontents
	
	\mainmatter

	\chapter{Introduction}

The main motivation behind research is to enhance mankind ability to predict and keep under control natural and artificial processes.
To this purpose, mathematical models have revealed to be a very compelling instrument. A mathematical model is a simplified representation of a phenomenon through several meaningful, quantitative parameters evolving with analytical laws. Once some faithful evolution equations are established, the role of mathematics is to provide as much information as possible on the solutions, even if often only qualitative properties can be derived. That is, mathematics does not study the reality, but the language in which we read it.  

On the other side, given a model, people often find the mathematical challenges interesting in themselves. It is natural that some
questions on the mathematical tools arise, or that variations of the model are proposed and discussed. 
This way, knowledge of mathematics is expanded, and more equipment is available to write new models.

At present, the problem of climate change and environment anthropization is a great concern for humankind. In order to activate effective countermeasures against biodiversity loss, it is important to understand as deeply as possible what conditions would entail such event. These conditions depend on quantitative and qualitative properties of the environment where the species lives, on a population's resilience to changes, but also on its interaction  with other species sharing the same habitat. We still know too little about the effects that these elements and their alteration have on the survival chances of species.

This thesis is far from giving a solution to these dreadful problems but aims to give a contribution to the field of evolution equations and systems with possible application to population dynamics.  

\subsubsection{Topics and aims of the thesis}

The thesis consists of three parts, each treating a different problem.

In the first part, corresponding to Chapter \ref{ch1}, we start from a reaction-diffusion model in a periodic environment with a fast diffusion line. The aim is to find conditions entailing survival or extinction of the population and to understand the influence of the line and the environment on the dynamics. Our analysis permits a comparison with the scenario where the fast diffusion line is not present for the general case of a medium with heterogeneity in one direction.
The content of Chapter \ref{ch1} is reflects the content of the paper \cite{periodic} by the author of this thesis.

The second part, contained in Chapter \ref{ch2}, is consecrated to a model of aggressive, asymmetric competition between two populations, derived from a Lotka-Volterra system. 
The presence of the aggression term naturally leads to a control problem, where a population tries to prevail on the other using an appropriate strategy. 
Hence, once the dynamics of the system is understood, we investigate conditions for the victory of the aggressive population, which quite surprisingly is not always possible. 
Moreover it is found that, depending on the initial condition, either a bang-bang or a constant strategy leads to the desired scenario.
Chapter \ref{ch2} corresponds to the paper \cite{wars} by Serena Dipierro, Luca Rossi, Enrico Valdinoci and the author of this thesis.

The last part of this thesis deals with a more abstract and general problem; we investigate asymptotic behaviour for a class of evolution equations with both fractional and classical time derivatives.
Our setting consists of an homogeneous evolution equation working on a bounded set. 
The framework comprehends both real and complex, local and nonlocal diffusion operators, and allow us to evaluate the impact of time derivatives on the decay of solutions.  
Depending on the type of time derivative, polynomial or exponential decays are entailed.
The results of Chapter \ref{ch3} are presented in the paper \cite{decay} in collaboration with Enrico Valdinoci and the note \cite{matrix} in collaboration with Serena Dipierro and Enrico Valdinoci. 

\subsubsection{Organisation of the manuscript}

In this introductory chapter, we make the reader familiar with the problems we investigate and the framework they are enclosed in. 
Following the historical path, 
we start by a general introduction that then branches in three sections corresponding to the precise research niches of our problems. In each section, after an overview of the state of the art of the topic, we introduce the corresponding problem in details and provide precise statements of our results.

As mentioned before, the rest of the manuscript consists of three chapters, corresponding respectively and in the same order to the topics we introduce in this introduction. 
Each chapter is meant to be a self-standing script.

%
%
%

\section{General historic background}

For apparent reasons of population control and resource organisation, one of the first themes for which modelisation has been used is population dynamics. 
The first example in this sense was written by Leonardo Fibonacci in \emph{Liber Abaci} and treats the size of a population of rabbits. 
Fibonacci supposed that each couple of rabbits that are older than one month gives birth to another couple of rabbits; calling $u_n$ the size of the population at the $n-$th month, under the previous hypothesis one deduces that
\begin{equation*}
	u_{n+2}=u_{n+1}+ u_n. 
\end{equation*}
Staring with $u_0=1$, it can be deduced that $u_n$ has an exponential behaviour \cite{bacaer2011short}. 
This deduction corresponds to the reality only as long as the food is abundant for all the individuals; moreover, the relation is involved and not easy to treat.

Another discrete model was proposed by Euler in the treatise \emph{Introduction to the Analysis of the Infinite}, published in 1748 \cite{bacaer2011short}. He assumed the annual growth rate to be a fixed quantity $\alpha>0$. Then, calling $P_n$ the size of the population at the year $n$, one has that
\begin{equation*}
	P_{n+1}=(1+\alpha)P_n,
\end{equation*} 
so one derives, calling $P_0$ the population at the initial time,
\begin{equation*}
	P_n=(1+\alpha)^{n} P_0.
\end{equation*}
The sequence $\{P_n\}_{n\in\N}$ is called a geometric sequence, and its behaviour is again exponential.
Thanks to these formulae, Euler treated some problems linked to the growth of urban population and he investigated the reliability of the biblical story of the Float. However, his model involve many computations that were hard to perform before the introduction of computers.

Thomas Malthus, in his work \emph{Essay on the Principle of Population} \cite{malthus}, used a simpler relation to represent the evolution of a population size; he supposed the growth of a population to be proportional to its size, that is, the growth rate to be a fixed constant, $a>0$. 
Moreover, as simplification, he assumed the size of the population to evolve in a continuous fashion with respect to time.
With these hypothesis, the evolution of $u$ follows the law
\begin{equation}\label{eq:malthus}
 u'(t)= a u(t) \quad \text{for} \ t\geq 0.
\end{equation}
The solutions to equation \eqref{eq:malthus} are exponentials, in accordance with the result of Fibonacci.

Again in \cite{malthus}, Malthus pointed out that the growth of a population is limited by the quantity of resources. This idea was taken into the equation by Verhulst \cite{verhulst1845loi}. He considered the number $k>0$ of individuals that the environment can support indefinitely with the available resources; this is called \emph{carrying capacity} of the environment. Then, he corrected Malthus's equation \eqref{eq:malthus} with the following:
\begin{equation}\label{logistic}
		  u'(t)= a u(t)\left( 1-\frac{u(t)}{k}  \right) \quad \text{for} \ t\geq 0.
\end{equation} 
Equation \eqref{logistic} presents two equilibria: $u=0$, that is repulsive, and $u=k$, which is attractive. 
In fact, for all $u(t)<k$, one has $u'(t)>0$, while for $u(t)>k$, it holds $u'(t)<0$; in both cases, the solution tends to get closer to the value $k$.
This means that, independently of the starting condition, as long as the initial datum is positive, the population size evolves approaching the value $k$, which is the maximum number of individuals that the environment can sustain.
The logistic model is much more realistic that the previous estimates. It is considered the precursor of interesting mathematical branches, including Lotka-Volterra systems and reaction-diffusion equations.

\section{The road-field model in a periodic medium}

\subsection{Reaction diffusion equations in the literature}

One important feature that is not taken into account in the logistic equation is dependence on space. 
The first effect to take into account for a space structured model is the fact that a population is subject to \emph{dispersion}. 
This is a result of the free movement for animals and of the dispersion of seeds for plants. 
The first hypothesis in the literature was to consider the individuals to move with random brownian walk, as particles of a gas.
Without taking account reproduction, calling $u(t, x)$ the size of a population and considering it in continuous dependence on time, 
the dispersal would follow the well-known heat equation
\begin{equation}\label{eq1}
	\partial_t u - \Delta u=0.
\end{equation}
Note that when speaking of population denisties and sizes, we only consider nonnegative solutions.

The first mathematicians who added a reaction term to equation \eqref{eq1} were 
Fisher \cite{fisher} and Kolmogorov, Petrovsky and Piskunov \cite{KPP}. 
They considered a function $u(t,x)$ representing the concentration of an advantageous gene in a population; it was supposed that the population lives in a one-dimensional environment and that the individuals move randomly. Taking these hypothesis, once the gene was introduced, it spreads according to the equation
\begin{equation}\label{eq:KPP}
\partial_t u - \partial_{xx}^2 u = f(u)
\end{equation}
where $f$ is a function such that
\begin{equation}\label{0225}
	f(0)=f(1)=0,
\end{equation}
moreover it is \emph{monostable}, that is,
\begin{equation*}
	 \quad f(u)>0 \ \text{for} \ u\in(0,1),
\end{equation*}
 and respects the condition called \emph{KPP} hypothesis 
\begin{equation}\label{0024}
 	f(u) < f'(0)u.
\end{equation}
The function $f$ represents the birth-death rate of individuals carrying the gene. The fact that $f(0)=0$ is a very natural assumption: if no individuals are present, no new individual is generated. On the other hand, the choice $f(1)=0$ suggests a saturation at the size $u=1$.
The hypothesis \eqref{0024} reflects the fact that the growth rate decreases as the size of the population grows, as it is the case for the logistic equation \eqref{logistic}. 
Actually, Fisher supposed $f(u)=au(1-u)$ for $a>0$, which is exactly the nonlinearity proposed by Verhulst, while Kolmogorov, Petrovsky and Piskunov selected $f(u)=a u(1-u)^2$.

For a large class of initial data, among which the Heaviside functions, the solutions to \eqref{eq:KPP} asymptotically converge to a function of the shape
\begin{equation}\label{sol}
u(t,x )=U(z) \quad \text{for} \ z=x+ct.
\end{equation}
Solutions of the form \eqref{sol} are called \emph{travelling waves} and the quantity $c$ is called \emph{speed of propagation} of the travelling wave. The travelling wave found in \cite{fisher} and \cite{KPP} has speed corresponding to $c_{KPP}=2\sqrt{ f'(0)}$; actually, a travelling wave exists for all $c\geq c_{KPP}$ and $c_{KPP}$ corresponds to the minimal speed.

The main questions addressed in \cite{fisher} and \cite{KPP} have been later asked for larger and larger class of nonlinearites. These questions concerns the existence of stationary solutions, the existence of travelling fronts and the asymptotic speed of propagation for the Cauchy problem.

For the sake of completeness, we must here name other two important settings.
In \cite{fife} and in \cite{weinberger2} for the multidimensional case, Fife and McLeod and Aronson and Weinberger treated equation \eqref{eq:KPP} in the case of a function $f$ satisfying the hypothesis \eqref{0225} and such that there exists a value $\theta$ for which
\begin{equation}\label{0044}
	  f(u)<0 \quad \text{if} \ u\in(0,\theta), \qquad f(u)>0 \quad \text{if} \ u\in(\theta, 1).
\end{equation}
A function satisfying \eqref{0044} is called \emph{bistable}, from the fact that the related equation has two attractive states, $0$ and $1$. This type of nonlinearity is particularly interesting because it embodies an important phenomenon in population dynamics, called the \emph{Allee effect} from the name of the scientist who discover it in the '30s.  It happens that in social animals, aggregation increases the survival rate of individuals; therefore, when the size of a population is under a certain threshold, the growth rate is negative; when the group size passes the threshold, the growth rate becomes positive.

A third important setting is the \emph{combustion} case, in which there exists a quantity $\theta\in(0,1)$ such that
\begin{equation*}
	f(u)=0 \quad \text{if} \ u\in[0,\theta], \qquad f(u)>0 \quad \text{if} \ u\in(\theta, 1).
\end{equation*}
This type of nonlinearity is used for ignition models, where to activate the combustion process the  temperature must pass a threshold.

As a matter of fact, Aronson and Weinberger investigated the equation
\begin{equation}\label{aw}
	\partial_t u- \Delta u= f(u) \quad \text{for} \ x \in\R^n
\end{equation} 
and asked under which conditions on the function $f$, other than \eqref{0225}, and on the initial datum $u_0$ one has \emph{invasion} or \emph{spreading}, that is,
\begin{equation*}
 u(t,x) \overset{t\to+\infty}{\longrightarrow} 1 \quad \text{locally uniformly in} \ x.	
\end{equation*}
The opposite behaviour is called \emph{extinction}, and it occurs when
\begin{equation*}
u(t,x) \overset{t\to+\infty}{\rightarrow} 0 \quad \text{uniformly in} \ x.	
\end{equation*}
We point out that for extinction a uniform convergence is required, otherwise, in some scenarios, one could have a positive mass escaping further and further in space as $t$ goes to infinity.
The authors found that for a compactly supported initial datum which is ``sufficiently large'' (depending on the nonlinearity), invasion occurs if and only if 
\begin{equation*}
	\int_0^1 f(x)dx>0.
\end{equation*}  
Let us give more details on the minimal requirements for the initial datum.  
In the monostable case, it is sufficient for $u_0$ to be greater than a positive constant in $(0,1)$ in a large enough ball. Moreover, if $f'(0)>0$, then  all solutions issued from a non zero, non negative initial datum  converges to 1 as $t$ goes to infinity; this is called \emph{hair trigger effect}.
In the bistable and monostable cases, the positive constant is necessarily greater than the threshold $\theta$.

Equation \eqref{eq:KPP} was the first example of a whole class of PDEs, the reaction-diffusion equations. From the initial works \cite{fisher, KPP, fife, weinberger2}, the literature on reaction-diffusion equations and the study on travelling waves have flourished. What is present here is a circumscribed niche, which is handy to provide context to our work.

\subsubsection{Reaction-diffusion equations in periodic media}

One of the other natural applications of equations \eqref{eq:KPP} and \eqref{aw} is of course population dynamics. Skellam \cite{skellam} was one of the firsts to study the effects of random dispersion on a population subject to  the malthusian law, after noticing that the framework given by \cite{fisher} and \cite{KPP} could be adapted to this problem.

In the optic of studying the survival and the distribution of a population in space, a homogeneous environment is not satisfying and one expects the growth of the population to vary according to the habitat conditions. 
On the other hand, from a mathematical point of view, heterogeneity in the nonlinearity creates great difficulties. Many new techniques were required to overcome these obstacles.

A first analysis was carried out by Shigesada, Kawasaki and Teramoto 
in \cite{shigesada1986traveling, shigesada1997biological}.
The authors observed that natural environments are a mosaic of different habitats, such as forests, meadows, brush, cultivated fields and villages. 
This led them to consider an environment which consists of two periodically alternating homogeneous habitats, one favourable, $E^+$, and one unfavourable, $E^-$, for the considered species. 
The heterogeneity of the living conditions is reflected by the birth-death rate, which they chose to be
\begin{equation*}
	f(x, u)= \left\{  
	\begin{array}{ll}
	u(\mu^+ -u), & \text{in} \ E^+, \\
	u(\mu^- -u), & \text{in} \ E^-,
	\end{array}	
	\right.
\end{equation*}
for some $\mu^+>\mu^-$.
Moreover, they also consider possibly varying diffusivity, hence they took
\begin{equation*}
A(x)= \left\{  
\begin{array}{ll}
A^+, & \text{in} \ E^+, \\
A^-, & \text{in} \ E^-.
\end{array}	
\right.
\end{equation*}
This is due to the observation of increased speed in unfavourable environments; hence we expect $A^+<A^-$ for a real population.
Then, the authors studied in \cite{shigesada1986traveling} the equation
\begin{equation}\label{0325}
	\partial_t u - \nabla \cdot (A(x) \nabla u) = f(x,u) \quad \text{for} \ x\in\R^n.
\end{equation}
This is known as the \emph{patch model}; they investigated long time behaviour, convergence to travelling fronts and propagation speeds.
Actually, since $u=1$ is no longer an equilibrium for equation \eqref{0325}, we have to modify our definition for species survival; from now on, we intend that  \emph{persistence} occurs if $u(x,t)$ approaches a non null stationary solution locally uniformly as $t$ tends to infinity. 
 
By making use of numerical simulations, it was found that the stability of the trivial solution $u=0$ plays a key role in determining if the population survives or not. 
It was already known (see \cite{coddington1955theory}) that a negative or positive sign of the principal eigenvalue resulting from the linearisation around $u=0$ entails respectively stability or instability of the $0$ solution. 
In \cite{shigesada1986traveling}, it was shown numerically that the stability of the trivial solution entails extinction, while its instability causes persistence of the population.
The authors also studied the sign of the eigenvalue depending on the values of $L$, the measures of $E^+$ and $E^-$ and the values of the parameters; this was possible because of the simplicity of the framework. 

Equation \eqref{0325} was later considered in \cite{kinezaki2003modeling} and \cite{bhroques} for general $A(x)$ and $f(x,u)$ depending on $x$ in a continuous fashion and perdiodically of period $L$ for some $L\in\R^n$.
In this second article, Berestycki, Hamel and Roques comprehended that the extinction or persistence of the population depends on the sign of a periodic eigenvalue $\lambda_p(-\mathcal{L}', \R^n)$, that is the unique real number such that the problem
\begin{equation}\label{sys:L_RN_p}
\left\{
\begin{array}{ll}
\mathcal{L'}(\psi) + \lambda \psi = 0, & x\in\R^n, \\
\psi> 0, &  x\in\R^n, \\
|| \psi ||_{\infty}=1, \\
\psi \ \text{is periodic in $x$ of periods $L$},
\end{array}
\right.
\end{equation}
where $\mathcal{L'}$ is given by
\begin{equation*}\label{def:mathcal_L'}
\mathcal{L'}(\psi):=  \nabla \cdot(A(x) \nabla  \psi)  + f_u(x,0)\psi,
\end{equation*}
has a solution $\psi_p\in W_{loc}^{2, 3}(\R^n)$. It was proved that when
$\lambda_p(-\mathcal{L}', \R^n)\geq 0$ extinction occurs. On the other hand,
when $\lambda_p(-\mathcal{L}', \R^n)<0$ there is persistence; moreover, there exists a unique stationary solution to \eqref{0325}, that is periodic of period $L$, and attracts all the solutions starting from a non negative, non zero bounded initial datum.

The studies on the patch model \cite{shigesada1986traveling, shigesada1997biological} and the ones on periodic media \cite{kinezaki2003modeling, bhroques} evidenced also the effect of fragmentation on the survival chances of a population. It was found that $\lambda_p(-\mathcal{L}', \R^n)$ decreases as the homogeneity increases, that is, a species has better survival chances when the environment is less fragmented.

\subsubsection{The case of a changing climate}

A new aspect that one may consider while studying ecological problems is a changing climate. If the environment changes in time, so does the fitness of a population. In this paragraph, we are going to analyse the difficulties produced by the new type on nonlinearity and how it has been overcome. 

A 1-dimensional model for population persistence under climate change was first proposed by Berestycki, Diekmann, Nagelkerke and Zegeling in
\cite{berestycki2009can}. 
The authors first imagined that a population lives in a favourable region enclosed into disadvantageous environment. 
Assuming that a global warming is in place, and that the population lives in the Boreal Emisphere,
the authors imagined that the favourable region moves to the north, so that for every favourable area lost in the South, an equivalent favourable area is gained in the North. 
The resulting equation is 
\begin{equation*}
\partial_t u -  \partial_{xx}^2 u=f(x-ct,u) \quad \text{for} \ x\in\R.
\end{equation*}
Later, in \cite{br2}, Berestycki and Rossi presented a model for climate change in $\R^n$ and for a larger class of nonlinearites; they dealt with equation 
\begin{equation}\label{1709}
\partial_t u -  \Delta u=f(x-ct e,u) \quad \text{for} \ x\in\R^n,
\end{equation}
with $e$ a direction in $\mathbb{S}^{n-1}$ and $f: \R^n\times \R^+ \to \R$. 
The authors focused on solutions
in the form of a travelling waves $u(x,t)=U(x-cte)$ which solve the equation
\begin{equation}\label{eq:cc}
\partial_t U -  \Delta U- c\,e\cdot \nabla U=f(x,U) \quad \text{for} \ x\in\R^n.
\end{equation}
This second equation is more treatable: in fact, the dependence in time of the nonlinearity, which poses a lot of problems, is transformed into a transport term; now, the equation has a nonlinearity depending only on space, and techniques for this type of heterogeneity are more familiar. 

The main question is if the population keeps pace with the shifting climate, that is, if a large enough group is able to migrate with the same speed of the climate. The answer to this question is positive if a solution to \eqref{eq:cc} exists; as happened for the periodic equation \eqref{0325}, this depends on the sign of the principal eigenvalue coming from the linearisation in $0$.

\subsubsection{The road-field model}

Spatial heterogeneity in natural environments may be the consequence not only of the diversity of the habitats, but also of the presence of obstacles or fast diffusion channels that affects the fitness and the mobility of individuals. 

In recent years, humans activity has caused drastic changes in the environment, causing different species to become invasive in areas they were not present \cite{shigesada1997biological}. In the case of the Processionary pine tree caterpillar, the diffusion in France has been even faster than anticipated. It has been observed that the insect was incidentally transported by humans from town to town, and from these settlements it spread in the surroundings \cite{robinet2012human}.
This in not the only example of ecological diffusion acceleration by fast diffusion lines. In Western Canadian Forest, GPS observations on wolves proved  that the animals exploit seismic lines, that are straight roads used by the oil companies to test reservoirs, to move faster and therefore to increase their probability of meeting a prey \cite{mckenzie2012linear}. 

Roads play a strong role also in the spreading of epidemics. The ``black death'' plague in the 14th century was one of the most devastating epidemics known in Europe. It is known that the plague was transported by animals and humans along the commercial trade line of the silk road, and from that  spread all over Europe. More recently, a similar effect has been conjectured for the COVID-19 infection. By tracing the spreading in Northen Italy in early March 2020, it was found that the diffusion occurred first along highways and then spread in the surrounding territory \cite{gatto2020spread}. 

Inspired by this behaviour, Berestycki, Roquejoffre and Rossi proposed in \cite{brr} a model of spreading in an environment presenting a fast diffusion channel. As a simplification, they considered the channel to be a straight line in $\R^2$, the $x$ axis $\R\times \{ y=0 \}$.
Their idea was to split the population into two groups; the first one, of density $u$, occupies the one dimensional environment $\R\times \{ y=0 \}$ representing the road, and the second one, of density $v$, occupies the surrounding territory; by symmetry, they considered just one half of the plan, thus $\Omega:=\{(x, y) \in\R^2 : y>0 \}$, which they called ``the field''. These two groups continuously exchange along the road: a fraction $\nu>0$ of the population in $\Omega$ at $y=0$ passes in the road, and a fraction $\mu>0$ of the population in the road passes in the field. 
The diffusivity is different in the two environments; its values are $D$ on the road and $d$ on the field, both positive. 
Moreover, it is supposed that population reproduces only in the field and that the environment is homogeneous; the corresponding function $f$ is required to satisfy \eqref{0225}, $f'(0)>0$ and a stronger version of the KPP hypothesis,  that is
\begin{equation*}
	v \mapsto \frac{f(v)}{v} \quad  \text{is decreasing}.
\end{equation*}
The resulting system, called \emph{road-field model}, is
\begin{equation}\label{sys:rf} 
\left\{
\begin{array}{ll}
\partial_t u(x,t) - D \partial_{xx}^2 u (x,t) = \nu v (x,0,t) - \mu u(x,t), &  x\in \R, \ t > 0, \\
\partial_t v(x,y,t) - d \Delta v (x,y,t)= f(v), & (x,y) \in \Omega, \ t>0, \\
-d \partial_y v(x,0,t) = -\nu v(x,0,t) + \mu u(x,t), & x \in \R, \ t>0.
\end{array} \right.
\end{equation}
The authors of \cite{brr} found that invasion occurs for any non negative, non zero initial datum, so the hair trigger effect holds; solutions converge to the unique steady state $\left(\frac{\nu}{\mu},1  \right)$. Moreover, they studied spreading speeds and found that it is enhanced by the presence of the road.
In a second paper \cite{BRR2}, the same authors investigated system \eqref{sys:rf} with a transport term and a reaction term on the line. 

Many variations of the road-field model were proposed. 
In \cite{pauthier2015uniform, pauthier2016influence}, the system was modified by introducing nonlocal exchanges between the road and the field. 
The case of a general nonlocal diffusion has been treated in \cite{berestycki2015effect, berestycki2014speed}. 
Different geometric settings have also been considered; in  \cite{rossi2017effect}, the model was extended in higher dimensions. 
For a complete list, we refer to the chapter in \cite{tellini2019comparison} by Tellini.

Treating system \eqref{sys:rf} poses some difficulties because of the interaction between functions living in different dimensions and the unusual boundary condition. Adding some heterogeneity in space increases the difficulties. This is why very few studies of this type have carried on so far, a part from an article by Giletti, Monsaingeon and Zhou \cite{giletti2015kpp}, where the authors considered the case of exchanges terms depending periodically on $x$.

Recently, Berestycki, Ducasse and Rossi introduced in \cite{romain} a new generalised principal eigenvalue fitting road-field models for a possibly heterogeneous reaction term.
Hence, they considered the system
\begin{equation*}
\left\{
\begin{array}{ll}
\partial_t u(x,t) - D \partial_{xx}^2 u (x,t) -c \partial_x u(t,x)= \nu v (x,0,t) - \mu u(x,t), &  x\in \R, \ t > 0, \\
\partial_t v(x,y,t) - d \Delta v (x,y,t)-c \partial_x u(t,x)= f(x,y,v), & (x,y) \in \Omega, \ t>0, \\
-d \partial_y v(x,0,t) = -\nu v(x,0,t) + \mu u(x,t), & x \in \R, \ t>0.
\end{array} \right.
\end{equation*}
Calling
\begin{equation*}\label{sys:operators}
\left\{
\begin{array}{l}
\mathcal{R}(\phi, \psi):=D \phi''+c \phi'+\nu {\psi}|_{y=0}-\mu \phi, \\
\mathcal{L}(\psi):= d\Delta \psi +c \partial_x \psi -f_v(x,y,0)\psi, \\
B(\phi, \psi):=d \partial_y {\psi}|_{y=0}+\mu \phi- \nu {\psi}|_{y=0},
\end{array}
\right.
\end{equation*}
this eigenvalue is defined as 
\begin{equation}\label{def:lambda1_S_Omega}
\begin{split}
\lambda_1( \Omega)=\sup \{ \lambda \in \R \ : \ \exists (\phi, \psi)\geq (0,0), \ (\phi, \psi) \not\equiv(0,0), \ \text{such that} \\ \mathcal{L}(\psi) + \lambda \psi \leq 0 \ \text{in} \ \Omega, \ \mathcal{R}(\phi, \psi) +\lambda \phi \leq 0  
\ \text{and} \ B(\phi, \psi)\leq 0 \ \text{in} \ \R \},
\end{split}
\end{equation}
with $(\phi, \psi)$ belonging to $W_{loc}^{2,3}(\R)\times W_{loc}^{2,3}(\overline{\Omega})$. Together with the definition, many interesting properties and bounds were studied.

Thanks to that, the same authors were able to investigate the case of 
a favourable ecological niche, possibly facing climate change, in \cite{econiches}. It was proven that the sign of $\lambda_1( \Omega)$ characterises
the extinction or the persistence of the population; moreover, comparing the results with  the ones found for the model without the road, in the absence of climate change a deleterious effect of the road on the survival chances was found. 
On the other hand, if the ecological niche shifts, the road has in some cases a positive effect on the persistence.

\subsection{A KPP model with a fast diffusion line in a periodic medium}

We are now ready to introduce in details the first problem dealt with in this thesis. 
We are going to investigate a road-field model in a periodic medium. 
This problem combines the interests of studying the effect of a fast diffusion line with the one of treating a heterogeneous nonlinearity, that, as we pointed out before, reflects a natural territory in a more realistic way than a homogeneous term. From a technical point of view, it also combines the difficulties of the two settings.

\subsubsection{The model}

We have already presented the road-field model. 
In our problem, we treat a road-field system with possible climate change and with a reaction term depending on the spatial variable $x$; in particular, we will focus on the case of periodic dependence. 
There is no dependence in the variable $y$, the heterogeneity in that direction is only due to the presence of the road.

Keeping the notation used so far, the system we investigate reads
\begin{equation}\label{sys:fieldroad}
\left\{
\begin{array}{lr}
\partial_t u-D  u '' -c u' - \nu  v|_{y=0} + \mu u= 0,   & x\in \R,  \\
\partial_t v -d \Delta v-c\partial_x v  =f(x,v),  & (x, y)\in \Omega, \\
-d  \partial_y{v}|_{y=0} + \nu v|_{y=0} -\mu u=0, & x\in\R.
\end{array}
\right.
\end{equation}
Recall that $D$, $d$, $\nu$, $\mu$ are positive constants and $c\geq 0$.
The function  $f:\R\times \R_{\geq 0}\to \R $
is always supposed to be $\mathcal{C}^{1}$ in $x$, locally in $v$, and Lipschitz in $v$, uniformly in $x$; moreover we suppose that the value $v=0$ is an equilibrium, that is
\begin{equation}\label{hyp:0}
f(x,0)=0, \quad \text{for all} \ x\in \R,
\end{equation}
and that 		
\begin{equation}\label{hyp:M}
\exists M>0 \ \text{such that} \ f(x, v)<0 \quad \text{for all} \ v>M \ \text{and all} \ x\in \R, 
\end{equation}
which indicates that there is a saturation level.
We will derive some inequalities on the generalised principal eigenvalue of \eqref{sys:fieldroad} for the general case of $f$ respecting these hypothesis and 
$c$ possibly nonzero.

The characterisation of extinction or persistence of the species is addressed in the case of $c=0$ and $f$ a periodic function, reflecting the periodicity of the environment in which the population diffuses, as we require with the forthcoming hypothesis. 
We will analyse the  case of a KPP nonlinearity, that is, we require that 
\begin{equation}\label{hyp:KPP}
\frac{f(x,s_2)}{s_2}< \frac{f(x,s_1)}{s_1} \quad  \text{for all} \ s_2>s_1>0 \  \text{and all} \ x\in\R.
\end{equation}
Then, we suppose that there exists $\ell> 0$ such that
\begin{equation}\label{hyp:per}
f(x+\ell, s)=f(x,s) \quad \text{for all} \ s >0 \  \text{and all} \ x\in\R.
\end{equation}


To study the effect of the line of fast diffusion, we will compare the behaviour of \eqref{sys:fieldroad} to the one of the system
\begin{equation}\label{sys:symmetric}
\left\{
\begin{array}{ll}
v_t-d\Delta v - c\partial_x v= f(x,v),  & (x,y)\in\Omega,\\
-\partial_y v|_{y=0} =0, & x\in\R,
\end{array}
\right.
\end{equation}
whose solution is a function $v(x,y)$ that can be extended by symmetry to the whole plane.
It is natural to consider system \eqref{sys:symmetric} as the counterpart of system \eqref{sys:fieldroad} in the case without the road, since it presents the same geometry, including the same boundary condition, exception made for the exchange terms that are in place for the case of a fast diffusion channel.

\subsection{Our results}

We are now ready to present the main results of this part of the thesis.

\paragraph{The case of a periodic $f(x,v)$.}

Here, we consider the case of a nonlinearity that respects the KPP hypothesis and is periodic in the direction of the road. Moreover, we consider $c=0$.

We begin by the following result on the long time behaviour for solutions of system \eqref{sys:fieldroad}. As already seen for similar problems, the key point lies in the stability of the $0$ solution. This is linked to the sign of the generalised principal eigenvalue for the road-field model, that we have defined in \eqref{def:lambda1_S_Omega}. With this notation, we have the following:

\begin{theorem} \label{thm:char}
	Let $f$ satisfy \eqref{hyp:0}-\eqref{hyp:per} and $c=0$.
	Then the following holds:
	\begin{enumerate}
		\item if $\lambda_1( \Omega)\geq 0$, then  extinction occurs.
		\item if $\lambda_1(\Omega)<0$, then persistence occurs and the  positive stationary solution $(u_{\infty}, v_{\infty})$ is unique and periodic in $x$.
	\end{enumerate} 
\end{theorem}

Next, we compare the behaviour of solutions of the system \eqref{sys:fieldroad} with the ones of \eqref{sys:symmetric}, or, equivalently, after  extension by symmetry to the whole plane, of
\begin{equation}\label{eq:bhroques}
	\partial_t v - d \Delta v = f(x,v), \quad \text{for} \ (x,y)\in\R^2.
\end{equation} 
Recalling the results of \cite{bhroques}, we know that the persistence or extinction of a population for a periodic equation in the whole $\R^2$ depends on the sign of the periodic eigenvalue $\lambda_p(-\mathcal{L}, \R^2)$, that was defined in \eqref{sys:L_RN_p} for a general case.
We obtain the following:

\begin{theorem}\label{thm:comparison}
	Assume $f$ fulfils hypotheses \eqref{hyp:0}-\eqref{hyp:per} and let $c=0$. Then:
	\begin{enumerate}
		\item if  $\lambda_p(-\mathcal{L}, \R^2)<0$, then $\lambda_1( \Omega)<0$, that is, if persistence occurs  for the system ``without the road'' \eqref{eq:bhroques}, then it occurs also for system ``with the road'' \eqref{sys:fieldroad}.
		\item if $\lambda_p(-\mathcal{L}, \R^2)\geq 0$, then $\lambda_1( \Omega)\geq 0$, that is, if extinction occurs for the system ``without the road'' \eqref{eq:bhroques}, then it occurs also for system ``with the road'' \eqref{sys:fieldroad}.
	\end{enumerate}	
\end{theorem}

Theorem \ref{thm:comparison} asserts that the road has no negative impact on the survival chances of the population in the case of a medium depending periodically on with respect to variable in the direction of the road.  

We recall the fact that fragmentation lowers the survival possibilities of a species (see \cite{bhroques, shigesada1986traveling}); also, even if we are not in the framework of an ecological niche, we remember from \cite{econiches} the fact that a road has a negative impact in the setting without climate change. For those reasons, the result in Theorem \ref{thm:comparison} may be somehow unexpected.
However, despite the fact that no reproduction takes place on the road, in the case of periodic media the presence of the fast diffusion channel does not interfere with the long time behaviour of the population, which depends only on the  environment of a periodicity cell.
As seen in \cite{bhroques}, where the dependence of persistence on the amplitude of fragmentation was studied, if the favourable zones are sufficiently large, the population will eventually spread in all of them; the presence of the road does not cause loss of favourable environment and consequently of persistence chances.
However, we expect the spreading speed to be influenced by the presence of the road, as it has been already proven in the case of homogeneous environment.

We point out that Theorem \eqref{thm:char} completes and is in accordance with the results on long time behaviour found in \cite{brr} for a homogeneous reaction function, which we can consider as a particular case of periodicity, satisfying a positive KPP request (thanks to the hypothesis $f'(0)>0$). In \cite{brr}, Theorem 4.1 states the convergence of any positive solution to the unique positive stationary solution of the system.  Since it is well known that for the homogeneous case it holds $\lambda_1(-\mathcal{L}, \R^2)=- f'(0)$, the hypothesis gives that $\lambda_1(-\mathcal{L}, \R^2)<0$ and, as a consequence of Theorem \ref{thm:char}, that persistence occurs. 
Instead if $f'(0)\leq0$, then we would be in the first case of Theorem \ref{thm:char}, yielding extinction of the population.

\paragraph{Effects of amplitude of heterogeneity}

One may ask if the presence of a road may alter the complex interaction between more favourable and less favourable zones; in particular, one could wonder if this could penalise the persistence, since it was shown that populations prefer a less fragmented environment. Nevertheless, owing from Theorem \ref{thm:comparison} that the road has no effect on the survival chances of the species, we can recover all the results on the effect of fragmentation.

Take a parameter $\alpha>0$ and consider system \eqref{sys:fieldroad} with nonlinearity
\begin{equation}\label{1421}
\tilde{f}(x,v)=\alpha f(x,v).
\end{equation}
To highlight the dependence on $\alpha$, we will call $\lambda_1(\Omega, \alpha)$ the  generalised principal eigenvalue defined in \eqref{def:lambda1_S_Omega} with nonlinearity $\tilde{f}$.
As a direct consequence of our Theorem \ref{thm:comparison} and of Theorem 2.12 in \cite{bhroques}, we have the following result on the amplitude of heterogeneity:

\begin{corollary}
	Assume $\tilde{f}$ is defined as in \eqref{1421}, $f$ satisfies \eqref{hyp:0}-\eqref{hyp:per}, and $c=0$. Then:
	\begin{enumerate}
		\item if $ \int_{0}^{\ell} f_v(x,0)>0$, or if $ \int_{0}^{\ell} f_v(x,0)=0$ and $f\neq 0$, then for all $\alpha >0$ we have $\lambda_1(\Omega, \alpha  )<0$.
		\item if $ \int_{0}^{\ell} f_v(x,0)<0$, then  $\lambda_1(\Omega, \alpha )>0$ for $\alpha$ small enough; if moreover there exists $x_0\in[0,\ell]$ such that $f_v(x_0,0)>0$, then for all $\alpha$ large enough $\lambda_1(\Omega, \alpha )<0$.
	\end{enumerate}
\end{corollary}

\paragraph{A climate change setting for a general $f(x,v)$.}

We consider now a general nonlinearity that depends on the variable in the direction of the road. We stress the fact that we do not suppose any periodicity, but the case of a periodic $f$ is a particular case of this setting. Moreover, the following results are done in the general framework of a possible climate change, so the parameter $c$ may be different from $0$.

Comparison between the systems with and without the road, in the general case, are done through comparison between $\lambda_1(\Omega)$ and the generalised principal eigenvalue of system \eqref{sys:symmetric}, given by
\begin{equation}\label{lambda:L_Omega}
\begin{split}
\lambda_1(-\mathcal{L}, \Omega)=\sup \{ \lambda \in \R \ : \ \exists \psi \geq 0, \psi \not\equiv 0 \ \text{such that} \\ 
\mathcal{L}(\psi) + \lambda \psi \leq 0 \ \text{on} \ \Omega, \ -\partial_y \psi|_{y=0}\leq 0 \ \text{on} \ \R \}
\end{split}
\end{equation}
for $\psi\in W_{loc}^{2,3}(\Omega)$. With this notation, we have the following:

\begin{theorem}\label{thm:ineq}
	Assume $\lambda_1(-\mathcal{L}, \R^2)$ as in \eqref{lambda:L_Omega} and $\lambda_1(\Omega)$ as in \eqref{def:lambda1_S_Omega}; then $\lambda_1(-\mathcal{L}, \R^2) \geq \lambda_1(\Omega)$.
\end{theorem}

In the special case $c=0$, some information on the relations between $\lambda_1(-\mathcal{L}, \R^2)$ and $\lambda_1(\Omega)$  was already available in \cite{econiches}: Proposition 3.1 gives that if $\lambda_1(-\mathcal{L}, \R^2)\geq 0$ then $\lambda_1(\Omega)\geq 0$. Thanks to that and Theorem \ref{thm:ineq}, the following result holds:

\begin{corollary}\label{thm:ineq2}
	If $c=0$, we have $\lambda_1(-\mathcal{L}, \R^2)<0$ if and only if $\lambda_1(\Omega)<0$.
\end{corollary}

As already pointed out in \cite{romain}, even for $c=0$ it is not true that $\lambda_1(-\mathcal{L}, \R^2) =\lambda_1(\Omega)$. In fact, it has been found that $\lambda_1(\Omega) \leq \mu$, while playing with $f$ one can have $\lambda_1(-\mathcal{L}, \R^2)$ as large as desired. However, the fact that they have the same sign reveals that they are profoundly linked.

\subsubsection{Perspectives}

The next problem to tackle for system \eqref{sys:fieldroad} in a periodic medium regards  the existence of travelling fronts and the study of their speed in all the direction of the plane. 
We point out that, with respect to the classical case, there are great difficulties linked to the anisotropy of the space, due both to the road and to the periodicity of the medium. 
An acceleration effect due to the presence of the road is expected to be found when $D>d$; however, the repercussions of the periodicity of the medium on the spreading speed in a general direction is hard to predict.

We also mention that it would be nice to extend the current results to the case of heterogeneous exchange terms, periodic in $x$, as already treated in \cite{giletti2015kpp}. The key point for attacking that problem is in the generalisation of the definition of $\lambda_1(\Omega)$ for non homogeneous coefficients.

\section{A new model for aggressive competition}

\subsection{Lotka-Volterra models: a literature overview}

Another issue that is overlooked in the logistic equation is the interaction of a species with the other ones living in the same environment. In the '20s, Lotka \cite{lotka} and Volterra \cite{volterra} observed independently some curios transitory oscillations in the concentration of chemicals during a reaction and in the population sizes of fishes. 

They formulated the following model; let $u$ be the quantity of a species of plants present in the environment and $v$ the size of a population of herbivores. It is supposed that the plants have a constant growth rate at all times, $a>0$. The herbivorous feed exclusively on the observed plant and have limitless appetite. 
The consumption of plants eaten by the animals is supposed to depend on the probability of meeting of the two, represented by $uv$; the actual loss of the plants is $-buv$ and the gain for the herbivores is $duv$ with $b>d>0$, owning the fact that some plants could be torn but not consumed. 
Moreover, as in the malthusian equation, the increase of a population is suppose to depend on its size. 
It is also supposed that the environment conditions are stable and than no mutation in the behaviour of the two species is possible.

 Then, the model reads
\begin{equation}\label{model:lv}
\left\{
\begin{array}{llr}
\dot{u}&= au-buv, & {\mbox{ for }}t>0,\\
\dot{v}&= -cv+duv, & {\mbox{ for }}t>0.
\end{array}
\right.
\end{equation}
This system has two equilibria, $(0,0)$ and $\left(  \frac{c}{d},\frac{a}{b}  \right)$. If the initial datum is any point of positive coordinates distinct from the equilibrium, the population sizes oscillate in time, running on a closed curve on the phase portrait.

\subsubsection{Competitive Lotka-Volterra models}

Since the pioneer works, many studies on the interaction between populations were carried out. In particular, after the studies of Gause \cite{gause1984ecology}, another model has been employed to investigate the dynamics between two populations in competition, that is, exploiting at least partly the same resources.
We propose here its construction using the example of two population of squirrels, the grey one and the red one, following the work in \cite{okubo1989spatial}.
These two species, one of the two recently introduced in Britain, both inhabit hardwood forests and rely on the same resources to live. Keeping in mind the derivation of the logistic equation, we realize that the resource term in this scenario depends on the size of both population.
Moreover, we take into consideration the fact that, due to the social organisation and sometimes the segregation between competing species, the presence of individuals of the rival group may obstruct food collection; if this is the case, there is an additional decrease of the available resources for both population.
Adding these corrections to the logistic of both groups,  the Lotka-Volterra competitive system reads
\begin{equation}\label{lv}
\begin{cases}
\dot{u}=a_u  u\left(1-\displaystyle
\frac{u+\alpha_{uv} v}{k_u}  \right), & t>0,\\
\dot{v}
=a_v v\left(1- \displaystyle\frac{v+\alpha_{vu} u}{k_v}  \right), & t>0,
\end{cases}
\end{equation}
where~$a_u$,~$a_v$,~$\alpha_{uv}$,~$\alpha_{vu}$,~$k_u$ and~$k_v$ are nonnegative real numbers.
The coefficients $a_u$ and $a_v$ are the intrinsic growth rates of the two population; $k_u$ and $k_v$ represent the carrying capacities of the environment for the two groups.
The coefficients~$\alpha_{uv}$ and~$\alpha_{vu}$ represent the competition between individuals of different species, and indeed they appear multiplied by the term~$uv$, which represents a probability of meeting. 
Taking the example of the squirrels, we expect that $\alpha_{uv}, \alpha_{vu} >1$. 
However, for other couple of populations relying on only partially overlapping food sets, one could have also $\alpha_{uv}, \alpha_{vu} \leq 1$. If finally the first population feeds on a subset of the resources of the second one, it would be $\alpha_{uv} \geq1$ and $\alpha_{vu} <1$.

For the sake of completeness, we recall that in the case of species mutually benefiting from the presence of the other, which is not part of the competitive framework, the dynamics prescribes negative values for $\alpha_{uv}$ and $\alpha_{vu}$.

The dynamics of system \eqref{lv} depends indeed on the values of the interspecific competition terms:
if~$\alpha_{uv}<1<\alpha_{vu}$, then the first species~$u$
has an advantage over the second one~$v$
and will eventually prevail; if~$\alpha_{uv}, \ \alpha_{vu} >1$, then the first population that penetrates the environment (that is, the one that has a greater size at the initial time) will persist while the other will extinguish;  if~$\alpha_{uv}, \ \alpha_{vu} <1$, there exists an attractive coexistence steady state.
The fact that, if two populations' ecological niches completely overlap, then one of the two species gets extinct, is exactly the statement of the Gause principle, a well-established law in ecology.  

The Lotka-Volterra models of ODEs have been extended in many ways and its applications range from  technology substitution to business competition. In the stochastic analysis community, system \eqref{model:lv} with additioned noise terms has been largely studied \cite{noise}. 
Another branch were these systems have been of huge influence is, of course, reaction-diffusion equations.
In the next paragraph we are spending some words on the results for diffusive Lotka-Volterra competitive systems.

\subsubsection{Competitive Lotka-Volterra model with diffusion}

In the interaction between different population, as already happens in the dynamic of a single species, spatial organisation plays an important role. 
A great literature has been devoted to the competitive Lotka-Volterra system with diffusion, that is, up to a rescaling,
\begin{equation}\label{model:diffusion}
\begin{cases}
\partial u - \Delta u=  u\left(1-u-\alpha_{uv} v  \right), & x\in\R^n, \  t>0,\\
\partial v - d\Delta v
=a v\left(1- v- \alpha_{vu} u \right), & x\in\R^n, \  t>0,
\end{cases}
\end{equation}
for some $d$, $a$, $\alpha_{uv}$, $\alpha_{vu}$  positive constants.
Richer dynamics and more questions naturally arise for system \eqref{model:diffusion} but, unsurprisingly, the study of these  involves many more difficulties. 
Just to give some flavour of those, we provide some details on the study of the speed of propagation of travelling waves connecting the steady states  $\left(1, 0\right)$ and $\left( 0, 1\right)$, to which many studies have been consecrated. From the works of Lewis, Li and Weinberger \cite{lewis2002spreading, li2005spreading} the minimal speed of propagation of a monotonic wave is called $c_{LLW}$. Even it the simplest case stated in \eqref{model:diffusion}, the exact value of $c_{LLW}$ is still not known. Calling $c_{KPP}$ the minimal speed of diffusion for the second equation under the assumption $u\equiv 0$, it holds that $c_{LLW}\geq c_{KPP}$, but the inequality may be strict depending on the parameters \cite{huang2010problem, huang2011non}.

Nevertheless, system \eqref{model:diffusion} is one of the simplest among many possibilities; in the literature one finds systems considering
nonlocal diffusion \cite{pan2012invasion},
free boundary \cite{du2017semi},
cross diffusion  \cite{lou2004limiting},
and many other variations.

\subsection{A model of Lotka-Volterra type for aggressive competition and analysis of the strategies}

Among the several models dealing with the
dynamics of biological systems, the case of populations in open hostility
seems to be rather unexplored.
Our model considers the case of
two populations competing for the same resource with one aggressive population
which attacks the other:
concretely, one may think of
a situation in which
two populations live together in the same territory and share the same
environmental resources,
till one population wants to prevail and try to overwhelm the other.
We consider this situation as a ``civil war'', since the two populations share
land and resources.

\subsubsection{The model}

We now describe in further detail our model of conflict between the
two populations and the attack strategies pursued by the aggressive population.
Given the lack of reliable data related to civil wars, 
the equations were derived by deduction
from the principles of population dynamics.
Our idea is to modify the Lotka-Volterra competitive system for two populations with
density~$u$ and~$v$,
adding to the usual competition for resources the fact that both populations suffer some losses as an outcome of the attacks. 
The key point in our analysis
is that the clashes do not depend on the chance of meeting of the two populations, given by the quantity~$uv$, as it happens in many other works in the literature, but they are sought by the first population and
depend only on the size~$u$ of the first population and on its level of aggressiveness~$a$.
The resulting model is
\begin{equation}\label{model}
\left\{
\begin{array}{llr}
\dot{u}&= u(1-u-v) - acu, & {\mbox{ for }}t>0,\\
\dot{v}&= \rho v(1-u-v) -au, & {\mbox{ for }}t>0,
\end{array}
\right.
\end{equation}
where~$a$,~$c$ and~$\rho$ are nonnegative real numbers. Here, the coefficient~$\rho$ models the fitness of the second population with respect to the first one. The parameter~$c$ here stands for the quotient of endured per inflicted  damages for the first population.

\subsubsection{Behaviour of solutions}

We denote by~$(u(t), v(t))$ a solution of~\eqref{model} starting from a point~$(u(0),v(0))\in [0,1] \times [0,1]$.
We will also refer to the \textit{orbit} of~$(u(0), v(0))$  as the collection of
points~$(u(t), v(t))$ for~$t\in \R$, thus both positive and negative times, while the \textit{trajectory} is the collection of points~$(u(t), v(t))$ for~$t\geq0$.

From the equations in~\eqref{model}, one promptly sees that $v=0$ is not an equilibrium, hence,
$v$ can reach the value~$0$ and even negative values in finite time. 
From a modelling point of view, negative values of $v$ are not acceptable, being it a population density. 
However, we will suppose that the dynamics stops when the value~$v=0$ is reached for the first time.
At this point, the conflict ends with the victory of the first population~$u$, that can continue its evolution with a classical Lotka-Volterra equation of the form
\begin{equation*}
\dot{u}= u (1- u)
\end{equation*}
and that would certainly fall into the attractive equilibrium~$u=1$.

In order to state our first result on the dynamics of the system~\eqref{model},
we first observe that, in a real-world situation, the value of~$a$ would probably be non-constant and discontinuous, so we allow this coefficient to take values in the class~$\mathcal{A}$
defined as follows:
\begin{equation}\begin{split}\label{DEFA}
\mathcal{A}&\; :=
\big\{a: [0, +\infty) \to [0, +\infty) {\mbox{ s.t.~$a$ is continuous}}\\
&\qquad \qquad {\mbox{except at most at a finite number of points}}\big\}.\end{split}\end{equation}
A \emph{solution related to a strategy~$a(t)\in \mathcal{A}$} is a pair~$(u(t), v(t)) \in C_0 (0,+\infty)\times C_0 (0,+\infty)$, which is~$C^1$ outside the discontinuity points of~$a(t)$ and
solves system~\eqref{model}.
Moreover, once the initial datum is imposed, the solution is assumed to be
continuous at~$t=0$. 

In this setting, we establish the existence of the solutions to problem~\eqref{model}
and we classify their behavior with respect to the possible exit from the domain~$[0,1]\times[0,1]$. Given~$(u(0), v(0))\in [0,1] \times [0,1]$ and~$a(t)\in\mathcal{A}$, two scenarios are possible for a solution~$(u(t),v(t))$ with~$a=a(t)$
of system~\eqref{model} starting at~$(u(0), v(0))$:

	\begin{enumerate}
		\item[(1)] The solution~$(u(t), v(t))$ issued from~$(u(0), v(0))$ belongs to~$ [0,1]\times (0,1]$ for all~$t\geq 0$.
		\item[(2)] There exists~$T\geq0$ such that the solution~$(u(t), v(t))$ issued from~$(u(0), v(0))$ exists unique for all~$t\leq T$, and~$v(T)=0$ and~$u(T)>0$.
	\end{enumerate} 

\noindent As a consequence, we can define the 
the \textit{stopping time} of the solution~$(u(t), v(t))$ as
\begin{equation}\label{def:T_s}
T_s (u(0), v(0)) = 
\left\{ 
\begin{array}{ll}
+\infty & \text{if situation (1) occurs}, \\
T & \text{if situation (2) occurs}.
\end{array}
\right.
\end{equation} 
{F}rom now on, we will implicitly  consider solutions~$(u(t),v(t))$ only for~$t\leq T_s(u(0), v(0))$.

We call \emph{victory of the first population} the scenario where~$T_s < +\infty$, that corresponds to the case where $v(T_s)=0$ and~$u(T_s)>0$. 
On the other hand,  we call \emph{victory of the second population} the scenario where~$(u,v)$ tends to $(0,1)$ as~$t$ tends to infinity.

Now we are going to analyze the dynamics of~\eqref{model} with a particular focus on possible strategies. To do this, we now define the \textit{basins of attraction}.
The first one is the basin of attraction of the point~$(0,1)$, that is
\begin{equation}\label{DEFB}\begin{split}
\mathcal{B}&\;:= \Big\{ (u(0),v(0))\in [0,1]\times[0,1] \;{\mbox{ s.t. }}\;\\
&\qquad\qquad T_s (u(0), v(0)) = +\infty, \ (u(t),v(t)) \overset{t\to\infty}{\longrightarrow} (0,1)  \Big\},
\end{split}\end{equation}
namely the set of the initial points for which the first population gets extinct (in infinite time) and the second one survives.
The other one is
\begin{equation}\label{DEFE}
\mathcal{E}:= \left\{ (u(0),v(0))\in ([0,1]\times[0,1])\setminus(0,0) \;{\mbox{ s.t. }}\;  T_s(u(0),v(0))< + \infty \right\},
\end{equation}
namely
the set of initial points for which we have the victory of the first population and the extinction of the second one. 


\subsubsection{A control problem}

In a rational, or at least well-organised population, one may expect that the parameter $a$, representing aggressiveness, is subject to control; we are suggesting that a population, by performing premeditated attacks, may control its strategy in the conflict and  would be able to choose the most appropriate one.

{F}rom now on, we may refer to the parameter~$a$ as the \textit{strategy}, that may also depend on time, and we will say that it is \textit{winning} if it leads to victory of the first population. 
We also notice that, with this choice, \eqref{model} is a \emph{control-affine} system.

The main problems that
we deal with are:
\begin{enumerate}
	\item The characterization of the {\em initial conditions for which there exists a winning strategy}.
	\item The {\em success of the constant strategies}, compared to all  possible strategies.
	\item The {\em construction of a winning strategy} for a given initial datum.
	\item The {\em existence of a single winning strategy independently of the initial datum}.
	\item The {\em existence of a winning strategy minimizing duration of the war}.
\end{enumerate}

The first question is a problem of \emph{target reachability} for a control-affine system.
The second point regards the choice of a suitable functional space where to choose the strategy. 
We also construct an actual winning strategy when victory is possible, answering the third and fourth question. 
The last question results to be an optimisation problem.

\subsection{Our results}

\subsubsection{Dynamics for a constant strategy}

The first step towards the understanding of the dynamics of the system
in~\eqref{model} is
is to analyze the behavior of the system for constant coefficients. 

To this end, we introduce some notation.
Following the terminology on pages~9-10 in~\cite{MR1056699},
we say that an equilibrium point (or fixed point) of the dynamics
is a (hyperbolic) {\em sink}
if all the eigenvalues of the linearized map have strictly
negative real parts, a (hyperbolic) {\em source}
if all the eigenvalues of the linearized map have strictly
positive real parts, and a (hyperbolic) {\em saddle}
if some of the eigenvalues of the linearized map have strictly
positive real parts
and some have negative real parts
(since in this problem we work in dimension~$2$,
saddles correspond to linearized maps with one
eigenvalue with
strictly positive real part
and one eigenvalue with
strictly negative real part).
We also recall that
sinks are asymptotically stable (and sources are
asymptotically stable for the reversed-time dynamics), see e.g. Theorem 1.1.1
in~\cite{MR1056699}.

With this terminology, we state the following theorem:

\begin{theorem}[Dynamics of system~\eqref{model}] \label{thm:dyn}
	For~$a > 0$, $c>0$ and~${\rho}> 0$ the system~\eqref{model} has the following features:
	\begin{itemize}
		\item[(i)] When~$0<ac<1$, the system has 3 equilibria:~$(0,0)$ is a source,~$(0,1)$ is 
		a sink, and 
		\begin{equation}\label{usvs}
		(u_s, v_s):= \left( \frac{1-ac}{1+{\rho}c} {\rho}c, \frac{1-ac}{1+{\rho}c} \right) \in (0,1)\times (0,1)
		\end{equation}
		is a saddle.
		
		\item[(ii)] When~$ac>1$, the system has 2 equilibria:~$(0,1)$ is a sink and~$(0,0)$ is a saddle.  
		\item[(iii)] When~$ac=1$, the system has 2 equilibria:~$(0,1)$ is a sink and~$(0,0)$
		corresponds to a strictly positive eigenvalue and a null one.
		\item[(iv)] We have 
		\begin{equation} \label{fml:division}
		[0,1]\times [0,1] = \mathcal{B} \cup \mathcal{E} \cup \mathcal{M}
		\end{equation}
		where~$\mathcal{B}~$ and~$\mathcal{E}$ are defined in~\eqref{DEFB}
		and~\eqref{DEFE}, respectively, and~$\mathcal{M}$ is a smooth curve.
		\item[(v)] The trajectories starting in~$\mathcal{M}$ tend to~$(u_s,v_s)$ if~$0<ac<1$,
		and to~$(0,0)$ if~$ac\ge1$ as~$t$ goes to~$+\infty$.
	\end{itemize}
\end{theorem}

More precisely, one can say that
the curve~$\mathcal{M}$ in Theorem~\ref{thm:dyn} is the stable manifold of the saddle
point~$(u_s,v_s)$ when~$0<ac<1$, and of the 
saddle point~$(0,0)$ when~$ac>1$. The case~$ac=1$ needs a special treatment,
due to the degeneracy of one eigenvalue, and in this case the curve~$\mathcal{M}$
corresponds to the center manifold of~$(0,0)$, and an ad-hoc argument will
be exploited
to show that also in this degenerate case orbits that start in~$\mathcal{M}$
are asymptotic in the future to~$(0,0)$.
As a matter of fact,~$\mathcal{M}$
acts as a dividing wall between the two basins of attraction $\mathcal{B}$ and $\mathcal{E}$, as described in~(iv)
of Theorem~\ref{thm:dyn}.

From a modelling point of view, Theorem~\ref{thm:dyn} shows that, also for our model, the Gause principle of exclusion is respected; that is, in general, two competing populations cannot coexist in the same territory, see e.g.~\cite{fath2018encyclopedia}. 

\begin{figure}[h] 
	\begin{subfigure}{.5\textwidth}
		\centering
		\includegraphics[width=.8\linewidth]{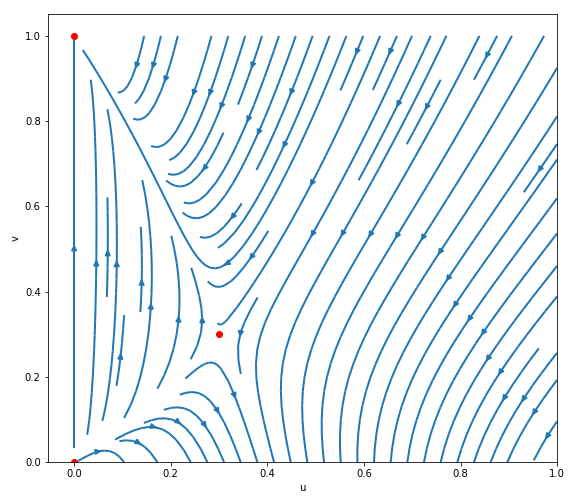}
		\caption{$a=0.8$,~$c=0.5$,~$\rho=2$}
	\end{subfigure}%
	\begin{subfigure}{.5\textwidth}
		\centering
		\includegraphics[width=.8\linewidth]{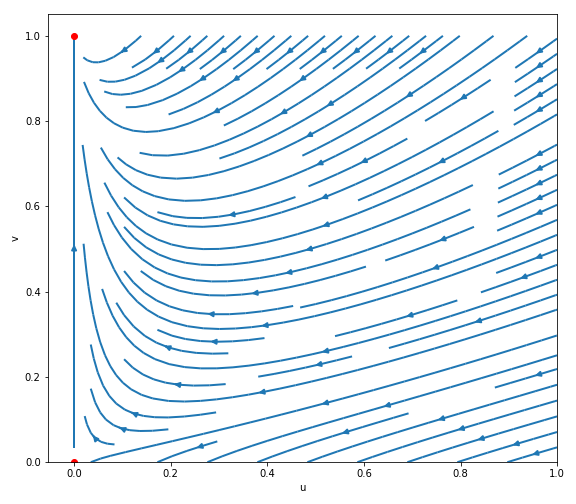}
		\caption{$a=0.8$,~$c=3$,~$\rho=2$}
	\end{subfigure}
	\caption{\it The figures show a phase portrait for the indicated values of the coefficients. In blue, the orbits of the points. The red dots represent the equilibria. The images are realised with Python.}
	\label{fig:dyn}
\end{figure}

One peculiar feature of our system is that, if the aggressiveness is too strong, the equilibrium~$(0,0)$ changes its ``stability'' properties, passing from a source (as in (i) of
Theorem~\ref{thm:dyn})
to a saddle point (as in (ii) of
Theorem~\ref{thm:dyn}). This shows that the war may have self-destructive outcomes, therefore it is important for the first population to analyze the situation in order to choose a proper level of aggressiveness. 
Figure~\ref{fig:dyn} shows one example of dynamics for each case.

The dedicated chapter contains further results on the dependence of $\mathcal{B}$ and $\mathcal{E}$ on the parameters $\rho$ and $c$. The parameter $a$, a part from having a more intricate influence on the system, may be interpreted not as a biological constant but rather as a choice of the first population.
Therefore, we perform a deeper analysis, whose result are presented in the next paragraph.

\subsubsection{Dynamics for variable strategies and optimisation results}

We now introduce some terminology.
Recalling~\eqref{DEFA},
for any~$\mathcal{T}\subseteq \mathcal{A}$, we set
\begin{equation}\label{DEFNU}
\mathcal{V}_{\mathcal{T}}:= \underset{a(\cdot)\in \mathcal{T}}{\bigcup} \mathcal{E}(a(\cdot)),
\end{equation}
where~$\mathcal{E}(a(\cdot))$ denotes the set of initial data~$(u_0,v_0)$
such that~$T_s(u_0,v_0)< +\infty$, when the coefficient~$a$ in~\eqref{model} is replaced by the function~$a(t)$.
Namely,~$\mathcal{V}_{\mathcal{T}}$ represents the set of initial conditions for which~$u$ is able to win by choosing a suitable strategy in~$\mathcal{T}$; we call~$\mathcal{V}_{\mathcal{T}}$ the \emph{victory set} with admissible strategies in~$\mathcal{T}$.
We also say that~$a(\cdot)$ is a \emph{winning strategy} for the point~$(u_0,v_0)$
if~$(u_0,v_0)\in \mathcal{E}(a(\cdot) )$.

Moreover, we will call
\begin{equation}\label{u0v0}
(u_s^0, v_s^0):= \left(\frac{\rho c}{1+\rho c}, \frac{1}{1+\rho c}\right).
\end{equation}
Notice that~$(u_s^0, v_s^0)$ is the limit point as~$a$ tends to~$0$ of the sequence of saddle points~$\{(u_s^a, v_s^a)\}_{a>0}$
defined in~\eqref{usvs}.
\medskip

With this notation,
the first question that we address is for which initial configurations it is possible for the population~$u$
to have a winning strategy, that is, to characterize the victory set. For this, we allow the strategy to take all the values in~$[0, +\infty)$.
In this setting, we have the following result:

\begin{theorem}\label{thm:Vbound}
	\begin{itemize}
		\item[(i)] For~$\rho=1$, we have that  
		\begin{equation}\label{Vbound1}\begin{split}
		\mathcal{V}_{\mathcal{A}} = \,&\Big\{ (u,v)\in[0,1] \times [0,1] \;
		{\mbox{ s.t. }}\;  v-\frac{u}{c}<0 \; {\mbox{ if }} u\in[0,c]\\
		&\qquad\qquad\qquad {\mbox{ and }}\; v\le1 \; {\mbox{ if }} u\in(c,1]\Big\},
		\end{split}\end{equation}
		with the convention that the last line in~\eqref{Vbound1} is not present if~$c\ge1$.
		\item[(ii)] 
		For~$\rho<1$, we have that
		\begin{equation}\label{bound:rho<1}
		\begin{split}
		\mathcal{V}_{\mathcal{A}} &\;= \Bigg\{ (u,v)\in[0,1] \times [0,1] \;{\mbox{ s.t. }}\;
		v< \gamma_0(u) \ \text{if} \ u\in [0, u_s^0], \\ 
		&\qquad\qquad\qquad\qquad
		v< \frac{u}{c} + \frac{1-\rho}{1+\rho c} \ \text{if} \ u\in \left[u_s^0, 
		\frac{\rho c(c+1)}{1+\rho c}\right]\\
		&\qquad\qquad\qquad\qquad
		{\mbox{and }}\; v\le1\ \text{if} \ u\in \left(
		\frac{\rho c(c+1)}{1+\rho c},1\right]
		\Bigg\},
		\end{split}
		\end{equation}
		where 
		\begin{equation*}
		\gamma_0(u):= \frac{u^{\rho}}{\rho c(u_s^0)^{\rho-1}},
		\end{equation*}
		and we use the convention that the last line
		in~\eqref{bound:rho<1} is not present if~$ \frac{\rho c(c+1)}{1+\rho c}\ge1$.
		\item[(iii)] For~$\rho>1$, we have that
		\begin{equation}\label{bound:rho>1}
		\begin{split}
		\mathcal{V}_{\mathcal{A}} &\;= \Bigg\{ (u,v)\in[0,1] \times [0,1]\;
		{\mbox{ s.t. }}\; v< \frac{u}{c} \ \text{if} \ u\in [0, u_{\infty}],\\&\qquad
		\qquad\qquad\qquad 
		v< \zeta(u)  \ \text{if} \ u\in\left(u_{\infty}, \frac{c}{(c+1)^{\frac{\rho-1}\rho}}\right] \\&\qquad
		\qquad\qquad\qquad 
		{\mbox{and }}\; v\le 1
		\ \text{if} \ u\in\left(\frac{c}{(c+1)^{\frac{\rho-1}\rho}},1\right] 
		\Bigg\},
		\end{split}
		\end{equation}
		where
		\begin{equation}\label{ZETADEF}
		u_{\infty}:= \frac{c}{c+1}
		\quad {\mbox{ and }}\quad \zeta (u):= \frac{u^{\rho}}{c \, u_{\infty}^{\rho-1}} .
		\end{equation}   
		and we use the convention that the last line
		in~\eqref{bound:rho>1} is not present if~$  \frac{c}{(c+1)^{\frac{\rho-1}\rho}}\ge1$.
	\end{itemize}
\end{theorem}

Theorem \ref{thm:Vbound} implies that the problem is not \emph{controllable}, that is, for some initial conditions the first population is not able to reach its target.

\medskip

In practice, 
constant strategies could be certainly easier to implement and
it is therefore natural to investigate whether or not
it suffices to restrict the control to constant strategies
without altering the possibility of victory.
The next result addresses this problem by showing that when~$\rho=1$
constant strategies are as good as all strategies,
but instead when $\rho\ne 1$ victory cannot be achieved by only
exploiting constant strategies:

\begin{theorem}\label{thm:W}
	Let $\mathcal{K}\subset \mathcal{A}$ be the set of constant functions. Then the following holds:
	\begin{itemize}
		\item[(i)] For~$\rho= 1$, we have that~$ \mathcal{V}_{\mathcal{A}}=\mathcal{V}_{\mathcal{K}}=\mathcal{E}(a)$ for all $a>0$;
		\item[(ii)] For~$\rho\neq 1$, we have that~$\mathcal{V}_{\mathcal{K}} \subsetneq \mathcal{V}_{\mathcal{A}}$.
	\end{itemize}	
\end{theorem} 

The result of Theorem~\ref{thm:W}, part~(i),
reveals a special rigidity of the case~$\rho=1$
in which the victory depends only on the initial conditions, but it is independent of the strategy~$a(t)$.
Instead, as stated in
Theorem~\ref{thm:W}, part~(ii),
for~$\rho \neq 1$ the choice of~$a(t)$ plays a crucial role in determining which population is going to win and constant strategies do not exhaust all the
possible winning scenarios.
We stress that~$\rho=1$ plays also
a special role in the biological interpretation of the model, since in this case the two
populations have the same fitness to the environmental resource, and hence, in a sense,
they are indistinguishable, up to the possible aggressive behavior of the first population.

Next, we show that for all points in the set~$\mathcal{V}_{\mathcal{A}}$ we can choose an appropriate piecewise constant strategy with at most one discontinuity; functions with these properties are called Heaviside functions. 

\begin{theorem}\label{thm:H}
	There holds that~$\mathcal{V}_{\mathcal{A}} = \mathcal{V}_{\mathcal{H}}$, where~$\mathcal{H}$ is the set of Heaviside functions.
\end{theorem}

The proof of Theorem~\ref{thm:H} solves also the third question mentioned in the introduction. As a matter of fact, it proves that for each point we either have a constant winning strategy or
a winning strategy of type 
\begin{equation*}
a(t) = \left\{
\begin{array}{lr}
a_1  &{\mbox{ if }} t<T ,\\
a_2  &{\mbox{ if }} t\geq T,
\end{array}
\right.
\end{equation*}
for some~$T\in(0,T_s)$, and
for suitable values~$a_1$,~$a_2 \in (0,+\infty)$ such that one is very small and the other one very large, the order depending on~$\rho$. 
The construction that
we give also puts in light the fact that the choice of the strategy depends on the initial datum, answering also our fourth question. 

It is interesting to observe that the winning strategy that switches abruptly from a small to a large value
could be considered, in the optimal control terminology, as a ``bang-bang'' strategy.
Even in a target reachability problem, the structure predicted by Pontryagin's Maximum Principle is brought in light: the bounds of the set~$\mathcal{V}_{\mathcal{A}}$, as
given in Theorem~\ref{thm:Vbound}, depend on the bounds that
we impose on the strategy, that are,~$a \in[0,+\infty)$.

It is natural to consider also the case
in which the level of aggressiveness 
is constrained between a minimal and maximal threshold,
which corresponds to the setting~$a\in[m,M]$ for suitable~$M\geq m\geq 0$, with the hypothesis that $M>0$.
In this setting, we denote by~$\mathcal{A}_{m,M}$ the class of piecewise continuous strategies~$a(\cdot)$
in~${\mathcal{A}}$ such that~$
m\leq a(t)\leq M$ for all~$t>0$ and we let
\begin{equation}\label{SPE}
\mathcal{V}_{m,M}:=\mathcal{V}_{\mathcal{A}_{m,M}}=\underset{{a(\cdot)\in \mathcal{A}}\atop{m\leq a(t)\leq M}
}{\bigcup} \mathcal{E}(a(\cdot))=
\underset{{a(\cdot)\in \mathcal{A}}_{m,M}
}{\bigcup} \mathcal{E}(a(\cdot)).\end{equation}
Then we have the following:

\begin{theorem}\label{thm:limit}
	Let $M$ and $m$ be two real numbers such that $M\geq m\geq 0$ and $M>0$. Then, for $\rho\neq 1$ we have the strict inclusion  $\mathcal{V}_{{m,M}}\subsetneq \mathcal{V}_{\mathcal{A}}$.
\end{theorem}

Notice that for $\rho=1$, Theorem \ref{thm:W} gives instead that $\mathcal{V}_{{m,M}}= \mathcal{V}_{\mathcal{A}}$
and we think that this is a nice feature, outlining a special role played by the parameter~$\rho$
(roughly speaking, when~$\rho=1$ constant strategies suffice
to detect all possible winning configurations, thanks to 
Theorem \ref{thm:W}, while when~$\rho\ne1$ non-constant strategies are necessary to detect
all winning configurations).

\paragraph{Time minimizing strategy.}
Once established that it is possible to win starting in a certain initial condition, we are interested in knowing which of the possible strategies is best to choose. One condition that may be taken into account is the duration of the war. Now, this question can be written as a minimization problem with a proper functional to minimize and therefore the classical Pontryagin theory applies. 

To state our next result, 
we recall the setting in~\eqref{SPE} and define
\begin{equation*}
\mathcal{S}(u_0, v_0) := \Big\{ a(\cdot)\in \mathcal{A}_{m,M}
\;\mbox{ s.t. }\; (u_0, v_0) \in \mathcal{E}(a(\cdot))  \Big\},
\end{equation*}
that is the set of all bounded strategies for which the trajectory starting at~$(u_0, v_0)$ leads to the victory of the first population.
To each~$a(\cdot)\in\mathcal{S}(u_0, v_0)$ we associate the stopping time defined in~\eqref{def:T_s}, and we express its dependence on~$a(\cdot)$ by writing~$T_s(a(\cdot))$.
In this setting, we provide the following statement concerning the strategy leading
to the quickest possible victory for the first population:

\begin{theorem}\label{thm:min}
	Given a point~$(u_0, v_0)\in \mathcal{V}_{m,M}$, there exists a winning strategy~$\tilde{a}(t)\in
	\mathcal{S}(u_0, v_0)$, and a trajectory~$(\tilde{u}(t), \tilde{v}(t) )$ associated with~$\tilde{a}(t)$,
	for~$t\in[0,T]$,
	with~$(\tilde{u}(0), \tilde{v}(0) )=(u_0,v_0)$, where~$T$ is given by
	\begin{equation*}
	T = \underset{a(\cdot)\in\mathcal{S}}{\min} T_s(a(\cdot)).
	\end{equation*}
	Moreover,
	\begin{equation*}
	\tilde{a}(t)\in \left\{m, \ M, \    a_s(t) \right\},
	\end{equation*}	  
	where
	\begin{equation}\label{KSM94rt3rjjjdfe}
	{a}_s(t) := \dfrac{(1-\tilde{u}(t)-\tilde{v}(t))[\tilde{u}(t) \, (2c+1-\rho c)+\rho c]}{\tilde{u}(t) \, 2c(c+1)}.
	\end{equation}	
\end{theorem}

\begin{figure}[h] 
	\begin{subfigure}{.5\textwidth}
		\centering
		\includegraphics[width=.9\linewidth]{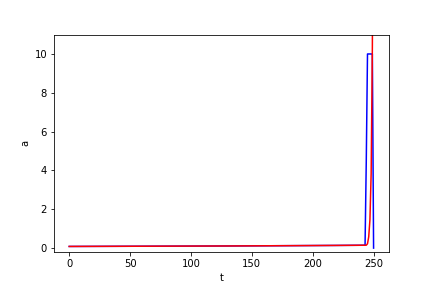}
	\end{subfigure}%
	\begin{subfigure}{.5\textwidth}
		\centering
		\includegraphics[width=.9\linewidth]{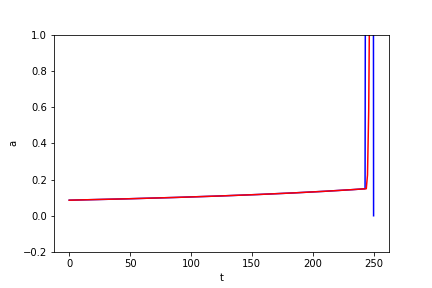}
	\end{subfigure}
	\caption{\it The figure shows the result of a numerical simulation searching a minimizing time strategy~$\tilde{a}(t)$ for the problem starting in~$(0.5, 0.1875)$ for the
		parameters~$\rho=0.5$,~$c=4.0$,~$m=0$ and~$M=10$. In blue, the value
		found for~$\tilde{a}(t)$; in red, the value of~$a_s(t)$ for the corresponding trajectory~$(u(t), v(t))$. As one can observe,~$\tilde{a}(t)\equiv a_s(t)$ in a long trait.
		The simulation was done using AMPL-Ipopt on the server NEOS and pictures have been made with Python. 
	}
	\label{fig:min}
\end{figure}

The surprising fact given by Theorem~\ref{thm:min}
is that the
minimizing strategy is not only of bang-bang type, but it may assume some values along a \emph{singular arc}, given by~$a_s(t)$.
This possibility is realized in some concrete cases, as we verified by running some numerical simulations, whose results can be visualized in Figure~\ref{fig:min}. 

\subsubsection{Perspectives}

The system of ODEs is the cornerstone for the study of the reaction-diffusion system
\begin{equation}\label{model2}
\left\{
\begin{array}{lr}
\partial_t u- \partial_{xx}^2 u = u(1-u-v) - acu, & {\mbox{ for }}  x\in\R, \ t>0,\\
\partial_t v- d\partial_{xx}^2 v = \rho v(1-u-v) -a \int_{\R} u, & {\mbox{ for }} x\in\R, \ t>0,
\end{array}
\right.
\end{equation}
for some $d>0$. 
We expect solutions to this system to have very interesting behaviours. It is possible that the second population reaches the value $0$ in only some points of the domain, giving an example of the interesting phenomenon known as dead-core, see e.g. \cite{souplet}.

\section{Evolution equations with classical and fractional derivatives}

\subsection{Fractional derivatives in evolution equations}

The idea of fractional calculus first appears in the discussions between Leibniz and  De l'Hospital (see \cite{Ross77}); namely, given the classical derivative $\frac{d^n f(x)}{dx^n}$ for $n\in\N$, it is quite natural to ask if it is possible to define a generalisation of this operator but with non entire order, thus $\frac{d^{\alpha} f(x)}{dx^{\alpha}}$ with $\alpha\in\R$. 

However, it is only in the last few decades that a good number of mathematicians started to work on fractional calculus. One of the reasons of this interest is the fact that fractional derivative can help in the modelization of processes with memory or of diffusion phenomena with spreading behaviour different from the one prescribed by Brownian motion. This has applications in charge carrier transport in semicondutors, nuclear magnetic resonance diffusometry, porous systems, dynamics in polymeric systems (see \cite{metzler2000random} and the references therein).

Here, we introduce some fractional derivatives and operators and justify their meaning and relations with the previous material.
Providing an overview on the state of art over existence, regularity and behaviour of evolution equations dealing with fractional operators is far from our purposes, due to the complexity of the topic. For that, we refer to \cite{arendt1992vector, ciraolo2018rigidity, clement1979abstract, pruss2013evolutionary, sanchez2013long, zacher2009weak} and the reference therein.

\subsubsection{The Caputo derivative}

There are several way of defining a fractional derivative. We choose to work with \emph{Caputo derivative}, 
that was first proposed in a geology model by Caputo in \cite{caputo1967linear}.
The Caputo derivative of order $\alpha\in(0,1)$ is defined by
\begin{equation*}
D_t^{\alpha} f(t):= \frac{1}{\Gamma(1-\alpha)} \int_{0}^{t} \frac{\dot{f}(\tau)}{(t-\tau)^{\alpha}} d\tau
\end{equation*}
where $\Gamma$ is Euler's Gamma-function and $\dot{f}$ is the classical derivative of $f$. 
For simplicity, we will omit the constant and work with 
\begin{equation*} \label{def:caputo}
\partial_t^{\alpha} f(t):= \Gamma(1-\alpha)	D_t^{\alpha} f(t).
\end{equation*}
Notice that $\partial_t^{\alpha} f(t)$ is defined for all $f\in\mathcal{C}^1([0, t])$ such that $\dot{f}\in L^1(0,t)$. It is also possible to define Caputo derivatives of higher order, thus with $m-1<\alpha<m$ for $m\in\N$, by the following:
\begin{equation*}
D_t^{\alpha} f(t):= \frac{1}{\Gamma(m -\alpha)} \int_{0}^{t} \frac{{f}^{(m)}(\tau)}{(t-\tau)^{1+\alpha-m}} d\tau.
\end{equation*}

The Caputo derivative describes a process ``with memory'', in the sense that the history of the process is encoded in the derivative, though old events ``count less'' than recent ones, since their value is counted with a smaller weight. Due to this memory property, this operator has
found several applications in models of hydrogeology, heat transmission, percolation.

The Caputo derivative is considered to be an operator of ``fractional'' order, as opposed to the entire order, that is proper to the classical derivatives. 
The value $\alpha$ corresponds to the order of the derivative:
indeed, for a power of order $r>0$, it holds $\partial_t^{\alpha} t^{r}=C t^{r-\alpha}$ for some constant $C>0$. 

\medskip

Among the other types of fractional derivatives, it is worth mentioning the Riemann-Liouville derivative because of its large diffusion. The Riemann-Liouville derivative is defined by
\begin{equation}
\mathcal{D}_t^{\alpha} f(t):= \frac{1}{\Gamma(1-\alpha)} \frac{d}{dt} \int_{0}^{t} \frac{f(\tau)}{(t-\tau)^{\alpha}} d\tau,
\end{equation}
and making the calculations one can show that it differs from the Caputo derivative by the term $\frac{f(0)}{t^{\alpha}}$. One of the reasons of the popularity of the Riemann-Liouville derivative is that its limit for $\alpha\to1$ coincides with the classical derivative. 

\paragraph{Evolution equations with Caputo time derivative.} 

Classical partial differential equations are often divided in three groups depending on the order of the time derivative: elliptic, parabolic, hyperbolic. Nevertheless, even in the Preface of \emph{Partial Differential Equations} by Evans \cite{evans}, the author states that this subdivision is fictive and ``creates the false impression that there is some kind of general and useful classification scheme available''.
This subdivision is supposed to put together object with similar behaviours and classic theory results are usually meant for one of these clusters. 

An evolution equation with the Caputo time derivative, for example
\begin{equation}\label{cap1}
	\partial_t^{\alpha} u - \Delta u =0, 
\end{equation}
is not part of any of this groups. 
Thus, many results that one may want to use are not available and must be recovered.

However, relations with the classical objects are present.
Because of some behaviour similarities, evolution equations with Caputo time derivative have often been compared to parabolic ones \cite{metzler2000random}.
Recently, in \cite{dipierro2018simple}, Dipierro and Valdinoci inspected a model of transmission in neural structures and derive equation \eqref{cap1} from basic principles. 
Doing so, they realized that it can be seen as a superposition of several hyperbolic equations acting with delay. 
Despite this, the behaviour of solutions of \eqref{cap1} is not similar to the one of wave equations: in fact, in opposition to hyperbolic equations, \eqref{cap1} has a regularising effect on initial data.

\subsubsection{The fractional Laplacian}

An operator that has been very popular in recent years is the fractional Laplacian, which is considered in some sense the fractional counterpart of the classic homonym operator.  
Given $s\in(0,1)$, we define the fractional Laplacian as
\begin{equation}\label{0302}
	- \left( \Delta\right)^s u(x) := \text{P.V.} \int_{\R^n} \frac{u(x)-u(y)}{|x-y|^{n+2s}} dy,
\end{equation}
where ``P.V.'' stands for Principal Value.
Curiously, many equivalent definitions of the fractional Laplacian are possible: a part the one in \eqref{0302}, in \cite{abatangelo2019getting, kwasnicki2017ten} one can find a very exhaustive list together with the most important properties of the operator known so far.

One of the reason for the popularity of the fractional Laplacian is its connection with L\'{e}vy processes. We now give the idea behind the derivation of an evolution equation containing the fractional Laplacian from a discrete L\'{e}vy process, which is similar to the derivation of the heat equation from a Brownian movement. Consider  an infinite grid $h \Z^n$, for some step $h>0$, and a discrete evolution of time step $\tau=h^{2s}$. Imagine to put a particle at the origin. At each time step $t\in \tau \N$, the particle can jump to any vertex of the grill different from its actual position, with a probability that depends on the length of the jump;  namely, if the particle is at the position $hk$, with $k\in\Z^n$, the probability to jump into the position $hj$, if $j\neq k$, is
\begin{equation*}
	P_h(k,j)=\frac{C}{|k-j|^{n+2s}},
\end{equation*} 
with $C$ a normalisation constant. Then, we call $u(t,x)$ the probability of finding the particle in $x\in h\Z^n$ at the time $t\in\tau\N$. The function $u(t,x)$ evolves according to the probability of the jumps; for example, the probability of finding the particle in the origin at some time $t+\tau$ is
\begin{equation*}
	u(t+\tau, 0)= \sum_{{j\in\Z^n\setminus 0}} P_h(0,j) u(t, j)= \sum_{{j\in\Z^n\setminus 0}} \frac{C }{|j|^{n+2s}}u(t,j)
\end{equation*}
By taking the limit as $h$ tends to $0$, 
and by performing suitable manipulations, 
the last equality becomes the evolution equation
\begin{equation*}
	\partial_t u - (\Delta)^s u=0.
\end{equation*} 
For all the details of the proof, we refer to \cite{jacob2005pseudo}.

Satellite-based measures of animal movement performed in the last years have shown that L\'{e}vy process are a better approximation of animal movement than Brownian motion. Some examples are provided by honey bees displacements and by movement of marine predators when prey is scarce \cite{bartumeus2009behavioral, humphries2010environmental, reynolds2009levy}. In general, it appears that L\'{e}vy-flights are more fitting hunt strategies than Brownian walks \cite{lomholt2008levy}.
For this reason, the fractional Laplacian has been introduced in population dynamics model, see \cite{bates2007existence, coville2007non,  genieys2006pattern} and the reference therein. However, the technical difficulties of dealing with such delicate operators have not been totally overcome.

\subsection{Decay estimates for evolution equations with classical and fractional derivatives}

Among the many open questions for fractional operators, we choose to study decay estimates of a class of evolution equations with possibly nonlocal or nonlinear diffusion operators.
In particular, we are going to study the decay in time of the Lebesgue norm of solutions to a Cauchy problem in a bounded domain. 
We present some general results that apply to a wide class of evolution equations, namely all the ones involving a diffusion operator that is satisfying a certain ellipticity property, involving an ``energy functional'' that suits for both local and non local operators, possibly complex. The time derivative may be of two types: purely classical or a linear combination of a classical derivative and a Caputo derivative.

\subsubsection{The problem}

We now set the problem. 
Let $\lambda_1, \lambda_2 \geq 0$ be fixed positive numbers. 
We suppose, for concreteness,
that $$\lambda_1 + \lambda_2=1,$$
but up to a rescaling of the operator we can take $\lambda_1, \lambda_2$
any nonnegative numbers with positive sum. 
Let $\Omega \subset \R^n$ be a
bounded open set and let $u_0\in L^{\infty}(\R^n)$ such that $\text{supp} \,u_0 \subset \Omega$. 
Consider the Cauchy problem
\begin{equation} \label{sys:generalform}
\left\{ \begin{array}{lr}
(\lambda_1 \partial_t^{\alpha} + \lambda_2 \partial_t) [u] + \mathcal{N}[u]=0, & {\mbox{for all }}x\in \Omega, \ t>0, \\
u(x,t)=0, & {\mbox{for all }}x\in \R^n \setminus \Omega , \ t>0, \\
u(x,0)=u_0(x), & {\mbox{for all }}x\in \R^n ,
\end{array} \right.
\end{equation}
where $\mathcal{N}$ is an operator, possibly involving fractional derivatives. 

We underline that we consider \emph{smooth} ($C^1$ often, $C^2$ if also the second derivative appears) and \emph{bounded} solutions of the problem \eqref{sys:generalform}. In fact, we want to avoid convergence problems with the integrals that appear in the statements and in the proofs. However, for certain operators, weaker hypothesis may be taken. 

Let us recall that for a complex valued function $v:\Omega\to\C$ the Lebesgue norm is
\begin{equation*}
\Vert v \Vert_{L^s(\Omega)} = \left( \int_{\Omega} |v(x)|^s \; dx \right)^{\frac{1}{s}}
\end{equation*}
for any $s\in[1, +\infty)$. Also, we call $\Re \{ z\}$ the real part of $z\in\C$.
The main assumption we take is the following: there exist $\gamma \in (0,+\infty) $ and $C\in (0,+\infty)$ such that
\begin{equation} \label{cond:complexstr}
\Vert u(\cdot,t) \Vert_{L^{s}(\Omega) }^{s-1+\gamma} \leq C \int_{\Omega} |u(x,t)|^{s-2} \Re \{ \bar{u}(x,t)\mathcal{N} [u](x,t)\} \; dx.
\end{equation}
The constants $\gamma$ and $C$ and their dependence from the parameters of the problem may vary from case to case. The righthandside of the equation may be seen as an energy functional linked to the diffusion operator.
This inequality implies, essentially, that the operator $\mathcal{N}$ is not too degenerate and the energy of the solution should control a certain power of the solution itself; here $\gamma$ plays the role of the degree of ellipticity. 
The inequality \eqref{cond:complexstr} strongly depends  on the validity of a Sobolev inequality for the solutions of the evolution equation.
To get an intuition of the roles of the factors, take the case of the Laplacian with $s=2$; integrating by parts on the righthandside one obtains $\Vert \nabla u (\cdot, t) \Vert_{L^{2}(\Omega) }^2$, thus the energy, which controls the $L^2$ norm of the solution by the Gagliardo-Nirenberg-Sobolev inequality.
In our setting, the structural inequality in~\eqref{cond:complexstr}
will be the cornerstone to obtain general energy estimates,
which, combined with appropriate barriers, in turn
produce time-decay estimates.

\subsection{Our Results}

Extending the method of \cite{SD.EV.VV}, we obtain a power-law decay in time
of the $L^s$ norm with $s\geq 1$. 
Also, for the case of classical time-derivatives,
we obtain exponential decays in time. The difference between
polynomial and exponential decays in time is thus related to
the possible presence of a fractional derivative in the operator involving the time variable.

\subsubsection{Decay estimate theorems}

First, we present this result for the more general setting, hence for a linear combination of classical and Caputo time derivative. We have the following:

\begin{theorem} \label{thm:complex}
	Let $u$ be a solution of the Cauchy problem \eqref{sys:generalform}, with $\mathcal{N}$ possibly complex
	valued. Suppose that there exist $s\in[1, +\infty)$, $\gamma\in(0,+\infty)$ and $C\in(0,+\infty)$ such that $u$ satisfies \eqref{cond:complexstr}.
	Then 
	\begin{equation} \label{claim1gen}
	(\lambda_1\partial_t^{\alpha} + \lambda_2\partial_t) \Vert u(\cdot,t) \Vert_{L^{s}(\Omega) } \leq -\dfrac{\Vert u(\cdot,t) \Vert_{L^{s}(\Omega) }^{\gamma}}{C},
	\qquad{\mbox{ for all }}t>0,\end{equation}
	where $C$ and $\gamma$ are the constants appearing in~\eqref{cond:complexstr}. 
	Furthermore,
	\begin{equation} \label{claim2gen}
	\Vert u(\cdot,t) \Vert_{L^{s}(\Omega) } \leq
	\dfrac{C_*}{1+t^{\frac{\alpha}{\gamma}}},\qquad{\mbox{ for all }}t>0,	
	\end{equation}
	for some~$C_*>0$, depending only on $C$, $\gamma$, $\alpha$
	and $\Vert u_0(\cdot) \Vert_{L^{s}(\R^n)}$.
\end{theorem}

A polynomial decay is a nice piece of information on the solution and we can expect this to be the best decay we can get for some fractional evolution equations \cite{gobbino, biler3al}. However, there is also evidence that for classical settings better decays can be achieved. In fact, the following theorem holds:

\begin{theorem}\label{thm:classic}
	Let $u$ be a solution of the Cauchy problem \eqref{sys:generalform} with only classical derivative (that is, $\lambda_1=0$) and $\mathcal{N}$ possibly complex
	valued. Suppose that there exist $s\in[1, +\infty)$, $\gamma\in(0,+\infty)$ and $C\in(0,+\infty)$ such that $u$ satisfies \eqref{cond:complexstr}.
	Then, for some~$C_*>0$, depending only on the constants~$C$ and~$\gamma$
	in~\eqref{cond:complexstr}, 
	and on~$\Vert u_0(\cdot) \Vert_{L^{s}(\R^n)}$, we have that:
	\begin{itemize}
		\item[1.]	if $0<\gamma \leq 1$ the solution $u$ satisfies
		\begin{equation} \label{claim3}
		\Vert u(\cdot,t) \Vert_{L^{s}(\Omega) } \leq
		C_* \, e^{-\frac{t}{C}},\qquad{\mbox{for all }}t>0;	
		\end{equation}
		\item[2.] if $ \gamma>1$, the solution $u$ satisfies
		\begin{equation} \label{claim4}
		\Vert u(\cdot,t) \Vert_{L^{s}(\Omega) } \leq
		\dfrac{C_*}{1+t^{\frac{1}{\gamma-1}}},\qquad{\mbox{for all }}t>0.	
		\end{equation}
	\end{itemize} 
\end{theorem}

As we will see in the proofs
of these two theorems, the idea is to find a supersolution of \eqref{claim1gen} and use a comparison principle in order to estimate the decay of the solution $u$. For the case of mixed derivatives, Vergara and Zacher \cite{VZ15} find both a supersolution and a subsolution decaying as $t^{-\frac{\alpha}{\gamma}}$. But, when $\alpha \to 1$, the subsolution tends to 0. 
On the other hand, the classical equation $\partial_t e =- e^{\gamma}$ has some exponential supersolutions.
This allows possibly better decays, which are in fact proven. 

We point out that the case of an evolution equation with only Caputo time derivative, i.e. of $\lambda_2=0$, was treated in \cite{SD.EV.VV}. The authors find in this case that the supersolution is still asymptotic to $t^{-\frac{\alpha}{\gamma}}$ and the decay is of polynomial type. 

It is interesting to notice the presence of a ``decoupling'' effect: for evolution equations with classical time derivative and fractional space derivative (take for example the fractional Laplacian, $-(\Delta)^{\sigma}u$, $\sigma\in(0,1)$, see \cite{SD.EV.VV}), the space derivative does not asymptotically interfere with the time derivative; 
thus the polynomial decay, typical of fractional derivatives, does not appear, leaving place for the exponential decay given by the classical time derivative. 
An example of this behaviour is found in \cite{MR3703556}, where a model inspired to atoms dislocation was studied.

\subsubsection{Applications} \label{ch1:app}

What makes Theorems \ref{thm:complex} and \ref{thm:classic} interesting is the fact that they may be applied to a wide range of equations. Indeed, the only hypothesis required in order to apply the theorems is the validity of the inequality \eqref{cond:complexstr} for suitable parameters $C$ and $\gamma$. In \cite{SD.EV.VV} and in our work, \eqref{cond:complexstr} was verified for many operators, that we are listing here together with some references on the origins of these operators:
\begin{itemize}[label=$\bullet$]
	\item the classic and fractional Laplacian, \cite{nonlocal},
	\item the classic and fractional $p$-Laplacian, \cite{CHEN20183065}, 
	\item the doubly nonlinear equation, \cite{RAVIART1970299}
	\item the classic and fractional porous medium equations, \cite{vazquez2007porous, porous} and \cite{MR2737788},
	\item the classic and fractional mean curvature equation, \cite{caffarelli2010nonlocal}
	\item the classic and fractional Kirchhoff equation, \cite{MR0002699} and \cite{kirchhoff},
	\item the classic and fractional magnetic operator, \cite{MR0142894} and \cite{groundstates}.
\end{itemize}
The list is not supposed to be exhaustive; in fact, the aim is only to provide some example of operators satisfying \eqref{cond:complexstr} and to encourage other mathematicians looking for some decay estimates to attempt with operators they are struggling with.






\chapter{\label{ch1} A Fisher-KPP model with a fast diffusion line in periodic media} 

In this chapter, we treat a model of population dynamics in a periodic environment presenting a fast diffusion line. The
``road-field'' model, introduced in \cite{brr}, is a system of coupled reaction-diffusion equations set in domains of different dimensions. Here, we consider for the first time the case of a reaction term depending on a spatial variable in a periodic fashion, which is of great interest for both its mathematical difficult and for its applications. We derive necessary and sufficient conditions for the survival of the species in terms of the sign of a suitable generalised principal eigenvalue, defined recently in \cite{romain}. Moreover, we compare the long time behaviour of a population in the same environment without the fast diffusion line, finding  that this element has no impact on the survival chances. This chapter corresponds to the paper \cite{periodic}.

\section{Setting and main results}

This chapter investigates some effects of a fast diffusion line in an ecological dynamics problem.
Various examples in the literature showed that, in the presence of roads or trade lines, some species or infections spread faster along these lines, and then diffuse in the surroundings. This was observed in the case of the Processionary caterpillar, whose spreading in France and Europe has been accelerated by accidental human transport \cite{robinet2012human}. Another striking proof was given in  \cite{gatto2020spread}, where the authors point out that the COVID-19 epidemics in Northern Italy at the beginning of 2020 diffused faster along the highways.  

A model for biological diffusion in a homogeneous medium presenting a fast diffusion line was proposed by Berestycki, Roquejoffre and Rossi in \cite{brr}, and since then is called the \emph{road-field model}. The authors proved an acceleration effect due to the road on the spreading speed of an invading species. 
Since then, a growing number of articles treated variations of the same system, investigating in particular the effect of different type of diffusion or different geometries \cite{berestycki2014speed, berestycki2015effect, rossi2017effect}.

However, natural environments are usually far from being homogeneous and, more often than not, territories are a composition of different habitats. Living conditions and heterogeneity play a strong impact on the survival chances of a species and on the equilibria at which the population can settle. 

Road-field models on heterogeneous environments have been  little studied so far, being more complex to treat. One of the few example is the paper \cite{giletti2015kpp} for periodic exchange terms between the population on the road and the one in the field. Recently, Berestycki, Ducasse and Rossi introduced a notion of generalised principal eigenvalue for the road-field system in \cite{romain} and, thanks to it, they were able to treat the case of an ecological niche facing climate change in \cite{econiches}. 

Here, we propose an analysis of the asymptotic behaviour of an invasive population under the assumption of spatial periodicity of the reaction term. Of course, under this hypothesis we can investigate deeper the dependence of the population on a natural-like environment and the effects of the road in this balance.
Under which conditions does the population survive in a periodic medium? And does the road play some role on the survival chances of a species, perturbing the environment and scattering the individuals, or rather permitting them to reach advantageous zones more easily?
These are the questions we are going to tackle.

\subsection{The model}

In this chapter, we study the reaction-diffusion model regulating the dynamics of a population living in a periodic environment with a fast diffusion channel. The equivalent of this model for homogeneous media was first introduced by Berestycki, Roquejoffre and Rossi in \cite{brr}. Consider the half plane $\Omega:=\R\times \R^+$, where we mean $\R^+=(0, +\infty)$.
The proposed model imposes the diffusion of a species in $\Omega$ and prescribes that on $\partial \Omega=\R\times \{ y=0\}$ the population diffuses at a different speed. 
We call $v(x,t)$ the density of population for $(x,y)\in\Omega$, hence on the ``field'', and $u(x)$ the density of population for $x\in\R$, i.e. on the ``road''; moreover, we take $D$, $d$, $\nu$, $\mu$ positive constants and $c\geq 0$. Then, the system we analyse reads
\begin{equation}\label{ch1sys:fieldroad}
\left\{
\begin{array}{lr}
\partial_t u-D  u '' -c u' - \nu  v|_{y=0} + \mu u= 0,   & x\in \R,  \\
\partial_t v -d \Delta v-c\partial_x v  =f(x,v),  & (x, y)\in \Omega, \\
-d  \partial_y{v}|_{y=0} + \nu v|_{y=0} -\mu u=0, & x\in\R.
\end{array}
\right.
\end{equation}
In $\Omega$, the population evolves with a net birth-death rate represented by $f$, that depends on the variable $x$. 
This embodies the heterogeneity of the media: in fact, environments are typically not uniform and some zone are more favourable than others. 
There is no dependence in the variable $y$, since the presence of the road itself creates enough heterogeneity in that direction.
The function  $f:\R\times \R_{\geq 0}\to \R $
is always supposed to be $\mathcal{C}^{1}$ in $x$, locally in $v$, and Lipschitz in $v$, uniformly in $x$; moreover, we suppose that the value $v=0$ is an equilibrium, that is
\begin{equation}\label{ch1hyp:0}
f(x,0)=0, \quad \text{for all} \ x\in \R,
\end{equation}
and that 		
\begin{equation}\label{ch1hyp:M}
\exists M>0 \ \text{such that} \ f(x, v)<0 \quad \text{for all} \ v>M \ \text{and all} \ x\in \R. 
\end{equation}
We will derive some inequalities on the generalised principal eigenvalue of \eqref{ch1sys:fieldroad} for the general case of $f$ respecting these hypothesis and 
$c$ possibly nonzero.

The characterisation of extinction or persistence of the species is performed for the case of $c=0$ and $f$ a periodic function, reflecting the periodicity of the environment in which the population diffuses. 
We will analyse the  case of a KPP nonlinearity, that is, we require that 
\begin{equation}\label{ch1hyp:KPP}
\frac{f(x,s_2)}{s_2}< \frac{f(x,s_1)}{s_1} \quad  \text{for all} \ s_2>s_1>0 \  \text{and all} \ x\in\R.
\end{equation}
Then, we suppose that there exists $\ell> 0$ such that
\begin{equation}\label{ch1hyp:per}
f(x+\ell, s)=f(x,s) \quad \text{for all} \ s >0 \  \text{and all} \ x\in\R.
\end{equation}


To study the effect of the line of fast diffusion, we will compare the behaviour of \eqref{ch1sys:fieldroad} to the one of the system
\begin{equation}\label{ch1sys:symmetric}
\left\{
\begin{array}{ll}
v_t-d\Delta v - c\partial_x v= f(x,v),  & (x,y)\in\Omega,\\
-\partial_y v|_{y=0} =0, & x\in\R,
\end{array}
\right.
\end{equation}
whose solution is a function $v(x,y)$ that can be extended by symmetry to the whole plane, thanks to the Neumann border condition.
It is natural to consider system \eqref{ch1sys:symmetric} as the counterpart of system \eqref{ch1sys:fieldroad} in the case without the road, since it presents the same geometry, including the same boundary condition exception made for the exchange terms that are in place for the case of a fast diffusion channel.

\subsection{State of the art}

We present here the background that led us consider system \eqref{ch1sys:fieldroad} and some useful results that are known in the community.

The study of reaction-diffusion equations started with the works by Fisher \cite{fisher} and by Kolmogorov, Petrowskii and Piskunov \cite{KPP}, who modelled the spacial diffusion of an advantageous gene in a population living in a one-dimensional environment through the equation
\begin{equation}\label{ch1eq:KPP}
\partial_t v -d \, \partial_{xx}^2 v = f(v) 
\end{equation}
for $x\in\R$ and $t\geq 0$. For \eqref{ch1eq:KPP}, it is supposed that $d>0$ and $f\geq 0$ is a $\mathcal{C}^1$ function satisfying $f(0)=f(1)=0$ and the KPP hypothesis $f(v)\leq f'(0)v$ for $v\in[0,1]$. The first example was a nonlinearity of logistic type, so $f(v)= av(1-v)$ for some $a>0$.
It was shown any solution $v$ issued from a nonnegative initial datum $v_0$ converges to 1 as $t$ goes to infinity, locally uniformly in space; this long time behaviour is called \emph{invasion}. 
The generalisation in higher dimension of equation \eqref{ch1eq:KPP} was then used to study the spatial diffusion of animals, plants, bacteria and epidemics \cite{skellam, okubo1980diffusion}.

A vast literature has been originated from the pioneer works, studying various aspects of the homogeneous equation \eqref{ch1eq:KPP}, in particular concerning the \emph{travelling fronts}. These are solutions of the form $v(t,x)= V(x \cdot e +ct)$ with $V:\R\to[0,1]$, for $e$ a direction, the \emph{direction of propagation},  and $c$ the \emph{speed of propagation of the travelling front}. Other than this, researchers have investigated the \emph{asymptotic speed of propagation} at which level sets of a solution starting from $v_0$ expands. These topics arose already in \cite{fisher} and \cite{KPP}, and their investigation was continued in many interesting articles, among which \cite{fife} and \cite{weinberger2}.

The correspondence of the theoretical results with actual data as seen in \cite{skellam} was encouraging, however it was clear that natural environments, even at macroscopic levels, were not well represented by a homogeneous medium, due to the alternation of forests, cultivated fields, plains, scrubs and many other habitats, as well as roads, rivers and other barriers \cite{kinezaki2003modeling}. It was necessary to look at more sophisticated features, as the effects of inhomogeneity, fragmentation, barriers and fast diffusion channels, and on the top of that, climate change.

A first analysis was carried out 
in \cite{shigesada1986traveling, shigesada1997biological} and 
in \cite{kinezaki2003modeling} for the so-called the \emph{patch model}. The authors considered a periodic mosaic of two different homogeneous habitats, one favorable and one unfavorable for the invading species. 
In \cite{kinezaki2003modeling}, the authors studied the long time behaviour of the population starting from any nonnegative initial datum. For further convenience, let us give the following definition:
\begin{definition}\label{ch1def:pers_extin}
	For the equation of \eqref{ch1eq:KPP} or the system \eqref{ch1sys:fieldroad}, we say that 
	\begin{enumerate}
		\item \emph{extinction} occurs if any solution starting from a non negative bounded initial datum converges to $0$ or to $(0,0)$  uniformly as $t$ goes to infinity.
		\item \emph{persistence} occurs if any solution starting from a non negative, non zero, bounded initial datum converges to a positive  stationary solution locally uniformly as $t$ goes to infinity.
	\end{enumerate}
\end{definition}

In \cite{kinezaki2003modeling}, it was empirically shown that the stability of the trivial solution $v=0$ determines the long time behaviour of the solutions. A solid mathematical framework for a general periodic environment was given in \cite{bhroques}. There, the authors considered the equation
\begin{equation}\label{ch1eq:bhroques}
\partial_t v - \nabla \cdot (A(x)\cdot \nabla v) = f(x, v)
\end{equation}
for $x\in \R^N$ and $t\geq 0$. 
The diffusion matrix $A(x)$ is supposed to be $\mathcal{C}^{1, \alpha}$, uniformly elliptic and periodic; however, for our interest we can suppose $A(x)=d\, I_{N}$, where $I_N$ is the identity matrix. The nonlinearity $f: \R^N \times \R_{\geq0} \to \R$ is supposed to be $\mathcal{C}^{1}$ in $x$, locally in $v$, and Lipshitz in $v$, uniformly in $x$, respecting hypothesis \eqref{ch1hyp:0}-\eqref{ch1hyp:KPP} and such that for some $L=(L_1, \dots, L_N)$, with $L_i\geq 0$, it holds
\begin{equation}\label{ch1hyp:per'}
f(x+L,s)=f(x,s) \quad \text{for all} \  s\geq 0 \ \text{and all} \ x\in\R^N.
\end{equation}

The criterion for persistence or extinction is given via a notion of periodic eigenvalue, that is the unique number $\lambda_p(-\mathcal{L}, \R^N)$
such that there exists a solution $\psi\in W_{loc}^{2, p}(\R^N)$ to the system
\begin{equation}\label{ch1sys:L_RN_p}
\left\{
\begin{array}{ll}
\mathcal{L'}(\psi) + \lambda \psi = 0, & x\in\R^N, \\
\psi> 0, &  x\in\R^N, \\
|| \psi ||_{\infty}=1, \\
\psi \ \text{is periodic in $x$ of periods $L$},
\end{array}
\right.
\end{equation}
where $\mathcal{L'}$ is given by
\begin{equation}\label{ch1def:mathcal_L'}
\mathcal{L'}(\psi):=  d \Delta \psi  + f_v(x,0)\psi.
\end{equation}
We point out that the existence and uniqueness of $\lambda_p(-\mathcal{L}, \R^N)$ is guaranteed by Krein-Rutman theory. The long time behaviour result in \cite{bhroques} is the following:

\begin{theorem}[Theorem 2.6 in \cite{bhroques}]\label{ch1thm:2.6inbhroques}
	Assume $f$ satisfies \eqref{ch1hyp:0}-\eqref{ch1hyp:KPP} and \eqref{ch1hyp:per'}. Then: 
	\begin{enumerate}
		\item If $\lambda_p(-\mathcal{L'}, \R^N)<0$, persistence occurs for \eqref{ch1eq:bhroques}.
		\item If $\lambda_p(-\mathcal{L'}, \R^N)\geq 0$, extinction occurs for \eqref{ch1eq:bhroques}.
	\end{enumerate}
\end{theorem}

To prove Theorem \ref{ch1thm:2.6inbhroques}, the authors performed an analysis of $\lambda_p(-\mathcal{L}, \R^N)$, proving that it coincide with the limit of eigenvalues for a sequence of domains invading $\R^N$, so that it coincides with the generalised principal eigenvalue of the system ``without the road'' \eqref{ch1sys:symmetric}. Nowadays, that and many other properties of this eigenvalue can be found as part of a broader framework in \cite{br}. In Section \ref{ch1s:eigenvalues}, we will provide further comments on it.

Another important fact highlighted both in the series in \cite{shigesada1986traveling, shigesada1997biological, kinezaki2003modeling} and in \cite{bhroques} is that the presence of multiple small unfavourable zones gives less chances of survival than one large one, the surface being equal. 

\medskip

A new difficulty that one may consider while studying ecological problems is, sadly, the issue of a changing climate. A 1-dimensional model in this sense was first proposed in
\cite{berestycki2009can} and \cite{potapov2004climate}, and was later treated in higher dimension in \cite{br2}. 
The authors first imagined that a population lives in a favourable region enclosed into a disadvantageous environment; due to the climate change, the favourable zone starts to move in one direction, but keeps the same surface. The resulting equation is 
\begin{equation}\label{ch11709}
\partial_t v -  \Delta v=f(x-ct e,v) \quad \text{for} \ x\in\R^N,
\end{equation}
with $e$ a direction in $\mathbb{S}^{N-1}$ and $f: \R^N\times \R_{\geq 0} \to \R$. It was observed that a solution to \eqref{ch11709} in the form of a travelling wave $v(x,t)=V(x-cte)$ solves the equation
\begin{equation}\label{ch1eq:cc}
\partial_t V -  \Delta V- c\,e\cdot \nabla V=f(x,V) \quad \text{for} \ x\in\R^N,
\end{equation}
which is more treatable. The main question is if the population keeps pace with the shifting climate, that is, if the species is able to migrate with the same speed of the climate. The answer to this question is positive if a solution to \eqref{ch1eq:cc} exists; this depends on the value of $c$. We point out that already in \cite{br2} the authors considered the general case of a possible periodic $f(x,v)$.

\medskip

As mentioned before, another feature worth investigation is the effect of fast diffusion channels on the survival and the spreading of species. In fact, the propagation of invasive species as well as epidemics is influenced by the presence of roads \cite{robinet2012human, gatto2020spread}.
This observations led Berestycki, Roquejoffre and Rossi to propose a model for ecological diffusion in the presence of a fast diffusion channel in \cite{brr}, the so-called \emph{road-field model}. The field is modelled with the halfplane $\Omega=\R \times \R_{+}$ and the line with the $x$ axis; the main idea is to use two different variables for modelling the density of population along the line, $u$, and on the half plane, $v$. The system reads
\begin{equation*} 
\left\{
\begin{array}{ll}
\partial_t u(x,t) - D \partial_{xx}^2 u (x,t) = \nu v (x,0,t) - \mu u(x,t), &  x\in \R, t > 0, \\
\partial_t v(x,y,t) - d \Delta v (x,y,t)= f(v), & (x,y) \in \Omega, t>0, \\
-d \partial_y v(x,0,t) = -\nu v(x,0,t) + \mu u(x,t), & x \in \R, t>0, 
\end{array} \right.
\end{equation*}
for $D$, $d$, $\nu$, $\mu$ positive constants; moreover, $f\in \mathcal{C}^1$ was supposed to satisfy
\begin{equation*}
f(0)=f(1)=0, \quad 0< f(s) < f'(0)s \ \text{for} \ s \in (0,1), \quad f(s)<0 \ \text{for}  \ s>1.
\end{equation*}
The three equations describe, respectively, the dynamic on the line, the dynamic on the half plane and the exchanges of population between the line and the half plane. On the line, the diffusion is faster than in  $\Omega $ if $D>d$. 
In \cite{brr}, the authors identify the unique positive stationary solution $\left(\frac{1}{\mu}, 1 \right)$ and prove persistence of the population. 
Moreover, they show that the presence of the line increases the spreading speed. Another version of the model with a reaction term for the line was presented by the same authors in \cite{BRR2}, while many variation of the models were proposed by other authors: with nonlocal exchanges in the direction of the road \cite{pauthier2015uniform, pauthier2016influence}, with nonlocal diffusion \cite{berestycki2015effect, berestycki2014speed}, and with different geometric settings \cite{rossi2017effect}. For a complete list, we refer to \cite{tellini2019comparison}.

The case of heterogeneous media for systems of road-field type has been so far not much treated, due to its difficulties. A first road-field model with exchange terms that are periodic in the direction of the road was proposed in \cite{giletti2015kpp}. There, the authors 
recovered the results of persistence and of acceleration on the propagation speed due to the road known in the homogeneous case;
they also studied the spreading of solution with exponentially decaying initial data and calculated their speeds.

Recently, Berestycki, Ducasse and Rossi introduced in \cite{romain} a new generalised principal eigenvalue fitting road-field system for possibly heterogeneous reaction term; here, we give its definition directly for the system \eqref{ch1sys:fieldroad}.  
Calling
\begin{equation}\label{ch1sys:operators}
\left\{
\begin{array}{l}
\mathcal{R}(\phi, \psi):=D \phi''+c \phi'+\nu {\psi}|_{y=0}-\mu \phi, \\
\mathcal{L}(\psi):= d\Delta \psi +c \partial_x \psi -f_v(x,0)\psi, \\
B(\phi, \psi):=d \partial_y {\psi}|_{y=0}+\mu \phi- \nu {\psi}|_{y=0},
\end{array}
\right.
\end{equation}
this eigenvalue is defined as 
\begin{equation}\label{ch1def:lambda1_S_Omega}
\begin{split}
\lambda_1( \Omega)=\sup \{ \lambda \in \R \ : \ \exists (\phi, \psi)\geq (0,0), \ (\phi, \psi) \not\equiv(0,0), \ \text{such that} \\ \mathcal{L}(\psi) + \lambda \psi \leq 0 \ \text{in} \ \Omega, \ \mathcal{R}(\phi, \psi) +\lambda \phi \leq 0  
\ \text{and} \ B(\phi, \psi)\leq 0 \ \text{in} \ \R \},
\end{split}
\end{equation}
with $(\phi, \psi)$ belonging to $W_{loc}^{2,3}(\R)\times W_{loc}^{2,3}(\overline{\Omega})$. Together with the definition, many interesting properties and bounds were studied; we will recall some of them later.

Thanks to that, the same authors were able to investigate the case of 
a favourable ecological niche, possibly facing climate change, in \cite{econiches}. It was proven that the sign of $\lambda_1( \Omega)$ characterises
the extinction or the persistence of the population; moreover, comparing the results with  the ones found for the model without the road, a deleterious effect of the road on the survival chances is always found when there is no climate change. On the other hand, if the ecological niche shifts, the road has in some cases a positive effect on the persistence.

\subsection{Main results}

We are now ready to present the main results of this chapter.

\subsubsection{The case of a periodic $f(x,v)$}

Here, we consider the case of a nonlinearity that respects the KPP hypothesis and is periodic in the direction of the road. Moreover, here we always consider $c=0$.

We begin by the following result on long time behaviour for solutions of system \eqref{ch1sys:fieldroad}:

\begin{theorem} \label{ch1thm:char}
	Assume $f$ satisfy \eqref{ch1hyp:0}-\eqref{ch1hyp:per}, $c=0$ and let $\lambda_1(\Omega)$ be as in \eqref{ch1def:lambda1_S_Omega}.
	Then the following holds:
	\begin{enumerate}
		\item if $\lambda_1( \Omega)\geq 0$, then  extinction occurs.
		\item if $\lambda_1(\Omega)<0$, then persistence occurs and the  positive stationary solution $(u_{\infty}, v_{\infty})$ is unique and periodic in $x$.
	\end{enumerate} 
\end{theorem}

Now, we compare the behaviour of solutions to the system \eqref{ch1sys:fieldroad} with the ones of system \eqref{ch1sys:symmetric}. 
This allows us to highlight the effects of the fast diffusion channel on the survival chances of the population.
Actually, since solutions of \eqref{ch1sys:symmetric} can be extended by refection to the whole plane, we can make the comparison with equation \eqref{ch1eq:bhroques} for $A(x)=d I_2$ and $L=(\ell, 0)$.
The comparison is performed thanks to the generalised principal eigenvalue $\lambda_1(\Omega)$ for system \eqref{ch1sys:fieldroad} and the periodic eigenvalue $\lambda_p(-\mathcal{L}, \R^2)$, as defined in \eqref{ch1sys:L_RN_p}, for the operator $\mathcal{L}$ in dimension 2.
We obtain the following:

\begin{theorem}\label{ch1thm:comparison}
	Assume $f$ respects hypothesis \eqref{ch1hyp:0}-\eqref{ch1hyp:per}, $c=0$. Then:
	\begin{enumerate}
		\item if  $\lambda_p(-\mathcal{L}, \R^2)<0$, then $\lambda_1( \Omega)<0$, that is, if persistence occurs  for the system ``without the road'' \eqref{ch1eq:bhroques}, then it occurs also for system ``with the road'' \eqref{ch1sys:fieldroad}.
		\item if $\lambda_p(-\mathcal{L}, \R^2)\geq 0$, then $\lambda_1( \Omega)\geq 0$, that is, if extinction occurs for the system ``without the road''  \eqref{ch1eq:bhroques}, then it occurs also for system ``with the road'' \eqref{ch1sys:fieldroad}.
	\end{enumerate}	
\end{theorem}


Theorem \ref{ch1thm:comparison} says that the road has no negative impact on the survival chances of the population in the case of a periodic medium depending only on the variable in the direction of the road.  
This is surprising if compared to the results obtained in \cite{econiches} (precisely Theorem 1.5, part \emph{(ii)}), where the authors find that the existence of the road is deleterious in presence of an ecological niche, and even more counter-intuitive owing the fact that fragmentation of the environment lessens the survival chances of a population, as shown in \cite{bhroques}. This means that, in the case of periodic media, the presence of the fast diffusion channel does not interfere with the persistence of the population, which depends only on the  environment of a periodicity cell.
As seen in \cite{bhroques}, where the dependence of persistence on the amplitude of fragmentation was studied, if the favourable zones are sufficiently large, the population will eventually spread in all of them; the presence of the road does not cause loss of favourable environment and consequently of persistence chances.
However, we expect the spreading speed to be influenced by the presence of the road, as it has been already proven in the case of homogeneous environment.

We point out that Theorem \eqref{ch1thm:char} completes and is in accordance with the results on long time behaviour found in \cite{brr} for a homogeneous reaction term, which we can see as a particular case of periodicity, which respects positive KPP hypothesis (where the positivity is requested through $f'(0)>0$). In \cite{brr}, Theorem 4.1 states the convergence of any positive solution to the unique positive stationary solution of the system.  Since it is well known that for the homogeneous case it holds $\lambda_1(-\mathcal{L}, \R^2)=- f'(0)$, the positivity hypothesis gives that $\lambda_1(-\mathcal{L}, \R^2)<0$ and, as a consequence of Theorem \ref{ch1thm:ineq}, that the second case in our Theorem \ref{ch1thm:char} occurs. 
If instead we asked for $f'(0)\leq0$, then we would be in the first case of Theorem \ref{ch1thm:char}, yielding extinction of the population.

\paragraph{Effects of amplitude of heterogeneity.}
One may expect that the presence of a road may alter the complex interaction between more favourable and less favourable zones, in particular penalising the persistence, since it was shown that populations prefer a less fragmented environment. However, the road does not interfere with that; as a consequence, also for environments presenting fast diffusion channels, some results of the analysis on the effect of fragmentation performed in \cite{bhroques} holds.

Take a parameter $\alpha>0$ and consider system \eqref{ch1sys:fieldroad} with nonlinearity
\begin{equation}\label{ch11421}
\tilde{f}(x,v)=\alpha f(x,v).
\end{equation}
To highlight the dependence on $\alpha$, we will call $\lambda_1(\Omega, \alpha)$ the  generalised principal eigenvalue defined in \eqref{ch1def:lambda1_S_Omega} with nonlinearity $\tilde{f}$.
As a direct consequence of Theorem \eqref{ch1thm:comparison} and Theorem 2.12 in \cite{bhroques}, we have the following result on the amplitude of heterogeneity:

\begin{corollary}
	Assume $\tilde{f}$ is defined as in \eqref{ch11421}, $f$ satisfies \eqref{ch1hyp:0}-\eqref{ch1hyp:per}, and $c=0$. Then:
	\begin{enumerate}
		\item if $ \int_{0}^{\ell} f_v(x,0)>0$, or if $ \int_{0}^{\ell} f_v(x,0)=0$ and $f\not\equiv 0$, then for all $\alpha >0$ we have $\lambda_1(\Omega, \alpha  )<0$.
		\item if $ \int_{0}^{\ell} f_v(x,0)<0$, then  $\lambda_1(\Omega, \alpha )>0$ for $\alpha$ small enough; if moreover there exists $x_0\in[0,\ell]$ such that $f_v(x_0,0)>0$, then for all $\alpha$ large enough $\lambda_1(\Omega, \alpha )<0$.
	\end{enumerate}
\end{corollary}

This result describes with precision the fact that, to persist, a species must have a sufficiently large favourable zone available. If the territory is more advantageous than not, then the population persist. If however there environment is generally unfavourable, the population persists only if there are some contiguous advantageous zones large enough; if instead the advantageous zones are fragmented, even if there is unlimited favourable territory, the population will encounter extinction.

\subsubsection{A climate change setting for a general $f(x,v)$}

We consider now a general nonlinearity that depends on the spatial variable in the direction of the road. We stress the fact that we do not suppose any periodicity, but the case of a periodic $f$ is a particular case of this setting. Moreover, the following result is done in the general framework of a possible climate change, so the parameter $c$ may be different from $0$.

Comparison between the systems with and without the road, in the general case, are done through comparison between $\lambda_1(\Omega)$ and the generalised principal eigenvalue of system \eqref{ch1sys:symmetric}, given by
\begin{equation}\label{ch1lambda:L_Omega}
\begin{split}
\lambda_1(-\mathcal{L}, \Omega)=\sup \{ \lambda \in \R \ : \ \exists \psi \geq 0, \psi \not\equiv 0 \ \text{such that} \\ 
\mathcal{L}(\psi) + \lambda \psi \leq 0 \ \text{on} \ \Omega, \ -\partial_y \psi|_{y=0}\leq 0 \ \text{on} \ \R \}
\end{split}
\end{equation}
for $\psi\in W_{loc}^{2,3}(\Omega)$. With this notation, we have the following:

\begin{theorem}\label{ch1thm:ineq}
	Assume $\lambda_1(-\mathcal{L}, \R^2)$ as in \eqref{ch1lambda:L_Omega} and $\lambda_1(\Omega)$ as in \eqref{ch1def:lambda1_S_Omega}; then $\lambda_1(-\mathcal{L}, \R^2) \geq \lambda_1(\Omega)$.
\end{theorem}

In the special case $c=0$, some information on the relations between $\lambda_1(-\mathcal{L}, \R^2)$ and $\lambda_1(\Omega)$  was already available in \cite{econiches}: Proposition 3.1 yields that $\lambda_1(-\mathcal{L}, \R^2)\geq 0$ implies $\lambda_1(\Omega)\geq 0$. Thanks to that and Theorem \ref{ch1thm:ineq}, the following result holds:

\begin{corollary}\label{ch1thm:ineq2}
	If $c=0$, we have $\lambda_1(-\mathcal{L}, \R^2)<0$ if and only if $\lambda_1(\Omega)<0$.
\end{corollary}

As already pointed out in \cite{romain}, even for $c=0$ it is not true that $\lambda_1(-\mathcal{L}, \R^2) =\lambda_1(\Omega)$. In fact, it has been found that $\lambda_1(\Omega) \leq \mu$, while playing with $f$ one can have $\lambda_1(-\mathcal{L}, \R^2)$ as large as desired. However, the fact that the two eigenvalues have the same sign reveals that they are profoundly linked.


\subsection{Organisation of the chapter}

In Section \ref{ch1s:eigenvalues}, we recall and discuss the properties of the eigenvalues $\lambda_1(\Omega)$, $\lambda_1(-\mathcal{L}, \R^2)$ and $\lambda_p(-\mathcal{L}, \R^2)$ already known in the literature. 
Furthermore, a periodic eigenvalue for the system \eqref{ch1sys:fieldroad} will be defined; because of the presence of the road, the periodicity is present only in the $x$ direction. As a consequence, it is useful to define an analogous generalised eigenvalue for the system without the road \eqref{ch1sys:symmetric} with periodicity only in the direction of the road.

In Section \ref{ch1s:ordering}, one finds the proof of Theorem \ref{ch1thm:ineq} and Theorem \ref{ch1thm:comparison}. Moreover, the relations between the newly defined generalised periodic eigenvalues and the known ones are shown.

The last Section \ref{ch1s:lb} treats large time behaviour for solutions to \eqref{ch1sys:fieldroad} with $c=0$ and periodic $f$; this includes the proof of Theorem \ref{ch1thm:char}.

\section{Generalised principal eigenvalues and their properties} \label{ch1s:eigenvalues}

Both road-field models and reaction-diffusion equations in periodic media have been treated in several papers.
In this section, we introduce some useful objects and recall their properties. 
All along this section we will make repeated use of the operators $\mathcal{L}$, $\mathcal{R}$ and $B$, that were  defined in \eqref{ch1sys:operators}.

\subsection{Eigenvalues in periodic media}

Since $\mathcal{L}$ has periodic terms, it is natural to look for eigenfunctions that have the same property. However, to begin the discussion on the periodic eigenvalue for the operator $\mathcal{L}$ in $\R^2$, we consider  its counterpart in $\R$. 
We look for the unique number $\lambda_p(-\mathcal{L}, \R)\in\R$ such that there exists a function $\psi\in W_{loc}^{2, 3}(\R)$ solution to the problem
\begin{equation}\label{ch1sys:L_R_p}
\left\{
\begin{array}{ll}
d\psi''+f_v(x, 0)\psi + \lambda \psi = 0, & x\in\R, \\
\psi> 0, &  x\in\R, \\
|| \psi ||_{\infty}=1, \\
\psi \ \text{is periodic in $x$ of period $\ell$.}
\end{array}
\right.
\end{equation}
In \eqref{ch1sys:L_R_p}, the operator $\mathcal{L}$ has been replaced by an operator working on $\R$, namely the Laplacian has been substituted by a double derivative.
Notice that existence and uniqueness of the solution to \eqref{ch1sys:L_R_p}, that we call $(\lambda_p(-\mathcal{L}, \R), \psi_p)$, is guaranteed  by Krein-Rutman theory.

For the operator $\mathcal{L}$, since it has no dependence on the $y$ variable,
we have to introduce a fictive periodicity in order to be able to use the Krein-Rutman theory. Thus, fix $\ell'>0$ and consider the problem in $\R^2$ of finding the value $\lambda_p(-\mathcal{L}, \R^2)\in\R$ such that there exists a solution $\psi\in W_{loc}^{2, 3}(\R^2)$ to the system
\begin{equation}\label{ch1sys:L_R2_p}
\left\{
\begin{array}{ll}
\mathcal{L}(\psi) + \lambda \psi = 0, & (x,y)\in\R^2, \\
\psi> 0, &  (x,y)\in\R^2, \\
|| \psi ||_{\infty}=1, \\
\psi \ \text{is periodic in $x$ and $y$ of periods $\ell$ and $\ell'$}.
\end{array}
\right.
\end{equation}
Again we can use the Krein-Rutman Theorem to see that there exists a unique pair $(\lambda_p(-\mathcal{L}, \R^2), \psi_{\ell'})$ solving \eqref{ch1sys:L_R2_p}. Now, with a slight abuse of notation, we consider the function $\psi_p(x,y)$ as the extension in $\R^2$ of $\psi_p$ solution to \eqref{ch1sys:L_R_p}. We observe that the pair $(\lambda_p(-\mathcal{L}, \R), \psi_p)$ gives a solution to \eqref{ch1sys:L_R2_p}. Hence, by the uniqueness of the positive eigenfunction, we get 
\begin{equation}\label{ch1eq:-5}
\lambda_p(-\mathcal{L}, \R^2)=\lambda_p(-\mathcal{L}, \R) \quad \text{and} \quad \psi_p\equiv \psi_{\ell'}.
\end{equation}
This also implies that neither $\lambda_p(-\mathcal{L}, \R^2)$ nor $\psi_{\ell'}$ depend on the parameter $\ell'$ that was artificially introduced. From now on, we will use only $ \psi_p$.

The properties of the eigenvalue $\lambda_p(-\mathcal{L}, \R^2)$ were also studied in \cite{br}, where it is called $\lambda'$ and defined as
\begin{equation}\label{ch1def:lambdap_dim2}
\begin{split} 
\lambda_p(-\mathcal{L}, \R^2)= &\inf \{ \lambda \in \R \ :\  \exists \varphi\in \mathcal{C}^2(\R^2)\cap L^{\infty}(\R^2), \ \varphi>0,  \\ &\hspace{8em} \varphi \ \text{periodic in $x$ and $y$}, \ \mathcal{L}(\varphi)+\lambda \varphi \geq 0    \}.
\end{split}
\end{equation}
In particular, in Proposition 2.3 of \cite{br} it is stated that the value found with \eqref{ch1sys:L_R2_p} coincides with the one defined in \eqref{ch1def:lambdap_dim2}.

\subsection{Generalised principal eigenvalues for the system with and without the road and some properties}

In this section, we are going to treat eigenvalues that are well defined also for non periodic reaction functions. 

The generalised eigenvalue $\lambda_1( \Omega)$ for the system \eqref{ch1sys:fieldroad}, that we defined in \eqref{ch1def:lambda1_S_Omega}, was first introduced in \cite{romain}. Together with this, the authors also proved the interesting property that $\lambda_1( \Omega)$ coincides with the limit of principal
eigenvalues of the same system restricted to a sequence of invading domains. 
They use some half ball domains defined as follow for $R>0$:
\begin{equation}\label{ch11722}
\Omega_R:=B_R\cap \Omega \quad \text{and} \quad I_R:=(-R, R).
\end{equation}
Then we have the following characterisation for $\lambda_1( \Omega)$:
\begin{proposition}[Theorem 1.1 of \cite{romain}] \label{ch1prop:romain}
	For $R>0$, 
	there is a unique 
	$\lambda_1( \Omega_R) \in \R$
	and a unique (up to multiplication by a positive scalar) positive 
	$(u_R, v_R) \in W^{2,3}(I_R) \times W^{2,3} (\Omega_R)$ that satisfy the eigenvalue problem
	\begin{equation}\label{ch1sys:halfball}
	\left\{
	\begin{array}{ll}
	\mathcal{R}(\phi, \psi) +\lambda\phi = 0, & x\in I_R, \\
	\mathcal{L}(\psi)  + \lambda \psi = 0, &(x,y)\in \Omega_R, \\
	B(\phi, \psi)= 0, & x\in I_R, \\
	\psi =0, & (x,y)\in (\partial\Omega_R ) \setminus (I_R\times \{0\}) \\
	\phi(R)=\phi(-R)=0. &
	\end{array}
	\right.
	\end{equation}
	Moreover, 
	\begin{equation*}
	\lambda_1( \Omega_R) \underset{R\to +\infty}{\searrow} \lambda_1( \Omega).
	\end{equation*}
\end{proposition}


We also consider the principal eigenvalue on the truncated domains for the linear operator $\mathcal{L}(\psi)$. 
To do that, for any $R>0$ we call $B_R^P$ the ball of centre $P=(x_P,y_P)$ and radius $R$.  We define $\lambda_1(-\mathcal{L}, B_R^P)$ as the unique real number such that
the problem
\begin{equation}\label{ch1sys:L_BR}
\left\{
\begin{array}{ll}
\mathcal{L}(\psi_R) + \lambda_1(-\mathcal{L}, B_R^P) \psi_R = 0, & (x,y)\in B_R^P, \\
\psi_R=0, & (x,y)\in \partial B_R^P, \\ 
\psi_R >0, & (x,y)\in B_R^P
\end{array}
\right.
\end{equation}
admits a solution $\psi_R\in W^{2,3}(B_R^P)$.
The existence and uniqueness of such quantity and its eigenfunction is a well-known result derived via the Krein-Rutman theory. 
We also notice that, calling $B_R$ the ball with radius $R$ and center $O=(0,0)$, the pair $(\lambda_1(-\mathcal{L}, B_R), \psi_R)$ is also a solution to the problem
\begin{equation}\label{ch1sys:L_OmegaR}
\left\{
\begin{array}{ll}
\mathcal{L}(\psi) + \lambda \psi = 0, & (x,y)\in \Omega_R, \\
\partial_y \psi(x,0)= 0, & x\in I_R, \\
\psi=0, & (x,y)\in (\partial \Omega_R)\setminus (I_R\times \{0\}), \\
\psi >0, & (x,y)\in \Omega_R.
\end{array}
\right.
\end{equation}
The proof of that is very simple. 
If $(\lambda, \psi)$ is the unique solution to \eqref{ch1sys:L_OmegaR}, extending $\psi$ by symmetry in $B_R$ we get a solution to \eqref{ch1sys:L_BR}. By the uniqueness of the solution to \eqref{ch1sys:L_BR}, we get $\lambda=\lambda_1(-\mathcal{L}, B_R)$.

Similarly to what happens with $\lambda_1( \Omega_R)$, thanks to the fact $\lambda_1(-\mathcal{L}, B_R)$ solves \eqref{ch1sys:L_OmegaR}, we have that the sequence $\lambda_1(-\mathcal{L}, B_R)$ converges to the value $\lambda_1(-\mathcal{L}, \Omega)$, that was defined in \eqref{ch1lambda:L_Omega}. This was precisely stated in \cite{econiches} as:

\begin{proposition}[Proposition  2.4 of \cite{econiches}]\label{ch1prop:lim_B_R}
	We have that
	\begin{equation}
	\lambda_1(-\mathcal{L}, B_R) \underset{R\to +\infty}{\searrow} \lambda_1(-\mathcal{L}, \Omega),
	\end{equation}
\end{proposition}

Another notion of generalised eigenvalue analysed in \cite{br} is the quantity 
\begin{equation}\label{ch1lambda:L_R2}
\begin{split}
\lambda_1(-\mathcal{L}, \R^2)=\sup \{ \lambda \in \R \ : \ \exists \psi \geq 0, \psi \not\equiv 0 \ \text{such that}  \
\mathcal{L}(\psi) + \lambda \psi \leq 0 \ \text{a.e on} \ \R^2 \}
\end{split}
\end{equation}
for test functions $\psi \in W_{loc}^{2,3}(\R^2)$. As stated in Proposition 2.2 of \cite{br}, we have
\begin{equation*}
	\lambda_1(-\mathcal{L}, B_R) \underset{R\to +\infty}{\searrow} \lambda_1(-\mathcal{L}, \R^2).
\end{equation*}
By that and \eqref{ch1prop:lim_B_R}, we have 
\begin{equation*}
\lambda_1(-\mathcal{L}, \R^2)=\lambda_1(-\mathcal{L}, \Omega)
\end{equation*}
With this notation, we can report the following affirmations deriving from Theorem 1.7 in \cite{br} for the case of a periodic reaction function:

\begin{theorem}[Theorem 1.7 in \cite{br}]\label{ch1thm:1.7inbr}
	Suppose $f$ satisfies \eqref{ch1hyp:per}.
	The following holds:
	\begin{enumerate}
		\item It holds that $\lambda_p(-\mathcal{L}, \R^2)\leq \lambda_1(-\mathcal{L}, \Omega)$.
		\item If $\mathcal{L}$ is self-adjoint (i.e, if $c=0$), then $\lambda_p(-\mathcal{L}, \R^2)=\lambda_1(-\mathcal{L}, \Omega)$.
	\end{enumerate}
\end{theorem}

At last, we recall the following result on the signs of the eigenvalues for the systems with and without the road:
\begin{proposition}[Proposition 3.1 in \cite{econiches}]\label{ch1prop:ineq}
	It holds that
	\begin{equation*}
	\lambda_1(-\mathcal{L}, \Omega) \geq 0 \quad \Rightarrow \quad \lambda_1( \Omega) \geq 0.
	\end{equation*} 
\end{proposition}
This is the result that, in combination with Theorem \ref{ch1thm:ineq}, gives Corollary \ref{ch1thm:ineq2}.

\subsection{The periodic generalised principal eigenvalue for the road-field system}

We introduce here two new eigenvalues that will be useful in the following proofs. They are somehow of mixed type, in the sense that
they are periodic in $x$ but not in $y$; this derives from the fact that the domains in which they are defined are periodic in the variable $x$ and truncated in the variable $y$. Here, we require $f$ to be periodic as in hypothesis \eqref{ch1hyp:per}. 

Given $r>0$, let  $(\lambda_p(-\mathcal{L}, \R\times(-r, r)), \psi_{r})$ be the unique pair solving the eigenvalue problem
\begin{equation}\label{ch1sys:bary2}
\left\{
\begin{array}{ll}
\mathcal{L}(\psi_{r})  + \lambda \psi_{r} = 0, \qquad(x,y)\in \R \times (-r, r), \\
\psi_{r} (x, \pm r)=0, \qquad x\in \R, \\
||\psi_{r}||_{\infty}=1, \ \psi_{r} \ \text{is periodic in} \ x.
\end{array}
\right.
\end{equation}
The existence and uniqueness  of the solution to \eqref{ch1sys:bary2} derives once again from Krein-Rutman theory.

We point out that $\lambda_p(-\mathcal{L}, \R\times(-r,r))$ is decreasing in $r$ by inclusion of domains. So,
there exists a well defined value, that with a slight abuse of notation we call $\lambda_p (-\mathcal{L}, \Omega)$, such that
\begin{equation}\label{ch1eq:-2}
\lambda_p(-\mathcal{L}, \R\times(-r,r)) \underset{r\to+\infty}{\searrow} \lambda_p (-\mathcal{L}, \Omega).
\end{equation}

Given $r>0$, there exists a unique  value $\lambda_p( \R\times(0, r))\in\R$ such that the problem
\begin{equation}\label{ch1sys:r}
\left\{
\begin{array}{ll}
\mathcal{R}(\phi, \psi) +\lambda\phi = 0, \qquad x\in \R, \\
\mathcal{L}(\psi)  + \lambda \psi = 0, \qquad(x,y)\in \R \times (0, r), \\
B(\phi, \psi)= 0, \qquad x\in \R, \\
\psi (\cdot, r)=0, \\
\phi \ \text{and} \ \psi  \ \text{are periodic in} \ x,
\end{array}
\right.
\end{equation}
has a solution.
The proof of the existence can be derived by modifying for periodic functions the proof of the existence of $\lambda_1( \Omega_R)$ that is found in the Appendix of \cite{romain}.

Moreover, we define
\begin{equation*}
\begin{split}
\lambda_p ( \Omega)= \sup \{ \lambda \in \R \ : \ \exists (\phi,\psi)\geq (0,0), \  (\phi,\psi) \ \text{periodic in} \ x, \ \text{such that} \\
\mathcal{R}(\phi, \psi) +\lambda \phi \leq 0,  
\mathcal{L}(\psi) + \lambda \psi \leq 0, \ \text{and} \ B(\phi, \psi)\leq 0  \}
\end{split}
\end{equation*}
with test functions $(\phi,\psi) \in W_{loc}^{2,3}(\R)\times W_{loc}^{2,3}(\overline{\Omega})$. 

Then, we have:

\begin{proposition}\label{ch1prop:-3}
	Suppose $f$ satisfies \eqref{ch1hyp:per}.
	We have that
	\begin{equation}\label{ch1eq:-3}
	\lambda_p( \R\times(0,r)) \underset{r\to+\infty}{\searrow} \lambda_p ( \Omega).
	\end{equation}
	Moreover, there exists a couple $(u_p, v_p)\in W_{loc}^{2,3}(\R)\times W_{loc}^{2,3}(\overline{\Omega})$ of positive functions periodic in $x$ such that satisfy
	\begin{equation}\label{ch1sys:upvp}
	\left\{
	\begin{array}{ll}
	\mathcal{R}(u_p, v_p)+ \lambda_p ( \Omega)v_p=0, & x\in\R, \\
	\mathcal{L}v_p+ \lambda_p ( \Omega) v_p=0, & (x,y)\in\Omega, \\
	B(u_p, v_p)=0, & x\in\R.
	\end{array}
	\right.
	\end{equation}

\end{proposition}

\begin{proof}
	By inclusion of domains, one has that $\lambda_p( \R\times(0, r))$ is decreasing in $r$.
	Let us call
	\begin{equation*}
	\bar{\lambda}:=\underset{r\to \infty}{\lim} \lambda_p( \R\times(0, r)).
	\end{equation*}
	
	\emph{Step 1}.
	We now want to show that there exists a couple $(\bar{\phi}, \bar{\psi})>(0,0)$, with $\bar{\phi}\in W_{loc}^{2,3}(\R)$ and $\bar{\psi}\in W_{loc}^{2,3}(\overline{\Omega})$, periodic in $x$, that satisfy
	\begin{equation}\label{ch11957}
	\left\{
	\begin{array}{ll}
	\mathcal{R}(\bar{\phi}, \bar{\psi})+ \bar{\lambda} \bar{\phi}=0, & x\in\R, \\
	\mathcal{L}( \bar{\psi})+ \bar{\lambda} \bar{\psi}=0, & (x,y)\in\Omega, \\
	B(\bar{\phi}, \bar{\psi})=0, & x\in\R.
	\end{array}
	\right.
	\end{equation}

	Fix $M>0$.
	First, for all $r>M+2$ consider the periodic eigenfunctions $(\phi_r, \psi_r)$ related to $\lambda_p( \R\times(0,r))$.
	We normalize $(\phi_r, \psi_r)$ so that
	\begin{equation*}
	\phi_r(0)+ \psi_r(0,0)=1.
	\end{equation*}
	
	Then, from the Harnack estimate in Theorem 2.3 of \cite{romain}, there exists $C>0$ such that
	\begin{equation}\label{ch11809}
	\max \{  \underset{I_{M+1}}{\sup} \phi_r, \ \underset{\Omega_{M+1}}{\sup} \psi_r  \} \leq C \min \{  \underset{I_{M+1}}{\inf} \phi_r, \ \underset{\Omega_{M+1}}{\inf} \psi_r  \} \leq C,
	\end{equation}
	where the last inequality comes from the normalization.
	We can use the interior estimate for $\phi_r$  and get
	\begin{equation*}
	|| \phi_r ||_{W^{2,3}(I_M)} \leq C' ( || \phi_r ||_{L^{3}(I_{M+1})}+ || \psi_r ||_{L^{3}(\Omega_{M+1})}   )
	\end{equation*}
	for some $C'$ depending on $M$, $\mu$, $\nu$, and $D$.
	By that and \eqref{ch11809}, we get
	\begin{equation}\label{ch11810}
	|| \phi_r ||_{W^{2,3}(I_M)} \leq C
	\end{equation}
	for a possibly different $C$.
	
	For $\psi_r$, in order to have estimates up to the border $y=0$ of $\Omega_M$, we need to make a construction. Recall that, calling $L:= \mathcal{L}+ \lambda_p( \R\times(0,r))$, $\psi_r$ solves
	\begin{equation*}
	\left\{
	\begin{array}{ll}
	L \psi_r =0, &  (x,y) \in \Omega_{M+1}, \\
	-d \partial_y \psi_r |_{y=0}  + \nu \psi_r|_{y=0}= \mu \phi_r, & x\in I_{M+1}.
	\end{array}
	\right.
	\end{equation*}
	We call
	\begin{equation*}
	\tilde{\psi}_r:= \psi_r e^{-\frac{\nu}{d}y}
	\end{equation*} 
	and the conjugate operator
	\begin{equation*}
	\tilde{L}(w):= e^{-\frac{\nu}{d}y} L\left( e^{\frac{\nu}{d}y} w \right).
	\end{equation*}
	Now, we have
	\begin{equation*}
	\left\{
	\begin{array}{ll}
	\tilde{L}\tilde{\psi}_r =0, &  (x,y) \in \Omega_{M+1}, \\
	-d \partial_y \tilde{\psi}_r |_{y=0}  = \mu \phi_r, & x\in I_{M+1}.
	\end{array}
	\right.
	\end{equation*}
	Next, calling
	\begin{equation}
	w_r(x,y)=\tilde{\psi}_r(x,y)- \frac{d}{\mu} \phi_r(x) y,  
	\end{equation}
	we have that
	\begin{equation}\label{ch11919}
	\left\{
	\begin{array}{ll}
	\tilde{L}w_r = -\dfrac{d}{\mu}  \tilde{L}( \phi_r(x) y), &  (x,y) \in \Omega_{M+1}, \\
	\partial_y {w_r}|_{y=0}  = 0, & x\in I_{M+1}.
	\end{array}
	\right.
	\end{equation}
	Now we define  in the open ball $B_{M+1}$ the function
	\begin{equation}\label{ch11911}
	\bar{w}_r(x,y):=w_r(x, |y|),
	\end{equation}
	that is the extension of $w_r$ by reflection; thanks to the Neumann condition in \eqref{ch11919} and the fact that ${w}_r \in W^{2,3}(\Omega_{M+1})$, we get that $\bar{w}_r \in W^{2,3}(B_{M+1})$. Also, we define the function 
	\begin{equation}\label{ch11912}
	g(x,y)= \frac{d}{\mu}  \tilde{L}( \phi_r(x) |y|).
	\end{equation}
	We also take the operator
	\begin{equation}\label{ch11913}
	\bar{L}w := d \Delta w+ c \partial_x w + 2{\nu} \sigma(y)\partial_y w + \left( f_v(x,0)+ \lambda_p( \R\times(0,r)) + \frac{\nu^2}{d} \right) w 
	\end{equation}
	where $\sigma(y)$ is the sign function given by
	\begin{equation*}
	\sigma(y) := \left\{
	\begin{array}{ll}
	1 & \text{if} \ y\geq 0, \\
	-1 & \text{if} \ y<0.
	\end{array}
	\right.
	\end{equation*}
	Thanks to the definition \eqref{ch11911}, \eqref{ch11912} and \eqref{ch11913}, we get that $\bar{w}_r $ is a weak solution to the equation
	\begin{equation}\label{ch11926}
	- \bar{L} \bar{w}_r = g \quad \text{for} \ (x,y)\in B_{M+1}.
	\end{equation}
	Finally, we can apply the interior estimates and get
	\begin{equation*}
	|| \bar{w}_r  ||_{W^{2,3}(B_M)} \leq C' ( || \bar{w}_r  ||_{L^{\infty}(B_{M+1})}+ || g ||_{L^{3}(B_{M+1})})
	\end{equation*}
	for some $C'$ depending on $M$ and the coefficients of the equation \eqref{ch11926}. But using the definition of $\bar{w}_r $ and the fact that $g$ is controlled by the norm of $\phi_r$,  we get, for a possible different $C'$,
	\begin{equation*}
	|| \bar{w}_r  ||_{W^{2,3}(B_M)} \leq C' ( || \psi_r  ||_{L^{\infty}(\Omega_{M+1})}+|| \phi_r  ||_{L^{\infty}(I_{M+1})}+ || \phi_r ||_{W^{2,3}(I_{M+1})}).
	\end{equation*}
	Using \eqref{ch11810} and \eqref{ch11809},
	we finally have 
	\begin{equation*}
	|| \psi_r  ||_{W^{2,3}(\Omega_M)} \leq C.
	\end{equation*}
	Thanks to that and \eqref{ch11810}, we have that $(\phi_r, \psi_r)$ is uniformly bounded in $W^{2,3}(I_M)\times W^{2,3}(\Omega_M)$ for all $M>0$.
	Hence, up to a diagonal extraction, $(\phi_r, \psi_r)$ converge weakly in $W_{loc}^{2,3}(I_M)\times W_{loc}^{2,3}(\Omega_M)$ to some $(\bar{\phi}, \bar{\psi}) \in W_{loc}^{2,3}(I_M)\times W_{loc}^{2,3}(\Omega_M)$. By Morrey inequality, the convergence is strong in $\mathcal{C}_{loc}^{1, \alpha}(\R)\times \mathcal{C}_{loc}^{1, \alpha}(\overline{\Omega})$ for $\alpha<1/6$. 
	Moreover, $(\bar{\phi}, \bar{\psi})$ are periodic in $x$ since all of the $(\phi_r, \psi_r)$ are periodic.
	Then, taking the limit of the equations in \eqref{ch1sys:r}, we obtain that $(\bar{\phi}, \bar{\psi})$ satisfy \eqref{ch11957}, as wished.

	\emph{Step 2.} We now prove that 
	\begin{equation}\label{ch11402}
	\bar{\lambda} \leq \lambda_p( \Omega).
	\end{equation}
	
	Take $\bar{\lambda}$ and 
	its associated periodic eigenfunctions couple $(\bar{\phi}, \bar{\psi})$ obtained in Step 1. 
	By definition, $\lambda_p( \Omega)$ is the supremum of the set 
	\begin{equation}\label{ch11747}
	\begin{split}
	\mathcal{A}:=  \{ \lambda \in \R \ : \ \exists (\phi,\psi)\geq (0,0), \  (\phi,\psi) \ \text{periodic in} \ x, \ 
	\mathcal{R}(\phi, \psi) +\lambda \phi \leq 0,  \\
	\mathcal{L}(\psi) + \lambda \psi \leq 0, \ \text{and} \ B(\phi, \psi)\leq 0  \}.
	\end{split}
	\end{equation}

	Then, using $(\bar{\phi}, \bar{\psi})$ as test functions, we obtain that $\bar{\lambda}$ is in the set $\mathcal{A}$ given in \eqref{ch11747}. By the fact that $\lambda_p( \Omega)$ is the supremum of $\mathcal{A}$, we get \eqref{ch11402}, as wished.
	
	\emph{Step 3}. We show 
	\begin{equation}\label{ch12006}
	\lambda_p(\Omega) \leq \bar{\lambda}.
	\end{equation}

	Now, take any $\lambda\in\mathcal{A}$ and one of its associate couple $(\phi, \psi)$. Then, by inclusion of domains, one gets that for all $r>0$ it holds 
	\begin{equation*}
	\lambda \leq \lambda_p( \R\times (0,r)).
	\end{equation*}
	Hence, by taking the supremum on the left hand side and the infimum on the right one, we get \eqref{ch12006}. By this and \eqref{ch11402}, equality is proven. Moreover, defining $(u_p, v_p)\equiv(\bar{\phi}, \bar{\psi})$, by \eqref{ch11957}, we have the second statement of the proposition.
\end{proof}

\section{Ordering of the eigenvalues}\label{ch1s:ordering}

This section is dedicated to show some inequalities and relations between the aforementioned eigenvalues.

\subsection{Proof of Theorem \ref{ch1thm:ineq}}

We start by proving Theorem \ref{ch1thm:ineq}.
We stress that this is done for the general setting of $c$ possibly non zero and $f(x,v)$ which may not be periodic.

\begin{proof}[Proof of Theorem \ref{ch1thm:ineq}] \quad
	
	Let us start by proving the first part of the theorem.
	For all $R>0$, there exists $R'>0$ and a point $C\in\R^2$ such that $B_R(C) \subset \Omega_{R'}$: it is sufficient to take $R'=3R$ and $C=(0, \frac{2}{3}R)$. We want to prove that
	\begin{equation}\label{ch11533}
	\lambda_1(-\mathcal{L}, B_R) \geq \lambda_1( \Omega_{R'}).
	\end{equation}
	Suppose by the absurd that \eqref{ch11533} is not true.
	Consider $\psi_R$ the eigenfunction related to $\lambda_1(-\mathcal{L}, B_R)$ and $v_{R'}$ the eigenfunction in the couple $(u_{R'},v_{R'})$ related to $\lambda_1( \Omega_{R'})$. 
	Since $\inf_{B_{R}(C)} v_{R'} >0$, and both eigenfunctions are bounded, there exists
	\begin{equation*}
	\theta^* := \sup \{  \theta\geq 0 \ : \ v_{R'}>\theta \psi_R \ \text{in} \   B_R(C)  \} >0.
	\end{equation*}
	Since $\theta^*$ is a supremum, then there exists $(x^*,y^*)\in \overline{B_R(C)}$ such that $v_{R'}(x^*, y^*)= \theta^* \psi_R (x^*, y^*)$. 
	Then,  $(x^*,y^*)\in {B_R(C)}$ because $v_{R'}>0$ and $\psi_R=0$ in $\partial B_R(C)$.
	Calling $\rho=v_{R'}-\theta^* \psi_R$, in a neighbourhood of $(x^*,y^*)$ we have that
	\begin{equation}\label{ch11602}
	-d \Delta\rho- c \cdot \nabla \rho - f_v(x,0)\rho=\lambda_1(-\mathcal{L}, B_R)\rho + (\lambda_1( \Omega_{R'}) - \lambda_1(-\mathcal{L}, B_R)) v_{R'}.
	\end{equation}
	We know that $\rho(x^*,y^*)=0$ and that $\rho \geq 0$ in $B_R(C)$. Then $(x^*,y^*)$ is a minimum for $\rho$, so $\nabla \rho(x^*,y^*) =0$ and $\Delta \rho (x^*,y^*) \geq 0$.
	Thus, the lefthandside of \eqref{ch11602} is non positive. But by the absurd hypotesis we have $(\lambda_1( \Omega_{R'}) - \lambda_1(-\mathcal{L}, B_R)) v_{R'}>0$. This gives
	\begin{equation*}
	0 \geq -d \Delta\rho(x^*,y^*) = (\lambda_1( \Omega_{R'}) - \lambda_1(-\mathcal{L}, B_R)) v_{R'} (x^*,y^*)>0,
	\end{equation*}
	which
	is a contradiction. With that we obtain that \eqref{ch11533} is true.
	
	Notice that the eigenvalue $\lambda_1(-\mathcal{L}, B_R(C))=\lambda_1(-\mathcal{L}, B_R)$, where $B_R$ is the ball centred in $(0,0)$, because $f(x,v)$ does not depend on $y$, thus system \eqref{ch1sys:L_BR} on $B_R(C)$ and $B_R$ are the same. As a consequence, also their eigenfunctions coincide.
	
	Recall that both $\lambda_1(-\mathcal{L}, \R^2)$ and $\lambda_1( \Omega)$ are limits of eigenvalues on limited domains, by \eqref{ch1prop:lim_B_R}  and Proposition \ref{ch1prop:romain}.
	Now, since for all $R>0$ there exists $R'$ such that \eqref{ch11533} is true, then passing to the limit we find the required inequality.
	
	%
	%
	%
	
\end{proof}

\subsection{Further inequalities between the eigenvalues}

In this section, we collect some results on the ordering of periodic and generalised eigenvalues for both system \eqref{ch1sys:fieldroad} and eqaution \eqref{ch1eq:bhroques}.  
Here we require $f$ to be periodic as in \eqref{ch1hyp:per}.

This first result is the analogue of Theorem \eqref{ch1thm:1.7inbr} for the system \eqref{ch1sys:fieldroad}:

\begin{theorem}\label{ch1thm:eq_1_p}
	Suppose $f$ respects hypothesis \eqref{ch1hyp:per}. Then:
	\begin{enumerate}
		\item It holds that $\lambda_1( \Omega)\geq \lambda_p( \Omega)$.
		\item If moreover $c=0$, then we have $\lambda_1( \Omega)= \lambda_p( \Omega)$.
	\end{enumerate}
\end{theorem}

\begin{proof}
	1.
	By definition, $\lambda_p( \Omega)$ is the supremum of the set $\mathcal{A}$ given in \eqref{ch11747},
	while $\lambda_1( \Omega)$ is the supremum of the set
	\begin{equation*}
	\begin{split}
	\{ \lambda \in \R \ : \ \exists (\phi,\psi)\geq (0,0),  \ 
	\mathcal{R}(\phi, \psi) +\lambda \phi \leq 0,  \\
	\mathcal{L}(\psi) + \lambda \psi \leq 0, \ \text{and} \ B(\phi, \psi)\leq 0  \} \supseteq \mathcal{A}.
	\end{split}
	\end{equation*}
	By inclusion of sets, we have the desired inequality.
	
	\medskip
	
	2.
	We call $$\mathcal{H}_R:= H_0^1(I_R)\times H_0^1(\Omega_R \cup (I_R\cup \{0\}) ). $$
	For $(u,v)\in \mathcal{H}_R$, we define
	\begin{equation*}
	Q_R(u,v):= \frac{ \mu \int_{I_R} D |u'|^2 + \nu \int_{\Omega_R} (d|\nabla v|^2-f_v(x,0)v^2) +  \int_{I_R} (\mu u- \nu v|_{y=0})^2      }{\mu \int_{I_R}u^2 + \nu \int_{\Omega_R} v^2}.
	\end{equation*}

	Now we fix $r>0$ and we consider $\lambda_p( \R \times (0,r)  )$ ad its periodic eigenfunctions $(\phi_{r}, \psi_{r})$. We consider $\psi_{r}$ to be extended to $0$ in $\Omega \setminus (\R\times (0,r))$. This way we have $\psi_{r}\in H^1(\Omega_R \cup (I_R\cup \{0\}) )$.
	
	Then for all $R>1$ we choose a $\mathcal{C}^2(\overline{\Omega})$ function $Y_R:\overline{\Omega}\to [0,1]$ such that
	\begin{align*}
	Y_R(x,y)=1 & \qquad  \text{if} \ |(x,y)|<R-1; \\
	Y_R(x,y)=0 & \qquad  \text{if} \ |(x,y)|\geq R; \\ 
	|\nabla Y_R|^2 \leq C; & \hspace{5em} 
	\end{align*}
	where $C$ is a fixed constant independent of $R$. To simplify the notation later, we call $X_R(x):=Y_R(x,y)|_{y=0}$; we also have that $X_R\in\mathcal{C}^2(\R)$ and $|X_R''|\leq C$. We have that
	\begin{equation*}
	(\phi_{r} X_R, \psi_{r} Y_R) \in \mathcal{H}_R.
	\end{equation*}

	Now we want to show that for a suitable diverging sequence $\{R_n\}_{n\in\N}$ we have
	\begin{equation} \label{ch1Claim}
	Q_{R_n} (\phi_{r} X_{R_n}, \psi_{r} Y_{R_n}) \overset{n\to \infty}{\longrightarrow} \lambda_p( \R \times (0,r)  )).
	\end{equation}

	First, let us show a few useful rearrangements of the integrals that define $Q_R (\phi_{r} X_R, \psi_{r} Y_R)$. We have that
	\begin{align*}
	\int_{I_R}|(\phi_{r} X_R)'|^2 &= \int_{I_R} (\phi_{r} X_R)' \, \phi_{r} \,  X_R ' + \int_{I_R} (\phi_{r} X_R)'   \, \phi_{r} ' \, X_R, \\
	& = \int_{I_R} (\phi_{r} X_R)' \, \phi_{r} \,  X_R '  + \left[ (\phi_{r} X_R^2)  \, \phi_{r} ' \right]_{-R}^R-\int_{I_R}  (\phi_{r} X_R)  \, \left( \phi_{r} '' \, X_R + \phi_{r} ' \,  X_R' \right), \\
	&= \int_{I_R} \phi_{r}^2 \, |X_R '|^2 +  \left[ (\phi_{r} X_R^2)  \, \phi_{r} ' \right]_{-R}^R -\int_{I_R}  \phi_{r} '' \, \phi_{r} \, X_R^2 ,
	\end{align*}
	by having applied integration by parts on the second line and trivial computation in the others.
	Since $X_R(R)= X_R(-R)=0$ and $X_R '$ is supported only in $I_R \setminus I_{R-1}$, we get 
	\begin{equation}\label{ch1eq:parte2}
	\mu D\int_{I_R}|(\phi_{r} X_R)'|^2 = -\mu D\int_{I_R}  \phi_{r} '' \, \phi_{r} \, X_R^2  + \mu D \int_{I_R \setminus I_{R-1}} \phi_{r}^2 \, |X_R '|^2.
	\end{equation}	
	With similar computations we get 
	\begin{equation}\label{ch1eq:parte1}
	\int_{\Omega_R} d|\nabla (\psi_{r} \, Y_R)|^2 = - \int_{\Omega_R} d\Delta \psi_{r} \, \psi_{r} \, Y_R^2 - \int_{I_R} (d\partial_y \psi_{r}) \psi_{r} \, X_R^2 + \int_{\Omega_R \setminus {\Omega_{R-1} }} d|\nabla Y_R|^2 \psi_{r}^2.
	\end{equation}
	Then, we also have 
	\begin{equation} \label{ch1eq:parte3}
	\int_{I_R} (\mu \phi_{r} \, X_R - \nu \psi_{r} \, X_R)^2 = \int_{I_R} \mu \phi_{r} \, X_R^2 (\mu \phi_{r}- \nu \psi_{r}) - \int_{I_R} \nu \psi_{r} \, X_R^2 (\mu \phi_{r}- \nu \psi_{r}).
	\end{equation}
	We now recall that $(\phi_{r}, \psi_{r})$ is an eigenfunction for the problem \eqref{ch1sys:r}.
	Thanks to the third equation of \eqref{ch1sys:r}, the second term in \eqref{ch1eq:parte1} cancel out with the second term in \eqref{ch1eq:parte3}. Moreover we can sum the first term of \eqref{ch1eq:parte2} and the first term of \eqref{ch1eq:parte3} and get
	\begin{equation*}
	-\int_{I_R} \mu D \phi_{r} '' \, \phi_{r} \, X_R^2 + \int_{I_R} \mu \phi_{r} \, X_R^2 (\mu \phi_{r} - \nu \psi_{r}) = \int_{I_R} \mu \lambda_p( \R\times (0, r) ) \phi_{r}^2 \, X_R^2.
	\end{equation*}
	Moreover we have that 
	\begin{equation*}
	- \int_{\Omega_R} d\Delta \psi_{r} \, \psi_{r} \, Y_R^2 - \int_{\Omega_R} f_v(x,0) \psi_{r}^2 \, Y_R^2=  \int_{\Omega_R}  \lambda_p( \R\times (0, r) ) \psi_{r}^2 \, Y_R^2  .
	\end{equation*}
	So, if we call
	\begin{equation*}
	P_R := \frac{ \mu  \int_{I_R \setminus I_{R-1}} D\phi_{r}^2 \, |X_R '|^2 + \nu \int_{\Omega_R \setminus {\Omega_{R-1} }} d|\nabla Y_R|^2 \psi_{r}^2}{\mu \int_{I_R}(\phi_{r} X_R)^2 + \nu \int_{\Omega_R} (\psi_{r} Y_R)^2},
	\end{equation*}
	we have that 
	\begin{equation*}\label{ch10014}
	Q_R (\phi_{r} X_R, \psi_{r} Y_R) = \lambda_p( \R\times (0, r) ) + P_R.
	\end{equation*}
	Proving \eqref{ch1Claim} is equivalent to show that
	\begin{equation} \label{ch11604}
	P_{R_n} \overset{n\to \infty}{\longrightarrow} 0
	\end{equation} 
	for some diverging sequence $\{R_n\}_{n\in \N}$.
	Suppose by the absurd \eqref{ch11604} is not true. 
	First, by the fact that the derivatives of $X_R$ and $Y_R$ are bounded, for some positive constant $C$ we have that 
	\begin{equation*}
	0 \leq P_R \leq C \frac{ \mu  \int_{I_R \setminus I_{R-1}} \phi_{r}^2  + \nu \int_{\Omega_R \setminus {\Omega_{R-1} }} \psi_{r}^2}{\mu \int_{I_R}(\phi_{r} X_R)^2 + \nu \int_{\Omega_R} (\psi_{r} Y_R)^2}
	\end{equation*}
	By the absurd hypothesis,  we have that
	\begin{equation} \label{ch11652}
	\underset{R\to \infty}{\liminf} \, P_R = \xi >0.
	\end{equation}
	Now let us define for all $R\in \N$ the quantity
	\begin{equation*}
	\alpha_R:= \mu \int_{I_R\setminus I_{R-1}}\phi_{r} ^2 + \nu \int_{\Omega_R \setminus \Omega_{R-1}} \psi_{r}^2.
	\end{equation*}
	Since $\phi_r$ and $\psi_r$ are bounded from above, we have that for some constant $k$ depending on $r$, $\mu$, and $\nu$, we have
	\begin{equation}\label{ch1H}
	\alpha_R \leq k R.
	\end{equation}
	For $R\in \N$ one has
	\begin{equation*}
	\mu \int_{I_R}(\phi_{r} X_R)^2 + \nu \int_{\Omega_R} (\psi_{r} Y_R)^2 = \sum_{n=1}^{R-1} \alpha_n + \mu \int_{I_R \setminus I_{R-1}}(\phi_{r} X_R)^2 + \nu \int_{\Omega_R \setminus \Omega_{R-1}} (\psi_{r} Y_R)^2.
	\end{equation*}
	By comparison with \eqref{ch11652}, we have 
	\begin{equation*}
	\underset{R\to \infty}{\liminf} \, \frac{\alpha_R}{\sum_{n=1}^{R-1} \alpha_n} \geq  \underset{R\to \infty}{\liminf} \, \frac{ \alpha_R}{\sum_{n=1}^{R-1} \alpha_n + \mu \int_{I_R \setminus I_{R-1}}(\phi_{r} X_R)^2 + \nu \int_{\Omega_R \setminus \Omega_{R-1}} (\psi_{r} Y_R)^2} \geq \frac{\xi}{C},
	\end{equation*}
	so for $0<\varepsilon< \xi /C$ we have
	\begin{equation}\label{ch1G}
	\alpha_R > \varepsilon \sum_{n=1}^{R-1} \alpha_n
	\end{equation}
	Thanks to \eqref{ch1G} we perform now a chain of inequalities:
	\begin{equation*}
	\alpha_{R+1} > \varepsilon \sum_{n=1}^{R} \alpha_n = \varepsilon \left( \alpha_R + \sum_{n=1}^{R-1} \alpha_n \right) > \varepsilon(1+\varepsilon)\sum_{n=1}^{R-1} \alpha_n > \dots > (1+\varepsilon)^{R+1} \frac{\varepsilon \alpha_1}{(1+\varepsilon)^3} .
	\end{equation*}
	from with we derive that $\alpha_{R}$ diverges as an exponential, in contradiction with the inequality in \eqref{ch1H}. Hence we obtain that 
	\eqref{ch11604} is true, so \eqref{ch1Claim} is also valid.

	By Proposition 4.5 in \cite{romain}, we have that
	\begin{equation}\label{ch11207}
	\lambda_1( \Omega_R) = \underset{  \substack{(u,v)\in \mathcal{H}_R, \\ (u,v)\neq (0,0)}   }{\min} Q_R(u,v).
	\end{equation}
	Hence by \eqref{ch11207} we have that
	\begin{equation*}
	\lambda_1( \Omega_R) \leq Q_R (\phi_{r} X_R, \psi_{r} Y_R).
	\end{equation*}

	Since for all $r>0$ there exist $R>0$ so that \eqref{ch1Claim} holds, we have moreover that
	\begin{equation*}
	\lambda_1( \Omega) \leq \lambda_p( \R\times (0, r) ).
	\end{equation*}
	Then, recalling Proposition \ref{ch1prop:-3}, we get that 
	\begin{equation*}
	\lambda_1( \Omega) \leq \lambda_p( \Omega ).
	\end{equation*}
	Since the reverse inequality was already stated in the first part of this theorem, one has the thesis.
\end{proof}

At last, we prove this proposition of the bounds for $\lambda_p(-\mathcal{L},\Omega)$.

\begin{proposition}\label{ch1prop:Lp_Omega}
	Suppose $f$ satisfies \eqref{ch1hyp:per}.
	We have that
	\begin{equation*}
	\lambda_p(-\mathcal{L}, \R^2)
	\leq \lambda_p(-\mathcal{L}, \Omega) \leq 
	\lambda_1(-\mathcal{L}, \Omega) 
	\end{equation*}
	and if $c=0$ the equality holds.
\end{proposition}

\begin{proof}
	Consider any $r>0$ and take $\lambda_p(-\mathcal{L}, \R\times(-r,r))$ and its eigenfunction $\psi_r$ solving \eqref{ch1sys:bary2}, that is periodic in $x$.
	Then take $\lambda_p(-\mathcal{L}, \R^2)$ and its periodic eigenfunction $\psi_p$, that as we have seen in \eqref{ch1eq:-5} does not depend on $y$, therefore it is limited and has positive infimum, and solves \eqref{ch1sys:L_R2_p}. Then, $\lambda_p(-\mathcal{L}, \R\times(-r,r))$ and $\lambda_p(-\mathcal{L}, \R^2)$ are eigenvalues for the same equation in two domains with one containing the other; hence, one gets that 
	\begin{equation}\label{ch11432}
	\lambda_p(-\mathcal{L}, \R^2)\leq \lambda_p(-\mathcal{L}, \R\times(-r,r)).
	\end{equation}
	By using \eqref{ch1eq:-2}, from \eqref{ch11432} we have
	\begin{equation}\label{ch11429}
	\lambda_p(-\mathcal{L}, \R^2)\leq \lambda_p(-\mathcal{L}, \Omega).
	\end{equation}
	
	Given $R<r$, we can repeat the same argument for $\lambda_1(-\mathcal{L}, B_R)$ and $\lambda_p(-\mathcal{L}, \R\times(-r,r))$  and get
	\begin{equation}\label{ch11433}
	\lambda_p(-\mathcal{L}, \R\times(-r,r)) \leq \lambda_1(-\mathcal{L}, B_R).
	\end{equation}
	By \eqref{ch1eq:-2} and by \eqref{ch1prop:lim_B_R}, we get
	\begin{equation*}
	\lambda_p(-\mathcal{L}, \Omega) \leq \lambda_1(-\mathcal{L}, \Omega).
	\end{equation*}
	This and \eqref{ch11429} give the first statement of the proposition.
	
	If $c=0$, by the second part of Theorem \ref{ch1thm:1.7inbr} we get that $\lambda_p(-\mathcal{L}, \R^2) =\lambda_1(-\mathcal{L}, \Omega)$, hence we have
	\begin{equation*}
	\lambda_p(-\mathcal{L}, \R^2)= \lambda_p(-\mathcal{L}, \Omega) =\lambda_1(-\mathcal{L}, \Omega),
	\end{equation*}
	as wished.
\end{proof}

\subsection{Proof of Theorem \ref{ch1thm:comparison}}

Owing Theorems \ref{ch1thm:char} and \ref{ch1thm:2.6inbhroques} together with the estimates on the eigenvalues proved in the last subsection, we are ready to prove Theorem \ref{ch1thm:comparison}.

\begin{proof}[Proof of Theorem \ref{ch1thm:comparison}]
	By Theorem \ref{ch1thm:1.7inbr}, we have that $\lambda_1(-\mathcal{L}, \R^2)=\lambda_p(-\mathcal{L}, \R^2)$. Then, by Corollary \ref{ch1thm:ineq2}, if $\lambda_1(\Omega)<0$ then
	$\lambda_p(-\mathcal{L}, \R^2)<0$, and if $\lambda_1(\Omega)\geq 0$ then
	$\lambda_p(-\mathcal{L}, \R^2)\geq 0$.
	
	Observe that, when $c=0$, choosing $N=2$ and $L=(\ell, 0)$, the operator $\mathcal{L'}$ defined in \eqref{ch1def:mathcal_L'} coincides with $\mathcal{L}$.
	Then, the affirmations on the asymptotic behaviour of the solutions of the system with and without the road comes from the characterisations in Theorem \ref{ch1thm:char} and \ref{ch1thm:2.6inbhroques}.
	
\end{proof}

\section{Large time behaviour for a periodic medium and $c=0$}\label{ch1s:lb}

We start considering the long time behaviour of the solutions. As already stated in Theorem \ref{ch1thm:char}, the two possibilities for a population evolving through \eqref{ch1sys:fieldroad} are persistence and extinction. We treat these two case in separate sections.

Before starting our analysis, we recall a comparison principle  first appeared in \cite{brr} that is fundamental for treating 
system \eqref{ch1sys:fieldroad}. We recall that a generalised subsolution (respectively, supersolution) is the supremum (resp. infimum) of two subsolutions (resp. supersolutions).

\begin{proposition}[Proposition 3.2 of \cite{brr}]\label{ch1prop:comparison}
	Let $(\underline{u}, \underline{v})$ and $(\overline{u}, \overline{v})$ be respectively a generalised subsolution bounded from
	above and a generalised supersolution bounded from below of \eqref{ch1sys:fieldroad} satisfying $\underline{u} \leq \overline{u}$ and $\underline{v} \leq \overline{v}$
	at $t = 0$. Then, either $\underline{u} \leq \overline{u}$ and $\underline{v} \leq \overline{v}$ for all $t$, or there exists $T > 0$ such that
	$(\underline{u}, \underline{v}) \equiv (\overline{u}, \overline{v})$ for $t\leq T$.
\end{proposition}

The original proof is given for the case of $f$ homogeneous in space; however, it can be adapted with changes so small that we find it useless to repeat it.

Proposition \ref{ch1prop:comparison} gives us important informations on the behaviour at microscopic level. In fact, it asserts that if two pairs of population densities are ``ordered'' at an initial time, then the order is preserved during the evolution according to the equations in \eqref{ch1sys:fieldroad}.

\subsection{Persistence}\label{ch1ss:persistence}

The aim of this section is to prove the second part of Theorem \ref{ch1thm:char}.
First, we are going to show a Liouville type result, that is Theorem \ref{ch1lemma:pbss}, and then we will use that to derive the suited convergence.

We start with some technical lemmas.

\begin{lemma}\label{ch1lemma:converg}
	Let $(u,v)$ be a bounded stationary solution to \eqref{ch1sys:fieldroad} and let $\{(x_n, y_n) \}_{n\in\N}\subset \Omega$ be a sequence of points such that $\{ x_n\}_{n\in\N}$ modulo $\ell$ tends to some $x'\in[0,\ell]$. 
	Then:
	\begin{enumerate}
		\item if $\{ y_n\}_{n\in\N}$ is bounded,
		the sequence of function $\{(u_n, v_n)  \}_{n\in\N}$ defined as
		\begin{equation}\label{ch11648}
		(u_n(x), v_n(x, y))=(u(x+x_n), v(x+x_n, y))
		\end{equation}
		converges  up to a subsequence  to $(\tilde{u}, \tilde{v})$ in $\mathcal{C}_{loc}^2(\R\times\Omega)$ and $(\tilde{u}(x-x'), \tilde{v}(x-x',y)$ is a bounded stationary solution to \eqref{ch1sys:fieldroad}.
		\item if $\{ y_n\}_{n\in\N}$ is unbounded,
		the sequence of function $\{ v_n \}_{n\in\N}$ defined as
		\begin{equation}\label{ch11649}
		v_n(x, y)= v(x+x_n, y+y_n)
		\end{equation}
		converges  up to a subsequence  to $ \tilde{v}$ and $\tilde{v}(x-x', y)$  in $\mathcal{C}_{loc}^2(\R^2)$ is a bounded stationary solution to the second equation in  \eqref{ch1sys:fieldroad} in $\R^2$.
	\end{enumerate}
\end{lemma}

\begin{proof}
	
	Let us call $V=\max\{ \sup u, \sup v  \}$.   
	For all $n\in\N$, there exists $x_n'\in[0,\ell)$ such that $x_n-x_n'\in\ell \Z$.
	
	We start with the case of bounded $\{y_n\}_{n\in\N}$.
	By the periodicity of $f$,
	we have that $(u_n, v_n)$ defined in \eqref{ch11648} is a solution to
	\begin{equation*}
	\left\{
	\begin{array}{lr}
	-D  u '' -c  u' - \nu  v|_{y=0} + \mu u= 0,   & x\in \R,  \\
	v -d \Delta v -c \partial_x v  =f(x+ x_n',v),  & (x, y)\in \Omega, \\
	-d  \partial_y{v}|_{y=0} + \nu v|_{y=0} -\mu u=0, & x\in\R,
	\end{array}
	\right.
	\end{equation*}
	
	Fix $p\geq 1$ and three numbers $j>h>k>0$; we use 
	the notation in \eqref{ch11722} for the sets $I_R$ and $\Omega_R$ for $R= k,\ h, \ j$.
	By Agmon-Douglis-Nirenberg estimates (see for example Theorem 9.11 in \cite{GT}), we have 
	\begin{equation*}
	\norm{u_n}_{W^{2,p}( I_h)} \leq C \left(  \norm{u_n}_{L^p(  I_j)} +\norm{v_n(x,0)}_{L^p(  I_j)}     \right).
	\end{equation*}
	To find the same estimate for the norm of $v_n$, we have to make the same construction used in the proof of Proposition \ref{ch1prop:-3} to find the bound for $\psi_r$. In the same way, we get
	\begin{equation*}
	\begin{split}
	\norm{v_n}_{W^{2,p}(  \Omega_h)} \leq C \Big(  \norm{u_n}_{L^p(  I_j)}   +\norm{v_n}_{L^p(  \Omega_j)} + \norm{f}_{L^p( I_j \times (0, V) )}     \Big) .
	\end{split}
	\end{equation*}
	where the constant $C$, possibly varying in each inequality, depends on $\nu$, $\mu$, $d$, $D$, $h$ and $j$. 
	Using the boundedness of $u$ and $v$, for a possible different $C$ depending on $f$ we get
	\begin{align*}
	\norm{u_n}_{W^{2,p}( I_h)} &\leq C V, \\
	\norm{v_n}_{W^{2,p}(  \Omega_h)} &\leq C V.
	\end{align*}

	Then, we apply the general Sobolev inequalities (see \cite{evans}, Theorem 6 in 5.6) and get for some $\alpha$ depending on $p$, that
	\begin{align*}
	\norm{u_n}_{\mathcal{C}^{\alpha}(  I_h)} & \leq C \norm{u_n}_{W^{2,p}(  I_h)}  \leq CV, \\
	\norm{v_n}_{\mathcal{C}^{\alpha}(  \Omega_h)} &\leq C \norm{v_n}_{W^{2,p}(   \Omega_h)}  \leq CV.
	\end{align*}

	Now we can apply Schauder interior estimates for the oblique boundary condition (see for example Theorem 6.30 in \cite{GT}) and find that 
	\begin{align*}
	\norm{u_n}_{\mathcal{C}^{2,\alpha}(I_k)} &\leq C \left(  \norm{u_n}_{\mathcal{C}^{\alpha}(I_h)} +\norm{v_n(x,0)}_{\mathcal{C}^{\alpha}(I_h)}     \right) \leq CV, \\
	\norm{v_n}_{\mathcal{C}^{2,\alpha}(\Omega_k)} &\leq C \Big(  \norm{u_n}_{\mathcal{C}^{\alpha}(I_h)} 
	+\norm{v_n}_{\mathcal{C}^{\alpha}(\Omega_h)} + \norm{f}_{\mathcal{C}^{\alpha}(I_h \times[0,V])}     \Big) \leq CV.
	\end{align*}
	So the sequences $\{u_n\}_{n\in\N}$ and $\{v_n\}_{n\in\N}$ are bounded locally in space in $C^{2,\alpha}$. By compactness we can extract converging subsequences with limits $\tilde{u}(x)$ and $\tilde{v}(x,y)$. Moreover, since by hypothesis $x_n'\to x'$ as $n\to+\infty$, we have that $(\tilde{u}, \tilde{v})$ is a solution
	\begin{equation*}
	\left\{
	\begin{array}{lr}
	-D  u '' -c  u' - \nu  v|_{y=0} + \mu u= 0,   & x\in \R,  \\
	v -d \Delta v -c \partial_x v  =f(x+ x',v),  & (x, y)\in \Omega, \\
	-d  \partial_y{v}|_{y=0} + \nu v|_{y=0} -\mu u=0, & x\in\R,
	\end{array}
	\right.
	\end{equation*}
	This concludes the proof of the first statement.
	
	Now suppose that $\{  y_n  \}_{n\in\N}$ is unbounded and, up to a subsequence, we can suppose that 
	\begin{equation}\label{ch11827}
	y_n \overset{n\to\infty}{\longrightarrow} +\infty.
	\end{equation}
	Then, 
	the function defined in \eqref{ch11649} solves the equation 
	\begin{equation*}
	-d\Delta v_n -c \partial_x v_n = f(x+x_n', v) \quad \text{for} \ (x,y)\in\R\times(-y_n,0) 
	\end{equation*}
	with the boundary condition $-d\partial_y v_n(x, y_n)  + \nu v_n(x, -y_n)- \mu u(x+x_n)=0$. 
	Fix $p\geq 1$ and three numbers $j>h>k>0$; we denote by $B_R$ the open ball of $\R^2$ centred in $(0,0)$ and with radius $R$, and we will consider $R=j, \ h, \ k$. Notice that by \eqref{ch11827} there exists $N\in\N$ we have that $y_n>j$ for all $n\geq N$.
	Hence, applying the previous estimates to $v_n$ for all $n\geq N$, we find that
	\begin{equation*}
	\begin{split}
	\norm{v_n}_{W^{2,p}(  B_h)} \leq C \Big(  \norm{v_n}_{L^p(  B_j)} + \norm{f}_{L^p( I_j \times (0, V) )}     \Big) \leq CV
	\end{split}
	\end{equation*}
	and then that
	\begin{equation*}
	\norm{v_n}_{\mathcal{C}^{2,\alpha}(B_k)} \leq C \Big(  \norm{v_n}_{\mathcal{C}^{\alpha}(B_h)} + \norm{f}_{\mathcal{C}^{\alpha}(I_h \times[0,V])}     \Big) \leq CV.
	\end{equation*}
	So the sequence $\{v_n\}_{n\in\N}$ is bounded locally in space in $C^{2,\alpha}(\R^2)$ and by compactness we can extract converging subsequence with limit $\tilde{v}(x,y)$, that satisfy
	\begin{equation*}
	-d\Delta v_n -c \partial_x v_n = f(x+x', v) \quad \text{for} \ (x,y)\in\R^2,
	\end{equation*}
	which gives the claim.
	
\end{proof}

The second lemma is similar to the first one, but treats a shifting in time.

\begin{lemma}\label{ch1lemma:mono}
	Let $(u,v)$ be a bounded solution to \eqref{ch1sys:fieldroad} which is monotone in time and let $\{ t_n\}_{n\in\N}\subset \R_{\geq 0}$ be a diverging sequence. Then, the sequence $\{(u_n, v_n)\}_{n\in \N}$ defined by
	\begin{equation}\label{ch11840}
	(u_n(t,x), v_n(t,x, y))=(u(t+t_n,x), v(t+t_n,x, y))
	\end{equation}
	converges in $C_{loc}^{1,2,\alpha}$ to a couple of functions $(\tilde{u}, \tilde{v})$ which is a stationary solution to \eqref{ch1sys:fieldroad}.
\end{lemma}

\begin{proof}
	We call $V=\max \{ \sup u, \sup v   \}$.
	For every fixed $x\in \R$  we have that $u_n(t,x)$ is an monotone bounded sequence. Then, we can define a function $\tilde{u}(x)$ as 
	\begin{equation}\label{ch11720}
	\tilde{u}(x) = \underset{n\to\infty}{\lim} {u_n(t,x)}
	\end{equation}
	and $0\leq \tilde{u}(x)\leq U$. 
	Analogously, for all $(x,y)\in \Omega$  we can define
	\begin{equation}\label{ch11721}
	\tilde{v}(x,y) = \underset{n\to\infty}{\lim} {v_n(t,x,y)}
	\end{equation}
	and $0 \leq \tilde{v}(x,y)\leq V$. 
	
	Fix $p\geq 1$, $T>0$ and three numbers $k<h<j$; we use 
	the notation in \eqref{ch11722} for the sets $I_R$ and $\Omega_R$ for $R= k,\ h, \ j$. 
	For $S$ an open subset in $\R^N$, in this proof we denote the space of function with one weak derivative in time and two weak derivatives in space by $W_p^{1,2}(S)$. 
	By Agmon-Douglis-Nirenberg estimates we have 
	\begin{equation*}
	\norm{u_n}_{W^{1,2}_p(  I_h)} \leq C \left(  \norm{u_n}_{L^p((0,T)\times I_j)} +\norm{v_n(t,x,0)}_{L^p((0,T)\times I_j)}     \right) \leq CV.
	\end{equation*}
	To find the same estimate for the norm of $v_n$, we have to make the same construction used in the proof of Proposition \ref{ch1prop:-3} to find the bound for $\psi_r$. In the same way, we get
	\begin{equation*}
	\begin{split}
	\norm{v_n}_{W^{1,2}_p((0,T)\times \Omega_h)} \leq C \Big(  \norm{u_n}_{L^p((0,T)\times I_j)} 
	+\norm{v_n}_{L^p((0,T)\times \Omega_j)} + \norm{f}_{L^p( I_j \times (0, V) )}     \Big) \leq CV.
	\end{split}
	\end{equation*}
	where the constant $C$, possibly varying in each inequality, depends on $\nu$, $\mu$, $d$, $D$, $T$, $h$ and $j$. Then, we apply the general Sobolev inequalities (see \cite{evans}, Theorem 6 in 5.6) and get for some $\alpha$ depending on $p$, that
	\begin{align*}
	\norm{u_n}_{\mathcal{C}^{\alpha}((0,T)\times I_h)} & \leq C \norm{u_n}_{W^{1,2}_p((0,T)\times I_h)}  \leq CV, \\
	\norm{v_n}_{\mathcal{C}^{\alpha}((0,T)\times \Omega_h)} &\leq C \norm{v_n}_{W^{1,2}_p((0,T)\times  \Omega_h)}  \leq CV.
	\end{align*}
	Moreover, since for $n\in \N$ the functions $u_n$ and $v_n$ are just time translation of the same functions $\tilde{u}$ and $\tilde{v}$, we also have that
	\begin{align*}
	\norm{u_n}_{\mathcal{C}^{\alpha}((0,+\infty)\times I_h)} &\leq CV, \\
	\norm{v_n}_{\mathcal{C}^{\alpha}((0,+\infty)\times  \Omega_h)} &\leq CV.
	\end{align*} 
	Now we can apply Schauder interior estimates (see for example Theorem 10.1 in Chapter IV of \cite{ladyzhenskaya}) and find that 
	\begin{align*}
	\norm{u_n}_{\mathcal{C}^{1,2,\alpha}((0,+\infty)\times I_k)} &\leq C \left(  \norm{u_n}_{\mathcal{C}^{\alpha}((0,+\infty)\times I_h)} +\norm{v_n(t,x,0)}_{\mathcal{C}^{\alpha}((0,+\infty)\times I_h)}     \right) \leq CV, \\
	\norm{v_n}_{\mathcal{C}^{1,2,\alpha}((0,+\infty)\times \Omega_k)} &\leq C \Big(  \norm{u_n}_{\mathcal{C}^{\alpha}((0,+\infty)\times I_h)} + \\ &\hspace{5em} 
	+\norm{v_n}_{\mathcal{C}^{\alpha}((0,+\infty)\times \Omega_h)} + \norm{f}_{\mathcal{C}^{\alpha}(I_h \times[0,V])}     \Big) \leq CV.
	\end{align*}
	So the sequences $\{u_n\}_{n\in\N}$ and $\{v_n\}_{n\in\N}$ are bounded locally in space in $C^{1,2,\alpha}$. By compactness we can extract converging subsequences with limits $q(t,x)$ and $p(t,x,y)$ that satisfy system \eqref{ch1sys:fieldroad}. But as said in \eqref{ch11720} and \eqref{ch11721} the sequences
	$\{u_n\}$ and $\{v_n\}$ also converge punctually to $\tilde{u}$ and $\tilde{v}$, that are stationary functions. Then, the couple $(\tilde{u}, \tilde{v})$ is a positive bounded stationary solution of system \eqref{ch1sys:fieldroad}.
\end{proof}

The following lemma gives essentials information on the stationary solutions, on which the uniqueness result of Theorem \ref{ch1lemma:pbss} will rely on.

\begin{lemma}\label{ch1lemma:pos_inf}
	Suppose that $c=0$, $f$ satisfies \eqref{ch1hyp:0}-\eqref{ch1hyp:per} and that
	$\lambda_1( \Omega)<0$. Then, every stationary bounded solution $(u, v)\not\equiv(0,0)$ of system \eqref{ch1sys:fieldroad} respects
	\begin{equation*}
	\underset{\R}{\inf} \, u >0, \quad \underset{\Omega}{\inf} v>0.
	\end{equation*}
\end{lemma}

\begin{proof}
	
	\emph{Step 1: sliding in $x$.}
	If $\lambda_1( \Omega)<0$, thanks to Proposition \ref{ch1prop:romain} there exists $R>0$ such that $\lambda_1( \Omega_R)<0$. 
	Since $\lambda_1( \Omega_R)$ is monotonically decreasing in $R$, 
	we can suppose that $R> \ell$.
	By a slight abuse of notation, let us call $(u_R,v_R)$ the eigenfunctions associated with $\lambda_1( \Omega_R)<0$ extended to 0 in $\R \times \Omega \setminus (I_R\times \Omega_R)$.
	
	We claim that there exists $\varepsilon>0$ such that $\varepsilon(u_R,v_R)$ is a subsolution for system \eqref{ch1sys:fieldroad}.
	In fact, we have that
	\begin{equation*}
	\underset{v\to 0^+}{\lim} \dfrac{f(x,v)}{v} = f_v(x,0),
	\end{equation*}
	so for $\varepsilon$ small enough we have that
	\begin{equation}\label{ch12220}
	\dfrac{f(x,\varepsilon v_R)}{\varepsilon v_R} > f_v(x,0) + \lambda_1( \Omega_R).
	\end{equation}
	Then, 
	\begin{equation}\label{ch12221}
	\left\{
	\begin{array}{ll}
	-D \varepsilon u_R '' -c \varepsilon  u_R' - \nu \varepsilon v_R|_{y=0} +\mu \varepsilon u_R=\lambda_1( \Omega_R ) \varepsilon u_R \leq 0,   & x\in I_R,  \\
	-d \Delta \varepsilon v_R -c \partial_x \varepsilon v_R  =(f_v(x,0)+ \lambda_1( \Omega_R) \varepsilon v_R \leq f(x, \varepsilon v_R),  & (x, y)\in \Omega_R, \\
	-d  \varepsilon\partial_y{v_R}|_{y=0} + \nu \varepsilon v_R|_{y=0} -\varepsilon u_R=0, & x\in I_R,
	\end{array}
	\right.
	\end{equation}
	so $\varepsilon(u_R,v_R)$ is a subsolution. 
	
	Decreasing $\varepsilon$ if necessary, we have that $\varepsilon(u_R,v_R)<(u,v)$ because $u$ and $v$ are strictly positive in all points of the domain while $(u_R, v_R)$ has compact support. Now we translate $\varepsilon(u_R,v_R)$ in the variable $x$ by multiples of $\ell$; given $k\in \Z$, we call
	\begin{align*}
	u_{R, k}(x):= \varepsilon u_R(x-k\ell), \quad  & I_{R,k}=(k\ell-R,k\ell+R), \\
	v_{R, k}(x,y):=\varepsilon v_R(x-k\ell,y), \quad  & \Omega_{R,k}= B_R(k\ell, 0) \cap \Omega.
	\end{align*}
	The couple $(u_{R,k}, v_{R,k})$ is still a subsolution to system \eqref{ch1sys:fieldroad} because is a translation of a subsolution by multiple of the periodicity of the coefficients in the equations.
	Suppose by the absurd that there exists $k\in \Z$ such that $(u_{R,k}, v_{R,k})\not < (u,v)$. 
	Since $u$ and $v$ are strictly positive in all points of respectively $\R$ and $\Omega$, while $u_{R,k}$ and $v_{R,k}$ have compact support, by decreasing $\varepsilon$ if necessary, we have that $(u_{R,k}, v_{R,k})\leq (u,v)$ and either  there exists $\bar{x}\in I_{R,k}$ such that $ u_{R,k}(\bar{x})=u(\bar{x})$ or there exists $(\bar{x}, \bar{y})\in \overline{\Omega}_{R,k} $ such that $v_{R,k}(\bar{x}, \bar{y})=v(\bar{x}, \bar{y})$. 
	Then, by the Comparison Principle, we have that $(u_{R,k}, v_{R,k})\equiv (u,v)$, which is absurd because $u_{R,k}$ and $v_{R,k}$ are compactly supported.
	Therefore, we have
	\begin{equation}\label{ch11811}
	\begin{array}{rl}
	u(x) > \varepsilon u_R(x+k\ell), &\quad \forall x\in\R, \ \forall k\in\Z, \\
	v(x,y) > \varepsilon v_R(x+k\ell,y), &\quad \forall (x,y)\in\Omega,  \ \forall k\in\Z.
	\end{array}
	\end{equation}
	Fix $Y<\sqrt{R^2-\ell^2}$. Then, let us call
	\begin{equation*}
	\delta_Y := \min\{ \underset{[0, \ell]}{\min} \, \varepsilon u_R(x),  \underset{[0,\ell]\times[0,Y] }{\min} \varepsilon v_R(x,y) \}. 
	\end{equation*}
	Since $[0,\ell]\times(0,Y) \subset \Omega_R$ and $[0,\ell]\subset I_R$,  we have that $\delta_Y>0$.
	Then,  \eqref{ch11811} implies that 
	\begin{equation}\label{ch11749}
	\begin{array}{ll}
	u(x)>\delta_Y, & \text{for} \ x\in \R, \\
	v(x,y)>\delta_Y, & \text{for} \ x\in \R, \ y\in[0,Y].
	\end{array}
	\end{equation}

	\emph{Step 2: sliding in $y$.}
	Recall that by Corollary \ref{ch1thm:ineq2} we have $\lambda_1( \Omega)<0$ implies $\lambda_1(-\mathcal{L},\Omega) <0$ and by Proposition \eqref{ch1prop:-3} it holds  $\lambda_p(-\mathcal{L}, \R^2)\leq \lambda_1(-\mathcal{L},\Omega)<0$. By Proposition \ref{ch1prop:Lp_Omega}
	and by \eqref{ch1eq:-2} we have that for some $r>0$ we have $\lambda_p(-\mathcal{L},\R\times (-r,r))<0$. Then, let us call $v_r$ the eigenfunction related to $\lambda_p(-\mathcal{L},\R\times  (-r,r))$ extended to 0 outside its support; repeating the classic argument, one has that for some $\theta>0$ we have $\theta v_r$ extended to 0 outside $\R\times  (-r,r)$ is a subsolution for the second equation in system \eqref{ch1sys:fieldroad}.
	For all $h>0$, let us now call $\varphi_h(x,y):=v_r(x, y+h)$.
	Since $v_r$ is periodic in the variable $x$, we have that $v_r$ is uniformly bounded. 
	Now take $Y>2r$ and $h_0>r$ such that $Y>h_0+r$; by decreasing $\theta$ if necessary, we get that $\theta v_r < \delta_Y$. 
	Hence, we get 
	\begin{equation}\label{ch11756}
	\theta \varphi_{h_0}(x,y) < v(x,y) \quad \text{for} \ x\in\R, \ y\geq 0.
	\end{equation}
	
	Now define
	\begin{equation*}
	h^*= \sup \{ h\geq h_0 \ : \  \theta \varphi_h(x,y) < v(x,y) \ \text{for} \ x\in\R, \ y\in[h-r, h+r] \}.
	\end{equation*}
	By \eqref{ch11756}, we get that $h^* \geq h_0>r$.
	
	We now take $\tilde{y}<h^*+r$ and define
	\begin{equation*}
	\tilde{h} = \left\{
	\begin{array}{ll}
	\tilde{y}, & \text{if} \ \tilde{y}\leq  h^*, \\
	\dfrac{\tilde{y}+h^*-r}{2}, & \text{if} \ h^* < \tilde{y} < h^*+r. 
	\end{array}
	\right.
	\end{equation*}
	Then, $\tilde{h}<h^*$: if $\tilde{h}=\tilde{y}$ it is trivial, otherwise one observe that $\tilde{y}-r < h^*$.
	Also, $\tilde{y}\in(\tilde{h}-r, \tilde{h}+r )$; in fact, that is obvious if $\tilde{h}=\tilde{y}$, otherwise we have that $\tilde{y}< h^*+r$ and
	\begin{equation*}
	\tilde{h}-r < h^*-r < \tilde{y}  <  \dfrac{\tilde{y}+h^*+r}{2} 
	= \tilde{h}+r.
	\end{equation*}
	Then, since $v_r$ and therefore $\varphi_{\tilde{h}}$ are periodic in $x$, we have that
	\begin{equation} \label{ch11606}
	v(x, \tilde{y}) > \theta \varphi_{\tilde{h}}(x, \tilde{y}) >  \underset{[0,\ell]}{\min} \, \theta \varphi_{\tilde{h}}(x, \tilde{y})>0,
	\end{equation}
	so $v(x,y)>0$ for all $y<h^*+r$, $x\in\R$ and moreover
	\begin{equation}\label{ch11759}
	v(x,y) > \theta \, \underset{[0, \ell]}{\min} \, v_r(x, 0) >0 \quad \text{for} \ x\in\R, \ y\leq h^*.
	\end{equation}
	
	Suppose by absurd that $h^*<+\infty$. Then there exists a sequence $\{h_n\}_{n\in \N}$ and a sequence $\{(x_n, y_n)\}_{n\in\N}$ with $(x_n, y_n)\in \R\times[h_n-r, h_n+r]$, such that $$\underset{n\to+\infty}{\lim} h_n=h^* \quad \text{and} \quad \underset{n\to+\infty}{\lim}   {\theta\varphi_{h_n}(x_n, y_n)-v(x_n, y_n) }=0.$$
	Up to a subsequence, $\{ y_n\}_{n\in\N} \subset [0, h^*+r] $ converges to some $\bar{y}\in[h^*-r, h^*+r]$ while  $\{ x_n\}_{n\in\N}$ either converges to some $\bar{x}\in \R$ or goes to infinity.

	For all $n\in \N$ there exists $x_n'\in[0,\ell)$ and $k\in\Z$ such that
	\begin{equation}\label{ch12040}
	x_n= x_n' + k\ell.
	\end{equation}
	Up to a subsequence,
	\begin{equation}\label{ch11744}
	x_n' \overset{n\to\infty}{\longrightarrow} x'\in[0,\ell].
	\end{equation}
	Define 
	\begin{equation*}
	(u_n(x), v_n(x, y)):=(u(x+x_n), v(x+x_n, y)).
	\end{equation*}
	Then, by Lemma \ref{ch1lemma:converg} we have that $\{(u_n,v_n)\}_{n\in\N}$ converges to some $(\tilde{u}, \tilde{v})$ such that 
	\begin{equation}\label{ch11959}
	\mbox{$(\tilde{u}(x+x'), \tilde{v}(x+x', y)$ is a bounded stationary solution to \eqref{ch1sys:fieldroad}.}
	\end{equation}

	By \eqref{ch11606}, we have
	\begin{equation}\label{ch11955}
	\tilde{v}(x,\tilde{y}) > \underset{x\in[0,\ell]}{\min} \theta \varphi_{\tilde{h}}(x, \tilde{y})>0 \quad \text{for} \ \tilde{y}<h^*+r.
	\end{equation}  
	
	We notice that if $\tilde{v}(0, \bar{y})=0$, since $\tilde{v}\geq 0$ and \eqref{ch11959} holds,
	by the maximum principle  we get $\tilde{v}\equiv 0$ in $\Omega$. But since \eqref{ch11955} holds, this is not possible and instead $\tilde{v}(0, \bar{y})>0$. Hence, $0<\tilde{v}(0, \bar{y})= \theta \varphi_{h^*}(0, \bar{y})$, so
	\begin{equation}\label{ch12050}
	\bar{y}\neq h^*\pm r.
	\end{equation}

	We have that $\theta \varphi_{h_n}$ is a subsolution for $ \mathcal{L}$ in $\R\times (h_n-r, h_n+r)$, since it is a translation of a subsolution. Moreover, thanks to the periodicity of $\varphi_{h_n}$ and the definition of $x_n'$ in \eqref{ch12040}, we have 
	\begin{equation*}
	\varphi_{h_n}(x+x_n, y)=\varphi_{h_n}(x+x_n', y).
	\end{equation*}
	It follows that the sequence $\varphi_{h_n}(x+x_n, y)$ converges to $\varphi_{h^*}(x+{x}', y)$. 
	Then, $\theta \varphi_{h^*}(x+{x}', y)$ is a subsolution for the second equation in \eqref{ch1sys:fieldroad} in $\R\times(h^*-r, h^*+r)$ and  by \eqref{ch12050} it holds $(0, \bar{y})\in \R\times(h^*-r, h^*+r)\subset \Omega$. Hence, we can apply the comparison principle to $\tilde{v}(x,y)$ and $\theta \varphi_{h^*}(x+{x}', y)$: since $\tilde{v}(0,\bar{y})=\theta\varphi_{h^*}({x}', \bar{y})$, then $\tilde{v}(x,y)\equiv\theta \varphi_{h^*}(x+{x}', y) $ in all the points of $\R\times(h^*-r, h^*+r)$. But then by continuity $\tilde{v}(x, h^*-r)=\theta \varphi_{h^*}(x+{x}', h^*-r)=0 $, which is absurd for \eqref{ch11955}. 
	Hence $h^*=+\infty$. From that and \eqref{ch11749} we have statement of the lemma.
\end{proof}

Finally, we are ready to prove existence and uniqueness of a positive bounded stationary solution to \eqref{ch1sys:fieldroad}. The existence of such couple of function is crucial to get the persistence result of Theorem \ref{ch1thm:char}.

\begin{theorem}\label{ch1lemma:pbss}	
	Suppose that $c=0$, $f$ satisfies \eqref{ch1hyp:0}-\eqref{ch1hyp:per} and that
	$\lambda_1( \Omega)<0$. Then, the following holds:
	\begin{enumerate}
		\item There exists a unique positive bounded stationary solution $(u_{\infty}, v_{\infty})$ to system \eqref{ch1sys:fieldroad}.
		\item The functions $u_{\infty}$ and $v_{\infty}$ are periodic in the variable $x$ of period $\ell$.
	\end{enumerate}
\end{theorem}

\begin{proof}
	\emph{Step 1: construction of a subsolution.}
	
	Since $\lambda_1( \Omega)<0$, by Theorem \ref{ch1thm:eq_1_p} it holds that $\lambda_p( \Omega)<0$ and moreover by Proposition \ref{ch1eq:-3}  
	there exists $r>1$ such that $\lambda_p( \R\times(0, r))<0$. Let us call $(\phi_r, \psi_r)$ the eigenfunction related to $\lambda_p( \R\times (0,r))$.
	
	We have that
	\begin{equation*}
	\underset{v\to 0^+}{\lim} \dfrac{f(x,v)}{v} = f_v(x,0),
	\end{equation*}
	so there exists $\varepsilon>0$ such that
	\begin{equation*}
	\dfrac{f(x,\varepsilon\psi_r)}{\varepsilon\psi_r} > f_v(x,0) + \lambda_p( \R\times (0,r)).
	\end{equation*}
	Then, 
	\begin{equation}\label{ch11933}
	\left\{
	\begin{array}{ll}
	-D \varepsilon \phi_r'' -c \varepsilon\phi_r' - \nu \varepsilon \psi_r\rvert_{y=0} + \varepsilon \phi_r = \lambda_p ( \R\times(0,r) ) \varepsilon\phi_r < 0,  & x\in\R,  \\
	-d \Delta \varepsilon\psi_r -c \partial_x \varepsilon\psi_r   < f(x, \varepsilon \psi_r),  & (x, y)\in \R\times (0, r), \\
	-d  \varepsilon\partial_y{\psi_r}|_{y=0} + \nu \varepsilon\psi_r|_{y=0} -\varepsilon\phi_r=0, & x\in\R,
	\end{array}
	\right.
	\end{equation}
	so $\varepsilon(\phi_r, \psi_r)$ is a subsolution to system \eqref{ch1sys:fieldroad}.

	%
	%
	%
	
	Thanks to Corollary \ref{ch1thm:ineq2}, $\lambda_1( \Omega)<0$ implies $\lambda_1(-\mathcal{L}, \R^2)<0$; then  Proposition \ref{ch1prop:-3} implies that  $\lambda_p(-\mathcal{L}, \R^2)<0$. By \eqref{ch1eq:-5}, also $\lambda_p(-\mathcal{L}, \R)<0$.
	
	Consider the periodic positive eigenfunction $\psi_p(x)$ related to $\lambda_p(-\mathcal{L}, \R)$. 
	With a slight abuse of notation, we can extend $\psi_p(x)$ in all $\R^2$ by considering constant with respect to the variable $y$.
	Repeating the same arguments as before, we can prove that for some $\theta$ the function
	$\theta \psi_p(x)$ is a subsolution for the second equation of system \eqref{ch1sys:fieldroad} in $\R^2$.
	
	Consider $\delta>0$.
	We have that $\psi_p(x)$ is limited, therefore there exists $\varepsilon'\in(0, \theta)$ such that 
	\begin{equation}\label{ch11653}
	\underset{[0,\ell]}{\max} \, \varepsilon' \psi_p(x) <\delta < \underset{\substack{[0,\ell]\times[0,r-1]}}{\min} \varepsilon\psi_r(x,y)  .
	\end{equation}
	Then,  let us define the functions
	\begin{align*}
	\underline{u}(x) & :=  \varepsilon\phi_r(x), \\
	\underline{v}(x,y) &:= \max \{  \varepsilon\psi_r(x,y)  , \, \varepsilon' \psi_p(x)\}.
	\end{align*}
	
	By \eqref{ch11653}, for $y\in(0, r-1)$ it holds that $\underline{v}(x,y)=\varepsilon\psi_r(x,y)$. Hence, we get that $(\underline{u}, \underline{v})$ is a subsolution for the first and third equation of \eqref{ch1sys:fieldroad}. Moreover, since $\varepsilon\psi_r(x,y)$ and $\varepsilon' \psi_p(x)$ are both subsolution to the second equation to \eqref{ch1sys:fieldroad}, so  the maximum between them is a generalised subsolution. Thanks to that, we can conclude that $(\underline{u}, \underline{v})$ is a generalised subsolution for the system \eqref{ch1sys:fieldroad}.
	
	Since $\phi_r$ and $\psi_p$ are periodic in $x$ and independent of $y$, we get $$\underset{\R}{\inf} \, \underline{u}(x)>0 \quad \text{and} \quad  \underset{\Omega}{\inf} \, \underline{v}(x,y) >0.$$ 
	
	So, $(\underline{u}, \underline{v})$ is a generalised subsolution for the system \eqref{ch1sys:fieldroad}, with positive infimum, and by the periodicity of $\phi_r$, $\psi_r$ and $\psi_p$, it is periodic in $x$ with period $\ell$. 
	
	\emph{Step 2: construction of a stationary solution.}
	
	Take the generalised subsolution $(\underline{u}, \underline{v})$.	We want to show that the solution $(\tilde{u}(t,x), \tilde{v}(t,x,y))$ having $(\underline{u}(x), \underline{v}(x,y))$ as initial datum is increasing in time and converge to a stationary solution.

	By the fact that $(\underline{u}, \underline{v})$ is a subsolution, at we have $(\underline{u}, \underline{v}) \leq (\tilde{u}, \tilde{v})$ for all $t\geq 0$. Hence, for all $\tau>0$,
	let us consider the solution $(z,w)$ stating at $t=\tau$ from the initial datum $(\underline{u}(x), \underline{v}(x,y))$. Then, at $t=\tau$ we have that $(\tilde{u}(\tau,x), \tilde{v}(\tau,x,y))\geq (z(\tau,x), w(\tau,x,y))$. By the comparison principle \ref{ch1prop:comparison}, we have that for all $t\geq \tau$ it holds that $(\tilde{u}(t,x), \tilde{v}(t,x,y))\geq (z(t,x), w(t,x,y))$. By the arbitrariness of $\tau$, we get that $(\tilde{u}(t,x), \tilde{v}(t,x,y))$ is increasing in time.
	
	Moreover, consider 
	\begin{equation*}
	V:= \max \left\{ M, \sup \underline{v}, \frac{\mu}{\nu} \sup \underline{u} \right\}, \quad U:= \frac{\nu}{\mu} V,
	\end{equation*}
	where $M>0$ is the threshold value defined in \eqref{ch1hyp:M}. 
	One immediatly checks that $(U,V)$ is a supersolution for the system \eqref{ch1sys:fieldroad}.
	Also, we have that $(\underline{u}(x), \underline{v}(x,y)) \leq (U,V)$, so by the comparison principle \ref{ch1prop:comparison} it holds that 
	\begin{equation*}
	(\tilde{u}(t,x), \tilde{v}(t,x,y)) \leq (U,V) \quad \text{for all} \ t>0.
	\end{equation*}
	Hence, $(\tilde{u}(t,x), \tilde{v}(t,x,y))$ is limited. 
	
	Now consider an increasing diverging sequence $\{t_n\}_{n\in N}\subset \R^+$. Then, define
	\begin{equation*}
	u_n(t,x):=\tilde{u}(t+t_n, x), \quad v_n(t,x,y):=\tilde{v}(t+t_n, x, y),
	\end{equation*}
	that is a sequence of functions. By Lemma \ref{ch1lemma:mono},  $(u_n, v_n)$ converge in $\mathcal{C}_{loc}^{1,2,\alpha}$ to a stationary bounded solution to \eqref{ch1sys:fieldroad}, that we call $(u_{\infty}, v_{\infty})$. We point out that $(u_{\infty}, v_{\infty})\not\equiv(0,0)$ since 
	\begin{equation*}
	(u_{\infty}, v_{\infty}) \geq (\underline{u}, \underline{v}) > (0,0).
	\end{equation*}
	Moreover, both functions are periodic of period $\ell$ in the variable $x$ since the initial datum is.
	
	\emph{Step 3: uniqueness.}
	
	Suppose that there exists another positive bounded stationary solution $(q,p)$ to \eqref{ch1sys:fieldroad}. Then, define
	\begin{equation*}
	k^* := \sup \left\{  k>0 \ : \ u_{\infty}(x) > k q(x) \ \forall x\in\R, \ v_{\infty}(x,y)> kp(x,y) \ \forall(x,y)\in\Omega     \right\}.
	\end{equation*}
	Since by Lemma \ref{ch1lemma:pos_inf} the functions $u_{\infty}$ and $v_{\infty}$ have positive infimum and since $p$ and $q$ are bounded, we have that $k^*>0$.
	
	We claim that 
	\begin{equation}\label{ch12337}
	k^* \geq 1.
	\end{equation}	
	
	By the definition of $k^*$, one of the following must hold: there exists
	\begin{equation}\label{ch1caso1}
	\begin{split}
	\mbox{ either a sequence $\{x_n\}_{n\in\N}\subset \R$ such that $u_{\infty}(x_n) - k^* q(x_n) \overset{n\to\infty}{\longrightarrow} 0$,}
	\end{split}	
	\end{equation}
	\begin{equation}\label{ch1caso2}
	\mbox{or a sequence $\{(x_n, y_n)\}_{n\in\N}\subset \Omega$ such that $v_{\infty}(x_n,y_n) - k^* p(x_n, y_n)  \overset{n\to\infty}{\longrightarrow} 0$.}
	\end{equation}
	There exists a sequence $\{ {x}_n' \}_{n\in\N} \subset[0,\ell)$ such that 
	\begin{equation}\label{ch11807}
	x_n-{x}_n' \in \ell \Z \quad \text{for all} \ n\in\N.
	\end{equation}
	Then, up to extraction of a converging subsequence, we have that there exists $ x'\in \R$ such that ${x}_n' \overset{n\to \infty}{\longrightarrow}  x'$.
	One can see that the sequence of couples  
	\begin{equation*}
	(q_n(x), p_n(x,y) ):= (q(x+x_n), p(x+x_n, y) )
	\end{equation*}
	is a stationary solution for \eqref{ch1sys:fieldroad} with reaction function $f(x+ {x}_n', v)$. By Lemma \ref{ch1lemma:converg}, up to a subsequence, $(q_n, p_n)$ converges in $\mathcal{C}_{loc}^{2}$  to some $(q_{\infty}, p_{\infty})$, which is a stationary solution of  \eqref{ch1sys:fieldroad} with reaction function $f(x+  x', v)$. 
	We also notice that, thanks to the periodicity of $u_{\infty}$ and $v_{\infty}$, $(u_{\infty}(x+ x'), v_{\infty}(x+ x', y))$ is also a stationary solution of  \eqref{ch1sys:fieldroad} with reaction function $f(x+  x', v)$. 
	Define the function 
	\begin{equation}
	\begin{split}
	\alpha(x)&:=u_{\infty}(x+ x') - k^*q_{\infty}(x), \\
	\beta(x)&:= v_{\infty}(x+ x', y) - k^*p_{\infty}(x,y),
	\end{split}
	\end{equation}
	and notice that $\alpha (x)\geq 0$, $\beta(x, y)\geq 0$.
	
	Now suppose that \eqref{ch1caso1} holds.
	We have that
	\begin{equation*}
	\alpha(0)=u_{\infty}( x') - k^*q_{\infty}(0)=0.
	\end{equation*}
	Moreover, $\alpha(x)$ is a solution to the equation
	\begin{equation*}
	-D \alpha''-c\alpha' -\nu \beta|_{y=0} +\mu \alpha =0.
	\end{equation*}
	By the maximum principle, we have that since $\alpha(x)$ attains its minimum in the interior of the domain then $\alpha(x)\equiv \min \alpha =0$. Then, one would have $u_{\infty}(x+ x')\equiv k^* q_{\infty}(x)$ and by the comparison principle \ref{ch1prop:comparison} we have $v_{\infty}(x+ x', y)\equiv k^*p_{\infty}(x+ x', y)$. 
	Subtracting the second equation of system \eqref{ch1sys:fieldroad} for $p_{\infty}$ from the one for $v_{\infty} $ we get
	\begin{equation}\label{ch11948}
	0=f(x+ x', v_{\infty}(x+ x',y))-k^* f(x+ x',p_{\infty}(x,y)).
	\end{equation}
	If by the absurd $k^*<1$, by the KPP hypothesis \eqref{ch1hyp:KPP} we have $k^*f(x+  x', p_{\infty}(x,y)) <f(x+  x', k^* p_{\infty}(x,y)) = f(x+ x', v_{\infty}(x+ x', y))$ and the right hand side of \eqref{ch11948} has a sign, that is absurd since the left hand side is 0. We can conclude that if we are in the case of \eqref{ch1caso1}, then \eqref{ch12337} holds.

	Suppose instead that \eqref{ch1caso2} is true. If  $\{ y_n\}_{n\in\N}$  is bounded, we define $y_n \overset{n\to \N}{\longrightarrow} y'\in\R$.
	Then,
	\begin{equation}\label{ch12007}
	\beta(0, y')= v_{\infty}( x', y')- k^* p_{\infty}(0,y' )=0.
	\end{equation}
	If by the absurd $k^*<1$, then by the  Fisher-KPP hypothesis \eqref{ch1hyp:KPP}  we have	
	\begin{equation}\label{ch12001}
	\begin{split}
	-d \Delta \beta(x,y) -c \partial_x \beta(x,y) &= f(x+ x', v_{\infty}(x+ x',y))- k^* f(x+ x', p_{\infty}(x,y)) \\
	&> f(x+ x', v_{\infty}(x+ x',y))- f(x, k^*p_{\infty}(x,y)).
	\end{split}
	\end{equation}
	Since $f$ is locally Lipschitz continuous in the second variable, one infers from \eqref{ch12001} that there exists a bounded function $b(x)$ such that
	\begin{equation}\label{ch12009}
	-d \Delta \beta -c \partial_x \beta + b \beta >0.
	\end{equation}
	Since that, $\beta \geq 0$ and by \eqref{ch12007} $\beta(0, y')=0$, if $y'>0$ we apply the strong maximum principle and we have $\beta \equiv 0$.  
	If $y'=0$, we point out that by the fact that $v_{\infty}$ and $p_{\infty}$ are solution to \eqref{ch1sys:fieldroad} it holds 
	\begin{equation*}
	d \partial_y \beta(x,0) = \nu (v_{\infty}(x+ x',0)  - k^* p_{\infty}(x,0) )  - \nu (u_{\infty}(x+ x')- k^* q_{\infty}(x)) \leq  0
	\end{equation*}
	By that, the inequality in \eqref{ch12009},  $\beta \geq 0$, $\beta(0, y')=0$,  we can apply Hopf's lemma and get again $\beta \equiv 0$.
	Then for both $y'>0$ and $y'=0$, $v_{\infty}(x+ x', y)\equiv k^*p_{\infty}(x+ x', y)$ and \eqref{ch11948} holds, but we have already saw that this is absurd. So, in the case of \eqref{ch1caso2}, if $\{ y_n\}_{n\in\N}$ is bounded, \eqref{ch12337} is true.
	
	At last, if $\{ y_n\}_{n\in\N}$ is unbounded, we define
	\begin{align*}
	V_n(x,y)&:=v_{\infty}(x+x_n, y+y_n), \\
	P_n(x,y)&:=p(x+x_n, y+y_n).
	\end{align*}
	By Lemma \ref{ch1lemma:mono}, up to subsequences, $V_n$ and $P_n$ converge in $\mathcal{C}_{loc}^{2}$  to some functions $V_{\infty}$ and $P_{\infty}$ solving
	\begin{equation*}
	- d \Delta v - c\partial_x v = f(x+ x', v) \quad \text{for} \ (x,y)\in\R^2 .
	\end{equation*}
	Moreover, if we suppose $k^*< 1$, by the Fisher-KPP hypothesis \eqref{ch1hyp:KPP} we have that
	\begin{equation*}
	k^*f(x+ x', P_{\infty})< f(x+ x', k^* P_{\infty})
	\end{equation*}
	and consequently, calling $\gamma:=V_{\infty} - k^* P_{\infty}$, we get
	\begin{equation*}
	- d \Delta\gamma - c\partial_x \gamma > f(x+ x', V_{\infty}) -  f(x+ x', k^*P_{\infty}) .
	\end{equation*}
	Once again using the local Lipschitz boundedness of $f$ in the second variable, for some bounded function $b$ we have that
	\begin{equation}\label{ch12336}
	- d \Delta\gamma - c\partial_x \gamma + b \gamma >0.
	\end{equation}
	Also, we have that
	\begin{equation*}
	\gamma(0,0)=V_{\infty}(0,0) - k^* P_{\infty} (0,0) = \underset{n\to \infty}{\lim} v_{\infty}(x_n, y_n) - k^* p(x_n, y_n)=0.
	\end{equation*}
	Since that, $\gamma \geq 0$ and \eqref{ch12336}, we  can apply the strong maximum principle and we have $\gamma \equiv 0$ in $\R^2$. Then, $V_{\infty}\equiv k^* P_{\infty}$ and 
	\begin{equation}
	0=- d \Delta\gamma - c\partial_x \gamma = f(x+ x', k^*P_{\infty}) -  k^* f(x+ x', P_{\infty}) >0,
	\end{equation}
	which is absurd. Since this was the last case to rule out, we can conclude that \eqref{ch12337} holds.
	
	From \eqref{ch12337}, we  have that 
	\begin{equation}\label{ch12346}
	(u_{\infty}, v_{\infty}) \geq (q,p).
	\end{equation}

	Now, we can repeat all the argument exchanging the role of $(u_{\infty}, v_{\infty})$ and $(q,p)$. We find 
	\begin{equation*}
	h^* := \sup \left\{  h>0 \ : \  q(x) > h u_{\infty}(x)  \ \forall x\in\R, \  p(x,y) >h  v_{\infty}(x,y) \ \forall(x,y)\in\Omega     \right\} \geq 1.
	\end{equation*}
	and
	\begin{equation*}
	(q,p) \geq  (u_{\infty}, v_{\infty}).
	\end{equation*}
	By that and \eqref{ch12346}, we have that $(u_{\infty}, v_{\infty}) \equiv (q,p)$. Hence, the uniqueness is proven.
\end{proof}

Now we are ready to give a result on the persistence of the population.
%

\begin{proof}[Proof of Theorem \ref{ch1thm:char}, part 1]

	Since $\lambda_1( \Omega)<0$, by Proposition \eqref{ch1prop:romain}, we have that there exists $R>0$ such that $\lambda_1( \Omega_R)<0$. Let us consider $(u_R, v_R)$ the eigenfunctions related to $\lambda_1( \Omega_R)<0$; then, with the argument already used in  the proof of Lemma \ref{ch1lemma:pos_inf} (precisely, in \eqref{ch12220} and \eqref{ch12221}), there exists a value $\varepsilon>0$ such that $(\varepsilon u_R, \varepsilon v_R)$ is a subsolution to \eqref{ch1sys:fieldroad} in $\Omega_R$. 
	Observe also that $u_R(x)=0$ for $x\in\partial I_R$ and $v_R(x,y)=0$ for $(x,y)\in(\partial \Omega_R)\cap\Omega$. Then, we can extend $\varepsilon u_R$ and $ \varepsilon v_R$ outside respectively $I_R$ and $\Omega_R$, obtaining the generalised subsolution $(\varepsilon u_R, \varepsilon v_R)$.
	
	Let us consider the solution $(u,v)$ issued from $(u_0, v_0)$. 
	Then, by the strong parabolic principle we have that
	\begin{equation}\label{ch12225}
	u(1, x)>0\quad \text{and} \quad  v(1,x,y)>0.
	\end{equation}
	Recall that $(u_{\infty}, v_{\infty})$ is the unique stationary solution of \eqref{ch1sys:fieldroad}, and that by Lemma \eqref{ch1lemma:pos_inf} we have
	\begin{equation}\label{ch12300}
	u_{\infty}>0\quad \text{and} \quad  v_{\infty}>0.
	\end{equation}
	By that and \eqref{ch12225},  we have that
	\begin{equation*}
	\delta:=\min\{\underset{x\in I_R}{\min} \, u(1,x), \underset{x\in I_R}{\min} \, u_{\infty}(x), 	\underset{(x,y)\in \Omega_R}{\min}  v(1,x,y), \underset{(x,y)\in \Omega_R}{\min} v_{\infty}(x,y) \} >0.
	\end{equation*}
	Without loss of generality, we can suppose 
	\begin{equation}\label{ch12301}
	\varepsilon <\delta
	\end{equation}
	and thus by \eqref{ch12225}, \eqref{ch12300}, and \eqref{ch12301}, we have
	\begin{equation}\label{ch12306}
	\begin{split}
	u_{\infty}(x) &> \varepsilon u_R(x) \quad  \text{for all}  \ x\in \R, \\
	v_{\infty}(x,y) &> \varepsilon v_R(x,y) \quad   \text{for all} \ (x,y)\in \Omega.
	\end{split}
	\end{equation}

	Now, consider the solution $(\underline{u}, \underline{v})$ issued from $(\varepsilon u_R, \varepsilon v_R)$.
	We point out that, by the comparison principle, for all $t>0$ we have
	\begin{equation}\label{ch10017}
	(\underline{u}(t,x), \underline{v}(t,x,y) \leq ({u}(t+1,x), {v}(t+1,x,y)).
	\end{equation}
	By the standard argument already used in the proof of Theorem \ref{ch1lemma:pbss}, we have that $(\underline{u}, \underline{v})$ is increasing in time and by Lemma \ref{ch1lemma:mono} it converges in $\mathcal{C}_{loc}^2$ to a stationary function $(\underline{u_{\infty}}, \underline{v_{\infty}})$ as $t$ tends to infinity. Since $(\underline{u}, \underline{v})$ is increasing in time and $(\varepsilon u_R, \varepsilon v_R)\not \equiv (0,0)$, by the strong maximum principle we have $(\underline{u_{\infty}}, \underline{v_{\infty}}) > (0,0)$. By \eqref{ch12306}, we also have 
	\begin{equation*}
	(\underline{u_{\infty}}, \underline{v_{\infty}}) \leq ({u_{\infty}}, {v_{\infty}})
	\end{equation*}
	Then, by the uniqueness of the bounded positive stationary solution proved in Theorem \ref{ch1lemma:pbss}, we have $(\underline{u_{\infty}}, \underline{v_{\infty}}) \equiv ({u_{\infty}}, {v_{\infty}})$.  
	
	Next, take  
	\begin{equation}\label{ch10008}
	V:= \max \left\{ M, \sup v_0, \frac{\mu}{\nu} \sup u_0, \sup v_{\infty}, \frac{\mu}{\nu} \sup u_{\infty}  \right\}, \quad U:= \frac{\nu}{\mu} V,
	\end{equation}
	where $M>0$ is the threshold value defined in \eqref{ch1hyp:M}.
	Making use of the hypothesis \eqref{ch1hyp:M} on $f$, one easily check that $(U, V)$ is a supersolution for \eqref{ch1sys:fieldroad}. Let us call $(\overline{u}, \overline{v})$ the solution to \eqref{ch1sys:fieldroad} issued from $(U,V)$.
	By definition, $(U,V)\geq (u_0, v_0)$, hence by the comparison principle for all $t>0$ we have
	\begin{equation}\label{ch10012}
	( u(t,x), v(t, x,y)  ) \leq (\overline{u}(t,x), \overline{v}(t, x,y) ).
	\end{equation}
	Repeating the argument used in the proof of Theorem \ref{ch1lemma:pbss}, we observe that $(\overline{u}, \overline{v})$ is decreasing in time and by Lemma \ref{ch1lemma:mono} it converges in $\mathcal{C}_{loc}^2$ to a stationary function $(\overline{u_{\infty}}, \overline{v_{\infty}})$ as $t$ tends to infinity.
	We have $(\overline{u_{\infty}}, \overline{v_{\infty}}) \leq (U,V)$, so the stationary solution is bounded. 
	Moreover, since by the definition of $(U,V)$ in \eqref{ch10008} we have $ ( {u_{\infty}},  {v_{\infty}}) \leq (U, V) $, by the comparison principle \ref{ch1prop:comparison} we get
	\begin{equation*}
	( {u_{\infty}},  {v_{\infty}}) \leq (\overline{u_{\infty}}, \overline{v_{\infty}}).
	\end{equation*}
	Since $(\overline{u_{\infty}}, \overline{v_{\infty}})$ is a bounded positive stationary solution of \eqref{ch1sys:fieldroad}, by Theorem \ref{ch1lemma:pbss} we have that $( {u_{\infty}},  {v_{\infty}}) \equiv (\overline{u_{\infty}}, \overline{v_{\infty}})$.
	
	By the comparison principle \ref{ch1prop:comparison} and by \eqref{ch10017} and \eqref{ch10012}, for all $t>1$ we have
	\begin{equation*}
	\begin{split}
	\underline{u}(t-1, x) \leq u(t,x) \leq \overline{u}(t,x) \quad \text{for all} \ x\in\R, \\
	\underline{v}(t-1, x,y) \leq v(t,x,y) \leq \overline{v}(t,x,y) \quad \text{for all} \ (x,y)\in\Omega.
	\end{split}
	\end{equation*}
	Since both $(\underline{u}, \underline{v})$ and $(\overline{u}, \overline{v})$ converge to $( {u_{\infty}},  {v_{\infty}})$ locally as $t$ tends to infinity, by the sandwich theorem we have that $(u, v)$ also does. This is precisely the statement that we wanted to prove. 
\end{proof}

\subsection{Extinction}

The first step to prove extinction is to show that there is no positive bounded stationary solution to system \eqref{ch1sys:fieldroad}, that is, the only bounded stationary solution is $(0,0)$.

\begin{lemma}\label{ch1lemma:ext}
	Suppose $c=0$ and $f$ satisfy \eqref{ch1hyp:0}-\eqref{ch1hyp:per}.
	If $\lambda_1( \Omega)\geq 0$, then there is no positive bounded stationary solution to system \eqref{ch1sys:fieldroad}.
\end{lemma}

\begin{proof}
	 	\emph{Step 1: construction of a supersolution.}
	Observe that in this case, since $c=0$, by Theorem \ref{ch1thm:eq_1_p} it holds $\lambda_p(\Omega)=\lambda_1(\Omega)\geq 0$.
	We take the couple of eigenfunctions $(u_p, v_p)$ related to $\lambda_p(\Omega)$ as prescribed by Proposition \ref{ch1prop:-3}; recall that $(u_p, v_p)$ are periodic in $x$.
	Suppose $(q,p)$ is a positive bounded stationary solution to \eqref{ch1sys:fieldroad}. Then, there exists $\eta>0$ such that
	\begin{equation}\label{ch12347}
	q(0) > \eta u_p(0).
	\end{equation}
	We now choose a smooth function $\chi :  \R_{\geq 0} \to \R_{\geq 0}$  such that  $\chi(y)=0$ for $y\in[0,\ell]$, $\chi(y)=1$ for $y\in[ 2\ell, +\infty)$. 

	By \eqref{ch1eq:-5} and Theorem \ref{ch1thm:1.7inbr}, we have $\lambda_p(-\mathcal{L}, \R)=\lambda_p(-\mathcal{L}, \R^2)=\lambda_1(-\mathcal{L}, \R^2)$. By that, Theorem \ref{ch1thm:ineq2} and the fact that $\lambda_1(\Omega)\geq 0$, we get $\lambda_p(-\mathcal{L}, \R)\geq 0$.
	We call $\psi_p$ the eigenfunction related to $\lambda_p(-\mathcal{L}, \R)$ and, with a slight abuse of notation, we extend it to $\R^2$ by considering it constant with respect to the variable $y$.
	Take $\varepsilon>0$ to be fixed after, and define
	\begin{equation*}
	(\overline{u}(x),  \overline{v}(x,y)):= (\eta u_p(x),  \eta v_p(x,y) + \varepsilon \chi(y) \psi_p(x)).
	\end{equation*}
	Then, it holds that
	\begin{equation}\label{ch12011}
	\begin{split}
	- d \Delta \overline{v}  
	&= -d \left(\Delta \eta v_p +\varepsilon \chi''\psi_p + \varepsilon \chi \psi_p''   \right), \\
	&= \left( f_v(x,0)+\lambda_p( \Omega) \right) \eta v_p + ( f_v(x,0)+ \lambda_p(-\mathcal{L}, \R) ) \varepsilon \chi \psi_p - d \varepsilon \chi''\psi_p, \\
	& = f_v(x,0) \overline{v} + \lambda_p( \Omega) \eta v_p + \lambda_p(-\mathcal{L}, \R) \varepsilon \chi \psi_p - d \varepsilon \chi''\psi_p.
	\end{split}
	\end{equation}
	Using the KPP hypothesis \eqref{ch1hyp:KPP} and the boundedness of $\chi''$, for $\varepsilon$ small enough we have
	\begin{equation*}
		f_v(x,0) \overline{v} - d \varepsilon \chi''\psi_p > f(x, \overline{v}).
	\end{equation*}
	By that, \eqref{ch12011} and the non negativity of $\lambda_p(\Omega)$ and $\lambda_p(-\mathcal{L}, \R)$, we have 
	\begin{equation*}
		- d \Delta \overline{v} > f(x, \overline{v}).
	\end{equation*}
	This means that $\overline{v}$ is a supersolution for the second equation of \eqref{ch1sys:fieldroad}. 
	
	Since by definition for $y\leq \ell$ we have $\chi(y)=0$, it holds that
	\begin{equation}\label{ch12010}
	(\overline{u}(x),  \overline{v}(x,y))\equiv (u_p(x),   v_p(x,y)) \quad \text{for all} \ (x, y)\in \R\times (0, \ell). 
	\end{equation} 
	By the fact that $\lambda_p(\Omega  ) \geq 0$, it is easy to check that $(u_p(x),   v_p(x,y))$ is a supersolution for the first and third equation in \eqref{ch1sys:fieldroad}. By \eqref{ch12010}, the same holds for  $(\overline{u}(x),  \overline{v}(x,y))$.
	This, together with \eqref{ch12011}, gives that $(\overline{u}(x),  \overline{v}(x,y))$ is a supersolution to \eqref{ch1sys:fieldroad}.

	\emph{Step 2: construction of a bounded supersolution}
	Now we distinguish two cases.
	If $ v_p$ is bounded, then we take
	\begin{equation}\label{ch1case1super}
	(\tilde{u}, \tilde{v}):= (\bar{u}, \bar{v})
	\end{equation}
	Otherwise, we proceed as follows.
	Since in this other case $v_p$ is unbounded, and since it is periodic in $x$, this means there exists a sequence $\{(x_n, y_n)\}_{n\in\N}$ such that 
	\begin{equation}\label{ch11721b}
	v_p(x_n, y_n) \to \infty, \ y_n\to \infty \quad \text{as} \ n\to\infty.
	\end{equation}
	Now, consider 
	\begin{equation}\label{ch11418}
	V:= \max \left\{ \underset{[0,\ell]\times [0, 3\ell]}{\max} v_p +1, \ \underset{[0,\ell]}{\max} \, \frac{\nu}{\mu} u_p +1 ,  \ M       \right\},
	\end{equation}
	where $M$ is the quantity defined in \eqref{ch1hyp:M}.
	Take the set $S:=(-\ell, \ell)\times(-\ell, \ell)$ and the constant $C$ of the Harnack inequality (see Theorem 5 in Chapter 6.4 of \cite{evans}) on the set $S$ for the operator $L(\psi)=\mathcal{L}(\psi)+\lambda_1(\Omega)\psi$. Then, by \eqref{ch11721b}, for some $N\in\N$ we have  
	\begin{equation*}
	V \leq \frac{1}{C} v_p(x_N, y_N).
	\end{equation*}
	Then by using that and Harnack inequality on $v_p(x+x_N,y+y_N)$ in the set $S$, we get 
	\begin{equation*}
	V \leq \frac{1}{C} \, \underset{S}{\sup} \, v_p(x, y)
	\leq \underset{S}{\inf} \, v_p(x,y),
	\end{equation*}
	Then, using the periodicity of $v_p$, we get
	\begin{equation}\label{ch1comp1}
	V \leq v_p(x, y_N)  \quad \text{for all} \ x\in\R.
	\end{equation}
	Now, define
	\begin{equation}\label{ch1case2superv}
	\tilde{v}(x,y):= \left\{
	\begin{array}{ll}
	\min \{  V, \bar{v}(x,y) \} & \text{if} \ y \leq y_N, \\
	V & \text{if} \ y > y_N.
	\end{array}
	\right.
	\end{equation}
	Also, we define
	\begin{equation*}
	U := \frac{\nu}{\mu} V
	\end{equation*}
	and
	\begin{equation}\label{ch1case2superu}
	\tilde{u}:= \min\{ U, u_p  \}.
	\end{equation}
	By the definition of $V$ in \eqref{ch11418}, one readily checks that $(U,V)$ is a supersolution for system \eqref{ch1sys:fieldroad} and that 
	\begin{equation}\label{ch1comp2}
	\tilde{u} = u_p \quad \text{and} \quad \tilde{v}(x,0)=v_p(x,0).
	\end{equation}
	We point out that by the definition of $(\tilde{u}, \tilde{v})$, \eqref{ch1comp1} and \eqref{ch1comp2}, for any $(\underline{u}, \underline{v})$ subsolution to system \eqref{ch1sys:fieldroad}, we will be able to apply
	the generalised comparison principle, Proposition 3.3 appeared in \cite{brr}. 
	Moreover, $(\tilde{u}, \tilde{v})$ is bounded from above by $(U,V)$.

	By the fact that $(u_p,  v_p)$ is a couple of generalised periodic eigenfunctions to \eqref{ch1sys:upvp}, by the strong maximum principle we have that 
	\begin{equation}\label{ch12341}
	\begin{split}
	\tilde{u}(x) &\geq \underset{[0,\ell]}{\min} \, \eta u_p(x') >0 \quad \text{for} \ x\in\R, \\
	\tilde{v}(x,y) &\geq \underset{ [0,\ell]\times[0,2\ell]}{\min} \eta  v_p(x',y') >0 \quad \text{for} \ (x,y)\in\R\times[0, 2\ell], \\
	\tilde{v}(x,y) &\geq \min\{\underset{[0,\ell]}{\min} \, \varepsilon\psi_p(x'), V\} >0 \quad \text{for} \ (x,y)\in\R\times(2\ell, +\infty).
	\end{split}
	\end{equation}
	
	\emph{Step 3: comparison with the stationary solution.}
	Next, define
	\begin{equation*}
	k^*:= \inf \{ k\geq 0 \ : \ k( \tilde{u}(x),   \tilde{v}(x,y)) > (q,p) \ \text{for all} \ (x,y)\in\Omega       \}.
	\end{equation*}
	Since by \eqref{ch12341} we have that $ \tilde{u}(x)$ and $ \tilde{v}(x,y)$ are bounded away from $0$, and since $(q,p)$ is bounded by hypothesis, we get that $k^*<+\infty$.
	By \eqref{ch12347}, we have that
	\begin{equation}\label{ch12234}
	k^*>1.
	\end{equation}
	Then, either
	\begin{equation}\label{ch1case1}
	\mbox{there exists a sequence $\{x_n\}_{n\in\N}\subset \R$ such that $k^*  \tilde{u}(x_n) - q(x_n) \overset{n\to\infty}{\longrightarrow} 0$,}
	\end{equation}
	or
	\begin{equation}\label{ch1case2}
	\mbox{
		there exists a sequence $\{(x_n, y_n)\}_{n\in\N}\subset \Omega$ such that $k^*  \tilde{v}(x_n,y_n) -  p(x_n, y_n)  \overset{n\to\infty}{\longrightarrow} 0$.}
	\end{equation}

	As usual, for all $n\in\N$ we take $x_n'\in[0,\ell)$ such that $x_n-x_n'\in\ell \Z$. Up to a subsequence, $\{ x_n'\}_{n\in\N}$ is convergent and we call
	\begin{equation*}
	x'= \underset{n\to\infty}{\lim} x_n' \in[0,\ell].
	\end{equation*}

	\emph{Step 4: $\{y_n\}_{n\in\N}$ is bounded}.
	If $\{y_n\}_{n\in\N}$ is bounded, consider a converging subsequence and call $y'= \underset{n\to \infty}{\lim} y_n$.
	
	We define 
	\begin{equation*}
	(q_n(x), p_n(x,y)):=(q(x+x_n), p(x+x_n,y)). 
	\end{equation*}
	By Lemma \ref{ch1lemma:converg}, $(q_n, p_n)$ converges in $\mathcal{C}_{loc}^2$ to some $(q_{\infty}, p_{\infty})$ such that $(q_{\infty}(x-x'), p_{\infty}(x-x', y))$ solves \eqref{ch1sys:fieldroad}.
	Define the functions 
	\begin{align*}
	\alpha(x) &:= k^*  \tilde{u}(x)-q_{\infty}(x-x'), \\
	\beta(x,y)&: = \tilde{v}(x,y)- p_{\infty}(x-x', y).
	\end{align*}
	
	If we are in the case of \eqref{ch1case1}, then by the periodicity of $\tilde{u}$ we get
	\begin{equation*}
	\alpha(x')= k^*  \tilde{u}(x')-q_{\infty}(0)= \underset{n\to\infty}{\lim} (  k^*  \tilde{u}(x_n)- q(x_n)   )=0.
	\end{equation*}
	Moreover, by the definition of $k^*$, we have that $\alpha\geq 0$. Also, $\alpha$ satisfies
	\begin{equation*}
	-D \alpha ''  -\nu \beta|_{y=0}+ \nu \alpha \geq 0. 
	\end{equation*}
	Then, the strong maximum principle yields that, since $\alpha$ attains its minimum at $x=x'$, then $\alpha\equiv0$. Then, by the comparison principle 3.3 in \cite{brr} we have that $\beta\equiv 0$, hence 
	\begin{equation}\label{ch12226}
	0= -d \Delta \beta \geq k^*f(x,  \tilde{v}) - f (x,p_{\infty}(x-x',y)).
	\end{equation}
	By \eqref{ch12234}, we have that $k^*  \tilde{v} >  \tilde{v}$. Hence, by the Fischer-KPP hypothesis \eqref{ch1hyp:KPP}, we have that
	\begin{equation}\label{ch12249}
	\frac{f(x, k^*  \tilde{v})}{k^*  \tilde{v}} < \frac{f(x,  \tilde{v})}{ \tilde{v}}.
	\end{equation}
	Hence, again by the fact that $\beta\equiv 0$, we  have $p_{\infty}(x-x',y)\equiv  k^* \tilde{v}$; by that and by \eqref{ch12249}, it holds
	\begin{equation}\label{ch12257}
	k^* f(x,  \tilde{v})- f(x,p_{\infty}(x-x',y))=k^* f(x,  \tilde{v})- f(x, k^* \tilde{v})>0.
	\end{equation}
	But this is in contradiction with \eqref{ch12226}, hence this case cannot be possible. 
	
	If instead \eqref{ch1case2} holds, we get that 
	\begin{equation}\label{ch12256}
	\beta(x', y')= k^*  \tilde{v}(x', y')- p_{\infty}(0, y')= \underset{n\to \infty}{\lim} k^*  \tilde{v}(x_n, y_n)- p(x_n, y_n)=0.
	\end{equation}
	By the definition of $k^*$ we also have that $\beta \geq 0$. Moreover, we get that
	\begin{equation*}
	-d \Delta \beta \geq f(x,k^*  \tilde{v}) - f (x,p_{\infty}(x-x',y))
	\end{equation*}
	using the fact that $ \tilde{v}(x,y)$ is a supersolution, $p_{\infty}(x-x',y)$ is a solution, and \eqref{ch12249}.
	Since $f$ is Lipschitz in the second variable, uniformly with respect to the first one, there exists some function $b$ such that
	\begin{equation*}
	-d \Delta \beta  - b \beta \geq 0.
	\end{equation*}
	If $y'>0$, using the strong maximum principle and owing \eqref{ch12256}, we have that $\beta \equiv 0$. 
	If instead $y'=0$, recall that it also holds
	\begin{equation*}
	-d \partial_y \beta|_{y=0} \geq \mu \alpha -\nu \beta.
	\end{equation*}
	Hence, in $(x,y)=(x', y')$, we get that $\partial_y \beta(x',y') \leq 0$. By Hopf's lemma, we get again that $\beta\equiv 0$.
	
	But $\beta\equiv 0$ leads again to  \eqref{ch12226} and \eqref{ch12257}, giving an absurd, hence also this case is not possible.

	\emph{Step 5: $\{y_n\}_{n\in\N}$ is unbounded}.
	We are left with the case of $\{y_n\}_{n\in\N}$ unbounded.
	Up to a subsequence, we can suppose that $\{y_n\}_{n\in\N}$ is increasing.
	We define
	\begin{equation*}
	P_n(x,y):= p(x+x_n, y+y_n).
	\end{equation*}
	By Lemma \ref{ch1lemma:converg} we have that, up to a subsequence, $\{P_n\}_{n\in\N}$ converges in $\mathcal{C}_{loc}^{2,\alpha}(\R^2)$ to some function $P_{\infty}$  such that $P_{\infty}(x-x',y)$ is a solution  to the second equation in \eqref{ch1sys:fieldroad} in $\R^2$.

	Now we have two cases depending on how $(\tilde{u}, \tilde{v})$ was constructed. If $v_p$ is bounded, we have defined the supersolution as in \eqref{ch1case1super}. Then,
	by defining 
	\begin{equation*}
	v_n(x,y):=v_p(x+x_n, y+y_n)
	\end{equation*}
	and applying Lemma \ref{ch1lemma:converg}, we have that $v_n$ converges locally uniformly to  a bounded function $v_{p, \infty}$ such that $v_{p, \infty}(x-x',y)$ satisfies
	\begin{equation}\label{ch11630}
	-d\Delta v_{p, \infty}(x-x',y) = (f_v(x,0)+\lambda_1(\Omega))v_{p, \infty}(x-x',y).
	\end{equation}
	In this case, we define 
	\begin{equation*}
	v_{\infty}(x,y):=\eta v_{p, \infty}(x,y) + \varepsilon \psi_p(x+x').
	\end{equation*}
	We point out that $v_{\infty}(x-x',y)$ is a  periodic supersolution of the second equation in \eqref{ch1sys:fieldroad} by \eqref{ch11630} and \eqref{ch12011}.
	
	If instead $v_p$ is unbounded, by \eqref{ch1case2superu} for $y>y_N$ we have $\tilde{v}=V$. In this case, we choose
	\begin{equation*}
	v_{\infty}:=V.
	\end{equation*}
	By the definition of $V$ in \eqref{ch11418}, we have that $v_{\infty}$ is also a supersolution to \eqref{ch1sys:fieldroad}.
	
	We call $\gamma(x,y):=k^* {v}_{\infty}(x-x',y) - P_{\infty}(x-x',y)$.
	Hence, $\gamma(x,y)\geq 0$ and
	\begin{equation}\label{ch12330}
	\gamma(x', 0)=k^* {v}_{\infty}(0, 0) - P_{\infty}(0,0)= \underset{n\to\infty}{\lim} k^*  \tilde{v}(x_n,y_n) - p(x_n, y_n)=0. 
	\end{equation}
	Notice than that, since \eqref{ch12234} holds, from the Fisher-KPP hypothesis on $f$ \eqref{ch1hyp:KPP}, we get
	\begin{equation*}
	\frac{f(x,k^* {v}_{\infty} )}{k^* {v}_{\infty}} < \frac{f(x, {v}_{\infty} )}{ {v}_{\infty}}.		
	\end{equation*}
	Using that, the fact that $k^* {v}_{\infty}(x-x',y)$ is a supersolution, and the fact that $P_{\infty}(x-x',y)$ is a solution, we obtain
	\begin{equation}\label{ch12333}
	-d \Delta \gamma  > f(x,  k^*{v}_{\infty}(x-x',y)) - f (x,P_{\infty}(x-x',y)).
	\end{equation}
	Since $f$ is Lipschitz in the second variable, uniformly with respect to the first one, there exists some function $b$ such that
	\begin{equation*}
	-d \Delta \gamma  - b \gamma \geq 0.
	\end{equation*}
	Using the strong maximum principle for the case of positive functions, since \eqref{ch12330} holds, we have that $\gamma\equiv 0$. As a consequence, from  \eqref{ch12333} we have
	\begin{equation*}
	f(x,k^*  {v}_{\infty}) - f (x,P_{\infty}) <0.
	\end{equation*} 
	but it also holds that $k^*  {v}_{\infty} \equiv P_{\infty}$, hence we have an absurd.

	Having ruled out all the possible cases, we can conclude that there exists no bounded positive stationary solution $(q,p)$ to \eqref{ch1sys:fieldroad}.
\end{proof}

At last, we are ready to prove the first part of Theorem \ref{ch1thm:char}.

\begin{proof}[Proof of Theorem \ref{ch1thm:char}, part 1]
	Define $$V:=\max \left\{ M, \sup v_0, \frac{\mu}{\nu} \sup u_0  \right\} \quad \text{and} \quad U:= \frac{\nu}{\mu} V.$$
	It is easy to check that $(U, V)$ is a supersolution for \eqref{ch1sys:fieldroad}. Then take $(\overline{U}, \overline{V})$ to be the solution to \eqref{ch1sys:fieldroad} with initial datum $(U,V)$. 
	Notice that by the comparison principle 
	\begin{equation}\label{ch10013}
	(0,0)\leq (u(t,x), v(t,x,y)) \leq (\overline{U}(t,x), \overline{V}(t,x,y)) \quad \text{for all} \ t>0, \ (x,y)\in\Omega.
	\end{equation}
	
	Since $(U, V)$ is a supersolution, we have that
	\begin{equation}\label{ch12317}
	(\overline{U}, \overline{V}) \leq (U,V) \quad \text{for all} \ t\geq 0.
	\end{equation} 
	Consider $\tau>0$ and call $(\tilde{U}, \tilde{V})$ the solution staring with initial datum $(U, V)$ at $t=\tau$. By \eqref{ch12317} we have that $(\overline{U}(\tau,x), \overline{V}(\tau,x, y)) \leq (U,V)$, hence by the comparison principle \eqref{ch1prop:comparison} we have  $(\overline{U}, \overline{V}) \leq (\tilde{U}, \tilde{V})$. By the arbitrariness of $\tau$, we get that $(\overline{U}, \overline{V})$ is decreasing in $t$.
	
	By Lemma \ref{ch1lemma:mono},  $(\overline{U}, \overline{V})$ converges locally uniformly to a stationary solution $(q,p)$. But by Lemma \ref{ch1lemma:ext}, the only stationary solution is $(0,0)$.
	By that and \eqref{ch10013}, we have that $(u(t,x), v(t,x,y))$ converges locally uniformly to $(0,0)$ as $t$ goes to infinity.
	
	Moreover, since $(U,V)$ is constant in $x$, and \eqref{ch1sys:fieldroad} is periodic in x, $(\overline{U}, \overline{V})$ is periodic in $x$.
	Hence, the convergence is uniform in $x$.

	Now suppose by the absurd that the convergence is not uniform in $y$; hence there exists some $\varepsilon>0$ such that for infinitely many $t_n\geq 0$, with $\{t_n\}_{n\in\N}$ an increasing sequence, and $( x_n, y_n)\in\Omega$, it holds  
	\begin{equation}\label{ch12327}
	\overline{V}(t_n,  x_n, y_n) >\varepsilon.
	\end{equation}
	Since $\overline{V}$ is periodic in $x$, without loss of generality we can suppose $x_n\in[0,\ell]$ and that up to a subsequence $\{x_n\}_{n\in\N}$ converges to some $x'\in[0,\ell]$. If $\{y_n\}_{n\in\N}$ were bounded, by \eqref{ch12327} the local convergence to $0$ would be contradicted; hence $y_n$ is unbounded.
	
	Then, define the sequence of functions
	\begin{equation*}
	V_n(t,x,y)=\overline{V}(t, x+x_n, y+y_n).
	\end{equation*}
	By \eqref{ch12327}, we have that
	\begin{equation}\label{ch10109}
	V_n(t_n,0,0)>\varepsilon \quad \text{for all} \ n\in\N.
	\end{equation}
	
	Also, since $V_n$ is bounded, by arguments similar to the ones used in Lemma \ref{ch1lemma:converg} and Lemma \ref{ch1lemma:mono}
	one can prove that, up to a subsequence, 
	$\{V_n\}_{n\in\N}$ converges in $\mathcal{C}_{loc}^2(\R^2)$ to a function $\tilde{V}$  that   solves
	\begin{equation}\label{ch10108}
	\partial_t \tilde{V} - d\Delta \tilde{V}  = f(x+x', \tilde{V}).
	\end{equation}
	Also by \eqref{ch10109}, we have that
	\begin{equation}\label{ch10110}
	\tilde{V}(t_n,0,0 )>\varepsilon \quad \text{for all} \ n\in\N.
	\end{equation}
	Recall that by the fact that $\lambda_1(\Omega)\geq 0$, Corollary \ref{ch1thm:ineq2} and Theorem \ref{ch1thm:1.7inbr}, $\lambda_p(-\mathcal{L}, \R^2)\geq 0$. Then by Theorem \ref{ch1thm:2.6inbhroques} we have that every solution to \eqref{ch10108} converges uniformly to $0$. But this is in contradiction with \eqref{ch10110}, hence we have an absurd and we must refuse the existence such positive $\varepsilon$. So, the convergence of $\overline{V}$ to $0$ is uniform in space.
	As a consequence, the convergence of $(u(t,x), v(t,x,y))$ to $(0,0)$ is uniform in space.
\end{proof}





\chapter[A New Lotka-Volterra Competitive System]{\label{ch2}Civil Wars: A New Lotka-Volterra Competitive System and Analysis of Winning Strategies} 

	We introduce a new model in population dynamics that describes
	two species sharing the same environmental resources in a situation of open hostility.
	Our basic assumption is that one of the populations deliberately seek for hostility through "targeted attacks". Hence, the
	interaction among these populations is described not in terms
	of random encounters but via the
	strategic decisions of one population
	that can attack the other according to different levels of aggressiveness.
	
	One of the features that distinguishes this model from usual competitive systems is that it allows one of the population to go extinct {\em in finite time}.
	This leads to a non-variational model for the two populations at war, taking into
	account structural parameters such as the relative fit of the two populations with
	respect to the available resources and the effectiveness of the attack strikes of the
	aggressive population.
	
	The analysis that we perform is rigorous and focuses
	on the dynamical properties of the system, by detecting
	and describing all the possible equilibria and their basins of attraction.
	
	Moreover, we will analyze the strategies that may lead to the victory of the aggressive
	population, i.e. the choices of the aggressiveness parameter,
	in dependence of the structural constants of the system and possibly varying in time
	in order to optimize the efficacy of the attacks, which take to the extinction in finite time
	of the defensive population. 
	
	The model that we present is flexible enough to also include commercial competition models of companies
	using aggressive policies against the competitors (such as misleading advertising, or releasing computer
	viruses to set rival companies out of the market).

	This chapter corresponds to the paper \cite{wars} in collaboration with Serena Dipierro, Luca Rossi and Enrico Valdinoci.

\section{Introduction}

Among the several models dealing with the
dynamics of biological systems, the case of populations engaging into a mutual conflict
seems to be unexplored.
This chapter aims at laying the foundations of a new model describing
two populations competing for the same resource with one aggressive population
which may attack the other:
concretely, one may think of
a situation in which
two populations live together in the same territory and share the same
environmental resources,
till one population wants to prevail and try to kill the other.
We consider this situation as a ``civil war'', since the two populations share
land and resources; the two populations may be equally fit to the environment
(and, in this sense, they are ``indistinguishable'', up to the aggressive attitude of
one of the populations), or they can have a different compatibility to the resources
(in which case one may think that the conflict could be motivated by the different
accessibility to environmental resources).

Given the lack of reliable data related to civil wars, a foundation of
a solid
mathematical theory for this type of conflicts may only leverage on the deduction
of the model from first principles: we follow this approach to obtain
the description of the problem in terms of a system of two
ordinary differential equations, each describing the evolution in time
of the density
of one of the two populations.

The method of analysis that we adopt is a combination
of techniques from different fields, including ordinary differential equations,
dynamical systems and optimal control.

This viewpoint will allow us to rigorously investigate the model,
with a special focus on a number of mathematical features of
concrete interest, such as the possible extinction of one of the two populations
and the analysis of the strategies that lead to the victory of the aggressive population.

In particular, we will analyze the {\it dynamics of the system},
characterizing the equilibria and their features (including possible basins of attraction)
in terms of the different parameters
of the model (such as relative fitness to the environment, aggressiveness
and effectiveness of strikes). Also, we will study the initial configurations which
may lead to the victory of the aggressive population, also taking into account
different possible {\it strategies} to achieve the victory: roughly speaking,
we suppose that the aggressive population may adjust the parameter
describing the aggressiveness in order to either dim or 
exacerbate the conflict with the aim of destroying the second population
(of course, the war has a cost in terms of life
for both the populations,
hence the aggressive population must select the appropriate strategy in terms
of the structural parameters of the system). We will show that the initial data
allowing the victory of the aggressive population
does not exhaust the all space, namely {\it there exists initial configurations
	for which the aggressive population cannot make the other extinct}, regardless the strategy adopted during the conflict. 

Furthermore, {\it for identical populations with the same fit to the environment
	the constant strategies suffices} for the aggressive population to possibly
achieve the victory: namely, if an initial configuration admits a piecewise continuous in time
strategy that leads to the victory of the aggressive population,
then it also admits a constant in time strategy that reaches the same objective
(and of course, for the aggressive population, the possibility of focusing
only on constant strategies would entail concrete practical advantages).

Conversely, {\it for populations with different fit to the environment,
	the constant strategies do not exhaust all the winning strategies}:
that is, in this case, there are initial conditions which allow
the victory of the aggressive population only under the exploitation
of a strategy that is not constant in time.

In any case, we will also prove that {\em strategies with at most one jump
	discontinuity are sufficient} for the aggressive population:
namely, independently from the relative fit to the environment,
if an initial condition allows the aggressive population to reach the victory
through a piecewise continuous in time
strategy, then the same goal can be reached using a ``bang-bang''
strategy with at most one jump.

We will also discuss the {\em winning strategies that minimize
	the duration of the war}: in this case, we will show that
jump discontinuous strategies may be not sufficient and interpolating
arcs have to be taken into account.
\medskip

We now describe in further detail our model of conflict between the
two populations and the attack strategies pursued by the aggressive population.
Our idea is to modify the Lotka-Volterra competitive system for two populations with
density~$u$ and~$v$,
adding to the usual competition for resources the fact that both populations suffer some losses as an outcome of the attacks. 
The key point in our analysis
is that the clashes do not depend on the chance of meeting of the two populations, given by the quantity~$uv$, as it happens in many other works in the literature  (starting from the publications of Lotka and Volterra,~\cite{lotka} and~\cite{volterra}), but they are sought by the first population and
depend only on the size~$u$ of the first population and on its level of aggressiveness~$a$.
The resulting model is
\begin{equation}\label{ch2model}
\left\{
\begin{array}{llr}
\dot{u}&= u(1-u-v) - acu, & {\mbox{ for }}t>0,\\
\dot{v}&= \rho v(1-u-v) -au, & {\mbox{ for }}t>0,
\end{array}
\right.
\end{equation}
where~$a$,~$c$ and~$\rho$ are nonnegative real numbers. Here, the coefficient~$\rho$ models the fitness of the second population with respect of the first one when resources are abundant for both; it is linked with the exponential growth rate of the two species. The parameter~$c$ here stands for the quotient of endured per inflicted  damages for the first population. 
Deeper justifications to the model~\eqref{ch2model} will be given in
Subsection~\ref{ch2ss:derivation}.

Notice that the size of the second population~$v$ may become negative in finite time while the first population is still alive. The situation where~$v=0$ and~$u>0$ represents the extinction of the second population and the victory of the first one. 

To describe our results, for communication convenience (and in spite of our
personal fully pacifist believes)
we take the perspective of the first population,
that is, the aggressive one; the objective
of this population is to win the war, and, to achieve that, it can influence the system by tuning the parameter~$a$. 

{F}rom now on, we may refer to the parameter~$a$ as the \textit{strategy}, that may also depend on time, and we will say that it is \textit{winning} if it leads to victory of the first population. 

The main problems that
we deal with in this chapter are:
\begin{enumerate}
	\item The characterization of the {\em initial conditions for which there exists a winning strategy}.
	\item The {\em success of the constant strategies}, compared to all  possible strategies.
	\item The {\em construction of a winning strategy} for a given initial datum.
	\item The {\em existence of a single winning strategy independently of the initial datum}.
\end{enumerate}

We discuss all these topics in Subsection~\ref{ch2ss:strategy}, presenting
concrete answers to each of these problems.

Also, since to our knowledge this is the first time that system~\eqref{ch2model} is considered,
in Subsections~\ref{ch2ss:notation} and~\ref{ch2ss:dynamics} we will discuss the dynamics and some interesting results about the dependence of the basins of attraction on the other parameters. 

It would also be extremely interesting to add the space component to our model, by considering a system of reaction-diffusion equations. This will be the subject of a further work. 


\subsection{Motivations and derivation of the model} \label{ch2ss:derivation}

The classic Lotka-Volterra equations were first introduced for modelling population dynamics between animals~\cite{volterra} and then used to model other phenomena involving competition, for example in technology substitution~\cite{substitution}. The competitive Lotka-Volterra system concerns the sizes~$u_1(t)$ and~$u_2(t)$ of two species competing for the same resources. The system that
the couple~$(u_1(t), u_2(t))$ solves is
\begin{equation}\label{ch2lv}
\begin{cases}
\dot{u}_1=r_1  u_1\left(\sigma-\displaystyle
\frac{u_1+\alpha_{12} u_2}{k_1}  \right), & t>0,\\
\dot{u}_2
=r_2 u_2\left(\sigma- \displaystyle\frac{u_2+\alpha_{21} u_1}{k_2}  \right), & t>0,
\end{cases}
\end{equation}
where~$r_1$,~$r_2$,~$\sigma$,~$\alpha_{12}$,~$\alpha_{21}$,~$k_1$ and~$k_2$ are nonnegative real numbers.

Here, the coefficients~$\alpha_{12}$ and~$\alpha_{21}$ represent the competition between individuals of different species, and indeed they
appear multiplied by the term~$u_1 u_2$, which represents a probability of meeting. 

The coefficient~$r_i$ is the exponential growth rate of the~$i-$th population, that is, the reproduction rate that is observed when the resources are abundant. 
The parameters~$k_i$ are called carrying capacity and represent the number of individuals of the~$i-$th population that can be fed with the resources of the territory, that are quantified by~$\sigma$. 
It is however usual to rescale the system in order to reduce the number of parameters. 
In general,~$u_1$ and~$u_2$ are rescaled so that they vary in the interval~$[0,1]$,
thus describing densities of populations.

The behavior of the system depends substantially on the values of~$\alpha_{12}$ and~$\alpha_{21}$ with respect to the threshold given by the value~$1$, see e.g.~\cite{nonlinear}: if~$\alpha_{12}<1<\alpha_{21}$, then the first species~$u_1$
has an advantage over the second one~$u_2$
and will eventually prevail; if~$\alpha_{12}$ and~$\alpha_{21}$ are both strictly above or below the threshold, then the first population that penetrates the environment (that is, the one that has a greater size at the initial time) will persist while the other will extinguish. 

Some modification of the Lotka-Volterra model were made in stochastic analysis by adding a noise term of the form~$-f(t)u_i$ in the~$i-$th equation, finding some interesting phenomena of phase transition,
see e.g.~\cite{noise}. 

The ODE system
in~\eqref{ch2lv} is of course the cornerstone to study the case of two competitive populations that diffuse in space. Many different types of diffusion have been compared and one can find a huge literature on the topic, see ~\cite{mimura, crooks, massaccesi} for some examples and~\cite{murray2} for a more general overview. 
We point out that other dynamic systems presenting finite time extinction of one or more species have been generalised for heterogeneous environments, see for example  the model in~\cite{gaucel} for the predator-prey behaviour of cats and birds, that has been thereafter widely studied. 
\medskip

In this chapter,
we will focus not only on {\em
	basic competition for resources}, but also on {\em
	situations of open hostility}.
In social sciences, war models are in general little studied; indeed, the collection of data up to modern times is hard for the lack of reliable sources. 
Also, there is still much discussion about what factors are involved and how to quantify them: in general, the outcome of a war does not only depend on the availability of resources, but also on more subtle factors as the commitment of the population and the knowledge of the battlefield,
see e.g.~\cite{toft2005state}. 
Instead, the causes of war were investigated by the statistician L.F. Richardson, who proposed some models for predicting the beginning of a conflict,
see~\cite{richardson1960arms}. 

In addition to the human populations, behavior of hostility between groups of the same species has been observed in chimpanzee. Other species with complex social behaviors are able to coordinate attacks against groups of different species: ants versus termites, agouti versus snakes, small birds versus hawk and owls, see e.g.~\cite{animalwar}.  

The model that
we present here
is clearly a simplification of reality. Nevertheless,
we tried to capture some important features of conflicts
between rational and strategic populations, 
introducing in the
mathematical modeling the new idea that a conflict may be sought
and the parameters that influence its development may be
conveniently adjusted. 

Specifically, in our model, the
interactions between populations are not merely driven
by chance and the strategic decisions of the population
play a crucial role in the final outcome of the conflict, and
we consider this perspective as an interesting novelty in the mathematical
description of competitive environments.


At a technical level, our aim is to introduce a model for conflict between
two populations~$u$ and~$v$,
starting from the model when the two populations compete for food and modifying it to add the information about the clashes. 
We imagine that
each individual of the first population~$u$ decides to attack an individual of the second population with some probability~$a$ in a given period of time. We assume that hostilities take the form of ``duels'', that is, one-to-one fights. 
In each duel, the  individual of the first population has a probability~$\zeta_u$ of being killed and a probability~$\zeta_v$ of killing his or her opponent; 
notice that in some duel the two fighters might be both killed. 
Thus, after one time-period, the casualties for the first and second populations
are~$a\zeta_u u$ and~$a\zeta_v u$
respectively.
The same conclusions are found if we imagine that the first population forms an army to attack the second, which tries to resist by recruting an army of proportional size. At the end of each battle, a ratio of the total soldiers is dead, and this is again of the form~$a\zeta_u u$ for the first population and~$a\zeta_v u$ for the second one.

Another effect that
we take into account is the drop in  the fertility of the population  during wars. 
This seems due to the fact that families suffer some income loss during war time, because of a lowering of the average productivity and lacking salaries only partially compensated by the state; another reason possibly discouraging couples to have children is the increased chance of death of the parents during war. 
As pointed out in~\cite{fertility}, in some cases the number of lost births during wars are comparable to the number of casualties. 
However, it is not reasonable to think that this information should be included in the exponential growth rates~$r_u$  and~$r_v$, because the fertility drop really depends on the intensity of the war. For this reason, 
we introduce the parameters~$c_u\geq 0$ and~$c_v\geq 0$   that are to be multiplied by~$a u$ for both populations.

Moreover, for simplicity,
we also suppose that the clashes take place apart from inhabited zone, without having influence on the harvesting of resources. 
\medskip

Now we derive the system of equations from a microlocal analysis.
As in the Lotka-Volterra model, it is assumed that the change of the size of the population in an interval of time~$\Delta t$ is proportional to the size of the population~$u(t)$, that is
\begin{equation*}
u(t+\Delta t)-u(t) \approx u(t) f(u,v)
\end{equation*}
for some appropriate function~$f(u,v)$. In particular,~$f(u,v)$ should depend on resources that are available and reachable for the population. 
The maximum number of individuals that can be fed with all the resources of the environment is~$k$; taking into account all the individuals of the two populations, the available resources are 
\begin{equation*}	
k-u-v.	
\end{equation*}
Notice that we suppose here that each individual consumes the same amount of resources, independently of its belonging. In our
model, this assumption
is reasonable since
all the individuals belong to the same species. 
Also, the competition for the resources depends only on the number of individuals, independently on their identity.

Furthermore, our model is sufficiently
general to take into account the fact that the growth rate of the populations can be possibly different. In practice,
this possible difference
could be the outcome of
a cultural distinction, or it may be also due to some slight genetic differentiation, as it happened for Homo Sapiens and Neanderthal,
see~\cite{flores}.

Let us call~$r_u$ and~$r_v$ the fertility of the first and second populations respectively. The contribution
to the population growth rate  is  given by
$$ f(u,v) := r_u \left(1-\frac{u+v}{k}   \right),~$$ and these effects
can be comprised in a typical Lotka-Volterra system.

Instead, in our model, we also take into
account the possible death rate due to casualties.
In this way,
we obtain a term such as~$-a\zeta_u$ to be added to~$f(u,v)$. 
The fertility losses give another term~$-ac_u$ for the first population.
We also perform the same analysis for the second population, with the appropriate coefficients.

With these considerations, the system of
the equations that we obtain is
\begin{equation}\label{ch2model1}
\left\{
\begin{array}{llr}
\dot{u}&= r_u u\left(1- \dfrac{u+v}{k} \right) - a(c_u + \zeta_u)u, & t>0,\\
\dot{v}&=r_v v\left(1- \dfrac{v+u}{k} \right) - a(c_v + \zeta_v)u, & t>0.
\end{array}
\right.
\end{equation}
As usual in these kinds of models, we can rescale the variables and the coefficients  in order to
find an equivalent model with fewer parameters. Hence, we perform the changes of variables
\begin{equation}\label{ch2changeofvar}\begin{split}
&\tilde{u}(\tilde t)= \dfrac{u(t)}{k}, \quad \tilde{v}(\tilde t)=\dfrac{v(t)}{k}, 
\quad {\mbox{ where }}\tilde{t}= r_u t,\\
&
\tilde{a}= \dfrac{a(c_v+\zeta_v)}{r_u}, \quad \tilde{c}= \dfrac{c_u+\zeta_u}{c_v+\zeta_v} \quad {\mbox{ and }}\quad \rho= \frac{r_v }{r_u },
\end{split}
\end{equation}
and, dropping the tildas for the sake of readability, we finally get the system
in~\eqref{ch2model}.
We will also refer to it as the civil war model (CW).

{F}rom the change of variables in~\eqref{ch2changeofvar}, we notice in particular that~$a$ may now take values in~$[0,+\infty)$.

\medskip

The competitive Lotka-Volterra system is already used to study some market phenomena as technology substitution, see e.g.~\cite{substitution, bhargava1989generalized, watanabe2004substitution}, and our model aims at adding new features to such models.

Concretely, in the technological competition model, one can think that~$u$ and~$v$
represent the capitals of two computer companies. 
In this setting, to start with, one can suppose that the first company produces a very successful product, say computers with a certain operating system, in an infinite market, reinvesting a proportion~$r_u$ of the profits into the production of additional items, which are purchased by the market, and so on: in this way, one obtains a linear equation of the type~$\dot u=r_u u$, with exponentially growing solutions. 
The case in which the market is not infinite, but reaches a saturation threshold~$k$, would correspond to the equation
$$\dot u=r_u u\left(1-\frac{u}{k}\right).$$ 
Then, when a second computer company comes into the business, selling computers with a different operating system to the same market, one obtains the competitive system of equations
$$ \begin{cases}
\dot u=r_u u\displaystyle\left(1-\frac{u+v}{k}\right),\\
\dot v=r_v v\displaystyle\left(1-\frac{v+u}{k}\right).
\end{cases}$$
At this stage, the first company may decide to use an ``aggressive'' strategy consisting in 
spreading a virus attacking the other company's operating system,
with the aim of setting the other company out of the market (once the competition of the second company is removed, the first company can then exploit the market in a monopolistic regime).
To model this strategy, one can suppose that the first company 
invests a proportion of its capital in the project and
diffusion of the virus, according to a quantifying parameter~$a_u\ge0$, thus producing the equation
\begin{equation}\label{ch2MAR1} \dot u=r_u u\left(1-\frac{u+v}{k}\right)-a_u u.\end{equation}
This directly impacts the capital
of the second company proportionally to the virus
spread,
since the second company has to spend money to project and release antiviruses,
as well as to repay unsatisfied customers,
hence resulting in a second equation of the form
\begin{equation}\label{ch2MAR2} \dot v=r_v v\left(1-\frac{v+u}{k}\right)-a_v u.\end{equation}
The case~$a_u=a_v$ would correspond to an ``even'' effect in which the costs of producing the virus is in balance with
the damages that it causes.
It is also realistic to take into account the case~$a_u<a_v$
(e.g., the first company manages to
produce and diffuse the virus at low cost,
with high impact on the functionality of the operating
system of the second company) as well as the case~$a_u>a_v$ (e.g., the cost of producing and diffusing
the virus is high with respect to the damages caused).

We remark that equations~\eqref{ch2MAR1} and~\eqref{ch2MAR2} can be set into the form~\eqref{ch2model1}, thus showing the interesting versatility of our model also
in financial mathematics.


\subsection{Some notation and basic results on the dynamics of system~\eqref{ch2model}}\label{ch2ss:notation}

We denote by~$(u(t), v(t))$ a solution of~\eqref{ch2model} starting from a point~$(u(0),v(0))\in [0,1] \times [0,1]$.
We will also refer to the \textit{orbit} of~$(u(0), v(0))$  as the collection of
points~$(u(t), v(t))$ for~$t\in \R$, thus both positive and negative times, while the \textit{trajectory} is the collection of points~$(u(t), v(t))$ for~$t\geq0$.

As already mentioned in the discussion below formula~\eqref{ch2model},
$v$ can reach the value~$0$ and even negative values in finite time. However, we will suppose that the dynamics stops when the value~$v=0$ is reached for the first time.
At this point, the conflict ends with the victory of the first population~$u$, that can continue its evolution with a classical Lotka-Volterra equation of the form
\begin{equation*}
\dot{u}= u (1- u)
\end{equation*}
and that would certainly fall into the attractive equilibrium~$u=1$.
The only other possibility is that the solutions are constrained in the set~$[0,1]\times(0,1]$.

In order to state our first result on the dynamics of the system~\eqref{ch2model},
we first observe that, in a real-world situation, the value of~$a$ would probably be non-constant and discontinuous, so we allow this coefficient to take values in the class~$\mathcal{A}$
defined as follows:
\begin{equation}\begin{split}\label{ch2DEFA}
\mathcal{A}&\; :=
\big\{a: [0, +\infty) \to [0, +\infty) {\mbox{ s.t.~$a$ is continuous}}\\
&\qquad \qquad {\mbox{except at most at a finite number of points}}\big\}.\end{split}\end{equation}
A \emph{solution related to a strategy~$a(t)\in \mathcal{A}$} is a pair~$(u(t), v(t)) \in C_0 (0,+\infty)\times C_0 (0,+\infty)$, which is~$C^1$ outside the discontinuous points of~$a(t)$ and
solves system~\eqref{ch2model}.
Moreover, once the initial datum is imposed, the solution is assumed to be
continuous at~$t=0$. 

In this setting, we establish the existence of the solutions of problem~\eqref{ch2model}
and we classify their behavior with respect to the possible exit from the domain~$[0,1]\times[0,1]$:

\begin{proposition}\label{ch2prop:dyn}
	Let~$a(t)\in\mathcal{A}$.
	Given~$(u(0), v(0))\in [0,1] \times [0,1]$, there exists
	a solution~$(u(t),v(t))$ with~$a=a(t)$
	of system~\eqref{ch2model} starting at~$(u(0), v(0))$.
	
	Furthermore, one of the two following situations occurs:
	\begin{enumerate}
		\item[(1)] The solution~$(u(t), v(t))$ issued from~$(u(0), v(0))$ belongs to~$ [0,1]\times (0,1]$ for all~$t\geq 0$.
		\item[(2)] There exists~$T\geq0$ such that the solution~$(u(t), v(t))$ issued from~$(u(0), v(0))$ exists unique for all~$t\leq T$, and~$v(T)=0$ and~$u(T)>0$.
	\end{enumerate} 
\end{proposition}

As a consequence  of Proposition~\ref{ch2prop:dyn}, we can define the 
the \textit{stopping time} of the solution~$(u(t), v(t))$ as
\begin{equation}\label{ch2def:T_s}
T_s (u(0), v(0)) = 
\left\{ 
\begin{array}{ll}
+\infty & \text{if situation (1) occurs}, \\
T & \text{if situation (2) occurs}.
\end{array}
\right.
\end{equation} 
{F}rom now on, we will implicitly  consider solutions~$(u(t),v(t))$ only for~$t\leq T_s(u(0), v(0))$.

Now we are going to analyze the dynamics of~\eqref{ch2model} with a particular focus on possible strategies. To do this, we now define the \textit{basins of attraction}.
The first one is the basin of attraction of the point~$(0,1)$, that is
\begin{equation}\label{ch2DEFB}\begin{split}
\mathcal{B}&\;:= \Big\{ (u(0),v(0))\in [0,1]\times[0,1] \;{\mbox{ s.t. }}\;\\
&\qquad\qquad T_s (u(0), v(0)) = +\infty, \ (u(t),v(t)) \overset{t\to\infty}{\longrightarrow} (0,1)  \Big\},
\end{split}\end{equation}
namely the set of the initial points for which the first population gets extinct (in infinite time) and the second one survives.
The other one is
\begin{equation}\label{ch2DEFE}
\mathcal{E}:= \left\{ (u(0),v(0))\in ([0,1]\times[0,1])\setminus(0,0) \;{\mbox{ s.t. }}\;  T_s(u(0),v(0))< + \infty \right\},
\end{equation}
namely
the set of initial points for which we have the victory of the first population and the extinction of the second one. 

Of course, the sets~$\mathcal{B}$ and~$\mathcal{E}$ depend on the parameters~$a$,~$c$, and~$\rho$; we will express this dependence by writing~$\mathcal{B}(a,c,\rho)$
and~$\mathcal{E}(a,c,\rho)$ when it is needed, and omit it otherwise for the sake of readability. The dependence on parameters will be carefully studied
in Subsection~\ref{ch2ss:dependence}.

\subsection{Dynamics of system~\eqref{ch2model} for constant strategies}\label{ch2ss:dynamics}

The first step towards the understanding of the dynamics of the system
in~\eqref{ch2model} is
is to analyze the behavior of the system for constant coefficients. 

To this end, we introduce some notation.
Following the terminology on pages~9-10 in~\cite{MR1056699},
we say that an equilibrium point (or fixed point) of the dynamics
is a (hyperbolic) {\em sink}
if all the eigenvalues of the linearized map have strictly
negative real parts, a (hyperbolic) {\em source}
if all the eigenvalues of the linearized map have strictly
positive real parts, and a (hyperbolic) {\em saddle}
if some of the eigenvalues of the linearized map have strictly
positive real parts
and some have negative real parts
(since in this chapter we work in dimension~$2$,
saddles correspond to linearized maps with one
eigenvalue with
strictly positive real part
and one eigenvalue with
strictly negative real part).
We also recall that
sinks are asymptotically stable (and sources are
asymptotically stable for the reversed-time dynamics), see e.g. Theorem 1.1.1
in~\cite{MR1056699}.

With this terminology, we state the following theorem:

\begin{theorem}[Dynamics of system~\eqref{ch2model}] \label{ch2thm:dyn}
	For~$a > 0$ and~${\rho}> 0$ the system~\eqref{ch2model} has the following features:
	\begin{itemize}
		\item[(i)] When~$0<ac<1$, the system has 3 equilibria:~$(0,0)$ is a source,~$(0,1)$ is 
		a sink, and 
		\begin{equation}\label{ch2usvs}
		(u_s, v_s):= \left( \frac{1-ac}{1+{\rho}c} {\rho}c, \frac{1-ac}{1+{\rho}c} \right) \in (0,1)\times (0,1)
		\end{equation}
		is a saddle.
		
		\item[(ii)] When~$ac>1$, the system has 2 equilibria:~$(0,1)$ is a sink and~$(0,0)$ is a saddle.  
		\item[(iii)] When~$ac=1$, the system has 2 equilibria:~$(0,1)$ is a sink and~$(0,0)$
		corresponds to a strictly positive eigenvalue and a null one.
		\item[(iv)] We have 
		\begin{equation} \label{ch2fml:division}
		[0,1]\times [0,1] = \mathcal{B} \cup \mathcal{E} \cup \mathcal{M}
		\end{equation}
		where~$\mathcal{B}~$ and~$\mathcal{E}$ are defined in~\eqref{ch2DEFB}
		and~\eqref{ch2DEFE}, respectively, and~$\mathcal{M}$ is a smooth curve.
		\item[(v)] The trajectories starting in~$\mathcal{M}$ tend to~$(u_s,v_s)$ if~$0<ac<1$,
		and to~$(0,0)$ if~$ac\ge1$ as~$t$ goes to~$+\infty$.
	\end{itemize}
\end{theorem}

More precisely, one can say that
the curve~$\mathcal{M}$ in Theorem~\ref{ch2thm:dyn} is the stable manifold of the saddle
point~$(u_s,v_s)$ when~$0<ac<1$, and of the 
saddle point~$(0,0)$ when~$ac>1$. The case~$ac=1$ needs a special treatment,
due to the degeneracy of one eigenvalue, and in this case the curve~$\mathcal{M}$
corresponds to the center manifold of~$(0,0)$, and an ad-hoc argument will
be exploited
to show that also in this degenerate case orbits that start in~$\mathcal{M}$
are asymptotic in the future to~$(0,0)$.

As a matter of fact,~$\mathcal{M}$
acts as a dividing wall between the two basins of attraction, as described in~(iv)
of Theorem~\ref{ch2thm:dyn} and in the forthcoming Proposition~\ref{ch2prop:char}.

Moreover, in the forthcoming
Propositions~\ref{ch2lemma:M}
and~\ref{ch2M:p045}
we will show that~$\mathcal{M}$ can be written as the graph of a function. This is particularly useful because, by studying the properties of this function, we gain relevant pieces of information on the sets~$\mathcal{B}$ and~$\mathcal{E}$
in~\eqref{ch2DEFB} and~\eqref{ch2DEFE}.

We point out that in Theorem~\ref{ch2thm:dyn}
we find that the set of initial data~$[0,1]\times[0,1]$ splits into three part:
the set~$\mathcal{E}$, given in~\eqref{ch2DEFE}, made
of points going to the extinction of the second population in finite time; the set~$\mathcal{B}$, given in~\eqref{ch2DEFB}, which is the
basin of attraction of the equilibrium~$(0,1)$;
the set~$\mathcal{M}$, which is
a manifold of dimension~$1$ that separates~$\mathcal{B}$ from~$\mathcal{E}$. 

In particular,
Theorem~\ref{ch2thm:dyn} shows that, also for our model, the Gause principle of exclusion is respected; that is, in general, two competing populations cannot coexist in the same territory, see e.g.~\cite{fath2018encyclopedia}. 

\begin{figure}[h]
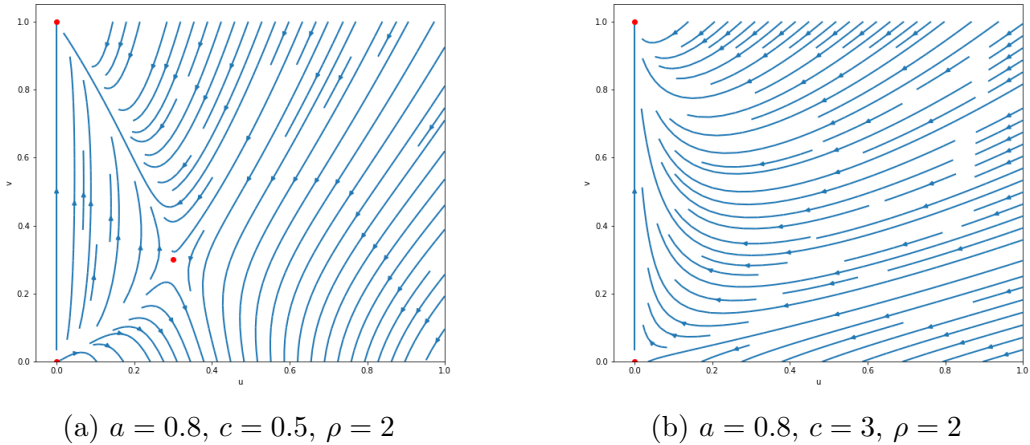
 
	\begin{subfigure}{.5\textwidth}
		\centering
		\includegraphics[width=.8\textwidth]{acless1dots}
		\caption{$a=0.8$,~$c=0.5$,~$\rho=2$}
	\end{subfigure}%
	\begin{subfigure}{.5\textwidth}
		\centering
		\includegraphics[width=.8\textwidth]{acgr1dots}
		\caption{$a=0.8$,~$c=3$,~$\rho=2$}
	\end{subfigure}
	\caption{\it The figures show a phase portrait for the indicated values of the coefficients. In blue, the orbits of the points. The red dots represent the equilibria.}
	\label{ch2fig:dyn}
\end{figure}

One peculiar feature of our system is that, if the aggressiveness is too strong, the equilibrium~$(0,0)$ changes its ``stability'' properties, passing from a source (as in (i) of
Theorem~\ref{ch2thm:dyn})
to a saddle point (as in (ii) of
Theorem~\ref{ch2thm:dyn}). This shows that the war may have self-destructive outcomes, therefore it is important for the first population to analyze the situation in order to choose a proper level of aggressiveness. 
Figure~\ref{ch2fig:dyn} shows one example of dynamics for each case.

%


\subsection{Dynamics of system~\eqref{ch2model} for variable strategies
	and optimal strategies for the first population}\label{ch2ss:strategy}

We now deal with the problem of choosing the strategy~$a$ such that the first population wins, that is a problem of \emph{target reachability} for a control-affine system. As we will see, the problem is not \emph{controllable}, meaning that, starting from a given initial point, it is not always possible to reach a given target.
\medskip

We now introduce some terminology, that we will use throughout the chapter.
Recalling~\eqref{ch2DEFA},
for any~$\mathcal{T}\subseteq \mathcal{A}$, we set
\begin{equation}\label{ch2DEFNU}
\mathcal{V}_{\mathcal{T}}:= \underset{a(\cdot)\in \mathcal{T}}{\bigcup} \mathcal{E}(a(\cdot)),
\end{equation}
where~$\mathcal{E}(a(\cdot))$ denotes the set of initial data~$(u_0,v_0)$
such that~$T_s(u_0,v_0)< +\infty$, when the coefficient~$a$ in~\eqref{ch2model} is replaced by the function~$a(t)$.

Namely,~$\mathcal{V}_{\mathcal{T}}$ represents the set of initial conditions for which~$u$ is able to win by choosing a suitable strategy in~$\mathcal{T}$; we call~$\mathcal{V}_{\mathcal{T}}$ the \emph{victory set} with admissible strategies in~$\mathcal{T}$.
We also say that~$a(\cdot)$ is a \emph{winning strategy} for the point~$(u_0,v_0)$
if~$(u_0,v_0)\in \mathcal{E}(a(\cdot) )$.

Moreover, we will call
\begin{equation}\label{ch2u0v0}
(u_s^0, v_s^0):= \left(\frac{\rho c}{1+\rho c}, \frac{1}{1+\rho c}\right).
\end{equation}
Notice that~$(u_s^0, v_s^0)$ is the limit point as~$a$ tends to~$0$ of the sequence of saddle points~$\{(u_s^a, v_s^a)\}_{a>0}$
defined in~\eqref{ch2usvs}.
\medskip

With this notation,
the first question that we address is for which initial configurations it is possible for the population~$u$
to have a winning strategy, that is, to characterize the victory set. For this, we allow the strategy to take all the values in~$[0, +\infty)$.
In this setting, we have the following result:

\begin{theorem}\label{ch2thm:Vbound}
	\begin{itemize}
		\item[(i)] For~$\rho=1$, we have that  
		\begin{equation}\label{ch2Vbound1}\begin{split}
		\mathcal{V}_{\mathcal{A}} = \,&\Big\{ (u,v)\in[0,1] \times [0,1] \;
		{\mbox{ s.t. }}\;  v-\frac{u}{c}<0 \; {\mbox{ if }} u\in[0,c]\\
		&\qquad\qquad\qquad {\mbox{ and }}\; v\le1 \; {\mbox{ if }} u\in(c,1]\Big\},
		\end{split}\end{equation}
		with the convention that the last line in~\eqref{ch2Vbound1} is not present if~$c\ge1$.
		\item[(ii)] 
		For~$\rho<1$, we have that
		\begin{equation}\label{ch2bound:rho<1}
		\begin{split}
		\mathcal{V}_{\mathcal{A}} &\;= \Bigg\{ (u,v)\in[0,1] \times [0,1] \;{\mbox{ s.t. }}\;
		v< \gamma_0(u) \ \text{if} \ u\in [0, u_s^0], \\ 
		&\qquad\qquad\qquad\qquad
		v< \frac{u}{c} + \frac{1-\rho}{1+\rho c} \ \text{if} \ u\in \left[u_s^0, 
		\frac{\rho c(c+1)}{1+\rho c}\right]\\
		&\qquad\qquad\qquad\qquad
		{\mbox{and }}\; v\le1\ \text{if} \ u\in \left(
		\frac{\rho c(c+1)}{1+\rho c},1\right]
		\Bigg\},
		\end{split}
		\end{equation}
		where 
		\begin{equation*}
		\gamma_0(u):= \frac{u^{\rho}}{\rho c(u_s^0)^{\rho-1}},
		\end{equation*}
		and we use the convention that the last line
		in~\eqref{ch2bound:rho<1} is not present if~$ \frac{\rho c(c+1)}{1+\rho c}\ge1$.
		\item[(iii)] For~$\rho>1$, we have that
		\begin{equation}\label{ch2bound:rho>1}
		\begin{split}
		\mathcal{V}_{\mathcal{A}} &\;= \Bigg\{ (u,v)\in[0,1] \times [0,1]\;
		{\mbox{ s.t. }}\; v< \frac{u}{c} \ \text{if} \ u\in [0, u_{\infty}],\\&\qquad
		\qquad\qquad\qquad 
		v< \zeta(u)  \ \text{if} \ u\in\left(u_{\infty}, \frac{c}{(c+1)^{\frac{\rho-1}\rho}}\right] \\&\qquad
		\qquad\qquad\qquad 
		{\mbox{and }}\; v\le 1
		\ \text{if} \ u\in\left(\frac{c}{(c+1)^{\frac{\rho-1}\rho}},1\right] 
		\Bigg\},
		\end{split}
		\end{equation}
		where
		\begin{equation}\label{ch2ZETADEF}
		u_{\infty}:= \frac{c}{c+1}
		\quad {\mbox{ and }}\quad \zeta (u):= \frac{u^{\rho}}{c \, u_{\infty}^{\rho-1}} .
		\end{equation}   
		and we use the convention that the last line
		in~\eqref{ch2bound:rho>1} is not present if~$  \frac{c}{(c+1)^{\frac{\rho-1}\rho}}\ge1$.
	\end{itemize}
\end{theorem}

In practice, 
constant strategies could be certainly easier to implement and
it is therefore natural to investigate whether or not
it suffices to restrict to constant strategies
without altering the possibility of victory.
The next result addresses this problem by showing that when~$\rho=1$
constant strategies are as good as all strategies,
but instead when $\rho\ne 1$ victory cannot be achieved by only
exploiting constant strategies:

\begin{theorem}\label{ch2thm:W}
	Let $\mathcal{K}\subset \mathcal{A}$ be the set of constant functions. Then the following holds:
	\begin{itemize}
		\item[(i)] For~$\rho= 1$, we have that~$ \mathcal{V}_{\mathcal{A}}=\mathcal{V}_{\mathcal{K}}=\mathcal{E}(a)$ for all $a>0$;
		\item[(ii)] For~$\rho\neq 1$, we have that~$\mathcal{V}_{\mathcal{K}} \subsetneq \mathcal{V}_{\mathcal{A}}$.
	\end{itemize}	
\end{theorem} 

The result of Theorem~\ref{ch2thm:W}, part~(i),
reveals a special rigidity of the case~$\rho=1$
in which, no matter which strategy~$u$ chooses,  the victory depends only on the initial conditions, but it is independent of the strategy~$a(t)$.
Instead, as stated in
Theorem~\ref{ch2thm:W}, part~(ii),
for~$\rho \neq 1$ the choice of~$a(t)$ plays a crucial role in determining which population is going to win and constant strategies do not exhaust all the
possible winning strategies.
We stress that~$\rho=1$ plays also
a special role in the biological interpretation of the model, since in this case the two
populations have the same fit to the environmental resource, and hence, in a sense,
they are indistinguishable, up to the possible aggressive behavior of the first population.

Next, we show that the set~$\mathcal{V}_{\mathcal{A}}$ can be recovered if we use piecewise constant functions with at most one discontinuity, that we call Heaviside functions. 

\begin{theorem}\label{ch2thm:H}
	There holds that~$\mathcal{V}_{\mathcal{A}} = \mathcal{V}_{\mathcal{H}}$, where~$\mathcal{H}$ is the set of Heaviside functions.
\end{theorem}

The proof of Theorem~\ref{ch2thm:H} solves also the third question mentioned in the Introduction. As a matter of fact, it proves that for each point we either have a constant winning strategy or
a winning strategy of type 
\begin{equation*}
a(t) = \left\{
\begin{array}{lr}
a_1  &{\mbox{ if }} t<T ,\\
a_2  &{\mbox{ if }} t\geq T,
\end{array}
\right.
\end{equation*}
for some~$T\in(0,T_s)$, and
for suitable values~$a_1$,~$a_2 \in (0,+\infty)$ such that one is very small and the other one very large, the order depending on~$\rho$. 
The construction that
we give also puts in light the fact that the choice of the strategy depends on the initial datum, answering also our fourth question. 

It is interesting to observe that the winning strategy that switches abruptly from a small to a large value
could be considered, in the optimal control terminology, as a ``bang-bang'' strategy.
Even in a target reachability problem, the structure predicted by Pontryagin's Maximum Principle is brought in light: the bounds of the set~$\mathcal{V}_{\mathcal{A}}$, as
given in Theorem~\ref{ch2thm:Vbound}, depend on the bounds that
we impose on the strategy, that are,~$a \in[0,+\infty)$.

It is natural to consider also the case
in which the level of aggressiveness 
is constrained between a minimal and maximal threshold,
which corresponds to the setting~$a\in[m,M]$ for suitable~$M\geq m\geq 0$, with the hypothesis that $M>0$.
In this setting, we denote by~$\mathcal{A}_{m,M}$ the class of piecewise continuous strategies~$a(\cdot)$
in~${\mathcal{A}}$ such that~$
m\leq a(t)\leq M$ for all~$t>0$ and we let
\begin{equation}\label{ch2SPE}
\mathcal{V}_{m,M}:=\mathcal{V}_{\mathcal{A}_{m,M}}=\underset{{a(\cdot)\in \mathcal{A}}\atop{m\leq a(t)\leq M}
}{\bigcup} \mathcal{E}(a(\cdot))=
\underset{{a(\cdot)\in \mathcal{A}}_{m,M}
}{\bigcup} \mathcal{E}(a(\cdot)).\end{equation}
Then we have the following:

\begin{theorem}\label{ch2thm:limit}
	Let $M$ and $m$ be two real numbers such that $M\geq m\geq 0$. Then, for $\rho\neq 1$ we have the strict inclusion  $\mathcal{V}_{{m,M}}\subsetneq \mathcal{V}_{\mathcal{A}}$.
\end{theorem}

Notice that for $\rho=1$, Theorem \ref{ch2thm:W} gives instead that $\mathcal{V}_{{m,M}}= \mathcal{V}_{\mathcal{A}}$
and we think that this is a nice feature, outlining a special role played by the parameter~$\rho$
(roughly speaking, when~$\rho=1$ constant strategies suffice
to detect all possible winning configurations, thanks to 
Theorem \ref{ch2thm:W}, while when~$\rho\ne1$ non-constant strategies are necessary to detect
all winning configurations).

\subsubsection{Time minimizing strategy}

Once established that it is possible to win starting in a certain initial condition, we are interested in knowing which of the possible strategies is best to choose. One condition that may be taken into account is the duration of the war. Now, this question can be written as a minimization problem with a proper functional to minimize and therefore the classical Pontryagin theory applies. 

To state our next result, 
we recall the setting in~\eqref{ch2SPE} and define
\begin{equation*}
\mathcal{S}(u_0, v_0) := \Big\{ a(\cdot)\in \mathcal{A}_{m,M}
\;\mbox{ s.t. }\; (u_0, v_0) \in \mathcal{E}(a(\cdot))  \Big\},
\end{equation*}
that is the set of all bounded strategies for which the trajectory starting at~$(u_0, v_0)$ leads to the victory of the first population.
To each~$a(\cdot)\in\mathcal{S}(u_0, v_0)$ we associate the stopping time defined in~\eqref{ch2def:T_s}, and we express its dependence on~$a(\cdot)$ by writing~$T_s(a(\cdot))$.
In this setting, we provide the following statement concerning the strategy leading
to the quickest possible victory for the first population:

\begin{theorem}\label{ch2thm:min}
	Given a point~$(u_0, v_0)\in \mathcal{V}_{m,M}$, there exists a winning strategy~$\tilde{a}(t)\in
	\mathcal{S}(u_0, v_0)$, and a trajectory~$(\tilde{u}(t), \tilde{v}(t) )$ associated with~$\tilde{a}(t)$,
	for~$t\in[0,T]$,
	with~$(\tilde{u}(0), \tilde{v}(0) )=(u_0,v_0)$, where~$T$ is given by
	\begin{equation*}
	T = \underset{a(\cdot)\in\mathcal{S}}{\min} T_s(a(\cdot)).
	\end{equation*}
	Moreover,
	\begin{equation*}
	\tilde{a}(t)\in \left\{m, \ M, \    a_s(t) \right\},
	\end{equation*}	  
	where
	\begin{equation}\label{ch2KSM94rt3rjjjdfe}
	{a}_s(t) := \dfrac{(1-\tilde{u}(t)-\tilde{v}(t))[\tilde{u}(t) \, (2c+1-\rho c)+\rho c]}{\tilde{u}(t) \, 2c(c+1)}.
	\end{equation}	
\end{theorem}

\begin{figure}[h] 
	\begin{subfigure}{.5\textwidth}
		\centering
		\includegraphics[width=.9\textwidth]{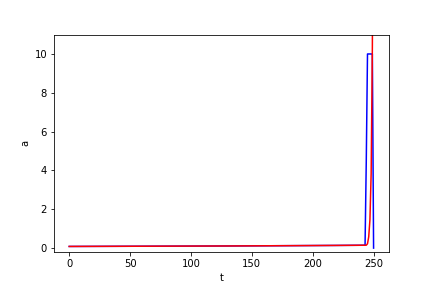}
	\end{subfigure}%
	\begin{subfigure}{.5\textwidth}
		\centering
		\includegraphics[width=.9\textwidth]{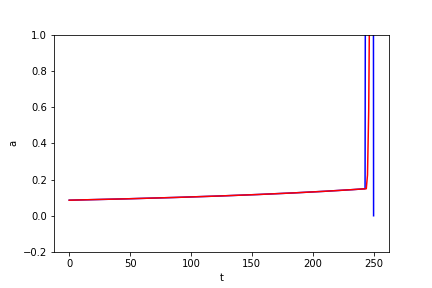}
	\end{subfigure}
	\caption{\it The figure shows the result of a numerical simulation searching a minimizing time strategy~$\tilde{a}(t)$ for the problem starting in~$(0.5, 0.1875)$ for the
		parameters~$\rho=0.5$,~$c=4.0$,~$m=0$ and~$M=10$. In blue, the value
		found for~$\tilde{a}(t)$; in red, the value of~$a_s(t)$ for the corresponding trajectory~$(u(t), v(t))$. As one can observe,~$\tilde{a}(t)\equiv a_s(t)$ in a long trait.
		The simulation was done using AMPL-Ipopt on the server NEOS and pictures have been made with Python. 
	}
	\label{ch2fig:min}
\end{figure}

The surprising fact given by Theorem~\ref{ch2thm:min}
is that the
minimizing strategy is not only of bang-bang type, but it may assume some values along a \emph{singular arc}, given by~$a_s(t)$.
This possibility is realized in some concrete cases, as we verified by running some numerical simulations, whose results can be visualized in Figure~\ref{ch2fig:min}. 

\subsection{Organization of the chapter}

In the forthcoming Section~\ref{ch2IKJM:plrg777}
we will exploit methods from ordinary differential equations and dynamical systems
to describe the equilibria of the system and their possible basins of attraction.
The dependence of the dynamics on the structural parameters,
such as fit to the environment, aggressiveness and efficacy of attacks,
is discussed in detail in Section~\ref{ch2ss:dependence}.

Section~\ref{ch2STRATE} is devoted to the analysis of the strategies
that allow the first population to eradicate the second one (this part needs
an original combination of methods from dynamical systems and optimal control theory).

\section{First results on the dynamics
	and proofs of
	Proposition~\ref{ch2prop:dyn} and Theorem~\ref{ch2thm:dyn}}\label{ch2IKJM:plrg777}

In this section we provide some useful results on the behavior of the solutions
of~\eqref{ch2model} and on the basin of attraction. In particular,
we provide the proofs
of Proposition~\ref{ch2prop:dyn} and Theorem~\ref{ch2thm:dyn}
and
we state a characterization of  the sets~$\mathcal{B}$
and~$\mathcal{E}$ given in~\eqref{ch2DEFB} and~\eqref{ch2DEFE}, respectively, see Propositions~\ref{ch2prop:char}.

This material will be extremely useful for the analysis of the strategy that we operate later. 

\smallskip

We start with some preliminary notation. Given a close set~$\mathcal{S} \subseteq [0,1]\times [0,1]$, we say that a trajectory~$(u(t),v(t))$ originated in~$\mathcal{S}$ \textit{exits} the set~$\mathcal{S}$ at some time~$T\geq 0$ if
\begin{itemize}
	\item~$(u(t),v(t))\in \mathcal{S}$ for~$t\leq T$, 
	\item~$(u(T),v(T))\in \partial \mathcal{S}$,
	\item for any vector ~$\nu$ normal to~$\partial \mathcal{S}$
	at the point~$(u(T),v(T))$, it holds that
	$$(\dot{u}(T), \dot{v}(T)) \cdot \nu >0.$$ 
\end{itemize}

Now, we prove Proposition~\ref{ch2prop:dyn}, which is fundamental to
the well-definition of our model:

\begin{proof}[Proof of Proposition~\ref{ch2prop:dyn}]
	We consider the function~$a(t)\in\mathcal{A}$, which is continuous except in a finite number of points~$0<t_1< \dots< t_n$.
	In all the intervals~$(0, t_1)$,~$(t_i, t_{i+1}]$, for~$i\in\{1,\cdots,n-1\}$, and~$(t_n, +\infty)$, the equations in~\eqref{ch2model}
	have smooth coefficients, and therefore a solution does exist.
	Now, it is sufficient to consider~$(u(t_i), v(t_i))$ as the initial datum for the dynamics
	in~$(t_i, t_{i+1}]$ to construct a solution~$(u(t), v(t))$ for all~$t>0$ satisfying
	system~\eqref{ch2model}.
	This is a rather classical result and we refer to~\cite{dynsyst} for more details. 
	
	Now, we prove that either the possibility in~(1) or the possibility in~(2) can occur.
	For this, by using the equation for~$v$ in~\eqref{ch2model},
	we notice that for~$v=1$ the inward pointing normal derivative  is
	\begin{equation*}
	-\dot{v}|_{v=1}=\left(-\rho v(1-u-v)+au\right)|_{v=1} = u(\rho+a) \geq 0.
	\end{equation*}
	This means that 
	no trajectory can exit~$[0,1]\times[0,1]$ on the edge~$v=1$. Similarly, using
	the equation for~$u$ in~\eqref{ch2model}, we see that
	for~$u=1$ the normal derivative inward pointing is
	\begin{equation*}
	-\dot{u}|_{u=1}=\left(-u(1-u-v)+acu\right)|_{u=1} = v +ac \geq 0,
	\end{equation*}
	and therefore
	no trajectory can exit~$[0,1]\times[0,1]$ on the edge~$u=1$.
	
	Moreover, it is easy to see that
	all points on the line~$u=0$ go to the equilibrium~$(0,1)$, thus trajectories do not cross the line~$u=0$. The only remaining possibilities are that the trajectories stay
	in~$[0,1]\times(0,1]$, that is possibility~(1), or they exit the square on the side~$v=0$,
	that is possibility~(2).
\end{proof}

Now, we give the proof of~(i), (ii) and~(iii) of Theorem~\ref{ch2thm:dyn}.

\begin{proof}[Proof of (i), (ii) and~(iii) of Theorem~\ref{ch2thm:dyn}]
	We first consider equilibria with first coordinate~$u=0$. In this case, 
	from the second equation in~\eqref{ch2model}, we have that
	the equilibria must satisfy~$\rho v(1-v)=0$, thus~$v=0$ or~$v=1$.
	As a consequence,~$(0,0)$ and~$(0,1)$ are two equilibria of the system.
	
	Now, we consider equilibria with first coordinate~$u>0$.
	Equilibria of this form must  satisfy~$\dot{u}=0$ with~$u\neq 0$, and therefore,
	from the first equation in~\eqref{ch2model},
	\begin{equation}\label{ch2curve:u'}
	1-u-v-ac=0.
	\end{equation}
	Moreover from the condition~$\dot{v}=0$ and the second equation in~\eqref{ch2model},
	we see that
	\begin{equation} \label{ch2curve:v'}
	\rho v(1-u-v)-au=0.
	\end{equation}
	Putting together~\eqref{ch2curve:u'} and~\eqref{ch2curve:v'}, we obtain that
	the intersection point must lie on the line~$\rho c v-u=0$.
	Since the equilibrium is at the intersection between two lines, it must be unique.
	One can easily verify that the values given in~\eqref{ch2usvs} satisfy~\eqref{ch2curve:u'}
	and~\eqref{ch2curve:v'}.
	
	{F}rom now on, we distinguish the three situations in~(i), (ii) and~(iii) of 
	Theorem~\ref{ch2thm:dyn}.
	\medskip
	
	{\emph{(i)}} If~$0<ac<1$, we have that the point~$(u_s,v_s)$ given in~\eqref{ch2usvs}
	lies in~$(0,1)\times(0,1)$. As a result, in this case the system has~$3$ equilibria,
	given by~$(0,0)$,~$(0,1)$ and~$(u_s,v_s)$.
	
	Now, we observe that the Jacobian of the system~\eqref{ch2model} is
	\begin{equation} \label{ch2Jmatrix}
	J(u,v)=	\begin{pmatrix}
	1-2u-v -ac & -u \\ 
	-\rho v-a & \rho (1-u-2v)
	\end{pmatrix}.
	\end{equation}
	At the point~$(0,0)$, the matrix has eigenvalues~$\rho >0$ and~$1-ac >0$, thus~$(0,0)$ is a source. At the point~$(0,1)$, the Jacobian~\eqref{ch2Jmatrix} has eigenvalues~$-ac <0$ and~$-\rho <0$, thus~$(0,1)$ is a sink. At the point~$(u_s,v_s)$, by exploiting the relations~\eqref{ch2curve:u'} and~\eqref{ch2curve:v'} we have that
	\begin{equation*}
	J(u_s,v_s)=\begin{pmatrix}
	-u_s  & -u_s \\ 
	-\rho v_s-a & \rho (ac-v_s)
	\end{pmatrix},
	\end{equation*}
	which, by the change of basis given by the matrix
	$$
	\begin{pmatrix}
	-\frac1{u_s}  & 0 \\ 
	-\frac1{u_s}\left[\left(\frac{u_s}{c}+a\right)\left(\frac{\rho c-c}{1+\rho c}\right)+ac	\right] & \frac{\rho c-c}{1+\rho c}
	\end{pmatrix},
	$$
	becomes	
	\begin{equation}\label{ch2degieerhfdj}
	\begin{pmatrix}
	1  & 1 \\ 
	ac & \rho ac
	\end{pmatrix}.	\end{equation}
	The characteristic polynomial of the matrix 
	in~\eqref{ch2degieerhfdj}	
	is~$\lambda^2-\lambda(1+\rho a c)+\rho a c-ac$, that has two real roots, as one can see by inspection. Hence,~$J(u_s, v_s)$ has two real eigenvalues. Moreover,
	the determinant of~$J(u_s, v_s)$ is~$-\rho ac u_s-au_s <0$, which implies
	that~$J(u_s, v_s)$ has one positive and one negative eigenvalues. These considerations
	give that~$(u_s, v_s)$ is a saddle point, as desired. This completes the proof
	of~(i) in Theorem~\ref{ch2thm:dyn}.
	\medskip
	
	{\emph{(ii) and (iii)}} We assume that~$ac\ge1$.
	We observe that the equilibrium described by the
	coordinates~$(u_s,v_s)$ in~\eqref{ch2usvs}
	coincides with~$(0,0)$ for~$ac=1$,
	and lies outside~$[0,1]\times[0,1]$ for~$ac>1$. 
	As a result, when~$ac\ge1$ the system has~$2$ equilibria, given by~$(0,0)$ and~$(0,1)$.
	
	Looking at the Jacobian in~\eqref{ch2Jmatrix},
	one sees that
	at the point~$(0,1)$, it has eigenvalues~$-ac <0$ and~$-\rho <0$,
	and therefore~$(0,1)$ is a sink when~$ac\ge1$.
	
	Furthermore, from~\eqref{ch2Jmatrix}
	one finds that
	if~$ac>1$ then~$J(0,0)$ has the positive eigenvalue~$\rho$ and the negative
	eigenvalue~$1-ac$, thus~$(0,0)$ is a saddle point.
	
	If instead~$ac=1$, then~$J(0,0)$
	has one positive eigenvalue and one null eigenvalue, as desired.
\end{proof}

To complete the proof of Theorem~\ref{ch2thm:dyn},
we will deal with the cases~$ac\neq1$ and~$ac=1$ separately. This analysis will be performed
in the forthcoming Sections~\ref{ch2sec:deg0} and~\ref{ch2sec:deg}.
The completion of the proof of
Theorem~\ref{ch2thm:dyn} will then be
given in Section~\ref{ch2sec:deg2}.

\subsection{Characterization of~$\mathcal{M}$ when~$ac\ne1$}
\label{ch2sec:deg0}

We consider here the case~$ac\ne1$. The case~$ac=1$ is degenerate
and it will be treated separately in Section~\ref{ch2sec:deg}.

We point out that
in the proof of~(i) and~(ii) in Theorem~\ref{ch2thm:dyn} we found a saddle point 
in both cases.
By the Stable Manifold Theorem (see for example~\cite{dynsyst}), the point~$(u_s, v_s)$ in~\eqref{ch2usvs} in the case~$0<ac<1$ and the
point~$(0,0)$ in the case~$ac> 1$ have a stable manifold and an unstable manifold.
These manifolds are unique,
they have dimension~$1$, and they are tangent to the eigenvectors of the linearized system.
We will denote by~$\mathcal{M}$ the stable manifold associated
with these saddle points.
Since we are interested in the dynamics in the square~$[0,1]\times[0,1]$, with a slight abuse of notation we will only consider the restriction of ~$\mathcal{M}$ in~$[0,1]\times[0,1]$.

In order to complete the
proof of Theorem~\ref{ch2thm:dyn}, we now analyze
some properties of~$\mathcal{M}$:

\begin{proposition} \label{ch2lemma:M}
	For $ac\ne1$ the set~$\mathcal{M}$ can be written as the graph of a unique increasing~${C}^2$ function~$\gamma:[0,u_{\mathcal{M}}] \to [0, v_{\mathcal{M}}]$ for some~$(u_{\mathcal{M}}, v_{\mathcal{M}}) \in
	\big(\{1\}\times[0,1]\big)\cup
	\big((0,1]\times\{1\}\big)$, such that~$\gamma(0)=0$,~$\gamma(u_{\mathcal{M}})=v_{\mathcal{M}}$ and
	\begin{itemize}
		\item if~$0<ac<1$,~$\gamma(u_s)=v_s$;
		\item if~$ac> 1$, in~$u=0$ the function~$\gamma$ is tangent to the
		line~$(\rho-1+ac)v-au=0$.  
	\end{itemize}
\end{proposition}

As a byproduct of the proof of Proposition~\ref{ch2lemma:M}, we also
obtain some useful information on the structure of the stable manifold and
the basins of attraction, that we summarize here below:

\begin{corollary}\label{ch2lemma:M1}
	Suppose that~$0<ac<1$. Then,
	the curves~\eqref{ch2curve:u'} and~\eqref{ch2curve:v'}, loci of the points such that~$\dot{u}=0$ and~$\dot{v}=0$ respectively, divide the square~$[0,1]\times[0,1]$ into four regions:
	\begin{equation}\begin{split}\label{ch2DEFA1234}
	\mathcal{A}_1 &\;:= \big\{ (u, v) \in [0,1]\times[0,1] \;{\mbox{ s.t }}\; \dot{u}\leq 0,\; \dot{v}\geq 0 \big\}, \\
	\mathcal{A}_2 &\;:= \big\{ (u, v) \in [0,1]\times[0,1] \;{\mbox{ s.t }}\; \dot{u}\leq 0,\; \dot{v}\leq 0 \big\}, \\
	\mathcal{A}_3 &\;:= \big\{ (u, v) \in [0,1]\times[0,1] \;{\mbox{ s.t }}\;\dot{u}\geq 0,\; \dot{v}\leq 0 \big\}, \\
	\mathcal{A}_4 &\;:= \big\{ (u, v) \in [0,1]\times[0,1] \;{\mbox{ s.t }}\;\dot{u}\geq 0, \;\dot{v}\geq 0 \big\}.
	\end{split}\end{equation}
	
	Furthermore, the sets~$\mathcal{A}_1\cup \mathcal{A}_4$
	and~$\mathcal{A}_2\cup\mathcal{A}_3$ are separated by the curve~$\dot{v}=0$, given by the graph of the continuous function
	\begin{equation}\label{ch2f:sigma}
	\sigma(v):= 1- \frac{\rho v^2+a}{\rho v+a},
	\end{equation}
	that satisfies~$\sigma(0)=0$,~$\sigma(1)=0$, and~$0<\sigma(v)<1$ for all~$v\in (0,1)$.
	
	In addition,
	\begin{equation}\label{ch2aggiunto}
	{\mbox{$\mathcal{M}\setminus \{(u_s,v_s) \}$
			is contained in~$\mathcal{A}_2\cup\mathcal{A}_4$,}}
	\end{equation}
	\begin{equation}\label{ch2primaBIS}
	(\mathcal{A}_3 \setminus \{ (0,0), (u_s, v_s) \} ) \subseteq \mathcal{E},
	\end{equation}
	and
	\begin{equation}\label{ch2prima2BIS}
	\mathcal{A}_1\setminus \{(u_s,v_s) \} \subset \mathcal{B},\end{equation}
	where the notation in~\eqref{ch2DEFB} and~\eqref{ch2DEFE} has been utilized.
\end{corollary}

\begin{figure}
	\centering
	\includegraphics[width=.4\textwidth]{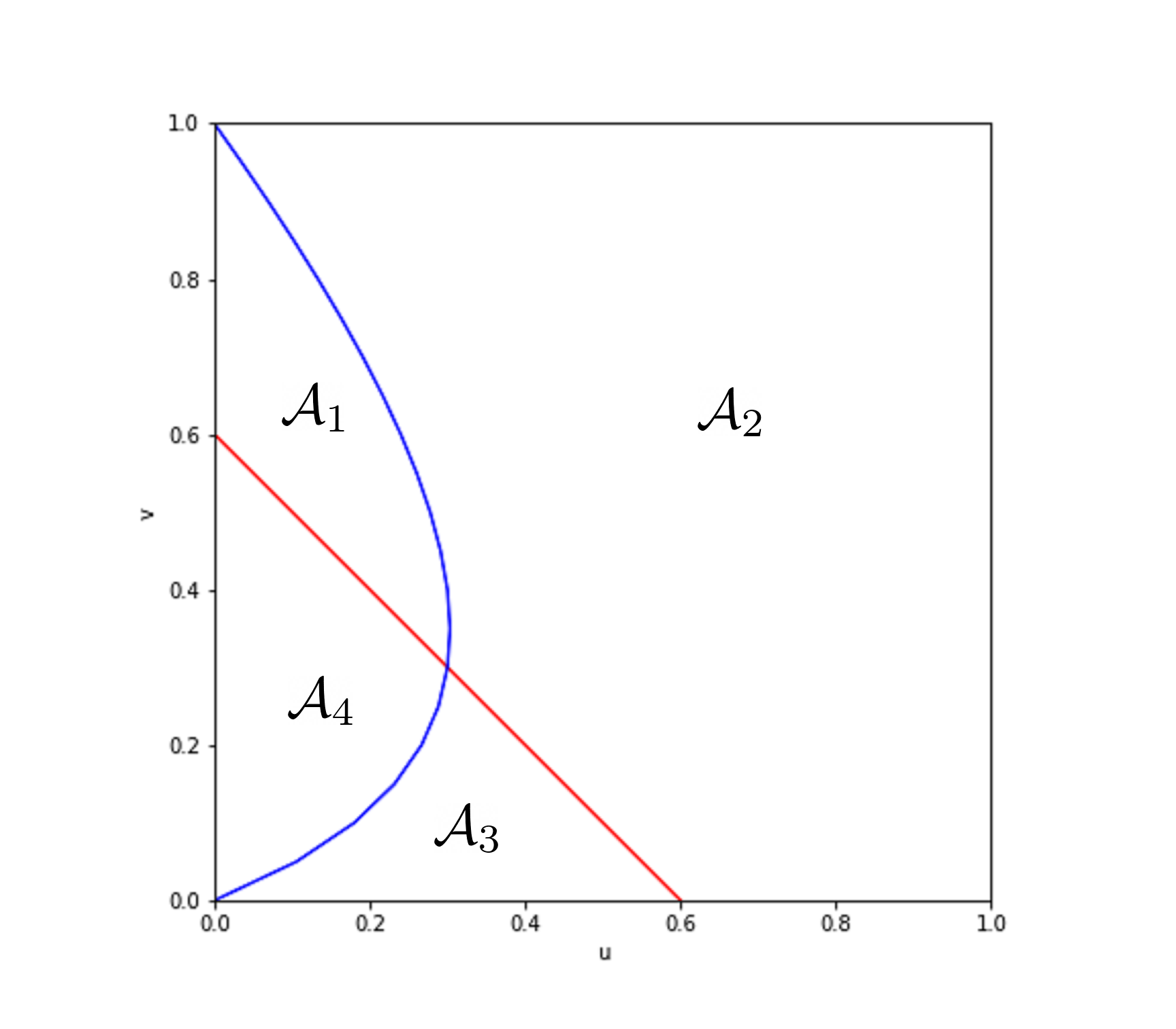}
	\caption{\it Partition of~$[0,1]\times[0,1]$ in the case~$a=0.8$,~$c=0.5$,~$\rho=2$,
		as given by~\eqref{ch2DEFA1234}.
		In red, the curve~$\dot{u}=0$. In blue, the curve~$\dot{v}=0$, parametrized
		by the function~$\sigma$ in~\eqref{ch2f:sigma}.}
	\label{ch2fig:zone}
\end{figure}

To visualize the statements in Corollary~\ref{ch2lemma:M1}, one can see Figure~\ref{ch2fig:zone}.

\begin{corollary}\label{ch2lemma:M2}
	Suppose that~$ac>1$.
	Then , we have that~$\dot{u}\leq 0$ in~$[0,1]\times [0,1]$,
	and the curve~\eqref{ch2curve:v'}
	divides the square~$[0,1]\times[0,1]$ into two regions:
	\begin{equation}\begin{split}\label{ch2DEFA12}
	\mathcal{A}_1 &\;:= \big\{ (u, v) \in [0,1]\times[0,1]\;{\mbox{ s.t. }}\; \dot{u}\leq 0,\; \dot{v}\geq 0 \big\}, \\
	\mathcal{A}_2 &\;:= \big\{ (u, v) \in [0,1]\times[0,1]\;{\mbox{ s.t. }}\; \dot{u}\leq 0, \;
	\dot{v}\leq 0 \big\}.
	\end{split}\end{equation}
	
	Furthermore, the sets~$\mathcal{A}_1$
	and~$\mathcal{A}_2$ are separated by the curve~$\dot{v}=0$, given by the graph of the continuous function~$\sigma$ given in~\eqref{ch2f:sigma}.
	
	In addition,
	\begin{equation}\label{ch2aggiun2}
	\mathcal{M}\subset \mathcal{A}_2.\end{equation} 
\end{corollary}

Proposition~\ref{ch2lemma:M} and Corollaries~\ref{ch2lemma:M1}
and~\ref{ch2lemma:M2}
are a bit technical, but provide fundamental information
to obtain a characterization of the sets~$\mathcal{E}$ and
$\mathcal{B}$, given in the forthcoming Proposition~\ref{ch2prop:char}.

We now provide the proof of Proposition~\ref{ch2lemma:M}
(and, as a byproduct, of Corollaries~\ref{ch2lemma:M1} and~\ref{ch2lemma:M2}).

\begin{proof}[Proof of Proposition~\ref{ch2lemma:M} and
	Corollaries~\ref{ch2lemma:M1} and~\ref{ch2lemma:M2}]
	We treat separately the cases~$0<ac<1$ and~$ac> 1$.
	We start with the case~$0<ac<1$, and divide the proof in three steps.
	\medskip
	
	\emph{Step 1: localizing~$\mathcal{M}$.}
	With the notation introduced in~\eqref{ch2DEFA1234},
	we prove that
	\begin{equation}\label{ch2dkoegjerig94768}\begin{split}&
	{\mbox{all trajectories starting in~$\mathcal{A}_3\setminus \{(0,0), (u_s,v_s) \}~$}}\\
	&{\mbox{exit the
			set~$\mathcal{A}_3$ on the side~$v=0$.}}\end{split}\end{equation} 
	To this aim, we first observe that
	\begin{equation}\label{ch2pouyi86}
	{\mbox{there are no cycles entirely contained
			in~$\mathcal{A}_3$,}}\end{equation}
	because~$\dot{u}$ and~$\dot{v}$ have a sign.
	Furthermore, 
	\begin{equation}\label{ch2pouyi862}{\mbox{there are no equilibria where a trajectory
			in the interior
			of~$\mathcal{A}_3$ can converge.}}\end{equation}
	Indeed,
	no point in~$\mathcal{A}_3$ with positive first coordinate can be mapped in~$(0,0)$ without exiting the set, because~$\dot{u}\geq 0$ in~$\mathcal{A}_3$. 
	Also, for all~$(u_0, v_0)\in \mathcal{A}_3 \setminus(u_s, v_s)$, we have that~$v_0<v_s$.
	On the other hand,~$\dot{v}\leq 0$ in~$\mathcal{A}_3$, so no trajectory that is entirely
	contained in~$\mathcal{A}_3$ can converge to~$(u_s, v_s)$. These observations prove~\eqref{ch2pouyi862}.
	
	As a consequence of~\eqref{ch2pouyi86},~\eqref{ch2pouyi862} and the
	Poincar\'e-Bendixson Theorem (see e.g.~\cite{TESCHL}), we have that
	all the trajectories in the interior of~$\mathcal{A}_3$ must exit the set at some time. 
	
	We remark that the side connecting~$(0,0)$ and~$(u_s, v_s)$ can be written as the
	of points belonging to
	$$\big\{ (u,v)\in [0,1]\times (0,v_s) \;{\mbox{ s.t. }}\;
	u=\sigma(v) \big\},$$
	where the function~$\sigma$ is defined in~\eqref{ch2f:sigma}.
	In this set, it holds that~$\dot{v}=0$ and~$\dot{u}>0$, thus the normal derivative pointing outward~$\mathcal{A}_3$ is negative, so the trajectories cannot go outside~$\mathcal{A}_3$ passing through this side. 
	
	Furthermore, on the side connecting~$(u_s, v_s)$ with~$(1-ac, 0)$, that lies on the straight line~$v=1-ac-u$, we have that~$\dot{u}= 0$  and~$\dot{v}<0$ for ~$(u,v)\neq (u_s,v_s)$, so also here the outer normal derivative is negative. Therefore, the trajectories cannot go outside~$\mathcal{A}_3$ passing through this side either. 
	
	These considerations complete the proof of~\eqref{ch2dkoegjerig94768}.
	Accordingly, recalling the definition of~$ \mathcal{E}$ in~\eqref{ch2DEFE}, we see
	that
	\begin{equation}\label{ch2prima}
	(\mathcal{A}_3 \setminus \{ (0,0), (u_s, v_s) \} ) \subseteq \mathcal{E}.
	\end{equation}
	In a similar way one can prove that all trajectories starting in~$\mathcal{A}_1\setminus \{(u_s,v_s) \}$ must converge to~$(0,1)$, which, recalling the definition of~$\mathcal{B}$
	in~\eqref{ch2DEFB}, implies that
	\begin{equation}\label{ch2prima2}
	\mathcal{A}_1\setminus \{(u_s,v_s) \} \subset \mathcal{B}.\end{equation}
	Thanks to~\eqref{ch2prima} and~\eqref{ch2prima2},
	we have that the stable manifold~$\mathcal{M}$ has no intersection with~$\mathcal{A}_1\setminus \{(u_s,v_s) \}~$ and~$\mathcal{A}_3\setminus \{(0,0),(u_s,v_s) \}~$,
	and therefore~$\mathcal{M}$ must lie
	in~$\mathcal{A}_2\cup \mathcal{A}_4$.
	
	Also, we know that~$\mathcal{M}$ is tangent to an eigenvector in~$(u_s, v_s)$,
	and we observe that
	\begin{equation}\label{ch2poui985004}
	{\mbox{$(1, -1)$ is not an eigenvector of the linearized system.}}\end{equation}
	Indeed, if~$(1, -1)$ were an eigenvector, then
	$$
	\begin{pmatrix}
	1-ac-2u_s-v_s & -u_s \\ 
	-\rho v_s-a & \rho-\rho u_s-2\rho v_s
	\end{pmatrix}\cdot
	\begin{pmatrix}
	1\\-1\end{pmatrix}=\lambda\begin{pmatrix}
	1\\-1\end{pmatrix},
	$$
	which implies that~$1-ac-a-\rho=(u_s+v_s)(1-\rho)$. Hence, recalling~\eqref{ch2usvs},
	we obtain that~$-a=\rho a c$, which is impossible. This
	establishes~\eqref{ch2poui985004}.
	
	In light of~\eqref{ch2poui985004}, we conclude that~$\mathcal{M}\setminus \{(u_s,v_s) \}$
	must have intersection with both~$\mathcal{A}_2$ and~$\mathcal{A}_4$.
	
	\medskip
	
	\emph{Step 2: defining~$\gamma(u)$.}
	Since~$\dot{u}> 0$ and~$\dot{v}>0$ in the interior of~$\mathcal{A}_4$,
	the portion of~$\mathcal{M}$
	in~$\mathcal{A}_4$ can be described globally as the graph of a
	monotone increasing smooth function~$\gamma_1:U\to[0,v_s]$,
	for a suitable interval~$U\subseteq[0,u_s]$ with~$u_s\in U$, and such that~$\gamma_1(u_s)=v_s$. 
	
	We stress that, for~$u>u_s$, the points~$(u,v)\in \mathcal{M}$ belong
	to~$\mathcal{A}_2$.
	
	Similarly, in the interior of~$\mathcal{A}_2$ we have that~$\dot{u}< 0$ and~$\dot{v}<0$.
	Therefore, we find that~$\mathcal{M}$ can be represented in~$\mathcal{A}_2$
	as the graph of a monotone increasing smooth function~$\gamma_2: V\to [v_s, 1]$, for a suitable interval~$V\subseteq[u_s,1]$ with~$u_s\in V$, and such that~$\gamma_2(u_s)=v_s$.  Notice that in the second case the trajectories and the parametrization run in opposite directions.
	
	Now, we define
	\begin{equation*}
	\gamma(u) := \begin{cases}
	\gamma_1(u)  & {\mbox{ if }}u\in U, \\
	\gamma_2(u)  &{\mbox{ if }} u\in V,
	\end{cases}
	\end{equation*}
	and we observe that it is an increasing smooth function locally
	parametrizing~$\mathcal{M}$ around~$(u_s,v_s)$ (thanks to the
	Stable Manifold Theorem).
	
	We point out that, in light of the
	Stable Manifold Theorem, the stable manifold~$\mathcal{M}$ is globally
	parametrized by
	an increasing smooth function on a set~$W\subset[0,1]$.
	\medskip
	
	\emph{Step 3:~$\gamma(0)=0$ and~$\gamma(u_{\mathcal{M}})=v_{\mathcal{M}}$
		for some~$(u_{\mathcal{M}},v_{\mathcal{M}})\in\partial\big([0,1]\times[0,1]\big)$.}
	We first prove that
	\begin{equation}\label{ch2r4gyghj}
	\gamma(0)=0.\end{equation}
	For this, we claim that
	\begin{equation}\label{ch2098ouitdbnb}
	{\mbox{orbits in the interior of~$\mathcal{A}_4$ do not come from
			outside~$\mathcal{A}_4$.}}
	\end{equation}
	Indeed, it is easy to see that points on the half axis~$\{u=0\}$ converge to~$(0,1)$,
	and therefore a trajectory cannot enter~$\mathcal{A}_4$ from this side. 
	
	As for the side connecting~$(0,0)$ to~$(u_s, v_s)$, here one has that~$\dot{u}\geq0$ and~$\dot{v}=0$, and so the inward pointing normal derivative is negative. Therefore,
	no trajectory can enter~$\mathcal{A}_4$ on this side.
	
	Moreover, on the side connecting~$(u_s, v_s)$ to~$(0, 1-ac)$ the inward pointing normal derivative is negative, because~$\dot{u}=0$ and~$\dot{v}\ge0$, thus we have that no
	trajectory can enter~$\mathcal{A}_4$ on this side either.
	These considerations prove~\eqref{ch2098ouitdbnb}.
	
	Furthermore, we have that
	\begin{equation}\label{ch2098ouitdbnb2}
	{\mbox{no cycles are allowed in~$\mathcal{A}_4$,}}
	\end{equation}
	because~$\dot{u}\ge0$ and~$\dot{v}\ge0$ in~$\mathcal{A}_4$.
	
	{F}rom~\eqref{ch2098ouitdbnb},~\eqref{ch2098ouitdbnb2} and the
	Poincar\'e-Bendixson Theorem (see e.g.~\cite{TESCHL}), we conclude that,
	given a point~$(\tilde u,\tilde v)\in\mathcal{M}$ in the
	interior of~$\mathcal{A}_4$, the~$\alpha$-limit set of~$(\tilde u,\tilde v)$,
	that we denote by~$\alpha_{(\tilde u,\tilde v)}$,
	can be 
	\begin{equation}\label{ch209765gjkd}\begin{split}&
	{\mbox{either an equilibrium or a union of (finitely many)}}\\
	&{\mbox{equilibria and non-closed orbits connecting these equilibria.}}\end{split}
	\end{equation}
	
	We stress that, being~$(\tilde u,\tilde v)$ in the
	interior of~$\mathcal{A}_4$, we have that
	\begin{equation}\label{ch28yfe993vcem}
	\tilde u<u_s.
	\end{equation}
	Now, we observe that
	\begin{equation}\label{ch2degfiewgh}
	{\mbox{$\alpha_{(\tilde u,\tilde v)}$ cannot contain the saddle point~$(u_s,v_s)$.}}\end{equation}
	Indeed, suppose by contradiction that~$\alpha_{(\tilde u,\tilde v)}$ does
	contain~$(u_s,v_s)$. Then,
	we denote by~$\phi_{(\tilde u,\tilde v)}(t)=\big(u_{(\tilde u,\tilde v)}(t),v_{(\tilde u,\tilde v)}(t)\big)$ the solution of~\eqref{ch2model}
	with~$\phi_{(\tilde u,\tilde v)}(0)=(\tilde u,\tilde v)$, and we have that
	there exists a sequence~$t_j\to-\infty$ such
	that~$\phi_{(\tilde u,\tilde v)}(t_j)$ converges to~$(u_s,v_s)$ as~$j\to+\infty$.
	In particular, in light of~\eqref{ch28yfe993vcem},
	there exists~$j_0$ sufficiently large
	such that
	$$ u_{(\tilde u,\tilde v)}(0)=\tilde u<u_{(\tilde u,\tilde v)}(t_{j_0}).$$
	Consequently, there exists~$t_\star\in(t_{j_0},0)$ such that~$\dot
	u_{(\tilde u,\tilde v)}(t_\star)<0$.
	
	As a result, it follows that~$\phi_{(\tilde u,\tilde v)}(t_\star)\not\in\mathcal{A}_4$.
	This, together with the fact that~$\phi_{(\tilde u,\tilde v)}(0)\in\mathcal{A}_4$,
	is in contradiction with~\eqref{ch2098ouitdbnb}, and the proof of~\eqref{ch2degfiewgh}
	is thereby complete.
	
	Thus, from~\eqref{ch209765gjkd} and~\eqref{ch2degfiewgh},
	we deduce that~$\alpha_{(\tilde u,\tilde v)}=\{(0,0)\}$.
	This gives that~$(0,0)$ lies on the stable manifold~$\mathcal{M}$,
	and therefore the proof of~\eqref{ch2r4gyghj} is complete.
	
	Now, we show that
	\begin{equation}\label{ch2ifregkjh0000}
	{\mbox{there exists~$(u_{\mathcal{M}},v_{\mathcal{M}})\in\partial\big([0,1]\times[0,1]\big)$
			such that~$\gamma(u_{\mathcal{M}})=v_{\mathcal{M}}$.}}
	\end{equation}
	To prove it, we first observe that
	\begin{equation}\label{ch2ifregkjh0000pre}
	{\mbox{orbits in~$\mathcal{A}_2$ converging to~$(u_s,v_s)$ come from
			outside~$\mathcal{A}_2$.}}
	\end{equation}
	Indeed, we suppose by contradiction that
	\begin{equation}\label{ch20696u833687}
	{\mbox{an orbit in~$\mathcal{A}_2$
			converging to~$(u_s,v_s)$ stays confined in~$\mathcal{A}_2$.}}\end{equation}
	We remark that, in this case, 
	\begin{equation}\label{ch2doeutoeru}
	{\mbox{an orbit in~$\mathcal{A}_2$ cannot be a cycle,}}
	\end{equation}
	because~$\dot{u}$
	and~$\dot{v}$ have a sign in~$\mathcal{A}_2$.
	Then, by the
	Poincar\'e-Bendixson Theorem (see e.g.~\cite{TESCHL}), we conclude that,
	given a point~$(\tilde u,\tilde v)\in\mathcal{M}$ in the
	interior of~$\mathcal{A}_2$, the~$\alpha$-limit set of~$(\tilde u,\tilde v)$,
	that we denote by~$\alpha_{(\tilde u,\tilde v)}$,
	can be either an equilibrium or a union of (finitely many)
	equilibria and non-closed orbits connecting these equilibria.
	We notice that the set~$\alpha_{(\tilde u,\tilde v)}$ cannot contain~$(0,1)$,
	since it is a stable equilibrium.
	We also claim that
	\begin{equation}\label{ch2qewytriyb}
	{\mbox{$\alpha_{(\tilde u,\tilde v)}$ cannot contain~$(u_s,v_s)$.}}\end{equation}
	Indeed, we
	suppose by contradiction that~$\alpha_{(\tilde u,\tilde v)}$ does
	contain~$(u_s,v_s)$. We observe that, since~$\dot{u}\le0$ in~$\mathcal{A}_2$,
	\begin{equation}\label{ch2koewtuyh}
	\tilde u>u_s.\end{equation}
	We denote by~$\phi_{(\tilde u,\tilde v)}(t)=\big(u_{(\tilde u,\tilde v)}(t),v_{(\tilde u,\tilde v)}(t)\big)$ the solution of~\eqref{ch2model}
	with~$\phi_{(\tilde u,\tilde v)}(0)=(\tilde u,\tilde v)$, and we have that
	there exists a sequence~$t_j\to-\infty$ such
	that~$\phi_{(\tilde u,\tilde v)}(t_j)$ converges to~$(u_s,v_s)$ as~$j\to+\infty$.
	In particular, in light of~\eqref{ch2koewtuyh},
	there exists~$j_0$ sufficiently large
	such that
	$$ u_{(\tilde u,\tilde v)}(0)=\tilde u>u_{(\tilde u,\tilde v)}(t_{j_0}).$$
	Consequently, there exists~$t_\star\in(t_{j_0},0)$ such that~$\dot
	u_{(\tilde u,\tilde v)}(t_\star)>0$.
	Accordingly, we have that~$\phi_{(\tilde u,\tilde v)}(t_\star)\not\in\mathcal{A}_2$.
	This and the fact that~$\phi_{(\tilde u,\tilde v)}(0)\in\mathcal{A}_2$ give
	a contradiction with~\eqref{ch20696u833687}, and therefore this
	establishes~\eqref{ch2qewytriyb}.
	
	These considerations complete the proof of~\eqref{ch2ifregkjh0000pre}.
	
	Now, we observe that
	the inward pointing normal derivative
	at every point in~$\mathcal{A}_2 \cap \mathcal{A}_3 \setminus\{(u_s, v_s)\}$
	is negative, since~$\dot{u}=0$ and~$\dot{v}\le0$. Hence, no trajectory can enter
	from this side.
	Also, 
	the inward pointing normal derivative
	at every point in~$\mathcal{A}_1 \cap \mathcal{A}_2 \setminus\{(u_s, v_s)\}$
	is negative, since~$\dot{u}\le0$ and~$\dot{v}=0$. Hence, no trajectory can enter
	from this side either.
	
	These observations and~\eqref{ch2ifregkjh0000pre} give the desired result
	in~\eqref{ch2ifregkjh0000}, and thus Proposition~\ref{ch2lemma:M}
	is established in the case~$ac<1$.
	
	\medskip
	
	Now we treat the case~$ac>1$, using the same ideas. In this setting,
	$\mathcal{M}$ is the stable manifold associated with the saddle point~$(0,0)$.
	We point out that, in this case, for all points in~$[0,1]\times [0,1]$ we have
	that~$\dot{u}\leq 0$. Hence, the curve of points satisfying~$\dot{v}=0$, that was also given in~\eqref{ch2curve:v'}, divides the square~$[0,1]\times[0,1]$ into two regions~$
	\mathcal{A}_1$ and~$\mathcal{A}_2$, defined in~\eqref{ch2DEFA12}.	
	
	Now, one can repeat verbatim the arguments in {\emph{Step 1}} with obvious modifications,
	to find that~$
	\mathcal{M}\subset \mathcal{A}_2$.
	
	Since the derivatives of~$u$ and~$v$ have a sign in~$\mathcal{A}_2$, and the
	set~$\mathcal{M}$ in this case is the trajectory of a point converging to~$(0,0)$, the set~$\mathcal{M}$ can be represented globally as the graph of a smooth increasing function~$\gamma: U\to [0,1]$ for a suitable interval~$U\subseteq[0,1]$ containing the origin.
	As a consequence, the condition~$\gamma(0)=0$ is trivially satisfied in this setting.
	The existence of a suitable~$(u_{\mathcal{M}},v_{\mathcal{M}})$ can be derived reasoning as in {\emph{Step 3}}
	with obvious modifications.
	
	Now, we prove that
	\begin{equation}\label{ch2ofriyty98579}
	{\mbox{at~$u=0$ the function~$\gamma$ is tangent
			to the line~$(\rho-1+ac)v-au=0$.}}\end{equation}
	For this, we recall~\eqref{ch2Jmatrix} and we see, by inspection, that	
	the Jacobian matrix~$J(0,0)$ has two eigenvectors, namely~$(0,1)$ and~~$(\rho-1+ac, a)$. The first one is tangent to the line~$u=0$, that is the unstable manifold of~$(0,0)$,
	as one can easily verify. Thus, the second eigenvector is the one tangent to~$\mathcal{M}$, as prescribed by the Stable Manifold Theorem (see e.g.~\cite{dynsyst}). Hence, in~$(0,0)$ the manifold~$\mathcal{M}$ is tangent to the line~$(\rho-1+ac)v-au=0$ and so is the function~$\gamma$ in~$u=0$. This proves~\eqref{ch2ofriyty98579}, and thus
	Proposition~\ref{ch2lemma:M}
	is established in the case~$ac>1$ as well.
\end{proof}	

\subsection{Characterization of~$\mathcal{M}$ when~$ac=1$}\label{ch2sec:deg}

Here we will prove the counterpart of Proposition~\ref{ch2lemma:M}
in the degenerate case~$ac=1$.

To this end, looking at the velocity fields,
we first observe that
\begin{equation}\label{ch2NOEX}
\begin{split}&
{\mbox{trajectories starting in~$(0,1)\times(-\infty,1)$ at time~$t=0$}}\\&{\mbox{remain in~$(0,1)\times(-\infty,1)$ for all time~$t>0$.}}\end{split}
\end{equation}
We also point out that
\begin{equation}\label{ch2CALR}
\begin{split}&
{\mbox{trajectories entering the region~${\mathcal{R}}:=
		\{u\in(0,1),\,u+v<0\}
		$}}\\&{\mbox{at some time~$t_0\in\R$}}\\&{\mbox{remain in that region for all time~$t>t_0$,}}\end{split}
\end{equation}
since~$\dot v=\rho v(1-u-v)-au=-\rho u-au<0$
along~$\{u\in(0,1),\,u+v=0\}$.

Also, by the
Center Manifold Theorem
(see e.g. Theorem~1 on page~16
of~\cite{MR635782} or pages 89-90 in~\cite{MR1031257}),
there exists a collection~$\mathcal{M}_0$
of invariant curves, which are all
tangent at the origin to the eigenvector corresponding to the null eigenvalue,
that is the straight line~$\rho v-au=0$. Then, we define~$\mathcal{M}:=
\mathcal{M}_0\cap ([0,1]\times[0,1])$ and we observe that this
intersection is nonvoid, given the tangency property of~$\mathcal{M}_0$
at the origin.

In what follows, for every~$t\in\R$, we denote by~$(u(t),v(t))=\phi_p(t)$ the orbit of~$p\in\mathcal{M}\setminus\{(0,0)\}$. We start by providing an observation
related to negative times:

\begin{lemma} \label{ch2NOCPA}
	If~$p\in\mathcal{M}\setminus\{(0,0)\}$
	then~$\phi_p(t)$ cannot approach the origin for negative values of~$t$.
\end{lemma}

\begin{proof}
	We argue by contradiction
	and denote by~$t_1,\dots,t_n,\dots$ a sequence of such negative values of~$t$, for which~$t_n\to-\infty$ and
	$$ \lim_{n\to+\infty}\phi_p(t_n)=(0,0).$$
	Up to a subsequence, we can also suppose that
	\begin{equation}\label{ch2bejv0565etP}
	u(t_{n+1})<u(t_n).
	\end{equation}
	In light of~\eqref{ch2CALR}, we have that, for all~$T\le0$,
	\begin{equation}\label{ch20okf3233}
	\phi_p(T)\not\in{\mathcal{R}}.
	\end{equation}
	Indeed, if~$\phi_p(T)\in{\mathcal{R}}$, we deduce from~\eqref{ch2CALR}
	that~$\phi_p(t)\in{\mathcal{R}}$ for all~$t\ge T$. In particular,
	we can take~$t=0\ge T$ and conclude that~$p=\phi_p(0)\in{\mathcal{R}}$,
	and this is in contradiction with the assumption that~$p\in{\mathcal{M}}\setminus\{(0,0)\}$.
	
	As a byproduct of~\eqref{ch20okf3233}, we obtain that, for all~$T\le0$,
	$$\phi_p(T)\in
	\{u\in(0,1),\,u+v\ge0\}\subseteq\{\dot u=-u(u+v)\le0\}.$$
	In particular
	$$ u(t_n)-u(t_{n+1})=\int_{t_{n+1}}^{t_n}\dot u(\tau)\,d\tau\le0,$$
	which is in contradiction with~\eqref{ch2bejv0565etP},
	and consequently we have established the desired result.
\end{proof}

Now we show that the~$\omega$-limit of any point lying
on the global center manifold coincides with the origin, according to the next result:

\begin{lemma}\label{ch2lerptj:lemma}
	If~$p\in\mathcal{M}$, then its~$\omega$-limit is~$(0,0)$.
\end{lemma}

\begin{proof}
	We observe that, for every~$t>0$,
	\begin{equation}\label{ch2T0z}
	\phi_p(t)\in[0,1]\times[0,1].
	\end{equation}
	Indeed, by~\eqref{ch2NOEX}, one
	sees that, for~$t>0$,
	$\phi_t(p)$ cannot cross~$\{0\}\times[0,1]$,
	$\{1\}\times[0,1]$ and~$[0,1]\times\{1\}$,
	hence the only possible escape side is given by~$[0,1]\times\{0\}$.
	Therefore, to prove~\eqref{ch2T0z}, we suppose, by contradiction,
	that there exists~$t_0\ge0$ such that~$\phi_p({t_0})\in[0,1]\times\{0\}$,
	that is~$v(t_0)=0$. Since~$(0,0)$ is an equilibrium, it follows that~$u(t_0)\ne0$.
	In particular,~$u(t_0)>0$ and accordingly~$\dot v(t_0)=-au(t_0)<0$.
	This means that~$v(t_0+\varepsilon)<0$ for all~$\varepsilon\in
	(0,\varepsilon_0)$ for a suitable~$\varepsilon_0>0$.
	Looking again at the velocity fields, this entails that~$\phi_p(t)\in(0,1)\times(-\infty,0)$
	for all~$t>\varepsilon_0$. Consequently,
	$\phi_p(t)$ cannot approach the
	straight line~$\rho v-au=0$ for~$t>\varepsilon_0$.
	
	This, combined with Lemma~\ref{ch2NOCPA}, says that the trajectory
	emanating from~$p$ can never
	approach the
	straight line~$\rho v-au=0$ at the origin, in contradiction with the definition
	of~${\mathcal{M}}$,
	and thus the proof of~\eqref{ch2T0z}
	is complete.
	
	{F}rom~\eqref{ch2T0z} 
	and the Poincar\'e-Bendixson Theorem (see e.g.~\cite{TESCHL}),
	we deduce that the~$\omega$-limit of~$p$
	can be either a cycle, or an equilibrium,
	or a union of (finitely many)
	equilibria and non-closed orbits connecting these equilibria.
	We observe that the~$\omega$-limit of~$p$ cannot be a cycle, since~$\dot u$
	has a sign in~$[0,1]\times[0,1]$. Moreover, it cannot contain the sink~$(0,1)$, due to
	Lemma~\ref{ch2NOCPA}. Hence, the only possibility is that
	the~$\omega$-limit of~$p$ coincides with~$(0,0)$,
	which is the desired result.
\end{proof}

As a consequence of Lemma~\ref{ch2lerptj:lemma} and the fact
that~$\dot u<0$ in~$(0,1]\times[0,1]$, we obtain the following statement:

\begin{corollary}\label{ch2N7u43566r}
	Every trajectory in~$\mathcal{M}$ has
	the form~$\{\phi_p(t),\,t\in\R\}$, with
	$$\lim_{t\to+\infty}\phi_p(t)=(0,0)$$
	and there exists~$t_p\in\R$ such that~$\phi_p(t_p)\in\big(\{1\}\times[0,1]\big)
	\cup\big([0,1]\times\{1\}\big)$.
\end{corollary}

The result in Corollary~\ref{ch2N7u43566r} can be sharpened
in view of the following statement (which can be seen as the counterpart
of Proposition~\ref{ch2lemma:M}
in the degenerate case~$ac=1$): namely, since the center manifold can in principle contain
many different trajectories
(see e.g. Figure~5.3 in~\cite{MR635782}), we provide a
tailor-made argument that excludes this possibility in the specific case
that we deal with.

\begin{proposition}\label{ch2M:p045}
	For $ac=1$	
	$\mathcal{M}$ contains one, and only one, trajectory,
	which is asymptotic to the origin as~$t\to+\infty$, and that can be written
	as a graph~$\gamma:[0,u_{\mathcal{M}}]\to[0,v_{\mathcal{M}}]$, for some~$(u_{\mathcal{M}},v_{\mathcal{M}})\in\big(\{1\}\times[0,1]\big)\cup
	\big((0,1]\times\{1\}\big)$, where~$\gamma$ is an increasing~$C^2$ function such that~$\gamma(0)=0$,
	$\gamma(u_{\mathcal{M}})=v_{\mathcal{M}}$ and the graph of~$\gamma$ at the origin is tangent to the line~$\rho v-au=0$.
\end{proposition}

\begin{proof} First of all, we show that
	\begin{equation}\label{ch2LIMSDD-0}
	{\mbox{$\mathcal{M}$ contains one, and only one, trajectory.}}\end{equation}
	Suppose, by contradiction, that~${\mathcal{M}}$
	contains two different orbits, that we denote by~${\mathcal{M}}_-$
	and~${\mathcal{M}}_+$.
	Using Corollary~\ref{ch2N7u43566r},
	we can suppose that~${\mathcal{M}}_+$
	lies above~${\mathcal{M}}_-$
	and
	\begin{equation}\label{ch2CQPSKD}
	\begin{split}&
	{\mbox{the region~${\mathcal{P}}\subset[0,1]\times[0,1]$ contained between~${\mathcal{M}}_+$
			and~${\mathcal{M}}_-$}}\\&{\mbox{lies in~$\{\dot u<0\}$.}}\end{split}
	\end{equation}
	Consequently, for every~$p\in{\mathcal{P}}$,
	it follows that
	\begin{equation}\label{ch29ikfjty} \lim_{t\to+\infty}\phi_p(t)=(0,0).\end{equation}
	In particular, we can take an open ball~$B\subset
	{\mathcal{P}}$ in the vicinity of the origin, denote by~$\mu(t)$ the Lebesgue measure of~${\mathcal{S}}(t):=\{\phi_p(t),\;
	p\in B\}$, and write that~$\mu(0)>0$
	and
	\begin{equation}\label{ch2LIMSDD} \lim_{t\to+\infty}\mu(t)=0.\end{equation}
	We point out that~${\mathcal{S}}(t)$
	lies in the vicinity of the origin for all~$t\ge0$, thanks to~\eqref{ch2CQPSKD}.
	As a consequence, for all~$t$,~$\tau>0$, changing variable
	$$ y:=\phi_{x}(\tau)=x+\int_0^\tau \frac{d\phi_x(\theta)}{d\theta}\,d\theta=
	x+\tau \frac{d\phi_x(0)}{dt}+O(\tau^2),
	$$
	we find that
	\begin{eqnarray*}
		\mu(t+\tau)&=&\int_{{\mathcal{S}}(t+\tau)}dy\\&=&
		\int_{{\mathcal{S}}(t)}\big|\det \big(D_x \phi_x(\tau)\big)\big|\,dx\\&=&
		\int_{{{\mathcal{S}}(t)}}\left|\det D_x \left(x+\tau \frac{d\phi_x(0)}{dt}
		+O(\tau^2)
		\right)\right|\,dx\\
		&=&
		\int_{{{\mathcal{S}}(t)}}\left( 1+\tau\,{\rm Tr}\left(D_x \frac{d\phi_x(0)}{dt}\right)
		+O(\tau^2)
		\right)\,dx
		\\&=&\mu(t)+\tau
		\int_{{{\mathcal{S}}(t)}} {\rm Tr}\left(D_x \frac{d\phi_x(0)}{dt}\right)\,dx
		+O(\tau^2)
		,
	\end{eqnarray*}
	where~${\rm Tr}$ denotes the trace of a~$(2\times2)$-matrix.
	
	As a consequence,
	\begin{equation}\label{ch2090987t4oyyorfg4}
	\frac{d\mu}{dt}(t)=
	\int_{{{\mathcal{S}}(t)}} {\rm Tr}\left(D_x \frac{d\phi_x(0)}{dt}\right)\,dx.\end{equation}
	Also, using the notation~$x=(u,v)$, we can write~\eqref{ch2model} when~$ac=1$ in the form
	$$ \frac{d\phi_x}{dt}(t)=
	\dot x(t)=\left(
	\begin{matrix}
	\dot u(t)\\
	\dot v(t)
	\end{matrix}\right)=
	\left(
	\begin{matrix}
	-u(t)(u(t)+v(t))\\
	\rho v(t)(1-u(t)-v(t))-au(t)
	\end{matrix}\right).
	$$
	Accordingly,
	\begin{eqnarray*} D_x \frac{d\phi_x(0)}{dt}&=&
		\left(
		\begin{matrix}
			-\partial_u\big(u(u+v)\big)&-\partial_v\big(u(u+v)\big)
			\\
			\partial_u\big(\rho v (1-u-v)-au\big)&\partial_v\big(\rho v (1-u-v)-au\big)
		\end{matrix}\right),
	\end{eqnarray*}
	whence
	\begin{equation}\label{ch20oj476ytgf9846}
	\begin{split}
	{\rm Tr}\left(D_x \frac{d\phi_x(0)}{dt}\right)\,&=-\partial_u\big(u(u+v)\big)
	+\partial_v\big(\rho v (1-u-v)-au\big)\\&=-2u-v+\rho(1-u-v)-\rho v\\&=\rho+O(|x|)
	\end{split}
	\end{equation}
	for~$x$ near the origin.
	
	As a result, recalling~\eqref{ch29ikfjty}, we can take~$t$ sufficiently large,
	such that~${{{\mathcal{S}}(t)}}$ lies in a neighborhood of the origin, exploit~\eqref{ch20oj476ytgf9846}
	to write that~${\rm Tr}\left(D_x \frac{d\phi_x(0)}{dt}\right)\ge\frac\rho2$
	and then~\eqref{ch2090987t4oyyorfg4} to conclude that
	$$ \frac{d\mu}{dt}(t)\ge \frac{\rho}2
	\int_{{{\mathcal{S}}(t)}} dx=\frac{\rho}{2}\,\mu(t).$$
	This implies that~$\mu(t)$ diverges (exponentially fast)
	as~$t\to+\infty$, which is in contradiction with~\eqref{ch2LIMSDD}.
	The proof of~\eqref{ch2LIMSDD-0}
	is thereby complete.
	
	Now, we check the other claims in the statement of Proposition~\ref{ch2M:p045}.
	The asymptotic property as~$t\to+\infty$ is a consequence of Corollary~\ref{ch2N7u43566r}.
	Also, the graphical property as well as the monotonicity
	property of the graph follow from the fact that~${\mathcal{M}}\subset\{\dot u<0\}$.
	The smoothness of the graph follows from the smoothness of the center manifold.
	The fact that~$\gamma(0)=0$ and
	$\gamma(u_{\mathcal{M}})=v_{\mathcal{M}}$ follow also from
	Corollary~\ref{ch2N7u43566r}. The tangency property at the origin is a consequence of
	the tangency property of the center manifold to the center eigenspace.
\end{proof}

As a byproduct of the proof of Proposition~\ref{ch2M:p045}
we also obtain the following information:

\begin{corollary}\label{ch2lemma:Mdeg}
	Suppose that~$ac=1$.
	Then , we have that~$\dot{u}\leq 0$ in~$[0,1]\times [0,1]$,
	and the curve~\eqref{ch2curve:v'}
	divides the square~$[0,1]\times[0,1]$ into two regions~$\mathcal{A}_1$
	and~$\mathcal{A}_2$,
	defined in~\eqref{ch2DEFA12}.
	
	Furthermore, the sets~$\mathcal{A}_1$
	and~$\mathcal{A}_2$ are separated by the curve~$\dot{v}=0$, given by the graph of the continuous function~$\sigma$
	given in~\eqref{ch2f:sigma}.
	
	In addition,
	\begin{equation}\label{ch2aggiun2BIS}
	\mathcal{M}\subset \mathcal{A}_2.\end{equation} 
\end{corollary}

\subsection{Completion of the proof of Theorem~\ref{ch2thm:dyn}}
\label{ch2sec:deg2}

We observe that, by the Stable Manifold Theorem
and the Center Manifold Theorem, the statement in~(v) of Theorem~\ref{ch2thm:dyn}
is obviously fulfilled. 

Hence, to complete the proof of Theorem~\ref{ch2thm:dyn},
it remains to show that the statement in~(iv) holds true.
To this aim, exploiting
the useful pieces of
information in Propositions~\ref{ch2lemma:M} and~\ref{ch2M:p045}, we first give a characterization of the sets~$\mathcal{E}$ and~$\mathcal{B}$:

\begin{proposition} \label{ch2prop:char}
	The following characterizations of the sets in~\eqref{ch2DEFB} and~\eqref{ch2DEFE}
	are true:
	\begin{equation}\label{ch2char:E}
	\begin{split}
	\mathcal{E}=\;&\Big\{ (u,v)\in [0,1]\times [0,1] \;{\mbox{ s.t. }}\; v<\gamma(u) \ \text{if} \ u\in[0,u_{\mathcal{M}}] \\
	&\qquad\qquad\qquad \qquad\qquad
	{\mbox{ and}} \;v\leq 1 \ \text{if} \ u\in(u_{\mathcal{M}}, 1]     \Big\},
	\end{split}
	\end{equation}
	and
	\begin{equation}\label{ch2char:B}
	\mathcal{B}=\Big\{ (u,v)\in [0,u_{\mathcal{M}}]\times [0,1]\;{\mbox{ s.t. }}\;  v>\gamma(u) \ \text{if} \ u\in[0,u_{\mathcal{M}}]  \Big\},
	\end{equation}
	for some~$(u_{\mathcal{M}}, v_{\mathcal{M}}) \in \partial \left( [0,1]\times [0,1] \right)$.
\end{proposition}

\begin{figure}[h] 
	\begin{subfigure}{.5\textwidth}
		\centering
		\includegraphics[width=.8\textwidth]{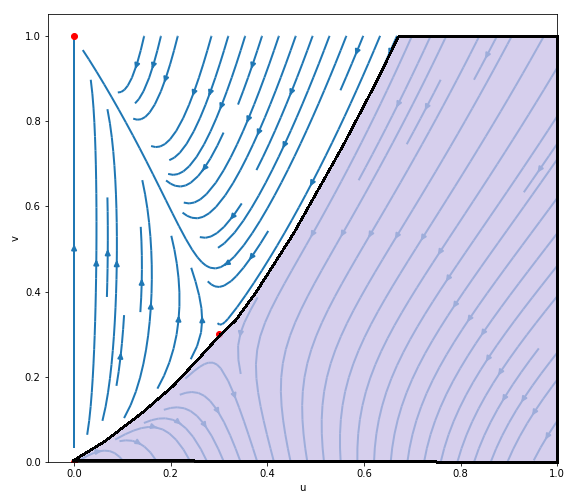}
		\caption{$a=0.8$,~$c=0.5$,~$\rho=2$}
	\end{subfigure}%
	\begin{subfigure}{.5\textwidth}
		\centering
		\includegraphics[width=.8\textwidth]{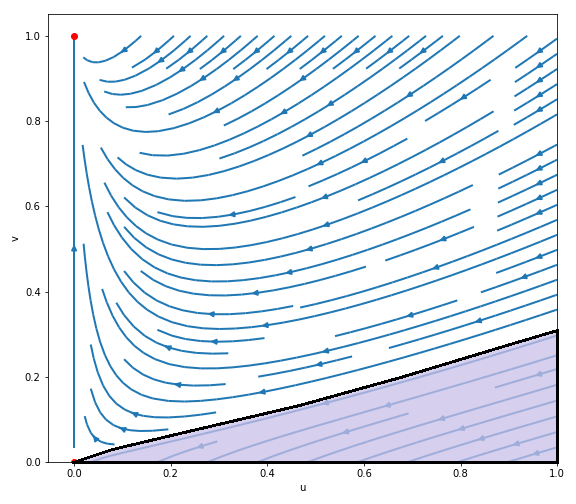}
		\caption{$a=0.8$,~$c=3$,~$\rho=2$}
	\end{subfigure}
	\caption{\it The figures show the phase portrait for the indicated values of the coefficients. In blue, the orbits of the points. The red dots show the equilibria. In violet, the set~$\mathcal{E}$.}
	\label{ch2fig:char}
\end{figure}

One can visualize the appearance of the set~$\mathcal{E}$ in~\eqref{ch2char:E}
in two particular cases in
Figure~\ref{ch2fig:char}.

\begin{proof}[Proof of Proposition~\ref{ch2prop:char}]
	We let~$\gamma$ be the parametrization of~$\mathcal{M}$,
	as given by Propositions~\ref{ch2lemma:M} (when~$ac\neq1$)
	and~\ref{ch2M:p045} (when~$ac=1$), and we
	consider the sets
	\begin{eqnarray*}
		&&\mathcal{X}:= \big\{ (u,v)\in [0,1]\times [0,1]\;{\mbox{ s.t. }}\; v < \gamma(u) \big\}\\
		{\mbox{and }}&&
		\mathcal{Y}:= \big\{ (u,v)\in [0,1]\times [0,1] \;{\mbox{ s.t. }}\;  v > \gamma(u) \big\}.
	\end{eqnarray*}
	The goal is to prove that~$\mathcal{X}\equiv\mathcal{E}$ and~$\mathcal{Y}\equiv\mathcal{B}$. We observe that, when~$u_{\mathcal{M}}=1$, then~$\mathcal{X} \cup \mathcal{Y} \cup \mathcal{M}=[0,1]\times[0,1]$.
	When instead~$u_{\mathcal{M}}\in(0,1)$, then~$\mathcal{X} \cup \mathcal{Y} \cup \mathcal{M}\cup\big((u_{\mathcal{M}},1]\times[0,1]\big)=[0,1]\times[0,1]$.
	Accordingly, if we show that 
	\begin{eqnarray}
	&&\mathcal{X}\cup\big((u_{\mathcal{M}},1]\times[0,1]\big)\subseteq\mathcal{E} \label{ch21627} \\
	{\mbox{and }}&&	\mathcal{Y}\subseteq\mathcal{B}, \label{ch21628}
	\end{eqnarray}
	we are done. 
	
	Hence, we now focus on the proof of~\eqref{ch21627}. 
	Namely, recalling~\eqref{ch2DEFE},
	we show that
	\begin{equation}\label{ch2alltra}
	{\mbox{all trajectories starting in~$\mathcal{X}$ exit the
			set on the side~$(0,1]\times \{0\}$.}}\end{equation} 
	For this, we first notice that, gathering together~\eqref{ch2primaBIS},~\eqref{ch2prima2BIS},~\eqref{ch2098ouitdbnb},
	\eqref{ch2098ouitdbnb2} and~\eqref{ch2doeutoeru}, we find that
	\begin{equation}\label{ch207960789djiewf2}
	{\mbox{no limit cycle exists in~$[0,1]\times[0,1]$}}\end{equation}
	(in the case~$0<ac<1$, and the same holds true in the case~$ac\ge1$ since~$\dot u$ has a sign). 
	
	In addition, 
	\begin{equation}\label{ch207960789djiewf}
	{\mbox{the~$\omega$-limit of any point in~$\mathcal{X}$
			cannot contain an equilibrium.}}\end{equation}
	Indeed, by Propositions~\ref{ch2lemma:M} (when~$(ac\neq1$)
	and~\ref{ch2M:p045} (when~$ac=1$),
	we have that~$\gamma(0)=0<1$, and therefore~$(0,1)\notin \overline{\mathcal{X}}$.
	Moreover, if~$ac<1$, 
	a trajectory in~$\mathcal{X}$ cannot converge to~$(u_s, v_s)$, since~$\mathcal{X}$ does not
	contain points of the stable manifold~$\mathcal{M}$, nor to~$(0,0)$,
	since this is a repulsive equilibrium and no trajectory converges here.
	If instead~$ac\geq 1$, then it cannot converge to~$(0,0)$, since~$\mathcal{X}$ does not
	contain points of~$\mathcal{M}$.
	These observations completes the proof of~\eqref{ch207960789djiewf}.
	
	{F}rom~\eqref{ch207960789djiewf2},~\eqref{ch207960789djiewf} and the 
	Poincar\'e-Bendixson Theorem (see e.g.~\cite{TESCHL}),
	we have that
	every trajectory starting in~$\mathcal{X}$ leaves the set (possibly in infinite time).
	
	If the trajectory leaves at~$t=+\infty$, then it converges to some equilibrium on~$\partial \mathcal{X}$, which is in contradiction with~\eqref{ch207960789djiewf}.
	
	As a consequence a trajectory in~$\mathcal{X}$ leaves the set in finite time.
	Suppose that a trajectory leaves~$\mathcal{X}$ at a point~$(u,v)\in \partial \mathcal{X}$; then either~$(u, v)\in \mathcal{M}$ or~$(u, v)\in \partial ( [0,1]\times[0,1] )$.
	The first possibility is impossible, otherwise the starting point of the trajectory would converge to~$(u_s, v_s)$. Hence, the only possibility is that the
	trajectory leaves~$\mathcal{X}$ at~$(u, v)\in \partial ( [0,1]\times[0,1] )$.
	By Proposition~\ref{ch2prop:dyn}
	this is possible only if~$u>0$ and~$v=0$,
	which proves~\eqref{ch2alltra}. As a consequence of~\eqref{ch2alltra} we obtain that
	\begin{equation}\label{ch2po0584yiugherghegk}
	\mathcal{X}\subseteq\mathcal{E}.\end{equation}
	
	We now claim that
	\begin{equation}\label{ch2po0584yiugherghegkbis}
	\big((u_{\mathcal{M}},1]\times[0,1]\big)\subseteq\mathcal{E}.
	\end{equation}
	To this end, we observe that there are neither cycles nor equlibria in~$(u_{\mathcal{M}},1]\times[0,1]$, and therefore
	we can use the Poincar\'e-Bendixson Theorem (see e.g.~\cite{TESCHL}) to conclude that any trajectory starting in~$(u_{\mathcal{M}},1]\times[0,1]$ must exit
	the set. Also, the inward normal velocity along the sides~$\{1\}\times
	(0,1]$ and~$(u_{\mathcal{M}},1)\times\{1\}$ is positive, and thus no trajectory can exit from these sides. Now, if a trajectory exits~$(u_{\mathcal{M}},1]\times[0,1]$
	from the side~$\{u_{\mathcal{M}}\}\times(0,1)$, then it enters the set~$\mathcal{X}$, and therefore~\eqref{ch2po0584yiugherghegkbis}
	is a consequence of~\eqref{ch2po0584yiugherghegk} in this case.
	If instead a trajectory exits~$(u_{\mathcal{M}},1]\times[0,1]$
	from the side~$(0,1)\times\{0\}$, then we directly obtain~\eqref{ch2po0584yiugherghegkbis}.
	
	{F}rom~\eqref{ch2po0584yiugherghegk} and~\eqref{ch2po0584yiugherghegkbis} we obtain~\eqref{ch21627},
	as desired.
	
	We now prove~\eqref{ch21628}, namely we show that
	\begin{equation}\label{ch2089ghvbdflpoiuytr}
	{\mbox{for all~$(u_0,v_0)\in \mathcal{Y}$ we have that~$(u(t), v(t))\to (0,1)$
			as~$t\to +\infty$.}}\end{equation}
	To this end, we observe that~$(u_s, v_s)$ (if~$0<ac<1$) and~$(0,0)$ are not in~$\mathcal{Y}$.
	Moreover,
	no trajectory starting in~$\mathcal{Y}$ converges to~$(u_s, v_s)$ (if~$0<ac<1$), nor to~$(0,0)$,
	since~$\mathcal{Y}$ does not contain points on~$\mathcal{M}$.
	
	In addition, recalling~\eqref{ch207960789djiewf2}, we have that
	there are no limit cycles in~$\mathcal{Y}$.
	As a consequence, by the
	Poincar\'e-Bendixson Theorem (see e.g.~\cite{TESCHL}),
	we have that
	every trajectory starting in~$\mathcal{Y}$ either go to~$(0,1)$ or it exits the set at some point of~$\partial \mathcal{Y}$.
	
	In the latter case,
	since no trajectory can cross~$\mathcal{M}$, the only possibility is that 
	the trajectory exits~$\mathcal{Y}$ at some
	point~$(u,v)\in\partial\big( [0,1]\times[0,1]\big)$. 
	We notice that,
	since~$\gamma$ is increasing, we have
	that~$\gamma(u)>0$ for all~$u>0$. As a consequence, 
	\begin{equation}\label{ch2ptouy988787}
	{\mbox{if~$(u,v)\in\overline{\mathcal{Y}}$,
			then~$v>\gamma(u)>0$ for all~$u>0$.}}\end{equation}
	Now, 
	thanks to Proposition~\ref{ch2prop:dyn}, 
	the only possibility that 
	a trajectory exits~$\mathcal{Y}$ at some
	point~$(u,v)\in\partial\big( [0,1]\times[0,1]\big)$ is for~$u>0$ and~$v=0$, 
	which would contradict~\eqref{ch2ptouy988787}.
	
	As a result, the only remaining possibility is that a trajectory
	in~$\mathcal{Y}$ converges to~$(0,1)$, which proves~\eqref{ch2089ghvbdflpoiuytr}.
	Hence, the proof of~\eqref{ch21628} is complete as well.
\end{proof}

With this, we are now able to complete the
proof of Theorem~\ref{ch2thm:dyn}:

\begin{proof}[Proof of (iv) of Theorem~\ref{ch2thm:dyn}]
	The statement in~(iv) of Theorem~\ref{ch2thm:dyn}
	is a direct consequence of the parametrization
	of the manifold~$\mathcal{M}$,
	as given by Proposition~\ref{ch2lemma:M} for~$ac\neq 1$ and by Proposition~\ref{ch2M:p045} for~$ac=1$,
	and the characterization of the sets~$\mathcal{B}$
	and~$\mathcal{E}$, as given by Proposition~\ref{ch2prop:char}.
\end{proof}

\section{Dependence of the dynamics
	on the parameters}\label{ch2ss:dependence}

In this section we discuss the dependence on the
parameters involved in the system~\eqref{ch2model}.

The dynamics of the system in~\eqref{ch2model} depends qualitatively
only on~$ac$, but of course the position of the saddle equilibrium
and the size and shape of the basins of attraction depend quantitatively
upon all the parameters. Here
we perform a deep analysis on each parameter separately.

We notice that the system in~\eqref{ch2model} does
not present a variational structure, due to the presence
of the
terms~$-acu$ in the first equation and~$-au$
in the second one, that are of first order in~$u$.
Thus, the classical methods of the calculus of variations cannot
be used and we have to make use of ad-hoc arguments,
of geometrical flavour.

\subsection{Dependence of the dynamics on the parameter~$c$}

We start by studying the dependence on~$c$, that represents the losses
(soldier death and missing births) caused
by the war for the first population with respect to the second one.
In the following proposition, we will express the dependence on~$c$
of the basin of attraction~$\mathcal{E}$
in~\eqref{ch2DEFE} by writing explicitly~$\mathcal{E}(c)$.

\begin{proposition}[Dependence of the dynamics on~$c$] \label{ch2prop:bhvc}
	With the notation in~\eqref{ch2DEFE}, we have that
	\begin{itemize}
		\item[(i)] If~$0< c_1 < c_2$, then 
		$\mathcal{E}(c_2) \subset \mathcal{E}(c_1)~$. 
		
		\item[(ii)] 	It holds that
		\begin{equation}\label{ch2242}
		\underset{c>0}{\bigcap} \, \mathcal{E}(c)= (0,1]\times \{0\}.
		\end{equation}
		
	\end{itemize}
\end{proposition}

We remark that the behavior for~$c$ sufficiently
small is given by~(i) of Theorem~\ref{ch2thm:dyn}: in this case,
there is a saddle point inside the domain~$[0,1]\times [0,1]$,
thus~$\mathcal{E}(c)\neq (0,1]\times [0,1]$. 
On the other hand, as~$c\to +\infty$, the set~$\mathcal{E}(c)$ gets
smaller and smaller until the first population has no
chances of victory if the second population has a positive size. 

The parameter~$c$ appears only in the first equation and it is
multiplied by~$-au$, that is always negative in the domain we
are interested in. Thus, the dependence on~$c$ is
independent of the other parameters. 
As one would expect, Proposition~\ref{ch2prop:bhvc}
tells us that
the greater the cost of the war for the first population,
the fewer possibilities of victory there are for it.

\begin{proof}[Proof of Proposition~\ref{ch2prop:bhvc}]
	{\emph{(i)}} We take~$c_2 > c_1 > 0$.
	According to Theorem~\ref{ch2thm:dyn},
	we denote by~$(u_s^{2}, v_s^{2})$
	the coexistence equilibrium for the parameter~${c_2}$ if~$a{c_2}<1$,
	otherwise we set~$(u_s^{2}, v_s^{2})=(0,0)$;
	similarly, we call~$(u_s^{1}, v_s^{1})$ the coexistence equilibrium
	for the parameter~${c_1}$ if~$a{c_1}<1$, and in the other cases we
	set~$(u_s^{1}, v_s^{1})=(0,0)$. 
	
	We observe that
	\begin{equation}\label{ch2poiuytrewq}
	v_s^{2} \leq v_s^{1}.\end{equation}
	Indeed, if~$a{c_2}<1$ then also~$a c_1<1$, and therefore,
	using the characterization in~\eqref{ch2usvs},
	\begin{equation*}
	\frac{\partial v_s}{\partial c} = \frac{-a(1+\rho {c}) -\rho 
		(1-a{c})  }{(1+\rho{c})^2} = \frac{-a-\rho}{(1+\rho{c})^2}  <0,
	\end{equation*}
	which implies~\eqref{ch2poiuytrewq} in this case.
	If instead~$a{c_2} \geq 1$ then the inequality in~\eqref{ch2poiuytrewq}
	is trivially satisfied, thanks to~(i), (ii) and~(iii)
	of Theorem~\ref{ch2thm:dyn}.
	
	Now, in the notation
	of Propositions~\ref{ch2lemma:M} (if~$ac\neq1$)
	and~\ref{ch2M:p045} (if~$ac=1$), thanks to the characterization
	in~\eqref{ch2char:E}, if we prove that 
	\begin{equation}\label{ch21907}
	\gamma_{c_1}(u)>\gamma_{c_2}(u)
	\quad \text{for any} \ u\in (0,\min \{ u_{\mathcal{M}}^{c_1}, u_{\mathcal{M}}^{c_2} \}],
	\end{equation}
	then the inclusion in~(i) is shown. Hence, we now
	focus on the proof of~\eqref{ch21907}.
	
	To this end, we observe that, since~$\gamma_{c_1}$
	is an increasing function, its inverse function~$f_{c_1}:[0, v_{\mathcal{M}}^{c_1}]\to
	[0,u_{\mathcal{M}}^{c_1}]$ is well defined and is increasing as well.
	In an analogue fashion, we define~$f_{c_2}(v)$ as the inverse
	of~$\gamma_{c_2}(u)$. We point out
	that the inequality in~\eqref{ch21907}
	holds true if 
	\begin{equation}\label{ch2eq:f{c_1}<f{c_2}}
	f_{c_1}(v)<f_{c_2}(v) \quad \text{for any} \ v\in (0,\min \{ v_{\mathcal{M}}^{c_1}, v_{\mathcal{M}}^{c_2} \}].
	\end{equation}
	Accordingly, we will show~\eqref{ch2eq:f{c_1}<f{c_2}}
	in three steps.
	
	First, in light of~\eqref{ch2poiuytrewq},
	we show that
	\begin{equation}\label{ch2yoiyjgsdn}
	{\mbox{the claim in~\eqref{ch2eq:f{c_1}<f{c_2}} is true
			in the interval~$[v_s^2, v_s^1]\cap(0,+\infty)$.}}\end{equation}
	For this, if~$ac_1\geq 1$, then also~$ac_2\ge1$, and therefore~$v_s^1=
	v_s^2=0$, thanks to~(ii) and~(iii)
	in Theorem~\ref{ch2thm:dyn}.
	Accordingly, in this case the interval~$[v_s^2, v_s^1]$
	coincides with the singleton~$\{ 0 \}$, and so there is nothing to prove.
	
	Otherwise, we
	recall that the curve~$u=\sigma(v)$, given in~\eqref{ch2f:sigma} and
	representing the points where~$\dot{v}=0$, is independent of~$c$.
	Moreover, thanks to formula~\eqref{ch2aggiunto} in
	Corollary~\ref{ch2lemma:M1} if~$ac<1$,
	formula~\eqref{ch2aggiun2} in
	Corollary~\ref{ch2lemma:M2}
	if~$ac>1$,
	and
	formula~\eqref{ch2aggiun2BIS} in Corollary~\ref{ch2lemma:Mdeg} if~$ac=1$
	(see also Figure~\ref{ch2fig:zone}), 
	we have that~$f_{c_1}(v)< \sigma(v)$ for~$v<v_s^1$
	and~$f_{c_2}(v)> \sigma(v)$ for~$v> v_s^2$, which
	proves~\eqref{ch2yoiyjgsdn} in the open
	interval~$(v_s^2,v_s^1)$.
	
	Moreover, it holds that
	\begin{equation}\label{ch2po0954757yuiewshfa}
	f_{c_1}(v_s^1)=\sigma(v_s^1)
	<f_{c_2}(v_s^1),\end{equation} and
	(if~$ac_2<1$, otherwise~$v_s^2=0$ and
	there is no need to perform this computation)
	\begin{equation}\label{ch2856sagd}
	f_{c_1}(v_s^2)<\sigma(v_s^2)=
	f_{c_2}(v_s^2).\end{equation}
	This completes the proof of~\eqref{ch2yoiyjgsdn}.
	
	Next we show that
	\begin{equation}\label{ch2yoiyjgsdn2}
	{\mbox{the claim in~\eqref{ch2eq:f{c_1}<f{c_2}} is true
			in the interval~$(0, v_s^{2})$.}}\end{equation}
	If~$ac_2\geq 1$, then~$v_s^2=0$, and so the claim in~\eqref{ch2yoiyjgsdn2}
	is trivial. Hence, we suppose that~$ac_2<1$ and we argue by contradiction,
	assuming
	that for some~$ v\in (0, v_s^{2})$
	it holds that~$f_{c_1}( v) \geq f_{c_2}( v)$.
	As a consequence, we can define
	$$ \bar v:=\sup\big\{v\in (0, v_s^{2})\;{\mbox{ s.t. }}\;
	f_{c_1}( v) \geq f_{c_2}( v)\big\}.
	$$
	We observe that, by continuity, we have that~$f_{c_1}(\bar{v}) 
	= f_{c_2}(\bar{v})$, and therefore, by~\eqref{ch2yoiyjgsdn}, we
	see that~$\bar{v} < v_s^{2}$.
	As a result, since~$f_{c_1}(v) <f_{c_2}(v)~$ for every~$
	v\in(\bar{v},v_s^{2}]$, then it holds that
	\begin{equation} \label{ch2eq:f{c_1}'<f{c_2}'}
	\frac{d f_{c_1}}{d v}(\bar{v}) < \frac{d f_{c_2}}{d v}(\bar{v}).
	\end{equation} 
	On the other hand, we can compute the derivatives by
	exploiting the fact that~$\gamma_{c_1}$ and~$\gamma_{c_2}$
	follow the flux for the system~\eqref{ch2model}. 
	Namely, setting~$\bar{u}:=f_{c_1}(\bar{v})$,
	we have that
	\begin{eqnarray*}
		&&\frac{d f_{c_1}}{d v}(\bar{v})
		= \frac{\dot{u}}{\dot{v}} (\bar{v})	
		= \frac{\bar{u} (1- \bar{u} - \bar{v}) - a{c_1} \bar{u}}{\rho \bar{v} (1- \bar{u} - \bar{v}) - a \bar{u}}
		\\
		{\mbox{and }}&&
		\frac{d f_{c_2}}{d v}(\bar{v})
		=  \frac{\dot{u}}{\dot{v}} (\bar{v})		
		= \frac{\bar{u} (1- \bar{u} - \bar{v}) - a{c_2} \bar{u}}{\rho\bar{v} (1- \bar{u} - \bar{v}) - a \bar{u}}
		.\end{eqnarray*} 
	Now, since~$\bar{v}\in[0,v_s^1)$, we have that~$
	\rho\bar{v} (1- \bar{u} - \bar{v}) - a \bar{u}>0$ (recall~\eqref{ch2aggiunto}
	and notice that~$(\bar u,\bar v)\in \mathcal{A}_4$).
	This and the fact that~${c_2}>{c_1}$ give that
	$$ \frac{d f_{c_1}}{d v}(\bar{v})> 
	\frac{d f_{c_2}}{d v}(\bar{v}),$$
	which is in contradiction with~\eqref{ch2eq:f{c_1}'<f{c_2}'},
	thus establishing~\eqref{ch2yoiyjgsdn2}.
	
	Now we prove that
	\begin{equation}\label{ch2yoiyjgsdn3}
	{\mbox{the claim in~\eqref{ch2eq:f{c_1}<f{c_2}} is true in the
			interval~$(v_s^{1},\min \{ u_{\mathcal{M}}^{c_1}, u_{\mathcal{M}}^{c_2} \}]$.}}\end{equation}
	Indeed, if~$ac_1< 1$, we argue towards a contradiction, supposing that
	there exists~$v>v_s^1$ such that~$f_{c_1}(v) \ge f_{c_2}(v)$. Hence, we can define
	$$\widehat v:= \inf \big\{ v>v_s^1\; {\mbox{ s.t. }}\;f_{c_1}(v) \ge f_{c_2}(v)
	\big\},$$
	and we deduce from~\eqref{ch2po0954757yuiewshfa} that~$\widehat v>v_s^1$.
	By continuity, we see that~$f_{c_1}(\widehat{v}) = f_{c_2}(\widehat{v})$.
	Therefore, since~$f_{c_1}(v) < f_{c_2}(v)~$ for any~$ v < \widehat{v}$, we conclude that
	\begin{equation} \label{ch2eq:f{c_1}'<f{c_2}'2BIS}
	\frac{d f_{c_1}}{dv}(\widehat{v}) > \frac{d f_{c_2}}{d v}(\widehat{v}).
	\end{equation} 
	On the other hand, setting~$\widehat{u}:=f_{c_1}(\widehat{v})$ and
	exploiting~\eqref{ch2model}, we get that
	\begin{align*}
	&\frac{d f_{c_1}}{d v}(\widehat{v})
	= \frac{\dot{u}}{\dot{v}} (\widehat{v})	
	= \frac{\widehat (1- \widehat{u} - \widehat{v}) - a{c_1}\widehat{u}}{\rho\widehat{v}
		(1- \widehat{u} -\widehat{v}) - a\widehat{u}}
	\\
	{\mbox{ and }}\qquad&	\frac{d f_{c_2}}{d v}(\widehat{v})
	=  \frac{\widehat{u}}{\widehat{v}} (\widehat{v})		
	= \frac{\rho\widehat{u} (1- \widehat{u} -\widehat{v}) - a{c_2} \widehat{u}}{
		\widehat{v} (1- \widehat{u} -\widehat{v}) - a \widehat{u}}.
	\end{align*} 
	Moreover, recalling~\eqref{ch2DEFA1234} and~\eqref{ch2aggiunto},
	we have that~$(f_{c_1}(\widehat{v}),\widehat{v})$
	and~$(f_{c_2}(\widehat{v}),\widehat{v})$ belong to the interior
	of~$\mathcal{A}_2$,
	and therefore~$\rho\widehat{v} (1- \widehat{u} -\widehat{v}) - a \widehat{u} <0$.
	This ad the
	fact that~${c_2}>{c_1}$ give that
	$$\frac{d f_{c_1}}{d v}(\widehat{v})<
	\frac{d f_{c_2}}{d v}(\widehat{v}),
	$$
	which is in contradiction with~\eqref{ch2eq:f{c_1}'<f{c_2}'2BIS}.
	This establishes~\eqref{ch2yoiyjgsdn3} in this case.
	
	If instead~$ac_1\ge1~$, then also~$ac_2\ge1$, and
	therefore we have that~$(u_s^{2}, v_s^{2})=(u_s^{1}, v_s^{1})=(0,0)$.
	In this setting, we use Propositions~\ref{ch2lemma:M}
	and~\ref{ch2M:p045}
	to say that at~$v=0$ the function~$f_{c_1}$ is tangent to the line~$(\rho-1+ac_1)v-au=0$, while~$f_{c_2}$ is tangent
	to~$(\rho-1+ac_2)v-au=0$. 
	Now, since
	\begin{equation*}
	\frac{\rho-1}{a} + c_1 < \frac{\rho-1}{a} + c_2,
	\end{equation*}
	we have that for positive~$v$ the second line is above the first one. 
	Also, thanks to the fact that~$f_{c_1}$ and ~$f_{c_2}$
	are tangent to these lines, we conclude that
	there exists~$\varepsilon>0$ 
	such that 
	\begin{equation}\label{ch22026}
	f_{c_1}(v)<f_{c_2}(v) \quad \text{for any } \ 	v<\varepsilon.
	\end{equation}
	
	Now, we suppose by contradiction that there exists some~$v> 0$
	such that~$f_{c_1}(v) \ge f_{c_2}(v)$. Hence, we can define
	$$\tilde v:= \inf \big\{ v>0\; {\mbox{ s.t. }}\;f_{c_1}(v) \ge f_{c_2}(v)
	\big\}.$$
	In light of~\eqref{ch22026}, we have that~$\tilde{v}\ge\varepsilon>0$.
	Moreover, by continuity, we see that~$f_{c_1}(\tilde{v}) = f_{c_2}(\tilde{v})$.
	Accordingly, since~$f_{c_1}(v) < f_{c_2}(v)~$ for any~$ v < \tilde{v}$, then it must be
	\begin{equation} \label{ch2eq:f{c_1}'<f{c_2}'2}
	\frac{d f_{c_1}}{dv}(\tilde{v}) > \frac{d f_{c_2}}{d v}(\tilde{v}).
	\end{equation} 
	On the other hand, setting~$\tilde{u}:=f_{c_1}(\tilde{v})$ and
	exploiting~\eqref{ch2model}, we see that
	\begin{align*}
	&\frac{d f_{c_1}}{d v}(\tilde{v})
	= \frac{\dot{u}}{\dot{v}} (\tilde{v})	
	= \frac{\tilde{u} (1- \tilde{u} - \tilde{v}) - a{c_1} \tilde{u}}{\rho\tilde{v} (1- \tilde{u} - \tilde{v}) - a \tilde{u}}
	\\
	{\mbox{ and }}\qquad&	\frac{d f_{c_2}}{d v}(\tilde{v})
	=  \frac{\tilde{u}}{\tilde{v}} (\tilde{v})		
	= \frac{\rho\tilde{u} (1- \tilde{u} - \tilde{v}) - a{c_2} \tilde{u}}{\tilde{v} (1- \tilde{u} - \tilde{v}) - a \tilde{u}}.
	\end{align*} 
	Now, thanks to~\eqref{ch2DEFA12} and~\eqref{ch2aggiun2},
	we have that~$(f_{c_1}(\tilde{v}),\tilde{v})$
	and~$(f_{c_2}(\tilde{v}),\tilde{v})$ belong to the interior
	of~$\mathcal{A}_2$,
	and therefore~$\rho\tilde{v} (1- \tilde{u} - \tilde{v}) - a \tilde{u} <0$. This ad the
	fact that~${c_2}>{c_1}$ give that
	$$\frac{d f_{c_1}}{d v}(\tilde{v})<
	\frac{d f_{c_2}}{d v}(\tilde{v}),
	$$
	which is in contradiction with~\eqref{ch2eq:f{c_1}'<f{c_2}'2}.
	This completes the proof of~\eqref{ch2yoiyjgsdn3}.
	
	Gathering together~\eqref{ch2yoiyjgsdn},~\eqref{ch2yoiyjgsdn2} and~\eqref{ch2yoiyjgsdn3},
	we obtain~\eqref{ch2eq:f{c_1}<f{c_2}}, as desired.
	\medskip
	
	{\emph{(ii)}} We first show that for all~$\varepsilon>0$ there exists~$c_{\varepsilon}>0$ such that for all~$c\ge c_{\varepsilon}$ it holds that
	\begin{equation} \label{ch2859}
	\mathcal{E}(c) \subset \big\{ 
	(u,v)\in [0,1]\times [0,1]\;{\mbox{ s.t. }}\; v < \varepsilon u    \big\}.
	\end{equation}
	The inclusion in~\eqref{ch2859} is also equivalent to 
	\begin{equation} \label{ch2255}
	\big\{ (u,v)\in [0,1]\times [0,1]\;{\mbox{ s.t. }}\; v > 
	\varepsilon u   \big\} \subset \mathcal{B}(c),
	\end{equation}
	and the strict inequality is justified by the fact that~$\mathcal{E}(c)$
	and~$\mathcal{B}(c)$ are separated by~$\mathcal{M}$, according
	to Proposition~\ref{ch2prop:char}.
	We now establish the inclusion in~\eqref{ch2255}.
	For this, let
	\begin{equation}\label{ch2255BIS}
	\mathcal{T}_{\varepsilon}:= \big\{ 
	(u,v)\in [0,1]\times [0,1] \;{\mbox{ s.t. }}\; v >  \varepsilon u    \big\}. \end{equation}
	Now, we can choose~$c$ large enough such that the condition~$ac\geq 1$
	is fulfilled. In this way, thanks to~(ii) and~(iii)
	of Theorem~\ref{ch2thm:dyn},
	the only equilibria are the points~$(0,0)$ and~$(0,1)$.
	
	Now, the component of the
	velocity in the inward normal direction to~$\mathcal{T}_{\varepsilon}$
	on the side~$\{v=\varepsilon u\}$ is given by
	\begin{eqnarray*}
		&&(\dot u,\dot v)\cdot \frac{(-\varepsilon,1)}{\sqrt{1+\varepsilon^2}}=
		\frac{\dot{v}-\varepsilon \dot{u}}{\sqrt{1+\varepsilon^2}}
		\\&&\qquad =\frac1{{\sqrt{1+\varepsilon^2}}}\big(
		\rho v(1-u-v) -au -\varepsilon u(1-u-v) + \varepsilon acu \big)\\
		&&\qquad=\frac1{{\sqrt{1+\varepsilon^2}}}\big[
		(\rho v-\varepsilon u)(1-u-v)  +  (\varepsilon c -1)au\big]\\
		&&\qquad=\frac1{{\sqrt{1+\varepsilon^2}}}\big[
		(\rho \varepsilon u-\varepsilon u)(1-u-\varepsilon u)  +  (\varepsilon c -1)au\big]	 ,
	\end{eqnarray*}
	that is positive for
	\begin{equation}\label{ch2possibly}
	c > c_{\varepsilon} := \frac{2\varepsilon(1+\rho) +a}{\varepsilon a}.
	\end{equation}
	This says that no trajectory in~$\mathcal{T}_{\varepsilon}$ can exit~$\mathcal{T}_{\varepsilon}$ from the side~$\{v=\varepsilon u\}$.
	
	The other parts of~$\partial \mathcal{T}_{\varepsilon}$ belong to~$\partial(
	(0,1)\times(0,1))$
	but not to~$[0,1]\times \{0 \}$.
	As a consequence, by Proposition~\ref{ch2prop:dyn},
	\begin{equation}\label{ch2po123097}
	{\mbox{every trajectory in~$\mathcal{T}_{\varepsilon}$ belongs to~$\mathcal{T}_{\varepsilon}$
			for all~$t\ge0$.}}\end{equation}	
	{F}rom this,~\eqref{ch207960789djiewf2} and the 
	Poincar\'e-Bendixson Theorem (see e.g.~\cite{TESCHL}), we conclude that
	the~$\omega$-limit of any trajectory starting in~$\mathcal{T}_{\varepsilon}$
	can be either an equilibrium or a union of (finitely many)
	equilibria and non-closed orbits connecting these equilibria.
	
	Now, we claim that, possibly taking~$c$ larger in~\eqref{ch2possibly},
	\begin{equation}\label{ch2po1230972}
	\mathcal{M}\subset \big([0,1]\times[0,1]\big)\setminus\mathcal{T}_{\varepsilon}.
	\end{equation}
	Indeed, suppose by contradiction that there exists~$(\tilde u,\tilde v)\in\mathcal{M}
	\cap\mathcal{T}_{\varepsilon}$. Then, in light of~\eqref{ch2po123097}, a trajectory passing
	through~$(\tilde u, \tilde v)$ and converging to~$(0,0)$ has to be entirely contained
	in~$\mathcal{T}_{\varepsilon}$.
	
	On the other hand, by Propositions~\ref{ch2lemma:M} and~\ref{ch2M:p045},
	we know that at~$u=0$ the manifold~$\mathcal{M}$ is tangent to the
	line~$(\rho-1+ac)v-au=0$.  
	Hence, if we choose~$c$ large enough such that
	$$ \frac{a}{\rho-1+ac}<\varepsilon,$$
	we obtain that this line is below the line~$v=\varepsilon u$, thus reaching
	a contradiction. This establishes~\eqref{ch2po1230972}.
	
	{F}rom~\eqref{ch2po1230972}, we deduce that,
	given~$(\tilde u,\tilde v)\in\mathcal{T}_{\varepsilon}$,
	and denoting~$\omega_{(\tilde u,\tilde v)}$ the~$\omega$-limit of~$(\tilde u,\tilde v)$,
	\begin{equation}\label{ch2podnjewbf215}
	\omega_{(\tilde u,\tilde v)}\neq \{(0,0)\},
	\end{equation}
	provided that~$c$ is taken large enough.
	
	Furthermore,~$\omega_{(\tilde u,\tilde v)}$ cannot consist of the two equilibria~$(0,0)$
	and~$(0,1)$ and non-closed orbits connecting these equilibria, since~$(0,1)$
	is a sink. As a consequence of this and~\eqref{ch2podnjewbf215},
	we obtain that~$\omega_{(\tilde u,\tilde v)}=\{(0,1)\}$
	for any~$(\tilde u,\tilde v)\in\mathcal{T}_{\varepsilon}$, provided that~$c$
	is large enough.
	
	Thus, recalling~\eqref{ch2DEFB} and~\eqref{ch2255BIS}, this proves~\eqref{ch2255},
	and therefore~\eqref{ch2859}.
	
	Now, using~\eqref{ch2859}, we see that for every~$\varepsilon>0$,
	$$\underset{c>0}{\bigcap}  \mathcal{E}(c)\subseteq \mathcal{E}(c_{\varepsilon})
	\subseteq \big\{
	(u,v)\in [0,1]\times [0,1]\; {\mbox{ s.t. }}\; v < \varepsilon u   \big \}.
	$$
	Accordingly,
	\begin{equation*}
	\underset{c>0}{\bigcap}  \mathcal{E}(c) \subseteq  \underset{\varepsilon>0}{\bigcap}  \big\{
	(u,v)\in [0,1]\times [0,1]\; {\mbox{ s.t. }}\; v < \varepsilon u   \big \} = (0,1] \times \{0\},
	\end{equation*}
	which implies~\eqref{ch2242}, as desired.
\end{proof}

\subsection{Dependence of the dynamics on the parameter~$\rho$}

Now we analyze the dependence of the dynamics on the parameter~$\rho$, that is the fitness of the second population~$v$ with respect to the fitness of the first one~$u$.

In the following proposition, we will make it explicit the dependence on~$\rho$ by
writing~$\mathcal{E}(\rho)$ and~$\mathcal{B}(\rho)$.

\begin{proposition}[Dependence of the dynamics on~$\rho$] \label{ch2prop:bhvA}
	With the notation in~\eqref{ch2DEFB} and~\eqref{ch2DEFE}, we have that
	\begin{itemize}
		\item[(i)]
		When~$\rho=0$, for any~$v \in [0,1]$ the point~$(0,v)$  is an equilibrium. If~$v\in(1-ac,1]$, then it corresponds to a strictly
		negative eigenvalue and a null one.
		If instead~$v\in[0,1-ac)$, then it corresponds to a strictly
		positive eigenvalue and a null one
		
		Moreover, 
		\begin{equation}\label{ch2first}
		\mathcal{B}(0)=  \varnothing,\end{equation} and for any~$\varepsilon<  ac/2~$ and
		any~$\delta< \varepsilon c/2$  we have that
		\begin{equation}\label{ch2first2}
		[0,1]\times [0,1-ac) \subseteq \mathcal{E}(0)  \subseteq \mathcal{T}_{\varepsilon, \delta}  ,
		\end{equation} where
		\begin{equation}\label{ch2TEPS}
		\mathcal{T}_{\varepsilon, \delta}:=\big\{ (u,v)\in[0,1] \times [0,1]\;
		{\mbox{ s.t. }}\; \delta v-\varepsilon u \leq \delta(1-\varepsilon) \big\}.
		\end{equation}
		\item[(ii)]
		For any~$\varepsilon<  ac/3~$ and any~$\delta< \varepsilon c/2$ it holds that 
		\begin{equation*}
		\underset{a>0}{\bigcap} \ \underset{{0<\rho<a/3}}{\bigcup} \mathcal{E}(\rho) \subseteq  \mathcal{T}_{\varepsilon, \delta} ,
		\end{equation*}
		where~$ \mathcal{T}_{\varepsilon, \delta}$ is defined in~\eqref{ch2TEPS}.
		\item[(iii)]
		It holds that 
		\begin{equation}\label{ch2ir4t4y4y}
		\underset{\omega>0}{\bigcap} \, \underset{{\rho>\omega}}{\bigcup} \mathcal{E}(\rho) =  (0,1] \times \{0\}.
		\end{equation}	
	\end{itemize}
\end{proposition}

We point out that
the case~$\rho =0$ is not comprehended in Theorem~\ref{ch2thm:dyn}.
As a matter of fact, the dynamics of this case is qualitatively very different from all the other cases. Indeed, for~$\rho =0$ the domain~$[0,1] \times [0,1]$ is not divided into~$\mathcal{E}$ and~$\mathcal{B}$, since more attractive equilibria appear on the line~$\{0\}\times(0,1)$. Thus, even if the second population cannot grow, it still has some chance of victory. 

As soon as~$\rho~$ is positive, on the line~$u=0$ only the equilibrium~$(0,1)$ survives, and it attracts all the points that were going to the line~$\{0\}\times(0,1)$ for~$\rho =0$.

When~$\rho \to +\infty$, the basin of attraction of~$(0,1)$ tends to invade the domain, thus the first population tends to have almost no chance of victory and the second population tends to win.
However, the dependence on the parameter~$\rho~$ is not monotone as one could think, at least not in~$[0,+\infty)\times[0,+\infty)$.

\begin{figure}[h] 
	
	\begin{subfigure}{0.5\textwidth} \label{ch2A3to01}
		\includegraphics[width=8cm, height=5cm]{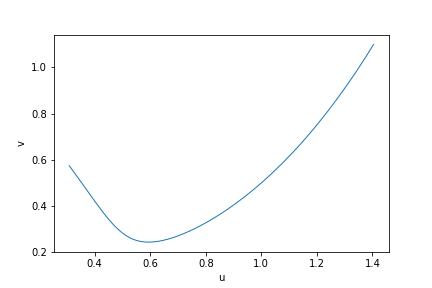} 
		\caption{~$a=0.2$,~$c=0.1$, and~$\rho=3$}
		\label{ch2fig:subim1}
	\end{subfigure}
	\begin{subfigure}{0.5\textwidth} \label{ch2A7notto01}
		\includegraphics[width=8cm, height=5cm]{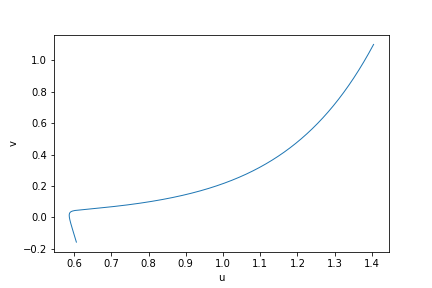}
		\caption{$a=0.2$,~$c=0.1$, and~$\rho=7$}
		\label{ch2fig:subim2}
	\end{subfigure}
	
	\caption{\it Figure (a) and Figure (b) show the trajectory starting from the point~$(u_0,v_0)=(1.4045, 1.1)$ for~$\rho=3$ and~$\rho=7$ respectively. For~$\rho=3$ the trajectory leads to the equilibrium~$(0,1)$, so~$(u_0,v_0)\notin \mathcal{E}(\rho=3)$, while for~$\rho=7$ the second population goes to extinction in finite time, so~$(u_0,v_0)\in \mathcal{E}(\rho=7)$. 
	}
	\label{ch2fig:trajrho}
\end{figure}

Indeed, by performing some simulation, one could find some values~${\rho}_1$ and~${\rho}_2$, with~$0<{\rho}_1 < {\rho}_2$, and a point~$(u^*, v^*)\in [0,+\infty)\times[0,+\infty)$ such that~$(u^*, v^*) \notin \mathcal{E}({\rho}_1)$ and~$(u^*, v^*) \in \mathcal{E}({\rho}_2)$, see~Figure~\ref{ch2fig:trajrho}. 

This means that, sometimes, a big value of fitness for the second population may lead to extinction while a small value brings to victory. This is counterintuitive, but can be easily explained: the parameter~$\rho~$ is multiplied by the term~$1-u-v$, that is negative past the counterdiagonal of the square~$[0,1]\times[0,1]$. So in the model~\eqref{ch2model}, as well as in any model of Lotka-Volterra type, the population that grows faster is also the one that suffers more the consequences of overpopulation. Moreover, the usual dynamics of Lotka-Volterra models is altered by the presence of the term~$-au$, and this leads to the lack of monotonicity that we observe.
\medskip

We now give the proof of Proposition~\ref{ch2prop:bhvA}:

\begin{proof}[Proof of Proposition~\ref{ch2prop:bhvA}]
	{\emph{(i)}}
	For~$\rho=0$, the equation~$\dot{v}=0$ collapses to~$u=0$. Since for~$u=0$ also the equation~$\dot{u}=0$ is satisfied,
	each point on the line~$u=0$ is an equilibrium. 
	
	Calculating the eigenvalues for the points~$(0, \tilde{v})$, with~$\tilde v\in[0,1]$,
	using the Jacobian matrix in~\eqref{ch2Jmatrix},
	one gets the values~$0$ and~$1-ac-\tilde{v}$.
	Accordingly, this entail that, if~$\tilde{v} < 1-ac$, the point~$(0, \tilde{v})$ corresponds to a strictly negative eigenvalue
	and a null one, while if~$\tilde{v}>1-ac$
	then~$(0, \tilde{v})$ corresponds to a strictly negative
	eigenvalue and a null one.
	These considerations proves the first statement in~(i).
	
	We notice also that in the whole square~$(0,1]\times[0,1]$ we have~$\dot{v}= -au < 0$, hence there is no trajectory that can go to~$(0,1)$, and there is no cycle.
	In particular this implies~\eqref{ch2first}.
	
	Now, we observe that
	on the side~$[0,1]\times \{1\}$ the inward normal derivative is 
	given by~$-\dot{v}=au$, which is nonnegative, and therefore
	a trajectory cannot exit the square on this side.
	Similarly, along the side~$\{1\}\times [0,1]$ the inward normal derivative is
	given by~$-\dot{u}=v+ac$, which is positive, hence
	a trajectory cannot exit the square on this side either.
	
	The side~$\{0\}\times[0,1]$ is made of equilibrium points at which the first population~$u$ is extinct, while on the side~$(0,1]\times\{0\}$ we have extinction of the population~$v$.
	Thus a trajectory either converges to one of the equilibria on the
	side~$\{0\}\times[0,1]$, or exits~$[0,1]\times[0,1]$ through the side~$(0,1]\times\{0\}$.
	
	In particular, since~$\{0\}\times [0,1-ac)$ consists of repulsive equilibria, we have that
	$$ [0,1]\times[0, 1-ac)\subseteq \mathcal{E}(0),~$$
	that is, trajectories starting in~$[0,1]\times[0, 1-ac)$ go to the extinction of~$v$.
	This proves the first inclusion in~\eqref{ch2first2}.
	
	To prove the second inclusion in~\eqref{ch2first2},
	we first show that
	\begin{equation}\label{ch2pot43 yb9 49y}
	{\mbox{points in~$\big([0,1]\times[0,1]\big)\setminus\mathcal{T}_{\varepsilon, \delta}$
			are mapped into~$\big([0,1]\times[0,1]\big)\setminus\mathcal{T}_{\varepsilon, \delta}$ itself.}}\end{equation}
	Indeed, 
	on the line~$\{\delta v-\varepsilon u = \delta(1-\varepsilon)\}$ we have that the inward-pointing normal derivative is given by
	\begin{equation}\begin{split}\label{ch2fagiano}
	&(\dot{u},\dot{v})\cdot\frac{(-\varepsilon,\delta)}{\sqrt{\varepsilon^2+\delta^2}}=
	\frac1{\sqrt{\varepsilon^2+\delta^2}}	\big(\delta \dot{v}- \varepsilon \dot{u}\big)\\
	&\qquad\qquad =\frac1{\sqrt{\varepsilon^2+\delta^2}}\big( -\delta a u - \varepsilon u(1-u-v) +\varepsilon ac u\big)\\ 
	&\qquad\qquad=
	\frac{ u}{\sqrt{\varepsilon^2+\delta^2}}\left[\varepsilon\left(-1+ac+u+\frac{\varepsilon}{\delta}u +1-\varepsilon\right)-\delta a\right] \\
	&\qquad\qquad
	=\frac{ 1}{\sqrt{\varepsilon^2+\delta^2}}\left[u^2\left( 1+\frac{\varepsilon}{\delta}\right) + u(\varepsilon ac -\delta a -\varepsilon^2)\right].
	\end{split}\end{equation}
	The first term is always positive; the second one is positive for the choice \begin{equation*}
	\delta < \frac{\varepsilon c}{2} \quad {\mbox{ and }}\quad\varepsilon< \frac{ac}{2}.
	\end{equation*} 
	Hence, under the assumption in~(i), on the line~$\{\delta v-\varepsilon u = \delta(1-\varepsilon)\}$ the inward-pointing normal derivative is positive, which implies
	that no trajectories in~$\big([0,1]\times[0,1]\big)\setminus\mathcal{T}_{\varepsilon, \delta}$ can exit
	from~$\big([0,1]\times[0,1]\big)\setminus\mathcal{T}_{\varepsilon, \delta}$. This establishes~\eqref{ch2pot43 yb9 49y}.
	
	As a consequence of~\eqref{ch2pot43 yb9 49y}, we obtain also the second
	inclusion~\eqref{ch2first2}, as desired.
	\medskip
	
	{\emph{(ii)}} 
	We claim that
	\begin{equation}\label{ch2poi8562eq2dvfkgjlkuykyuasdawre}
	\big([0,1]\times[0,1]\big)\setminus\mathcal{T}_{\varepsilon, \delta}
	\subseteq \mathcal{B}(\rho),\end{equation}
	for all~$0<\rho< a/3$.
	To this end, we observe that, in order to determine
	the sign of the inward pointing normal derivative on the side~$\{\delta v -\varepsilon u = \delta(1-\varepsilon)\}$, by~\eqref{ch2fagiano} we have to 
	check that~$\delta \dot{v}- \varepsilon\dot{u}\ge0$. In order to simplify the calculation, we use the change of coordinates~$x:=u$ and~$y:=1-v$.
	In this way, one needs to verify that~$\delta \dot{y}+\varepsilon \dot{x} \leq 0$
	on the line~$\{\delta y + \varepsilon x = \delta \varepsilon\}$. For this, we compute
	\begin{equation}
	\label{ch2cneq0}
	\begin{split}
	\delta \dot{y}+\varepsilon \dot{x} &= \delta \rho (y-1)(y-x)+\delta a x + \varepsilon x (y-x) - \varepsilon acx, \\
	& = -\delta \rho (1-y) y+x \big(
	\delta \rho (1-y) +\delta a + \varepsilon (y-x) -\varepsilon ac   \big) ,\\
	& = -\delta \rho (1-y) y +  x \big(  \delta \rho  -\delta \rho y 
	+\delta a + \varepsilon y-\varepsilon x -\varepsilon a c   \big)\\
	&\le  x \big(  \delta \rho  -\delta \rho y 
	+\delta a + \varepsilon y-\varepsilon x -\varepsilon a c   \big)  .
	\end{split}
	\end{equation}
	Now we choose~$\delta<\varepsilon c / 2$ and we
	recall that~$\rho < a/3$. Moreover, we notice
	that
	$$y= \varepsilon-\frac{\varepsilon}{\delta}x\le\varepsilon,
	$$ 
	and therefore~$\varepsilon y \leq \varepsilon^2$. Thus, we have that
	\begin{eqnarray*}
		-\delta \rho y + \delta \rho  +\delta a + \varepsilon y-\varepsilon x 
		-\varepsilon a c  \le \frac{\varepsilon ac }{6} + \frac{\varepsilon ac }{2} +
		\varepsilon^2
		-\varepsilon a c= \varepsilon\left( \frac{2}{3} ac  + \varepsilon
		- a c\right)
	\end{eqnarray*}
	that is negative for~$\varepsilon < ac/3$. Plugging this
	information into~\eqref{ch2cneq0}, we obtain
	that~$\delta \dot{y}+\varepsilon \dot{x}  \leq 0~$, as desired.
	
	This proves that
	trajectories in~$
	\big([0,1]\times[0,1]\big)\setminus\mathcal{T}_{\varepsilon, \delta}$
	cannot exit~$
	\big([0,1]\times[0,1]\big)\setminus\mathcal{T}_{\varepsilon, \delta}$.
	This, the fact that there are no cycles in~$[0,1]\times[0,1]$
	and the Poincar\'e-Bendixson Theorem (see e.g.~\cite{TESCHL})
	give that
	trajectories in~$\big([0,1]\times[0,1]\big)\setminus\mathcal{T}_{
		\varepsilon, \delta}$ converge to~$(0,1)$,
	that is the only equilibrium in~$\big([0,1]\times[0,1]\big)\setminus
	\mathcal{T}_{\varepsilon, \delta}$. Hence, 
	\eqref{ch2poi8562eq2dvfkgjlkuykyuasdawre}
	is established.
	
	{F}rom~\eqref{ch2poi8562eq2dvfkgjlkuykyuasdawre}
	we deduce that
	$$ \mathcal{E}(\rho)\subseteq\mathcal{T}_{\varepsilon, \delta}$$
	for all~$0<\rho<a/3$, which implies the desired result
	in~(ii).
	\medskip
	
	{\emph{(iii)}} 
	We consider~$\varepsilon_1>\varepsilon_2 >0$ to be
	taken sufficiently small in what follows,
	and we show that
	there exists~$R>0$, depending on~$\varepsilon_1$
	and~$\varepsilon_2$, such that for
	all~$\rho\geq R$ it holds that
	\begin{equation}\label{ch2qeruyjy8790}
	\mathcal{R}_{\varepsilon_1, \varepsilon_2}:=
	[0, 1-\varepsilon_1]\times [\varepsilon_2,1] \subseteq \mathcal{B}(\rho).\end{equation}
	For this, we first observe that 
	\begin{equation}\label{ch2po089egdgdkjfkghjighywrv58465v8}
	{\mbox{no trajectory starting
			in~$\mathcal{R}_{\varepsilon_1, \varepsilon_2}$
			can exit the set.}}\end{equation}
	Indeed, looking at the velocity fields on the sides~$\{0\}\times
	[\varepsilon_2, 1]$ and~$[0,1-\varepsilon_1]\times\{1\}$,
	one sees that no trajectory in~$\mathcal{R}_{\varepsilon_1, \varepsilon_2}$ can exit from these sides.
	
	Moreover,
	on the side~$\{1-\varepsilon_1\} \times [\varepsilon_2, 1]$, the normal inward derivative is
	\begin{equation*}
	-\dot{u}=-[u(1-u-v)-acu] = -(1-\varepsilon_1)(\varepsilon_1-v-ac),
	\end{equation*}
	and this is positive for~$\varepsilon_1\leq ac$ (which is fixed
	from now on).
	In addition, on the side~$[0,1-\varepsilon_1]\times\{ \varepsilon_2 \}$, the inward normal derivative is
	\begin{eqnarray*}
		&&	\dot{v}= [\rho v(1-u-v)-au] = 
		\rho \varepsilon_2(1-u-\varepsilon_2) - au\\&&\qquad
		\ge	\rho \varepsilon_2(\varepsilon_1-\varepsilon_2) - a(1-
		\varepsilon_1),
	\end{eqnarray*}
	and this is positive
	for 
	\begin{equation}\label{ch2rhodef}
	\rho > \frac{a(1-\varepsilon_1)}{\varepsilon_2(\varepsilon_1
		-\varepsilon_2)}=:R.\end{equation}
	These observations complete the proof
	of~\eqref{ch2po089egdgdkjfkghjighywrv58465v8}.
	
	{F}rom~\eqref{ch207960789djiewf2},
	\eqref{ch2po089egdgdkjfkghjighywrv58465v8}
	and the Poincar\'e-Bendixson Theorem (see e.g.~\cite{TESCHL}), we have that
	all the trajectories in the interior 
	of~$\mathcal{R}_{\varepsilon_1, \varepsilon_2}$
	must converge to
	either an equilibrium or a union of (finitely many)
	equilibria and non-closed orbits connecting these equilibria.
	
	In addition, we claim that, if~$0<ac<1$, recalling~\eqref{ch2usvs}
	and possibly enlarging~$\rho$ 
	in~\eqref{ch2rhodef},
	\begin{equation}\label{ch2possibly2}
	(u_s,v_s)\notin \mathcal{R}_{\varepsilon_1, \varepsilon_2}.
	\end{equation}
	Indeed, we have that~$u_s \to 1-ac$ and~$v_s \to 0$,
	as~$\rho \to +\infty$. Hence, we can choose~$\rho$ large enough
	such that the statement in~\eqref{ch2possibly2} is satisfied.
	
	As a consequence of~\eqref{ch2possibly2},
	we get that all the trajectories in the interior 
	of~$\mathcal{R}_{\varepsilon_1, \varepsilon_2}$
	must converge to
	the equilibrium~$(0,1)$,
	and this establishes~\eqref{ch2qeruyjy8790}.
	
	Accordingly,~\eqref{ch2qeruyjy8790} entails that, for~$\varepsilon_1>
	\varepsilon_2>0$ sufficiently small, there exists~$R>0$, depending on~$\varepsilon_1$
	and~$\varepsilon_2$, such that for
	all~$\rho\geq R$
	\begin{equation*}
	\mathcal{E}(\rho)
	\subset\big((0,1]\times[0,\varepsilon_2)\big)\cup
	\big( (1-\varepsilon_1,1]\times(\varepsilon_2,1]\big)
	.\end{equation*}
	This implies~\eqref{ch2ir4t4y4y}, as desired.
\end{proof}

\subsection{Dependence of the dynamics on the parameter~$a$}

The consequences of the lack of variational structure
become even more extreme when we observe the dependence
of the dynamics on
the parameter~$a$, that is the aggressiveness of the first population
towards the other. 
Throughout this section, we take~$\rho>0$ and~$c>0$,
and we perform our analysis
taking into account the limit cases~$a\to0$ and~$a\to+\infty$.
We start analyzing the dynamics of~\eqref{ch2model}
in the case~$a=0$.

\begin{proposition}[Dynamics of~\eqref{ch2model} when~$a=0$] \label{ch2prop:bhvaPRE}
	For~$a=0$
	the system~\eqref{ch2model} has the following
	features:
	\begin{itemize}
		\item[i)] The system has the equilibrium~$(0,0)$, which is a source,
		and a straight line of equilibria~$(u,1-u)$, for all~$u\in[0,1]$,
		which correspond to a strictly negative eigenvalue and a null one.
		\item[ii)] Given any~$(u(0), v(0))\in (0,1)\times(0,1)$ we have that
		\begin{equation}\label{ch2form}
		(u(t), v(t)) \to(\bar{u}, 1-\bar{u})\quad{\mbox{ as }}t\to+\infty, 
		\end{equation}
		where~$\bar{u}$ satisfies
		\begin{equation}\label{ch21650}
		\frac{v(0) }{u^{\rho}(0)}\bar{u}^{\rho} + \bar{u} -1=0.
		\end{equation}
		\item[iii)] The equilibrium~$(u_s^0, v_s^0)$ given in~\eqref{ch2u0v0}
		has a stable manifold, which can be written as the graph of an
		increasing smooth function~$\gamma_0:[0,u_{\mathcal{M}}^0]\to[0,v_{\mathcal{M}}^0]$,
		for some~$(u_{\mathcal{M}}^0,v_{\mathcal{M}}^0)\in\big(\{1\}\times[0,1]\big)\cup
		\big((0,1]\times\{1\}\big)$, such that~$\gamma_0(0)=0$,
		$\gamma_0(u_{\mathcal{M}}^0)=v_{\mathcal{M}}^0$.
		
		More precisely, 
		\begin{equation}\label{ch2def:gamma0}
		\gamma_0 (u):= \frac{v_s^0}{(u_s^0)^{\rho}} u^{\rho} \quad
		{\mbox{ and }}\quad
		u_{\mathcal{M}}^0:=\min \left\{1, \frac{u_s^0}{(v_s^0)^{\frac{1}{{\rho}}}}\right\},	\end{equation} 
		being~$(u_s^0,v_s^0)$ defined in~\eqref{ch2u0v0}.
	\end{itemize}
\end{proposition}

We point out that formula~\eqref{ch2form}
says that for~$a=0$ every point in the interior
of~$[0,1]\times[0,1]$ tends to a coexistence equilibrium.
The shape of the trajectories depends on~$\rho$, being
convex in the case~$\rho>1$,
a straight line in the case~$\rho=1$, and concave in the case~$\rho< 1$. This means that if the second population~$v$ is alive at the beginning, then it does not get extinct in finite time.

\begin{proof}[Proof of Proposition~\ref{ch2prop:bhvaPRE}]
	{\emph{(i)}}	
	For~$a=0$, we look for the equilibria of the system~\eqref{ch2model}
	by studying when~$\dot{u}=0$ and~$\dot{v}=0$. It is easy to see that
	the point~$(0,0)$ and all the points on the line~$u+v=1$ are the only equilibria.
	
	The Jacobian of the system (see~\eqref{ch2Jmatrix}, with~$a=0$) at the point~$(0,0)$
	has two positive eigenvalues,~$1$ and~$\rho~$, and thereofore~$(0,0)$ is a
	source.
	
	Furthermore,
	the characteristic polynomial at a point~$(\tilde{u}, \tilde{v})$ on the line~$u+v=1$
	is given by
	$$(\lambda+\tilde{u})(\lambda+\rho \tilde{v})-\rho \tilde{u}\tilde{v}
	=\lambda(\lambda+\tilde{u} +\rho \tilde{v}),$$
	and therefore, the eigenvalues are~$0$ and~$-\tilde{u} -\rho \tilde{v}<0$.
	\medskip
	
	{\emph{(ii)}}  We point out that
	when~$a=0$ 
	\begin{equation}\label{ch2intprim678}
	{\mbox{$\mu(t):=v(t)/u^{\rho} (t)$ is a prime integral for the system.}}
	\end{equation}
	Indeed,
	\begin{equation*}
	\mu'= \frac{\dot{v}u^\rho - {\rho} u^{{\rho}-1}  \dot{u} v }{u^{2{\rho}}}=u^{{\rho}-1} \frac{{\rho}uv(1-u-v)- {\rho} uv(1-u-v) }{u^{2{\rho}}}=0.
	\end{equation*}  
	As a result, the trajectory starting at a point~$(u(0),v(0) )\in(0,1)\times(0,1)$ lies on the curve
	\begin{equation}\label{ch2intprim}
	v(t)=\frac{v(0)}{ u^{\rho}(0)}\,  u^{\rho} (t).\end{equation}
	Moreover, the trajectory starting at~$(u(0),v(0) )$ is asymptotic as~$t\to+\infty$
	to an equilibrium on this curve. Since~$(0,0)$ is a source, the only possibility
	is that the trajectory starting at~$(u(0),v(0) )$ converges to an
	equlibrium~$(\bar{u}, \bar{v})$ 
	such that
	$\bar{v}=1-\bar{u}$. This entails that
	\begin{equation*}
	1-\bar{u} =\bar{v}=(v(0)/ u^{\rho}(0))  \bar{u}^{\rho},
	\end{equation*}
	which is exactly equation~\eqref{ch21650}.	
	
	\medskip
	
	{\emph{(iii)}} We observe that the point~$(u_s^0, v_s^0)$ given in~\eqref{ch2u0v0}
	lies on the straight line~$u+v=1$, and therefore, thanks to~(i) here, it is
	an equilibrium of the system~\eqref{ch2model}, which corresponds
	to a strictly negative eigenvalue~$-u_s^0-\rho v_s^0$ and a null one.
	
	Hence, by the Center Manifold Theorem
	(see e.g. Theorem~1 on page~16
	of~\cite{MR635782}), the point~$(u_s^0, v_s^0)$ 
	has a stable manifold, which has dimension~$1$ and
	is tangent to the eigenvector of the linearized system associated to the
	strictly negative eigenvalue~$-u_s^0-\rho v_s^0$.
	
	Also, the graphicality and the monotonicity
	properties follow from the strict sign of~$\dot{u}$
	and~$\dot{u}$. The smoothnes of the graphs follows from the smoothness
	of the center manifold.
	The fact that~$\gamma_0(0)=0$ is a consequence of the monotonicity
	property of~${u}$ and~${v}$, which ensures that the limit at~$t\to-\infty$
	exists, and the fact that this limit has to lie on the prime integral
	in~\eqref{ch2intprim}.
	The fact that~$\gamma_0(u_{\mathcal{M}}^0)=v_{\mathcal{M}}^0$
	follows from formula~\eqref{ch2form} and the monotonicity property.
	Formula~\eqref{ch2def:gamma0}
	follows from the fact that any trajectory has to lie on the prime integral
	in~\eqref{ch2intprim}.
\end{proof}


To state our next result concerning the dependence of the basin of
attraction~$\mathcal{E}$ defined in~\eqref{ch2DEFE} on the parameter~$a$,
we give some notation.
We will make it explicit the dependence of the sets~$\mathcal{E}$
and~$\mathcal{B}$
on the parameter~$a$, by writing
explicitly~$\mathcal{E}(a)$ and~$\mathcal{B}(a)$, and we will call
\begin{equation*}
\mathcal{E}_0:=\underset{a'>0}{\bigcap} \, \underset{a'>a>0}{\bigcup} \mathcal{E}(a)
\end{equation*}
and
\begin{equation}\label{ch2def:Einfty}
\mathcal{E}_{\infty}:=\underset{a'>0}{\bigcap} \, \underset{a>a'}{\bigcup} \mathcal{E}(a).
\end{equation}
In this setting, we have the following statements:

\begin{proposition}[Dependence of the dynamics on~$a$] \label{ch2prop:bhva}
	\quad
	\begin{itemize}
		\item[(i)] We have that
		\begin{equation} \label{ch2char:E0}
		\begin{split}
		&	 \big\{  (u,v)\in [0,1]\times [0,1]\;{\mbox{ s.t. }}\; 
		v < \gamma_{0}(u) \,\text{ if } \, u\in[0, u_{\mathcal{M}}^0]\\
		&\qquad\qquad\qquad\qquad{\mbox{and }}		\;
		v \leq 1  \, \text{ if }\,  u\in(u_{\mathcal{M}}^0, 1]  \big\}\\&
		\qquad \subseteq
		\mathcal{E}_0 \subseteq \\&\big\{  (u,v)\in [0,1]\times [0,1]\;{\mbox{ s.t. }}\; v \le \gamma_{0}(u) \,\text{ if } \, u\in[0, u_{\mathcal{M}}^0]\\
		&\qquad\qquad\qquad\qquad{\mbox{and }}		\;
		v \leq 1  \, \text{ if }\,  u\in(u_{\mathcal{M}}^0, 1]  \big\},
		\end{split}
		\end{equation}
		where~$\gamma_0$ and~$u_{\mathcal{M}}^0$ are given in~\eqref{ch2def:gamma0}.
		\item[(ii)] It holds that 
		\begin{equation}\label{ch2asdfgert019283}
		\mathcal{S}_c\subseteq
		\mathcal{E}_{\infty} \subseteq\overline{\mathcal{S}_c},
		\end{equation}
		where
		\begin{equation} \label{ch2def:S_c}
		\mathcal{S}_c:=\left\{ (u,v)\in[0,1] \times [0,1]\;
		{\mbox{ s.t. }}\; v-\frac{u}{c}<0 \right\}.
		\end{equation}	
	\end{itemize}
\end{proposition}

We point out that the set~$\mathcal{E}_0$ in~\eqref{ch2char:E0}
does not coincide with
the basin of attraction for the system~\eqref{ch2model} when~$a=0$.
Indeed, as already mentioned, formula~\eqref{ch2form}
in Proposition~\ref{ch2prop:bhvaPRE}
says that for~$a=0$ every point in the interior
of~$[0,1]\times[0,1]$ tends to a coexistence equilibrium and thus
if~$v(0)\neq0$ then~$v(t)$ does not get extinct in finite time.

Also, as~$a\to+\infty$, we have that the set~$\mathcal{E}_{\infty}$
is determined by~$\mathcal{S}_c$, defined in~\eqref{ch2def:S_c},
that depends only on the parameter~$c$. 

\medskip

The statement in~(i) of Proposition~\ref{ch2prop:bhva}
will be a direct consequence of the following result. 
Recalling the function~$\gamma$ introduced in
Propositions~\ref{ch2lemma:M} and~\ref{ch2M:p045}, 
we express here the dependence on the parameter
$a$ by writing~$\gamma_a$,~$u_a$,~$v_a$,
$u_s^a$,~$u_{\mathcal{M}}^a$.
We will also denote by~$\mathcal{M}^a$ the stable manifold
of the point~$(u_s, v_s)$ in~\eqref{ch2usvs}, and by~$\mathcal{M}^0$
the stable manifold
of the point~$(u_s^0, v_s^0)$ in~\eqref{ch2u0v0}.
The key lemma is the following:

\begin{lemma}\label{ch2lemma:conv_gamma}
	For all~$u\in[0,1]$, we have that~$\gamma_a(u) \to \gamma_0(u)$
	uniformly as~$a\to0$, where~$\gamma_0(u)$
	is the function defined in~\eqref{ch2def:gamma0}.
\end{lemma}

\begin{proof}
	Since we are dealing with the limit as~$a$ goes to zero,
	throughout this proof we will always assume that
	we are in the case~$ac<1$.
	
	Also, we denote by~$\phi_p^{a}(t)$ the flow at time~$t$
	of the point~$p\in[0,1]\times[0,1]$ associated with~\eqref{ch2model},
	and similarly by~$\phi_p^{(0)}(t)$ the flow at time~$t$
	of the point~$p$ associated with~\eqref{ch2model} when~$a=0$. 
	With a slight abuse of notation,
	we will also write~$\phi_p^{a}(t)=(u_a(t),v_a(t))$,
	with~$p=(u_a(0),v_a(0))$.
	
	Let us start by proving that
	\begin{equation}\label{ch2340}
	\mathcal{M}^a\cap\big([0,u_s^0]\times[0,v_s^0]\big)\to
	\mathcal{M}^0\cap\big([0,u_s^0]\times[0,v_s^0]\big) \quad {\mbox{ as }}a\to0.
	\end{equation}
	For this, we claim that, for every~$\varepsilon>0$,
	if
	\begin{equation}\label{ch2aqwzero}
	(u_a(0))^2+(v_a(0))^2 \ge\frac{\varepsilon^2}{4}
	\end{equation}
	and
	\begin{equation}\label{ch2zz}
	\big| (u_a(t),v_a(t) ) - (u_s^a, v_s^a)\big| > \frac{\varepsilon}{2},
	\end{equation}
	then
	\begin{equation}  \label{ch2z}
	|\dot{u}_a(t)|^2 +|\dot{v}_a(t)|^2 > \frac{\varepsilon^{4}}{C_0},
	\end{equation}
	for some~$C_0>0$, depending only on~$\rho$ and~$c$.
	
	Indeed, by~(v) of
	Theorem~\ref{ch2thm:dyn} and~\eqref{ch2zz},
	the trajectory~$(u_a(t), v_a(t))$ belongs to the set~$[0, u_s^a] 
	\times [0, v_s^a] \setminus B_{\frac{\varepsilon}{2}}(u_s^a, v_s^a)~$.
	
	Moreover, we claim that
	\begin{equation}\label{ch2123456poi}
	1-ac-u_a(t)-v_a(t)\ge \frac{\varepsilon \sqrt{2}}{4},
	\end{equation}
	for any~$t>0$ such that~\eqref{ch2zz} is satisfied.
	To prove this, we recall that~$(u_s^a, v_s^a)$ lies on the straight
	line~$\ell$ given by~$v=-u+1-ac$
	when~$0<ac<1$ (see~\eqref{ch2curve:u'}). 
	Clearly, there is no point of 
	the set~$[0, u_s^a] 
	\times [0, v_s^a] \setminus B_{\frac{\varepsilon}{2}}(u_s^a, v_s^a)~$
	lying on~$\ell$, and we notice that the points 
	in the set~$[0, u_s^a] 
	\times [0, v_s^a] \setminus B_{\frac{\varepsilon}{2}}(u_s^a, v_s^a)~$
	with minimal distance from~$\ell$ are given by~$p:=(u_s^a-\varepsilon/2,
	v_s^a)$ and~$q:=(u_s^a, v_s^a-\varepsilon/2)$.
	Also, the distance of the point~$p$ from the straight line~$\ell$
	is given by~$\frac{\varepsilon}2\cdot \tan\frac\pi4=
	\frac{\varepsilon \sqrt{2}}{4}$.
	Thus, the distance between~$(u_a(t),v_a(t) )$ 
	and the line~$\ell$ is greater than 
	$\frac{\varepsilon \sqrt{2}}{4}$,
	and this implies~\eqref{ch2123456poi}.
	
	As a consequence of~\eqref{ch2123456poi}, we obtain that
	\begin{equation}\label{ch2pggdeyw087968754}
	(\dot{u}_a(t))^2 = 
	\big(u_a(t)(1-ac-u_a(t)-v_a(t)) \big)^2 > 
	(u_a(t))^2\left(\frac{\varepsilon \sqrt{2}}{4}\right)^2
	\end{equation}
	and that
	\begin{equation}\begin{split}\label{ch2pggdeyw087968754BIS}
	(\dot{v}_a(t))^2\,=& 
	\big(\rho v_a(t)(1-u_a(t)-v_a(t))-au_a(t) \big)^2 \\
	\ge&
	\left(\rho v_a(t)\left(ac+\frac{\varepsilon \sqrt{2}}{4}\right)-au_a(t) \right)^2.
	\end{split}\end{equation}
	Now, if~$u_a(t)\ge\rho cv_a(t)$, then from~\eqref{ch2pggdeyw087968754}
	and~\eqref{ch2aqwzero}
	we obtain that
	\begin{eqnarray*}&&
		(\dot{u}_a(t))^2+(\dot{v}_a(t))^2 \ge 
		(\dot{u}_a(t))^2>
		(u_a(t))^2\left(\frac{\varepsilon \sqrt{2}}{4}\right)^2\\&&\qquad\qquad
		\ge \frac{(u_a(t))^2}2\left(\frac{\varepsilon \sqrt{2}}{4}\right)^2
		+\frac{(\rho cv_a(t))^2}2\left(\frac{\varepsilon \sqrt{2}}{4}\right)^2
		\\&&\qquad\qquad
		\ge \min \{1,\rho^2c^2\} \frac{\varepsilon^2}{16}
		\big((u_a(t))^2+(v_a(t))^2\big)\\
		&&\qquad\qquad
		\ge \min \{1,\rho^2c^2\} \frac{\varepsilon^2}{16}
		\big((u_a(0))^2+(v_a(0))^2\big)\\
		&&\qquad\qquad
		\ge \min \{1,\rho^2c^2\} \frac{\varepsilon^4}{64},
	\end{eqnarray*}
	which proves~\eqref{ch2z} in this case.
	
	If instead~$u_a(t)<\rho cv_a(t)$, we use~\eqref{ch2pggdeyw087968754BIS}
	to see that
	\begin{eqnarray*}&&
		(\dot{u}_a(t))^2+(\dot{v}_a(t))^2 \ge 
		(\dot{v}_a(t))^2\ge
		\left(\rho v_a(t)\left(ac+\frac{\varepsilon \sqrt{2}}{4}\right)-au_a(t) \right)^2
		\\&&\qquad\qquad =
		\left(\frac{\varepsilon \sqrt{2}\rho v_a(t)}{4}+a\big(\rho cv_a(t)-u_a(t)\big) \right)^2
		\ge
		\left(\frac{\varepsilon \sqrt{2}\rho v_a(t)}{4}\right)^2\\&&\qquad\qquad
		\ge
		\frac12\left(\frac{\varepsilon \sqrt{2}\rho v_a(t)}{4}\right)^2
		+\frac12\left(\frac{\varepsilon \sqrt{2} u_a(t)}{4c}\right)^2\\&&\qquad\qquad
		\ge \min \left\{\rho^2,\frac1{c^2}\right\}\frac{\varepsilon^2}{16}
		\big( (u_a(t))^2 +(v_a(t))^2\big)\\
		&&\qquad\qquad
		\ge \min \left\{\rho^2,\frac1{c^2}\right\}\frac{\varepsilon^2}{16}
		\big( (u_a(0))^2 +(v_a(0))^2\big)\\
		&&\qquad\qquad
		\ge \min \left\{\rho^2,\frac1{c^2}\right\}\frac{\varepsilon^4}{64},
	\end{eqnarray*}
	which completes the proof of~\eqref{ch2z}.
	
	Now, for any~$\eta>0$, we define
	$$\mathcal{P}_\eta:=\left\{(u,v)\in[0,1]\times[0,1]\;
	{\mbox{ s.t. }}\; v=\frac{v_s^0-\eta'}{(u_s^0+\eta')^\rho}u^\rho\;
	{\mbox{ with }} |\eta'|\le\eta
	\right\}.$$
	Given~$\varepsilon>0$, we define
	\begin{equation}\label{ch2mettiin}
	{\mbox{$\eta(\varepsilon)$ to be the smallest~$\eta$
			for which~$\mathcal{P}_\eta \supset B_{\varepsilon}(u_s^0,v_s^0)$.}}
	\end{equation}
	We remark that
	\begin{equation}\label{ch2ricorda}
	\lim_{\varepsilon\to0}\eta(\varepsilon)=0.
	\end{equation}
	Also, given~$\delta>0$, we define a tubular
	neighborhood~$\mathcal{U}_\delta$
	of~$\mathcal{M}^0$ as
	$$ \mathcal{U}_\delta :=\bigcup_{q\in\mathcal{M}^0}
	B_{\delta}(q).$$
	Furthermore, we define
	\begin{equation}\label{ch2lometto}
	{\mbox{$\delta(\varepsilon)$ the smallest~$\delta$ such that
			$\mathcal{U}_\delta\supset \mathcal{P}_{\eta(\varepsilon)}$.}}
	\end{equation}
	Recalling~\eqref{ch2ricorda}, we have that
	\begin{equation}\label{ch2ricorda2}
	\lim_{\varepsilon\to0}\delta(\varepsilon)=0.
	\end{equation}
	
	We remark that, as~$a\to0$, the point~$(u_s^a, v_s^a)$ in~\eqref{ch2usvs},
	which is a saddle point for the dynamics of~\eqref{ch2model}
	when~$ac<1$ (recall Theorem~\ref{ch2thm:dyn}),
	tends to the point~$(u_s^0,v_s^0)$ in~\eqref{ch2u0v0}, that belongs
	to the line~$v+u=1$, which is an equilibrium point for the dynamics of~\eqref{ch2model}
	when~$a=0$, according to Proposition~\ref{ch2prop:bhvaPRE}.
	
	As a consequence, for every~$\varepsilon>0$, there exists~$a_\varepsilon>0$
	such that if~$a\in(0,a_\varepsilon)$,
	\begin{equation}\label{ch2qwertyuisdfghjxcvbn}
	|(u_s^a,v_s^a)-(u_s^0,v_s^0)|\le\frac\varepsilon8.
	\end{equation}
	This gives that the intersection of~$\mathcal{M}^a$ with~$B_{\varepsilon/2}(u_s^0,v_s^0)$
	is nonempty.
	
	Furthermore, since~$\gamma_a(0)=0$, in light of Proposition~\ref{ch2lemma:M},
	we have that the intersection of~$\mathcal{M}^a$ with~$B_{\varepsilon/2}$
	is nonempty. Hence, there exists~$p_{\varepsilon,a}\in \mathcal{M}^a\cap
	\partial B_{\varepsilon/2}$.
	
	We also notice that
	\begin{equation}\label{ch2swqdbvsdjvksdv097654}
	\mathcal{M}^a=\phi_{p_{\varepsilon,a}}^{a}(\R).\end{equation}
	In addition, 
	\begin{equation}\label{ch2a12ewgerheh}
	\phi_{p_{\varepsilon,a}}^{a}\big((-\infty,0]\big)\subset B_{\varepsilon/2}.
	\end{equation}
	Also, since the origin belongs to~$\mathcal{M}^0$, we have that~$
	B_{\varepsilon/2}\subset \mathcal{U}_\varepsilon$. {F}rom
	this and~\eqref{ch2a12ewgerheh}, we deduce that
	\begin{equation}\label{ch2lsdgrdhtrjb yrweur748v6348900}
	\phi_{p_{\varepsilon,a}}^{a}\big((-\infty,0]\big)
	\subset \mathcal{U}_\varepsilon.\end{equation}
	
	Now, we let~$C_0$ be as in~\eqref{ch2z} and
	we claim that there exists~$t_{\varepsilon,a}\in(0,3\sqrt{C_0}\varepsilon^{-2})$
	such that
	\begin{equation}\label{ch2esisteuntempo}
	\phi_{p_{\varepsilon,a}}^{a}(t_{\varepsilon,a})\in\partial B_{3\varepsilon/4}
	(u_s^0,v_s^0).
	\end{equation}
	To check this, we argue by contradiction and we suppose that
	$$ \phi_{p_{\varepsilon,a}}^{a}\big((0,3\sqrt{C_0}\varepsilon^{-2})\big)
	\cap B_{3\varepsilon/4}(u_s^0,v_s^0)=\varnothing.$$
	Then, for every~$t\in(0,3\sqrt{C_0}\varepsilon^{-2})$,
	recalling also~\eqref{ch2qwertyuisdfghjxcvbn},
	$$ \big|\phi_{p_{\varepsilon,a}}^{a}(t)-(u_s^a,v_s^a)\big|\ge
	\big|\phi_{p_{\varepsilon,a}}^{a}(t)-(u_s^0,v_s^0)\big| -
	\big|(u_s^a,v_s^a)-(u_s^0,v_s^0)\big|\ge
	\frac{3\varepsilon}4-\frac\varepsilon8>\frac{\varepsilon}2,
	$$
	and consequently~\eqref{ch2zz} is satisfied for every~$t\in(0,3\sqrt{C_0}
	\varepsilon^{-2})$.
	
	Moreover,
	we observe that~$p_{\varepsilon,a}$ satisfies~\eqref{ch2aqwzero},
	and therefore, by~\eqref{ch2z},
	$$ | \dot{u}_a(t)|^2 +| \dot{v}_a(t)|^2
	> \frac{\varepsilon^{4}}{C_0},$$
	for all~$t\in(0,3\sqrt{C_0}\varepsilon^{-2})$,
	where we used the notation~${\phi}_{p_{\varepsilon,a}}^{a}(t)=
	(u_a(t),v_a(t))$, being~$p_{\varepsilon,a}=(u_a(0),v_a(0))$.
	As a result,
	$$ \big( \dot{u}_a(t)+  \dot{v}_a(t)\big)^2>\frac{\varepsilon^{4}}{C_0},$$
	and thus
	$$ \dot{u}_a(t)+  \dot{v}_a(t)>\frac{\varepsilon^{2}}{\sqrt{C_0}}.$$
	This leads to
	\begin{eqnarray*}
		&&u_a\left(\frac{3\sqrt{C_0}}{\varepsilon^2}\right)
		+v_a\left(\frac{3\sqrt{C_0}}{\varepsilon^2}\right)=u_a(0)+v_a(0)
		+\int_0^{\frac{3\sqrt{C_0}}{\varepsilon^2}}\big( \dot{u}_a(t)+  \dot{v}_a(t)\big)\,dt
		\\&&\qquad\quad
		\ge u_a(0)+v_a(0)+
		\int_0^{\frac{3\sqrt{C_0}}{\varepsilon^2}}\frac{\varepsilon^{2}}{\sqrt{C_0}}\,dt
		=u_a(0)+v_a(0) +3\ge 3,
	\end{eqnarray*}
	which forces the trajectory to exit the region~$[0,1]\times[0,1]$.
	This is against the assumption that~$p_{\varepsilon, a}\in\mathcal{M}^a$,
	and therefore the proof of~\eqref{ch2esisteuntempo} is complete.
	
	In light of~\eqref{ch2esisteuntempo}, we can set~$q_{\varepsilon,a}:=
	\phi_{p_{\varepsilon,a}}^{a}(t_{\varepsilon,a})$,
	and we deduce from~\eqref{ch2mettiin} that~$q_{\varepsilon,a}\in\mathcal{P}_{\eta(\varepsilon)}$.
	We also observe that the set~$\mathcal{P}_\eta$ is invariant for the
	flow with~$a=0$, thanks to~\eqref{ch2intprim678}. These observations
	give that~$\phi_{q_{\varepsilon,a}}^{0}(t)\in\mathcal{P}_{\eta(\varepsilon)}$
	for all~$t\in\R$. 
	
	As a result, using~\eqref{ch2lometto}, we conclude that
	\begin{equation}\label{ch2ASDFGHJtergyfhgj}
	\phi_{q_{\varepsilon,a}}^{0}(t)\in\mathcal{U}_{\delta(\varepsilon)}\quad
	{\mbox{ for all }} t\in\R.
	\end{equation}
	In addition, by the continuous dependence of the flow
	on the parameter~$a$ (see e.g. Section~2.4
	in~\cite{MR3791466},
	or Theorem~2.4.2 in~\cite{MR3186036}),
	$$ \big|\phi_{q_{\varepsilon,a}}^{0}(t)-\phi_{q_{\varepsilon,a}}^{a}(t)\big|
	<\varepsilon,$$
	for all~$t\in[-3\sqrt{C_0}\varepsilon^{-2},0]$, provided that~$a$
	is sufficiently small, possibly in dependence of~$\varepsilon$.
	This fact and~\eqref{ch2ASDFGHJtergyfhgj} entail that
	$$ \phi_{q_{\varepsilon,a}}^{a}(t)\in\mathcal{U}_{\delta(\varepsilon)+\varepsilon}\quad
	{\mbox{ for all }} t\in[-3\sqrt{C_0}\varepsilon^{-2},0].
	$$
	In particular, for all~$t\in[0,t_{\varepsilon,a}]$,
	\begin{equation}\label{ch2andatosu}
	\phi_{p_{\varepsilon,a}}^{a}(t)=\phi_{q_{\varepsilon,a}}^{a}(t-t_{\varepsilon,a})
	\in\mathcal{U}_{\delta(\varepsilon)+\varepsilon}.\end{equation}
	
	We now claim that
	for all~$t\ge t_{\varepsilon,a}$,
	\begin{equation}\label{ch2qwertyuiop}
	\phi_{p_{\varepsilon,a}}^{a}(t)\subset B_{\varepsilon}(u_s^a,v_s^a).
	\end{equation}
	Indeed, this is true when~$t=t_{\varepsilon,a}$ thanks to~\eqref{ch2qwertyuisdfghjxcvbn}
	and~\eqref{ch2esisteuntempo}. Hence, since
	the trajectory~$\phi_{p_{\varepsilon,a}}^{a}(t)$
	is contained in the domain where~$\dot{u}\ge0$ and~$\dot{v}\ge0$,
	thanks to~\eqref{ch2aggiunto}, we deduce that~\eqref{ch2qwertyuiop} holds true.
	
	{F}rom~\eqref{ch2qwertyuisdfghjxcvbn} and~\eqref{ch2qwertyuiop}, we conclude that
	$$ \phi_{p_{\varepsilon,a}}^{a}(t)\subset B_{2\varepsilon}(u_s^0,v_s^0),$$
	for all~$t\ge t_{\varepsilon,a}$.
	
	Using this,~\eqref{ch2lsdgrdhtrjb yrweur748v6348900}
	and~\eqref{ch2andatosu}, we obtain that
	$$ \phi_{p_{\varepsilon,a}}^{a}(\R)\subset\mathcal{U}_{\delta(\varepsilon)+
		2\varepsilon}.$$
	This and~\eqref{ch2ricorda2} give that~\eqref{ch2340} is satisfied, as desired.
	
	One can also show that
	\begin{equation}\label{ch2340BIS}
	\mathcal{M}^a\cap\big([u_s^0, u_{\mathcal{M}}^0]\times[v_s^0,v_{\mathcal{M}}^0]\big)\to
	\mathcal{M}^0\cap\big([u_s^0, u_{\mathcal{M}}^0]\times[v_s^0,v_{\mathcal{M}}^0]\big) \quad {\mbox{ as }}a\to0.
	\end{equation}
	The proof of~\eqref{ch2340BIS} is similar to that of~\eqref{ch2340},
	just replacing~$p_{\varepsilon,a}$ with~$(u_{\mathcal{M}}^a,v_{\mathcal{M}}^a)$ (in this case
	the analysis near the origin is simply omitted since the trajectory
	has only one limit point).
	
	With~\eqref{ch2340} and~\eqref{ch2340BIS}
	the proof of Lemma~\ref{ch2lemma:conv_gamma} is thereby complete.
\end{proof}

Now we are ready to give the proof of Proposition~\ref{ch2prop:bhva}:

\begin{proof}[Proof of Proposition~\ref{ch2prop:bhva}]
	{\emph{(i)}} We call~$\mathcal{G}$ the right-hand-side of~\eqref{ch2char:E0}, that is
	\begin{eqnarray*}
		\mathcal{G}&:=& \big\{  (u,v)\in [0,1]\times [0,1]\;{\mbox{ s.t. }}\; 
		v < \gamma_{0}(u) \,\text{ if } \, u\in[0, u_{\mathcal{M}}^0]\\
		&&\qquad\qquad\qquad\qquad{\mbox{and }}		\;
		v \leq 1  \, \text{ if }\,  u\in(u_{\mathcal{M}}^0, 1]  \big\},
	\end{eqnarray*}
	and we aim at proving that~$\mathcal{G}\subseteq\mathcal{E}_0
	\subseteq\overline{\mathcal{G}}$.
	
	For this, we observe that, by Lemma~\ref{ch2lemma:conv_gamma},
	$\gamma_a(u)$ converges to~$\gamma_{0}(u)$ pointwise as~$a\to0$.
	In particular,~$u_{\mathcal{M}}^a\to u_{\mathcal{M}}^0$ as~$a\to0$.
	
	Also, recalling~\eqref{ch2def:gamma0},
	we notice that if~$u_{\mathcal{M}}^0= u_s^0 / (v_s^0)^{\frac{1}{\rho}}<1$,
	then~$\gamma_0(u_{\mathcal{M}}^0)= 1$,
	otherwise if~$u_{\mathcal{M}}^0=1$
	then~$\gamma_0(u_{\mathcal{M}}^0)<1$, being~$\gamma_0(u)$ strictly
	monotone increasing.
	
	Furthermore, thanks to Proposition~\ref{ch2prop:char},
	we know that
	he set~$\mathcal{E}(a)$ is bounded from above by the graph of the function~$\gamma_a(u)$ for~$u\in [0, u_{\mathcal{M}}^a]$ and from the
	straight line~$v=1$ for~$u\in(u_{\mathcal{M}}^a, 1]$ (that is non empty for~$u_{\mathcal{M}}^a<1$). 
	
	Now we claim that, for all~$a'>0$,
	\begin{equation}\label{ch2gocont123}
	\mathcal{G} \subseteq \underset{0<a<a'}{\bigcup} \mathcal{E}(a).
	\end{equation}
	To show this, we take a point~$(u,v)\in\mathcal{G}$.
	Hence, in light of the considerations
	above, we have that~$(u,v)\in\mathcal{E}(a)$ for any~$a$ sufficiently small,
	which proves~\eqref{ch2gocont123}.
	
	{F}rom~\eqref{ch2gocont123}, we deduce that
	\begin{equation}\label{ch2gocont1232233}
	\mathcal{G} \subseteq \underset{a'>0}{\bigcap} \, \underset{0<a<a'}{\bigcup} \mathcal{E}(a).
	\end{equation}
	Now we show that
	\begin{equation}\label{ch2gocont12322}
	\underset{a'>0}{\bigcap} \,\underset{0<a<a'}{\bigcup} \mathcal{E}(a)\subseteq
	\overline{\mathcal{G}} .
	\end{equation}
	For this, we take 
	$$(\hat{u},\hat{v})\in  \underset{a'>0}{\bigcap} \, \underset{0<a<a'}{\bigcup} \mathcal{E}(a),$$
	then it must hold that for every~$a'>0$
	there exists~$a<a'$ such that~$(\hat{u},\hat{v})\in\mathcal{E}(a)$,
	namely~$\hat v < \gamma_{a}(\hat u)$ if~$\hat u\in[0, u_{\mathcal{M}}^a]$ and~$\hat v
	\leq 1$ if~$\hat u\in(u_{\mathcal{M}}^a, 1]$.  
	Thus, by the pointwise convergence,
	we have that~$  \hat{v} \le\gamma_0(\hat{u})~$ if~$\hat u\in[0, u_{\mathcal{M}}^0]$ and~$\hat v
	\leq 1$ if~$\hat u\in(u_{\mathcal{M}}^0, 1]$, which proves~\eqref{ch2gocont12322}.
	
	{F}rom~\eqref{ch2gocont1232233} and~\eqref{ch2gocont12322},
	we conclude that
	\begin{equation*}
	\mathcal{G}\subseteq
	\underset{a'>0}{\bigcap} \, \underset{0<a<a'}{\bigcup} \mathcal{E}(a) =\mathcal{E}_0\subseteq \overline{\mathcal{G}} ,
	\end{equation*}
	as desired.
	\medskip	
	
	{\emph{(ii)}} Since we deal with the limit case as~$a\to+\infty$, from now on
	we suppose from now on that~$ac>1$.
	We fix~$\varepsilon>0$ and we consider the set  
	\begin{equation*}
	\mathcal{S}_{\varepsilon^+} := 
	\left\{  (u,v)\in [0,1]\times[0,1]\;{\mbox{ s.t. }}\;
	v>u \left( \frac{1}{c}+\varepsilon \right)      \right\}.
	\end{equation*}
	We claim that
	\begin{equation}\label{ch2prova1}
	\mathcal{S}_{\varepsilon^+} \subseteq \mathcal{B}(a)
	\end{equation}
	for~$a$ big enough, possibly in dependence of~$\varepsilon$.
	For this, 
	we first analyze the component of the velocity in the inward normal directions
	along the boundary of~$\mathcal{S}_{\varepsilon^+}$.
	On the side~$\{0\}\times [0,1]$, 
	the trajectories cannot cross the boundary thanks to Proposition~\ref{ch2prop:dyn}, and the same happens for the sides~$[0,1]\times \{1\}$ and~$\{1\} \times [\varepsilon + 1/c, 1]$. 
	
	Hence, it remains to check the sign of the normal derivative along the
	side given by the straight line~$v-u(\varepsilon +1/c )=0$.
	We compute
	\begin{align*}& (\dot{u},\dot{v})\cdot\left(-\left(\varepsilon+\frac1c\right),
	1\right)=
	\dot{v}- \dot{u}\left(\varepsilon+ \frac{1}{c} \right) \\
	&\quad = {\rho}v(1-u-v)-au - \left(\varepsilon+ \frac{1}{c} \right)u(1-u-v) + \left(\varepsilon+ \frac{1}{c} \right) acu \\
	&\quad=  \Bigg[ {\rho}v - \left(\varepsilon+ \frac{1}{c} \right) u     \Bigg] (1-u-v) +\varepsilon ac u	.	
	\end{align*} 	
	Thus, by using that~$v-u(\varepsilon +1/c )=0$, we obtain that
	\begin{align*}
	(\dot{u},\dot{v})\cdot\left(-\left(\varepsilon+\frac1c\right),
	1\right) \geq u\left[a\varepsilon c + ({\rho}-1)(1-u-v)   \left( \varepsilon+ \frac{1}{c} \right)  \right]. 
	\end{align*}
	Notice that~$u\leq 1$ and~$|1-u-v|\leq 2$, and therefore
	$$ (\dot{u},\dot{v})\cdot\left(-\left(\varepsilon+\frac1c\right),
	1\right) \geq u\left[a\varepsilon c -2 ({\rho}+1) \left( \varepsilon+ \frac{1}{c} \right)  \right] .$$
	Accordingly, the normal velocity is positive for~$a \geq {a}_1$, where
	\begin{equation*}
	{a}_1:= 2({\rho}+1) \left( \varepsilon+ \frac{1}{c}  \right)\frac{1}{\varepsilon c}.
	\end{equation*}
	These considerations, together with the fact that there are no cycles
	in~$[0,1]\times[0,1]$ and the
	Poincar\'e-Bendixson Theorem (see e.g.~\cite{TESCHL})
	give that the~$\omega$-limit set of any trajectory starting
	in the interior of~$\mathcal{S}_{\varepsilon^+}$
	can be either an equilibrium or a union of (finitely many)
	equilibria and non-closed orbits connecting these equilibria.
	
	We remark that 
	\begin{equation}\label{ch2asdfgzxcv098t76re}
	{\mbox{the~$\omega$-limit set of any trajectory cannot
			be the equilibrium~$(0,0)$.}}\end{equation}
	Indeed, if the~$\omega$-limit of a trajectory
	were~$(0,0)$, then this trajectory must lie on the stable manifold of~$(0,0)$,
	and moreover it must be contained in~$\mathcal{S}_{\varepsilon^+}$,
	since no trajectory can exit~$\mathcal{S}_{\varepsilon^+}$.
	On the other hand, by Proposition~\ref{ch2lemma:M},
	we have that at~$u=0$ the stable manifold is tangent to the
	line
	$$v=\frac{a}{\rho-1+ac}u=\frac{1}{\frac{\rho-1}{a}+c}u.
	$$
	Now, if we take~$a$ sufficiently large, this line lies
	below the line~$v=u(1/c+\varepsilon)$, thus providing a contradiction.
	Hence, the proof of~\eqref{ch2asdfgzxcv098t76re} is complete.
	
	Accordingly, since~$(0,1)$ is a sink, the only possibility is that
	the~$\omega$-limit set of any trajectory starting
	in the interior of~$\mathcal{S}_{\varepsilon^+}$ is the equilibrium~$(0,1)$.
	Namely, we have established~\eqref{ch2prova1}.
	
	As a consequence of~\eqref{ch2prova1}, we deduce that for every~$\varepsilon>0$
	there exists~$a_{\varepsilon}>0$ such that
	\begin{equation}\label{ch2qwt5uktkjer464586897}
	\underset{a\ge a_\varepsilon}{\bigcup}  \mathcal{E}(a) \subseteq
	\left\{  (u,v)\in [0,1]\times[0,1]\;{\mbox{ s.t. }}\;
	v\le u \left( \frac{1}{c}+\varepsilon \right)      \right\}.\end{equation}
	In addition,
	\begin{equation*}\begin{split}
	&	\underset{\varepsilon >0 }{\bigcap} \left\{  (u,v)\in [0,1]\times[0,1]\;{\mbox{ s.t. }}\;
	v\le u \left( \frac{1}{c}+\varepsilon \right)      \right\}\\&\qquad
	=\left\{  (u,v)\in [0,1]\times[0,1]\;{\mbox{ s.t. }}\;
	v\le \frac{u}{c} \right\}=\overline{\mathcal{S}_c}.\end{split}
	\end{equation*}
	This and~\eqref{ch2qwt5uktkjer464586897} entail that 
	\begin{equation*}
	\underset{a'>0}{\bigcap} \, \underset{a>a'}{\bigcup}
	\mathcal{E}(a)\subseteq \overline{\mathcal{S}_c},
	\end{equation*}
	which implies the second inclusion in~\eqref{ch2asdfgert019283}.
	
	Now, to show the first inclusion in~\eqref{ch2asdfgert019283},
	for every~$\varepsilon\in(0,1/c)$ we consider the set
	\begin{equation*}
	\mathcal{S}_{\varepsilon^-} := \left\{  
	(u,v)\in [0,1]\times[0,1] \;{\mbox{ s.t. }}\; v<u \left( \frac{1}{c}-\varepsilon \right)      \right\}.
	\end{equation*}
	We claim that, for all~$\varepsilon\in(0,1/c)$,
	\begin{equation}\label{ch2chefus}
	\mathcal{S}_{\varepsilon^-} \subseteq \mathcal{E}_{\infty}.
	\end{equation}
	For this, we first show that if~$a$ is sufficiently large, possibly
	in dependence of~$\varepsilon$,
	\begin{equation}\label{ch2forse33}\begin{split}&
	{\mbox{every trajectory starting in the interior of~$\mathcal{S}_{\varepsilon^-}$}}\\
	&{\mbox{can exit~$\mathcal{S}_{\varepsilon^-}$ from the side~$[0,1]\times\{0\}$.}}
	\end{split}\end{equation}
	Indeed, 
	on the side~$\{1\}\times [0,1]$ the trajectory cannot exit the set,
	thanks to Proposition~\ref{ch2prop:dyn}.  
	On the side given by~$v-(-\varepsilon+1/c)u=0$, the component of the velocity
	in the direction of the
	outward normal is 
	\begin{eqnarray*}&&
		(\dot{u},\dot{v})\cdot\left(- \left( \frac{1}{c} -\varepsilon \right),1
		\right)=
		\dot{v} - \dot{u} \left( \frac{1}{c} -\varepsilon \right)\\
		&& \qquad= \rho v(1-u-v) -au - \left( \frac{1}{c} -\varepsilon \right)u(1-u-v)   + \left( \frac{1}{c} -\varepsilon \right)acu\\
		&&\qquad
		=u\left[\left( \frac{1}{c} -\varepsilon \right)(\rho-1)(1-u-v)    - 
		\varepsilon ac \right]\\
		&&\qquad \le u\left[2\left( \frac{1}{c} -\varepsilon \right)(\rho+1)  - 
		\varepsilon ac \right]
		,
	\end{eqnarray*}
	which is negative if~$a \geq {a}_2$, with
	\begin{equation*}
	{a}_2:= 2\left( \frac{1}{c} -\varepsilon \right)
	\left( \rho+1  \right) \frac{1}{\varepsilon c} .
	\end{equation*}
	Hence, if~$(u(0), v(0))\in \mathcal{S}_{\varepsilon^-}$, then
	either~$T_s(u(0), v(0)) <\infty$ or~$(u(t), v(t))\in \mathcal{S}_{\varepsilon^-}$
	for all~$t\geq 0$, where the notation in~\eqref{ch2def:T_s} has been used.
	We also notice that,
	for~$a>1/c$, the points~$(0,1)$ and~$(0,0)$ are the only equilibria
	of the system, and there are no cycles. 
	We have that~$(0,1) \notin \overline{\mathcal{S}_{\varepsilon^-}}$
	and~$(0,0) \in \overline{\mathcal{S}_{\varepsilon^-}}$, thus if
	\begin{equation}\label{ch2dweioterygvhsdjk}
	{\mbox{$(u(t), v(t))\in
			\mathcal{S}_{\varepsilon^-}$ for all~$t\geq 0$}}\end{equation} then
	\begin{equation}\label{ch2tendere}
	(u(t), v(t)) \to (0,0).
	\end{equation}
	On the other hand, by Proposition~\ref{ch2lemma:M},
	we have that at~$u=0$ the stable manifold is tangent to the
	line
	$$v=\frac{a}{\rho-1+ac}u=\frac{1}{\frac{\rho-1}{a}+c}u,
	$$
	and, if we take~$a$ large enough, this line lies
	above the line~$v=u(1/c-\varepsilon)$. This says that, for sufficiently large~$t$,
	the trajectory must lie outside~$ \mathcal{S}_{\varepsilon^-}$,
	and this is in contradiction with~\eqref{ch2dweioterygvhsdjk}.
	
	As a result of these considerations, we conclude 
	that if~$(u(0), v(0))\in \mathcal{S}_{\varepsilon^-}$ then~$T_s(u(0), v(0)) <\infty$,
	which implies~\eqref{ch2forse33}.
	
	As a consequence of~\eqref{ch2forse33}, we obtain that for every~$\varepsilon\in(0,1/c)$
	there exists~$a_\varepsilon>0$ such that
	$$\mathcal{S}_{\varepsilon^-}\subseteq\underset{a\ge a_\varepsilon}{\bigcap}
	\mathcal{E}(a).$$
	In particular for all~$\varepsilon\in(0,1/c)$ it holds that
	\begin{equation*}
	\mathcal{S}_{\varepsilon^-}\subseteq \underset{a'>0}{\bigcap} \,
	\underset{a> a'}{\bigcup} \mathcal{E}(a)=\mathcal{E}_\infty,
	\end{equation*}
	which proves~\eqref{ch2chefus}, as desired.
	
	Then, the first inclusion in~\eqref{ch2asdfgert019283} plainly follows
	from~\eqref{ch2chefus}.
\end{proof}

\section{Analysis of the strategies for the first population}\label{ch2STRATE}

The main theorems on the winning strategy have
been stated in Subsection~\ref{ch2ss:strategy}.
In particular, Theorem~\ref{ch2thm:Vbound} gives the characterization of the set
of points that have a winning strategy~$\mathcal{V}_{\mathcal{A}}$
in~\eqref{ch2DEFNU},
and Theorem~\ref{ch2thm:W} establishes the non equivalence of constant
and non-constant strategies when~$\rho\ne1$
(and their equivalence when~$\rho=1$).
Nonetheless, in Theorem~\ref{ch2thm:H} we state that
Heaviside functions are enough to construct a
winning strategy for every point in~$\mathcal{V}_{\mathcal{A}}$.

In the following subsections we will give the proofs of these results.

\subsection{Construction of winning non-constant strategies}\label{ch2explowrewwt}

We want to put in light the construction of non-constant winning
strategies for the points for which constant strategies fail.

For this, we recall the notation introduced in~\eqref{ch2u0v0},
\eqref{ch2ZETADEF} and~\eqref{ch2def:gamma0}, and we have the following statement:

\begin{proposition}\label{ch2prop:construction}
	Let~$M>1$. Then we have:
	\begin{itemize}
		\item[1.] For~$\rho<1$, let~$(u_0, v_0)$ be a point of the set
		\begin{equation}\label{ch2PPDEFA}
		\mathcal{P}:=\left\{ (u, v)\in [0,1]\times[0,1] \;{\mbox{ s.t. }}\;  u\in [u_s^0, 1], \ \gamma_{0}(u) \leq v < \frac{u}{c} + \frac{1-\rho}{1+\rho c}    \right\}.	\end{equation}
		Then there exist~$a^*>M$,~$a_*<\frac{1}{M}$,
		and~$T\ge0$, depending on~$(u_0, v_0)$,~$c$, and~$\rho$, such that
		the Heaviside strategy defined by
		\begin{equation}\label{ch2NSJmldsf965to}
		a(t) = \left\{
		\begin{array}{lr}
		a^*,  & {\mbox{ if }} t<T, \\
		a_*,  &  {\mbox{ if }} t\geq T,
		\end{array}
		\right.
		\end{equation}
		belongs to~$\mathcal{V}_{\mathcal{A}}$.
		\item[2.] For~$\rho>1$, let~$(u_0, v_0)$ be a point of the set
		\begin{equation}\label{ch2DEFQ}
		\mathcal{Q}:=\left\{ (u, v)\in [0,1]\times[0,1] \;{\mbox{ s.t. }}\;u\in [u_{\infty}, 1], \ \frac{u}{c}  \leq v < \zeta(u)   \right\}.
		\end{equation}
		Then there exist~$a^*>M$,~$a_*<\frac{1}{M}$, and~$T\ge0$, depending on~$(u_0, v_0)$,~$c$, and~$\rho$, such that the Heaviside strategy
		defined by
		\begin{equation*}
		a(t) = \left\{
		\begin{array}{lr}
		a_*,  &{\mbox{ if }} t<T, \\
		a^*,  &{\mbox{ if }} t\geq T,
		\end{array}
		\right.
		\end{equation*}
		belongs to~$\mathcal{V}_{\mathcal{A}}$.
	\end{itemize}
\end{proposition}

\begin{proof} We start by proving the first claim in Proposition~\ref{ch2prop:construction}.
	To this aim, we take~$(\bar{u}, \bar{v})\in \mathcal{P}$, and we observe that
	\begin{equation*}
	\bar{v}-\frac{\bar{u}}{c}< \frac{1-\rho}{1+\rho c} = v_s^0-\frac{u_s^0}{c}.
	\end{equation*}
	Therefore, there exists~$\xi>0$ such that 
	\begin{equation*}
	\xi < \frac{v_s^0- \bar{v}-\frac{1}{c}(u_s^0-\bar{u})}{\bar{u}-u_s^0}.
	\end{equation*}
	Hence, setting
	\begin{equation}\label{ch2VUESSE0}
	v_S:=\left(\frac{1}{c} -\xi \right)(u_s^0-\bar{u}) + \bar{v},
	\end{equation}	
	we see that
	\begin{equation}\label{ch2VUESSE}v_S<v_s^0.\end{equation}	
	Now, we want to show that there exists~$a^*>0$ such that, for any~$a>a^*$ and~$u>u_s^0$, we have that
	\begin{equation}\label{ch21632}
	\frac{\dot{v}}{\dot{u}} > \frac{1}{c}- \xi.
	\end{equation}
	To prove this, we first notice that
	\begin{equation}\label{ch2questab}
	{\mbox{if~$a>\displaystyle\frac2c$, then~$\dot{u}\le -u<0$.}}\end{equation}
	Moreover, we set~$$a_1:=\frac{1+\rho c}{4c} , ~$$
	and we claim that, 
	\begin{equation}\label{ch2questae}
	{\mbox{if~$a>a_1$ and~$u>u_s^0$, then~$\dot{v}<0$.}}\end{equation}
	Indeed, we recall that the function~$\sigma$
	defined in~\eqref{ch2f:sigma} represents the points
	in~$[0,1]\times[0,1]$ where~$\dot v=0$
	and separates the points where~$ \dot v>0$, which lie on the left of
	the curve described by~$\sigma$, from the points where~$ \dot v<0$, which lie on the right of
	the curve described by~$\sigma$.	
	
	Therefore, in order to show~\eqref{ch2questae}, it is sufficient to prove that
	the curve described by~$\sigma$
	is contained in~$\{u\le u_s^0\}$ whenever~$a>a_1$. For this, one computes that, if~$u=\sigma(v)$
	and~$a>a_1$, then
	\begin{eqnarray*}&&
		u-u_s^0=\sigma(v)-\frac{\rho c}{1+\rho c}=
		1-\frac{\rho v^2+a}{\rho v+a}-\frac{\rho c}{1+\rho c}\\&&\qquad=
		\frac{\rho v-\rho v^2}{\rho v+a}-\frac{\rho c}{1+\rho c}=\frac{\rho v(1-v)}{\rho v+a}-\frac{\rho c}{1+\rho c}
		\\&&\qquad\le\frac{\rho }{4(\rho v+a)}-\frac{\rho c}{1+\rho c}\le
		\frac{\rho }{4a}-\frac{\rho c}{1+\rho c}\\&&\qquad\le
		\frac{\rho }{4a_1}-\frac{\rho c}{1+\rho c}\le0.
	\end{eqnarray*}
	This completes the proof of~\eqref{ch2questae}.
	
	Now we define
	\begin{equation*}
	a_2:=\left( \rho+\frac{1}{c}+\xi  \right) \frac{2}{u_s^0 c \xi}.
	\end{equation*}
	and we claim that
	\begin{equation}\label{ch2questaX}
	{\mbox{if~$a>a_2$ and~$u>u_s^0$, then }}
	\dot{v} < \left( \frac{1}{c}- \xi \right) \dot{u}.
	\end{equation}
	Indeed, under the assumptions of~\eqref{ch2questaX},
	we deduce that
	\begin{eqnarray*}&&
		\dot{v} -\left( \frac{1}{c}- \xi \right) \dot{u}
		=\rho v(1-u-v)-au-\left( \frac{1}{c}- \xi \right)\Big(
		u(1-u-v)-acu
		\Big)\\
		&&\qquad=(1-u-v)\left(
		\rho v-\left( \frac{1}{c}- \xi \right)u\right)-ac \xi u
		\le 2\left(
		\rho v+ \frac{u}{c}+ \xi u\right)-ac \xi u\\&&\qquad< 2\left(
		\rho + \frac{1}{c}+ \xi \right)-a_2\,c \xi {u_s^0}=0,
	\end{eqnarray*}
	and this establishes the claim in~\eqref{ch2questaX}.
	
	Then, choosing
	$$ a^*:=
	\max\left\{\displaystyle\frac2c,a_1,a_2,M\right\},$$
	we can exploit~\eqref{ch2questab},~\eqref{ch2questae} and~\eqref{ch2questaX}
	to deduce~\eqref{ch21632}, as desired.
	
	Now we claim that, for any~$a>a^*$, there exists~$T\ge0$ 
	such that the trajectory~$(u(t), v(t))$ starting from~$(\bar{u}, \bar{v})$ satisfies
	\begin{equation}\label{ch2SM -kg}
	{\mbox{$u(T)=u_s^0$ and~$v(T)< v_S$.}}
	\end{equation}
	Indeed, we define~$T\ge0$ to be the first time for which~$u(T)=u_s^0$.
	This is a fair definition, since~$u(0)=\bar{u}\ge u_s^0$
	and~$\dot u$ is negative, and bounded away from zero
	till~$u\ge u_s^0$, thanks to~\eqref{ch2questab}.
	Then, we see that
	\begin{eqnarray*}&&
		v(T)=\bar{v}+\int_0^T \dot v(t)\,dt<
		\bar{v}+\int_0^T\left(
		\frac{1}{c}- \xi\right)\,\dot u(t)\,dt=
		\bar{v}+\left(
		\frac{1}{c}- \xi\right)(u(T)-u(0))\\&&\qquad\qquad=
		\bar{v}+\left(
		\frac{1}{c}- \xi\right)(u_s^0-\bar u)=v_S,
	\end{eqnarray*}
	thanks to~\eqref{ch2VUESSE0} and~\eqref{ch21632}, and this establishes~\eqref{ch2SM -kg}.
	
	Now we observe that
	$$ v(T)<v_S<v_s^0=\gamma_0(u_s^0)=\gamma_0(u(T))$$
	due to~\eqref{ch2VUESSE}
	and~\eqref{ch2SM -kg}
	
	As a result, recalling Lemma~\ref{ch2lemma:conv_gamma}, we can choose~$a_*<1/M$
	such that
	$$ v(T)<\gamma_{a_*}(u(T)).$$
	Accordingly, by Proposition~\ref{ch2prop:char},
	we obtain that~$(u(T),v(T))\in{\mathcal{E}}(a_*)$.
	Hence, applying the strategy
	in~\eqref{ch2NSJmldsf965to},
	we accomplish the desired result and complete the proof of the first claim
	in Proposition~\ref{ch2prop:construction}.\medskip
	
	Now we focus on the proof of the second claim in Proposition~\ref{ch2prop:construction}.
	For this, let \begin{equation}\label{ch28ygdw}
	(u_0,v_0)\in \mathcal{Q},\end{equation} and consider the trajectory~$(u_0(t),v_0(t))$
	starting from~$(u_0,v_0)$ for the strategy~$a=0$. In light of formula~\eqref{ch2form}
	of Proposition~\ref{ch2prop:bhvaPRE}, we have that
	\begin{equation}\label{ch2Ecijerrin}\begin{split}&
	{\mbox{the trajectory~$(u_0(t),v_0(t))$
			converges}}\\&{\mbox{to a point of the form~$(u_F, 1-u_F)$ as~$t\to+\infty$.}}\end{split}
	\end{equation}
	
	We define
	\begin{equation}\label{ch21713}
	v_F:=1-u_F, \quad	v_{\infty}:=1-u_{\infty}=\frac{1}{c+1},
	\end{equation} 
	where the last equality can be checked starting from the value of~$u_{\infty}$ given in~\eqref{ch2ZETADEF}.
	Using the definition of~$\zeta$ in~\eqref{ch2ZETADEF} and the information in~\eqref{ch2intprim678}, 
	we also notice that the curve given by~$v=\zeta(u)$ is a trajectory for~$a=0$. Moreover 
	\begin{equation*}
	\zeta(u_{\infty})= \frac{1}{c(u_{\infty})^{\rho-1}} u_{\infty}^{\rho}=\frac{c}{c(c+1)}=v_{\infty}
	\end{equation*}
	and, recalling~\eqref{ch21713} and formula~\eqref{ch2form} of Proposition~\ref{ch2prop:bhvaPRE}, we get that 
	the graph of~$\zeta$ is a trajectory for~$a=0$ that converges to~$(u_\infty,1-u_\infty)$ as~$t\to+\infty$.
	
	Also, by~\eqref{ch28ygdw}, we have that~$v_0 < \zeta(u_0)$. Thus, since by Cauchy's uniqueness result for ODEs, two orbits never intersect, we have that
	\begin{equation}\label{ch2JS145DD-0}
	{\mbox{the orbit~$(u_0(t),v_0(t))$ must lie below the graph of~$\zeta$.}}\end{equation}
	Since both~$(u_F, v_F)$ and~$(u_{\infty}, v_{\infty})$ belong to the line given by~$v=1-u$, from~\eqref{ch2JS145DD-0} we get that
	\begin{equation}\label{ch22304}
	u_{\infty} < u_F
	\end{equation}
	and 
	\begin{equation}\label{ch22305}
	v_{\infty} > v_F.
	\end{equation}
	Thanks to~\eqref{ch22304} and~\eqref{ch22305} and recalling the values of~$u_{\infty}$ from~\eqref{ch2ZETADEF} and of~$v_{\infty}$ from~\eqref{ch21713}, we get that
	\begin{equation}\label{ch21524}
	v_F < v_{\infty} = \frac{u_{\infty}}{c} < \frac{u_F}{c}.
	\end{equation}
	As a consequence, since the inequality in~\eqref{ch21524} is strict,
	we find that there exists~$T'>0$
	such that
	\begin{equation}\label{ch22334}
	v_0(T') < \frac{u_0(T')}{c}.
	\end{equation}
	Moreover, since~$\dot{u}<0$ for~$v>1-u$ and~$a=0$, we get that 
	$u_0(t)$ is decreasing in~$t$, and therefore~$u_F < u_0(T') < u_0$.
	
	
	By the strict inequality in~\eqref{ch22334}, and 
	claim~(ii) in Proposition~\ref{ch2prop:bhva},
	we have that~$(u_0(T'), v_0(T')) \in \mathcal{E}_{\infty}$, where~$\mathcal{E}_{\infty}$ is defined in~\eqref{ch2def:Einfty}.
	In particular, we have that~$(u_0(T'), v_0(T')) \in
	\underset{a>a'}{\bigcup} \mathcal{E}(a)$, for every~$a'>0$.
	Consequently, there exists~$a^*>M$ such that~$(u_0(T'),v_0(T'))\in\mathcal{E}(a^*)$. Therefore, applying the strategy
	\begin{equation*}
	a(t) = \left\{
	\begin{array}{lr}
	0,  & t<T', \\
	a^*,  & t\geq T,
	\end{array}
	\right.
	\end{equation*}
	we reach the victory.
\end{proof}

\subsection{Proof of Theorem~\ref{ch2thm:Vbound}}

To avoid repeating passages in the proofs of Theorems \ref{ch2thm:Vbound} and \ref{ch2thm:W}, we first state and prove the following lemma:

\begin{lemma}\label{ch2lemma:rho=1}
	If~$\rho=1$, then for all~$a>0$ we have~$\mathcal{E}(a)=\mathcal{S}_c$, where~$\mathcal{S}_c$ was defined in~\eqref{ch2def:S_c}.
\end{lemma}

\begin{proof}
	Let~$(u(t),v(t))$ be a trajectory starting at a point
	in~$[0,1]\times[0,1]$.
	For any~$a>0$, we consider the function
	$$\mu(t):= \frac{\displaystyle
		v\left(\frac{t}{a}\right)}{\displaystyle u\left(\frac{t}{a}\right)}.$$
	Notice that
	\begin{equation}\label{ch2-1-e}
	(u(0),v(0))\in\mathcal{E}(a)
	{\mbox{ if and only if there exists~$T>0$ such that }}\mu(T)=0.
	\end{equation}
	In addition, we observe that
	\begin{equation}	\label{ch2eq:A=1}\begin{split}
	\dot{\mu}(t) \,&= 
	\frac{\displaystyle\dot v\left(\frac{t}{a}\right)u\left(\frac{t}{a}\right)-
		v\left(\frac{t}{a}\right)\dot u\left(\frac{t}{a}\right)}{\displaystyle au^2\left(\frac{t}{a}\right)}
	\\&=
	\frac{\displaystyle -u^2 \left(\frac{t}{a}\right)+c u\left(\frac{t}{a}\right)v\left(\frac{t}{a}\right)}{\displaystyle u^2\left(\frac{t}{a}\right)}\\&=c\mu(t) -1.
	\end{split}
	\end{equation}
	The equation in~\eqref{ch2eq:A=1} is integrable and leads to
	\begin{equation*}
	\mu(t)=\frac{e^{ct} \left( c\mu(0)-1 \right) +1}{c}. 
	\end{equation*}
	{F}rom this and~\eqref{ch2-1-e}, we deduce that
	\begin{equation*}
	(u(0),v(0))\in\mathcal{E}(a)
	{\mbox{ if and only if }}
	c\mu(0)-1<0.
	\end{equation*}
	This leads to
	\begin{equation*}
	(u(0),v(0))\in\mathcal{E}(a)
	{\mbox{ if and only if }}\,
	\frac{v(0)}{u(0)}<\frac1c,
	\end{equation*}
	which, recalling the definition of~$\mathcal{S}_c$
	in~\eqref{ch2def:S_c}, ends the proof.
\end{proof}

Now we provide the proof of Theorem~\ref{ch2thm:Vbound}, exploiting the result obtained
in Section~\ref{ch2explowrewwt}.

\begin{proof}[Proof of Theorem~\ref{ch2thm:Vbound}] \quad \\
	{\em (i)} Let~$\rho=1$. For the sake of simplicity, we suppose that~$c\ge1$, and therefore the second line in~\eqref{ch2Vbound1} is not present
	(the proof of~\eqref{ch2Vbound1} when~$c<1$ is similar, but one has to take into account also the set~$(c,1]\times[0,1]$ and show that it is contained
	in~$\mathcal{V}_{\mathcal{A}}$ by checking the sign of the component
	of the velocity field in the normal direction).
	
	We claim that
	\begin{equation}\label{ch2Thns932}
	{\mbox{$\mathcal{V}_{\mathcal{A}}=\mathcal{S}_c$,}}\end{equation} where~$\mathcal{S}_c$ was
	defined in~\eqref{ch2def:S_c}
	(incidentally,~$\mathcal{S}_c$
	is precisely the right-hand-side of equation~\eqref{ch2Vbound1}).
	
	{F}rom Lemma \ref{ch2lemma:rho=1} we have that for $\rho=1$ and $a>0$ it holds $\mathcal{S}_c=\mathcal{E}(a)\subset \mathcal{V}_{\mathcal{A}}$. Thus, to show \eqref{ch2Thns932} we just need to check that
	\begin{equation}\label{ch2inturn}
	\mathcal{V}_{\mathcal{A}} \subseteq \mathcal{S}_c,\end{equation}
	which is equivalent to
	\begin{equation}\label{ch2dotdot}
	\mathcal{S}_c^C \subseteq \mathcal{V}_{\mathcal{A}}^C,
	\end{equation}
	where the superscript~$C$ denotes the complement of the set in the topology of~$[0,1]\times[0,1]$.
	
	First, by definition we have that
	\begin{equation}\label{ch21617}
	\mathcal{S}_c^C \cap ((0,1]\times\{0\})=\varnothing.
	\end{equation}
	Now, we analyze the behavior of the trajectories at~$\partial \mathcal{S}_c^C$. By Proposition~\ref{ch2prop:dyn}, no trajectory can exit~$\mathcal{S}_c^C$ from a point on~$\partial ([0,1]\times[0,1]) \setminus((0,1]\times\{0\})$. Moreover,
	$\partial \mathcal{S}_c^C \cap ((0,1]\times\{0\})=\varnothing$
	thanks to~\eqref{ch21617} and
	the fact that~$\mathcal{S}_c^C$ is closed in the topology of~$[0,1]\times[0,1]$.
	Hence,
	\begin{equation}\label{ch21648}
	{\mbox{no trajectory can exit~$\mathcal{S}_c^C$ from a point on~$\partial ([0,1]\times[0,1])$.}}
	\end{equation}
	Furthermore, it holds that~$$\partial \mathcal{S}_c^C \cap \big(
	(0,1)\times(0,1)\big)= \left\{ (u,v)\in (0,1)\times(0,1) \;{\mbox{ s.t. }}\; v=\frac{u}{c} \right\}.$$
	The velocity of a trajectory starting on the line~$v=\frac{u}{c}$  in the orthogonal direction pointing inward~$\mathcal{S}_c^C$ is 
	\begin{equation*}
	(\dot{u}, \dot{v})\cdot\frac{(-1,c)}{\sqrt{c^2+1}}=\frac{1}{\sqrt{c^2+1}} (cv-u)(1-u-v)=0,
	\end{equation*}
	the last equality coming from the fact that~$cv=u$ on~$\partial \mathcal{S}_c^C \cap\big( (0,1)\times(0,1)\big)$.
	This means that
	\begin{equation}\label{ch21649}
	{\mbox{no trajectory can exit~$\mathcal{S}_c^C$ from a point on the line~$v=\frac{u}{c}$.}}
	\end{equation}
	
	{F}rom~\eqref{ch21648} and~\eqref{ch21649}, we get that no trajectory exits~$\mathcal{S}_c^C$. Then, by~\eqref{ch21617}, no trajectory starting in~$\mathcal{S}_c^C$ can reach the set~$(0,1]\times\{0\}$, therefore~$\mathcal{S}_c^C \cap \mathcal{V}_{\mathcal{A}}= \varnothing$ and this implies that~\eqref{ch2dotdot} is true. 
	As a result, the proof of~\eqref{ch2inturn} is established and the proof is completed for $\rho=1$.

	\medskip
	
	{\em (ii)} Let~$\rho<1$. For the sake of simplicity, we suppose
	that~$\frac{\rho c(c+1)}{1+\rho c}\ge1$.
	Let~$\mathcal{Y}$ be the set in the right-hand-side of~\eqref{ch2bound:rho<1}, and
	\begin{equation}\label{ch2qwertyuiolkjhgf}
	\mathcal{F}_0:=
	\big\{  (u,v)\in [0,1]\times [0,1]\;{\mbox{ s.t. }}\; 
	v < \gamma_{0}(u) \,\text{ if } \, u\in[0, 1]  \big\}.\end{equation}
	Notice that 
	\begin{equation}\label{ch28ujff994-p-1}
	\mathcal{Y} = \mathcal{F}_0 \cup \mathcal{P},
	\end{equation}
	being~$\mathcal{P}$ the set defined in~\eqref{ch2PPDEFA}.
	
	Moreover,
	\begin{equation}\label{ch28ujff994-p-2}
	\mathcal{P}\subseteq \mathcal{V}_{\mathcal{A}},
	\end{equation}
	thanks to
	Proposition~\ref{ch2prop:construction}.
	
	We also claim that
	\begin{equation}	\label{ch28ujff994-p-3BIS}\mathcal{F}_0\subseteq \mathcal{V}_{\mathcal{K}},\end{equation}
	where~$\mathcal{K}$ is the set of constant functions.
	Indeed, if~$(u,v)\in\mathcal{F}_0$, we have that~$v < \gamma_{0}(u)$
	and consequently~$v < \gamma_{a}(u)$, as long as~$a$ is small enough,
	due to Lemma~\ref{ch2lemma:conv_gamma}.
	
	{F}rom this and Proposition~\ref{ch2prop:char}, we deduce that~$(u,v)$
	belongs to~${\mathcal{E}}(a)$, as long as~$a$ is small enough,
	and this proves~\eqref{ch28ujff994-p-3BIS}.
	
	{F}rom~\eqref{ch28ujff994-p-3BIS} and the fact that~$\mathcal{K}\subseteq\mathcal{A}$,
	we obtain that
	\begin{equation}	\label{ch28ujff994-p-3}\mathcal{F}_0\subseteq \mathcal{V}_{\mathcal{A}}.\end{equation}
	Then, as a consequence of~\eqref{ch28ujff994-p-1},
	\eqref{ch28ujff994-p-2} and~\eqref{ch28ujff994-p-3},
	we get that~$\mathcal{Y}\subseteq \mathcal{V}_{\mathcal{A}}$.
	
	Hence, we are left with proving that
	\begin{equation}\label{ch28iujdpp-1}
	\mathcal{V}_{\mathcal{A}} \subseteq \mathcal{Y}.\end{equation}
	For this, we show that
	\begin{equation}\label{ch28iujdpp-2}
	{\mbox{on~$\partial \mathcal{Y}\cap\big((0,1)\times(0,1)
			\big)$ the outward normal derivative is nonnegative.}}
	\end{equation}
	To prove this,
	we calculate the outward normal derivative on the part of~$\partial \mathcal{Y}$ lying on the graph of~$v=\gamma_0(u)$, that is
	\begin{equation*}
	\dot{v}-\frac{ u^{{\rho} -1}\dot{u}}{c(u^0_s)^{\rho-1} }={\rho} v(1-u-v)-a u -\frac{ u^{{\rho} }(1-u-v-ac)}{ c(u^0_s)^{\rho-1} }.
	\end{equation*}
	By substituting~$v=\gamma_0(u)=\frac{u^\rho}{\rho c(u_s^0)^{\rho-1}}$ we get 
	\begin{eqnarray*}&&
		\dot{v}-\frac{u^{{\rho} -1}\dot{u}}{ c(u^0_s)^{\rho-1} } =
		\frac{u^\rho}{ c(u_s^0)^{\rho-1}}(1-u-v)-a u -\frac{ u^{{\rho} }(1-u-v-ac)}{ c(u^0_s)^{\rho-1} }\\
		&&\qquad=-a u +\frac{ acu^{{\rho} }}{ c(u^0_s)^{\rho-1} }	
		=au^{\rho}  \left( - u^{1-{\rho} } + \frac{1}{(u^0_s)^{\rho-1} } \right)
		.\end{eqnarray*}
	As a result, since~${\rho} <1$, we have
	\begin{equation}\label{ch2ygfbv7r9yty4}
	\dot{v}-\frac{ u^{{\rho} -1}\dot{u}}{c(u^0_s)^{\rho-1} }\geq0 \quad \text{for} \ u\leq u_s^0.
	\end{equation}
	On the part of~$\partial \mathcal{Y}$ contained on the line~$v=\frac{u}{c} + \frac{1-{\rho} }{1+{\rho} c}$, the outward normal derivative is
	\begin{equation}\label{ch27undws8uf8v}\begin{split}&
	\dot{v}-\frac{\dot{u}}{c}= {\rho} v(1-u-v) -au -\frac{u(1-ac-u-v)}{c}=\left({\rho} v-\frac{u}{c}\right)(1-u-v)\\&\qquad\qquad=
	\left(
	\frac{\rho u}{c} + \frac{\rho(1-{\rho}) }{1+{\rho} c}-\frac{u}{c}\right)\left(1-u-
	\frac{u}{c}- \frac{1-{\rho} }{1+{\rho} c}\right)\\&
	\qquad\qquad=
	\left(
	\frac{(\rho-1) u}{c} + \frac{\rho(1-{\rho}) }{1+{\rho} c}\right)\left(-
	\frac{u(c+1)}{c}+ \frac{\rho(1+c) }{1+{\rho} c}\right)
	.\end{split}
	\end{equation}
	We also observe that, when~$u>u_s^0=\frac{\rho c}{1+\rho c}$, the condition~$\rho<1$ gives that
	\begin{eqnarray*}
		\frac{(\rho-1) u}{c} + \frac{\rho(1-{\rho}) }{1+{\rho} c}<
		\frac{\rho (\rho-1)}{1+\rho c}+ \frac{\rho(1-{\rho}) }{1+{\rho} c}=0
	\end{eqnarray*}
	and
	\begin{eqnarray*}-\frac{u(c+1)}{c}+ \frac{\rho(1+c) }{1+{\rho} c}<
		-\frac{\rho(c+1)}{1+\rho c}+ \frac{\rho(1+c) }{1+{\rho} c}=0.
	\end{eqnarray*}
	Therefore, when~$u>u_s^0$, we deduce from~\eqref{ch27undws8uf8v} that~$$
	\dot{v}-\frac{\dot{u}}{c}>0.$$
	Combining this and~\eqref{ch2ygfbv7r9yty4}, we obtain~\eqref{ch28iujdpp-2}, as desired.
	
	Now, by~\eqref{ch28iujdpp-2}, we have that, for any value of~$a$, no trajectory starting in~$\big([0,1]\times[0,1]
	\big)\setminus\mathcal{Y}$ can enter in~$\mathcal{Y}$,
	and in particular no trajectory starting
	in~$\big([0,1]\times[0,1]
	\big)\setminus\mathcal{Y}$
	can hit~$\{v=0\}$, which ends the proof of~\eqref{ch28iujdpp-1}.
	\medskip
	
	{\em (iii)} Let~$\rho>1$. For the sake of simplicity, we suppose
	that~$\frac{c}{(c+1)^\rho}\ge1$.
	Let~$\mathcal{X}$ be the right-hand-side of~\eqref{ch2bound:rho>1}.
	We observe that
	\begin{equation}\label{ch27hperpre923i5}
	\mathcal{X}= \mathcal{S}_{c} \cup \mathcal{Q},
	\end{equation}
	where~$\mathcal{S}_{c}$ was defined in~\eqref{ch2def:S_c} and~$\mathcal{Q}$ in~\eqref{ch2DEFQ}.
	Thanks to Proposition~\ref{ch2prop:bhva},
	one has that~$\mathcal{S}_{c}\subseteq \underset{a>a'}{\bigcup} \mathcal{E}(a)$, for every~$a'>0$, and therefore~$\mathcal{S}_{c}\subseteq\mathcal{V}_{\mathcal{A}}$. 
	Moreover, by the second claim in
	Proposition~\ref{ch2prop:construction}, one also has that~$\mathcal{Q}\subseteq \mathcal{V}_{\mathcal{A}}$. Hence, 
	\begin{equation}\label{ch21923}
	\mathcal{X}\subseteq \mathcal{V}_{\mathcal{A}}.
	\end{equation}
	Accordingly, to prove equality in~\eqref{ch21923} and thus complete
	the proof of~\eqref{ch2bound:rho>1}, we need to show that~$\mathcal{V}_{\mathcal{A}} \subseteq \mathcal{X}$.
	First, we prove that
	\begin{equation}\label{ch21229}
	(0,1]\times\{0\} \subseteq \mathcal{X}.
	\end{equation}
	Indeed, for~$u>0$ we have~$v=\frac{u}{c}>0$, therefore~$(u, 0)\in \mathcal{X}$ for~$u\in (0, u_{\infty}]$. Then,~$\zeta(u)$ is increasing in~$u$ since it is a positive power function, therefore~$v=\zeta(u)>0$ for~$u\in(u_{\infty}, 1]$, hence ~$(u, 0)\in \mathcal{X}$ for~$u\in ( u_{\infty}, 1]$. These observations prove~\eqref{ch21229}.
	
	We now prove that the component of the velocity field in the outward normal direction with respect to~$\mathcal{X}$  is nonnegative on
	\begin{multline*}
	\partial\mathcal{X}\cap \partial( \mathcal{X}^C)= \\
	\left\{ (u,v)\in(0,u_{\infty}]\times(0,1) \ : \ v=\frac{u}{c} \right\} \cup  \left\{ (u,v)\in(u_{\infty},1)\times(0,1) \: \ v=\zeta(u)  \right\}.
	\end{multline*}
	To this end. we observe that on the line~$v=\frac{u}{c}$, the outward normal derivative is
	\begin{equation}\label{ch21853}
	\dot{v}-\frac{1}{c}\dot{u}= \rho v(1-u-v)-au -\frac{u}{c}(1-ac-u-v)=(\rho v -\frac{u}{c})(1-u-v).
	\end{equation}
	The first term is positive because for~$\rho >1$ we have
	\begin{equation*}
	\rho v > v =\frac{u}{c}.
	\end{equation*}
	Moreover, for~$u\leq u_{\infty}$ we have that~$$1-u-v\geq1-u_{\infty}-\frac{u_{\infty}}{c}=0,$$
	thanks to~\eqref{ch2ZETADEF}.
	Thus, the left hand side of~\eqref{ch21853} is nonnegative,
	which proves that the component of the velocity field in the outward normal direction is nonnegative
	on~$\partial\mathcal{X}\cap\left\{v=\frac{u}{c} \right\}$.
	
	On the part of~$\partial \mathcal{X}$ lying in the graph of~$v=\zeta(u)$, the component of the velocity field in the outward normal direction is
	given by
	\begin{equation}\label{ch2po091326uthgbvfjf}
	\dot{v}-\frac{\rho u^{{\rho} -1}\dot{u}}{\rho c(u_{\infty})^{\rho-1} } ={\rho} v(1-u-v)-a u -\frac{{\rho} u^{{\rho} }}{\rho c(u_{\infty})^{\rho-1} }(1-u-v-ac).
	\end{equation}
	Now we substitute~$v=\zeta(u)= \frac{u^{{\rho} }}{\rho c(u_{\infty})^{\rho-1} }$ in~\eqref{ch2po091326uthgbvfjf} and we get 
	\begin{align*}
	\dot{v}-\frac{u^{{\rho} -1}\dot{u}}{ c(u_{\infty})^{\rho-1} } =	au  \left( - 1 + \frac{u^{\rho-1}}{(u_{\infty})^{\rho-1} } \right)
	\end{align*}
	which leads to
	\begin{equation*}
	\dot{v}-\frac{\rho u^{{\rho} -1}\dot{u}}{\rho c(u_{\infty})^{\rho-1} }>0 \quad \text{if} \ u>u_{\infty},
	\end{equation*}
	as desired. 
	
	As a consequence of these considerations, we find that no trajectory starting in~$\mathcal{X}^C$ can enter in~$\mathcal{X}$ and therefore hit~$\{v=0\}$, by~\eqref{ch21229}. Hence, we conclude that~$\mathcal{V}_{\mathcal{A}}\subseteq \mathcal{X}$, which, together with~\eqref{ch21923}, establishes~\eqref{ch2bound:rho>1}.
\end{proof}

\subsection{Proof of Theorem~\ref{ch2thm:W}}

In order to prove Theorem~\ref{ch2thm:W},
we will establish a geometrical lemma in order to
understand the reciprocal position of the function~$\gamma$,
as given by Propositions~\ref{ch2lemma:M} and~\ref{ch2M:p045},
and the straight line where the saddle equilibria lie. To emphasize the dependence of~$\gamma$
on the parameter~$a$ we will often use the notation~$\gamma=\gamma_a$.
Moreover, we recall the notation of the saddle points~$(u_s,v_s)$
defined in~\eqref{ch2usvs} and of the points~$(u_{\mathcal{M}},v_{\mathcal{M}})$ given by
Propositions~\ref{ch2lemma:M} and~\ref{ch2M:p045}, with the convention
that
\begin{equation}\label{ch2usvs2}
{ \mbox{$(u_s,v_s)=(0,0)$ if~$ac\ge1$,}  }
\end{equation}
and we state the following
result:

\begin{lemma}\label{ch2lemma:vett_tg}
	If~$\rho<1$, then 
	\begin{equation}\label{ch2gamma1>r}
	\frac{u}{{\rho}c} \leq \gamma_a(u) \quad \text{ for }  u\in[0, u_s]
	\end{equation}
	and
	\begin{equation}\label{ch2gamma<r}
	\gamma_a(u) \leq \frac{u}{{\rho}c} \quad \text{ for }  u\in[u_s, u_{\mathcal{M}}].
	\end{equation}
	If instead~${\rho}>1$, then
	\begin{equation}\label{ch2gamma1<r}
	\gamma_a(u) \leq \frac{u}{{\rho}c} \quad \text{ for }  u\in[0, u_s]
	\end{equation}
	and
	\begin{equation}\label{ch2gamma>r}
	\frac{u}{{\rho}c} \leq \gamma_a(u)  \quad \text{ for }  u\in[u_s, u_{\mathcal{M}}].
	\end{equation}
	Moreover equality holds in~\eqref{ch2gamma1>r} and~\eqref{ch2gamma1<r} 
	if and only if either~$u=u_s$ or~$u=0$. Also, strict inequality
	holds in~\eqref{ch2gamma<r}
	and~\eqref{ch2gamma>r} for~$u\in(u_s, u_{\mathcal{M}})$.
\end{lemma}

\begin{proof}
	We focus here on the proof of~\eqref{ch2gamma<r}, since
	the other inequalities are proven in a similar way. Moreover,
	we deal with the case~$ac<1$, being the case~$ac\ge1$ analogous with
	obvious modifications.
	
	We suppose by contradiction that~\eqref{ch2gamma<r} does not hold true.
	Namely, we assume that there exists~$\tilde{u}\in(u_s,u_{\mathcal{M}}]$ such that
	$$  \gamma_a(\tilde u) > \frac{\tilde u}{{\rho}c}.$$
	Since~$\gamma_a$ is continuous thanks to Propositions~\ref{ch2lemma:M},
	we have that
	$$  \gamma_a( u) > \frac{ u}{{\rho}c} \quad \mbox{in a neighborhood of~$\tilde u$. }$$
	Hence, we consider the largest open
	interval~$(u_1,u_2)\subset(u_s,u_{\mathcal{M}}]$ containing~$\tilde u$ and such that
	\begin{equation} \label{ch2g>r}
	\gamma_a(u) > \frac{u}{{\rho}c} \quad {\mbox{ for all }} u \in (u_1,u_2).
	\end{equation}
	Moreover, in light of~\eqref{ch2usvs}, we see that
	\begin{equation}\label{ch2togdfgheter}
	\gamma_a(u_s)=v_s= \frac{1-ac}{1+\rho c}=
	\frac{u_s}{\rho c}.\end{equation}
	Hence, by the continuity of~$\gamma_a$, we have that~$\gamma_a(u_1)=\frac{u_1}{\rho c}$
	and 
	\begin{equation}\label{ch2doesnotcon}
	{\mbox{either~$\gamma_a(u_2)=\displaystyle \frac{u_2}{\rho c}$ or~$u_2=u_{\mathcal{M}}$.}}\end{equation}
	
	Now, we consider the set
	\begin{equation*}
	\mathcal{T}:= \left\{  (u,v)\in [u_1,u_2]\times[0,1] \;
	{\mbox{ s.t. }}\; \frac{u}{{\rho}c} < v<  \gamma_a(u)      \right\},
	\end{equation*}
	that is non empty, thanks to~\eqref{ch2g>r}. 
	We claim that
	\begin{equation}\label{ch2lkjhgfds1234567}
	{\mbox{for all~$(u(0), v(0))\in \mathcal{T}$, the~$\omega$-limit
			of its trajectory is~$(u_s,v_s)$.}}\end{equation}
	To prove this, 
	we analyze the normal derivative on 
	\begin{equation*}\begin{split}
	&\partial \mathcal{T} =\mathcal{T}_1\cup\mathcal{T}_2\cup
	\mathcal{T}_3,\\
	{\mbox{where }}\quad &
	\mathcal{T}_1:=\big\{ (u, \gamma_a(u)) \;{\mbox{ with }} 
	u \in (u_1,u_2) \big\},\\
	&\mathcal{T}_2:= \left\{  \left( u, \frac{u}{{\rho}c} \right) \;{\mbox{ with }}
	u \in (u_1,u_2)  \right\}\\
	{\mbox{and }}\quad &
	\mathcal{T}_3:=\left\{ (u_2, v) \;{\mbox{ with }} 
	v \in \left( \frac{u_2}{{\rho}c},\min\{\gamma_a(u_2),1\}\right) \right\}
	,
	\end{split}\end{equation*}
	with the convention that~$\partial \mathcal{T}~$ does contain~$\mathcal{T}_3$
	only if the second possibility in~\eqref{ch2doesnotcon} occurs.
	
	We notice that the set~$\mathcal{T}_1$ is an orbit for the system, and
	thus the component of the velocity in the normal direction is null. 
	On~$\mathcal{T}_2$, we have that the sign of
	the component of the velocity in the inward normal direction is given by
	\begin{equation}\label{ch2der:r}
	\begin{split}&
	(\dot{u},\dot{v})\cdot\left(-\frac1{\rho c},1\right)
	=
	\dot{v} - \frac{1}{{\rho}c} \dot{u} = {\rho} v(1-u-v) 
	-au - \frac{u}{{\rho}c}(1-u-v) + \frac{au}{\rho} \\
	&\qquad\qquad = \frac{u}{c} \left( 1-u-\frac{u}{{\rho}c}   \right)
	\left( 1 - \frac{1}{{\rho}}   \right) -au\left( 1-\frac{1}{\rho}  
	\right)  \\
	&\qquad\qquad = \frac{u}{c} \left(1-\frac{1}{\rho} \right) \left(  
	1-u-\frac{u}{{\rho}c}  -ac \right) .
	\end{split}
	\end{equation}
	Notice that  for~$u \geq u_s$ we have that
	\begin{equation}\label{ch2acapo}
	1-u-v -ac  \leq0,
	\end{equation}
	thus the sign of last term in~\eqref{ch2der:r} depends only on
	the quantity~$1-\frac{1}{\rho}$.
	Consequently, since~${\rho}<1$ the sign of
	the component of the velocity in the inward normal direction is positive. 
	
	Furthermore, in the case in which the second possibility in~\eqref{ch2doesnotcon} occurs,
	we also check the sign of the component of the velocity in the inward normal direction
	along~$\mathcal{T}_3$. In this case, if~$\gamma_a(u_2)<1$ then~$u_2=1$,
	and therefore we find that
	$$(\dot{u},\dot{v})\cdot\left(-1 ,0 \right)=-\dot{u}=-u(1-u-v)+acu=
	v+ac,
	$$
	which is positive. If instead~$\gamma_a(u_2)=1$
	$$(\dot{u},\dot{v})\cdot\left(-1 ,0 \right)=-\dot{u}=-u(1-u-v)+acu=
	-u(1-ac-u-v)
	,
	$$
	which is positive, thanks to~\eqref{ch2acapo}.
	
	We also point out that there are no cycle in~$\mathcal{T}$, since~$\dot{u}$
	has a sign. These considerations and the
	Poincar\'e-Bendixson Theorem (see e.g.~\cite{TESCHL})
	give that the~$\omega$-limit set of~$(u(0),v(0))$
	can be either an equilibrium or a union of (finitely many)
	equilibria and non-closed orbits connecting these equilibria.
	Since~$(0,0)$ and~$(0,1)$ do not belong to the closure of~$\mathcal{T}$,
	in this case the only possibility is that the~$\omega$-limit is the equilibrium~$(u_s,v_s)$.
	Consequently, we have that~$u_1=u_s$, and that~\eqref{ch2lkjhgfds1234567}
	is satisfied.
	
	Accordingly, in light of~\eqref{ch2lkjhgfds1234567}, we have that the
	set~$\mathcal{T}$
	is contained in the stable manifold of~$(u_s,v_s)$, which is in contradiction
	with the definition of~$\mathcal{T}$.
	Hence,~\eqref{ch2gamma<r} is established, as desired.
	\medskip
	
	Now we show that strict inequality holds true in~\eqref{ch2gamma<r}
	if~$u\in(u_s,u_{\mathcal{M}})$.	
	To this end, we suppose by contradiction that
	there exists~$\bar{u}\in (u_s,u_{\mathcal{M}})$ such that
	\begin{equation}\label{ch2equality}
	\gamma_a(\bar{u})=\frac{\bar{u}}{{\rho}c}.
	\end{equation} 
	Now, since~\eqref{ch2gamma<r} holds true, we have that
	the line~$v-\frac{u}{{\rho}c}=0$ is tangent to the curve~$
	v=\gamma_a(u)$ at~$(\bar{u}, \gamma_a(\bar{u}))$,
	and therefore at this point the components of the velocity
	along the normal directions to the curve and to the line coincide. 
	On the other hand,
	the normal derivative at a point on
	the line has a sign, as computed in~\eqref{ch2der:r}, while the normal derivative
	to~$v=\gamma_a(u)$ is~$0$ because the curve is an orbit. 
	
	This, together with~\eqref{ch2togdfgheter}, proves that equality
	in~\eqref{ch2gamma<r} holds true if~$u=u_s$, but strict inequality holds true
	for all~$u\in(u_s,u_{\mathcal{M}})$,
	and thus
	the proof of Lemma~\ref{ch2lemma:vett_tg} is complete.
\end{proof}

For each~$a>0$, we define~$(u_d^a, v_d^a)\in [0,1]\times[0,1]$ as the unique intersection of the graph of~$\gamma_a$ with the line~$\{v=1-u\}$, that is the solution of the system
\begin{equation}\label{ch2ki87yh556g}
\left\{
\begin{array}{l}
v_d^a=\gamma_a(u_d^a),\\
v_d^a=1- u_d^a.
\end{array}
\right.
\end{equation}
We recall that the above intersection is
unique since the function~$\gamma_a$
is increasing. Also, by construction,
\begin{equation}\label{ch2CALM}
u_d^a\le u_{\mathcal{M}}.
\end{equation}
Now, recalling~\eqref{ch2usvs}  and making explicit the dependence on $a$ by writing $u_s^a$
(with the convention in~\eqref{ch2usvs2}), we give the following result:

\begin{lemma}\label{ch2lemma:ord}
	We have that:
	\begin{enumerate}
		\item For $\rho<1$, for all $a^*>0$ it holds that
		\begin{equation}\label{ch21304b}
		\gamma_a(u) \leq \gamma_{a^*}(u) \quad \text{for all} \ a > a^* \ \text{and for all} \ u\in[u_s^{a^*}, u_d^{a^*}].
		\end{equation}
		
		\item  For $\rho>1$, for all $a^*>0$ it holds that
		\begin{equation}\label{ch21819b}
		\gamma_a(u) \leq \gamma_{a^*}(u) \quad\text{for all} \ a < a^* \ \text{and for all} \ u\in[u_s^{a^*}, u_d^{a^*}].
		\end{equation}
	\end{enumerate}
\end{lemma}

\begin{proof}
	We claim that
	\begin{equation}\label{ch21306}
	u_s^{a^*} < u_d^{a^*}.
	\end{equation}
	Indeed, when~$a^* c\ge1$, we have that~$u_s^{a^*} =0< u_d^{a^*}$
	and thus~\eqref{ch21306} holds true. If instead~$
	a^* c<1$, by \eqref{ch2usvs} and~\eqref{ch2ki87yh556g} we have that
	\begin{equation} \label{ch28uhj76tuyg6446r6f6}\gamma_{a^*}(u_s^{a^*})+u_s^{a^*}=1-a^* c< 1=
	\gamma_{a^*}(u_d^{a^*})+u_d^{a^*}.\end{equation}
	Also, since~$\gamma_{a^*}$ is increasing,
	we have that the map~$r\mapsto \gamma_{a^*}(r)+r$
	is strictly increasing. Consequently, we deduce from~\eqref{ch28uhj76tuyg6446r6f6} that~\eqref{ch21306}
	holds true in this case as well.
	
	Now we suppose that $\rho<1$ and we prove~\eqref{ch21304b}.
	For this, we claim that, for every~$a^*>0$ and every~$a>a^*$,
	\begin{equation}\label{ch2xcvbn881300}\gamma_{a}(
	u_s^{a^*})\le\gamma_{a^*}( u_s^{a^*})\quad{\mbox{ with strict inequality when }}
	a^*\in\left(0,\frac{1}{c}\right).
	\end{equation}
	To check this, we distinguish two cases.
	If $a^*\in\left(0,\frac{1}{c}\right)$, then for all $a>a^*$
	\begin{equation}\label{ch21300}
	u_s^a=\max \left\{ 0, \rho c \frac{1-ac}{1+ \rho c} \right\} <  \rho c \frac{1-a^*c}{1+ \rho c}  =u_s^{a^*}.
	\end{equation}
	By \eqref{ch21300} and formula \eqref{ch2gamma<r} in Lemma \ref{ch2lemma:vett_tg}, we have that
	\begin{equation}\label{ch21647}
	\gamma_a(u_s^{a^*}) < \frac{u_s^{a^*}}{\rho c} = \gamma_{a^*}(u_s^{a^*}) \quad \text{for all} \ a> a^*.
	\end{equation}
	If instead $a^*\geq \frac{1}{c}$, then $u_s^{a^*}=0$ and for all $a>a^*$ we have $u_s^a=0$. As a consequence, 
	\begin{equation}\label{ch21647b}
	\gamma_{{a^*}}(u_s^{a^*})=\gamma_{{a}}(u_s^{a^*}) \quad \text{for all} \ a> a^* .
	\end{equation}
	The claim in \eqref{ch2xcvbn881300} thus follows from~\eqref{ch21647} and~\eqref{ch21647b}.
	
	Furthermore, by Propositions \ref{ch2lemma:M} and \ref{ch2M:p045},
	\begin{equation}\label{ch21641b}
	\gamma_a'(0)= \frac{a}{\rho+ac-1} < \frac{{a^*}}{\rho+{a^*}c-1}=\gamma_{a^*}'(0) \quad \text{for all} \ a> a^*\ge\frac1c.
	\end{equation}
	Moreover, for all $a\ge a^*$ and $u>u_s^{a^*}$  it holds that, when~$v=\gamma_{a^*}(u)$,
	\begin{equation}\label{ch21623b}
	-\big(
	acu-u(1-u-v)
	\big)=
	u(1-u-\gamma_{a^*}(u)- ac) < u(1-u_s^{a^*}-v_s^{a^*}-ac)\le 0.
	\end{equation}
	Now, we establish that
	\begin{equation}\label{ch2767675747372}
	u(\rho c  v-u)(1-u-v)(a-a^*) < 0 \ \ \text{for all} \ a > a^*, \   u
	\in(u_s^{a^*}, u_d^{a^*}), \  v=\gamma_{a^*}(u).
	\end{equation}
	Indeed,
	for the values of $a$, $u$ and $v$ as in \eqref{ch2767675747372} we have
	that~$v\le\gamma_{a^*}(u_d^{a^*})$ and hence
	\begin{equation}\label{ch2767675747372-2}
	(1-u-v)> (1-u_d^{a^*}-\gamma_{a^*}(u_d^{a^*}))=0.
	\end{equation}
	Moreover, by formula \eqref{ch2gamma<r} in Lemma \ref{ch2lemma:vett_tg},
	for $u\in(u_s^{a^*}, u_d^{a^*})$ and $v=\gamma_{a^*}(u)$ and  we have that $$\rho c v -u =
	\rho c \gamma_{a^*}(u) -u<
	0.$$
	{F}rom this and~\eqref{ch2767675747372-2}, we see that~\eqref{ch2767675747372} plainly follows, as desired.
	
	As a consequence of~\eqref{ch21623b} and~\eqref{ch2767675747372}, one deduces that,
	for all~$a > a^*$, $u\in(u_s^{a^*}, u_d^{a^*})$ and~$v=\gamma_{a^*}(u)$,
	\begin{equation}\label{ch21335}\begin{split}
	&
	\frac{au- \rho v(1-u-v)}{acu-u(1-u-v)} - \frac{a^* u- \rho v(1-u-v)}{a^* cu-u(1-u-v)} \\=\,&
	\frac{(a-a^*)c \rho uv(1-u-v)-(a-a^*) u^2(1-u-v)}{\big(a cu-u(1-u-v)\big)\big(a^* cu-u(1-u-v)\big)}
	\\=\,&
	\frac{(a-a^*)(1-u-v)u( c \rho v- u)}{\big(a cu-u(1-u-v)\big)\big(a^* cu-u(1-u-v)\big)}
	\\ \le\,&0.
	\end{split}
	\end{equation}
	Now, we define
	\begin{equation}\label{ch23456784jncdkc6knsbd vc83456789} {\mathcal{Z}}(u):=\gamma_a(u)-\gamma_{a^*}(u)\end{equation}
	and we claim that
	\begin{equation}\label{ch24jncdkc6knsbd vc8}
	{\mbox{if~$u_o\in(u_s^{a^*}, u_d^{a^*})$ is such that~${\mathcal{Z}}(u_o)=0$,
			then~${\mathcal{Z}}'(u_o)<0$.}}
	\end{equation}
	Indeed,
	since~$\gamma_a$ is a trajectory for~\eqref{ch2model},
	if~$(u_a(t),v_a(t))$ is a solution of~\eqref{ch2model}, we have that~$
	v_a(t)=\gamma_a(u_a(t))$, whence
	\begin{equation}\label{ch2989u:SMNDnb csn44}
	\begin{split}&\rho v_a(t)(1-u_a(t)-v_a(t)) -au_a(t)=
	\dot v_a(t)=\gamma_a'(u_a(t))\,\dot u_a(t)\\&\qquad=
	\gamma_a'(u_a(t))\big( u_a(t)(1-u_a(t)-v_a(t)) - acu_a(t)\big)
	.\end{split}\end{equation}
	Then, we let~$v_o:=\gamma_a(u_o)$
	and we notice that~$v_o$ coincides also with~$\gamma_{a^*}(u_o)$.
	Hence, we
	take trajectories of the system with parameter~$a$
	and~$a^*$ starting at~$(u_o,v_o)$,
	and by~\eqref{ch21335} we obtain that
	\begin{eqnarray*}0
		&>&\frac{au_o- \rho v(1-u_o-v_o)}{acu_o-u_o(1-u_o-v_o)}-
		\frac{a^*u_o- \rho v(1-u_o-v_o)}{a^*cu_o-u_o(1-u_o-v_o)}\\&=&
		\frac{au_a(0)- \rho v(1-u_a(0)-v_a(0))}{acu_a(0)-u(1-u_a(0)-v_a(0))}-
		\frac{a^*u_{a^*}(0)- \rho v(1-u_{a^*}(0)-v_a(0))}{a^*cu_{a^*}(0)-u(1
			-u_{a^*}(0)-v_{a^*}(0))}\\&
		=&\gamma'_a(u_a(0))-\gamma'_{a^*}(u_{a^*}(0))\\&=&
		\gamma'_a(u_o)-\gamma'_{a^*}(u_o),
	\end{eqnarray*}
	which establishes~\eqref{ch24jncdkc6knsbd vc8}.
	
	Now we claim that
	\begin{equation}\label{ch2xx124ff469}\begin{split}&
	{\mbox{there exists~$\underline{u}\in[u_s^{a^*}, u_d^{a^*}]$
			such that~${\mathcal{Z}}(\underline{u})<0$}}\\&{\mbox{and~${\mathcal{Z}}(u)\le0$ for every~$u\in[u_s^{a^*},\underline{u}]$.}}
	\end{split}\end{equation}
	Indeed, if~$a^*\in\left(0,\frac{1}{c}\right)$,
	we deduce from~\eqref{ch2xcvbn881300}
	that~${\mathcal{Z}}( u_s^{a^*})<0$
	and therefore~\eqref{ch2xx124ff469} holds true with~$\underline{u}:=
	u_s^{a^*}$. If instead~$a^*\ge\frac{1}{c}$, 
	we have that~$u_s^{a}=u_s^{a^*}=0$
	and we deduce from~\eqref{ch2xcvbn881300}
	and~\eqref{ch21641b} that~${\mathcal{Z}}(u_s^{a^*})=0$
	and~${\mathcal{Z}}'(u_s^{a^*})<0$, from which~\eqref{ch2xx124ff469}
	follows by choosing~$\underline{u}:=u_s^{a^*}+\epsilon$
	with~$\epsilon>0$ sufficiently small.
	
	Now we claim that
	\begin{equation}\label{ch2TBP-SP-EL-34}
	{\mathcal{Z}}(u)\le0\qquad{\mbox{for every }}u\in[u_s^{a^*}, u_d^{a^*}].
	\end{equation}
	To prove this, 
	in light of~\eqref{ch2xx124ff469}, it suffices to check that~${\mathcal{Z}}(u)\le0$
	for every~$u\in(\underline{u}, u_d^{a^*}]$.
	Suppose not. Then there exists~$u^\sharp\in(\underline{u}, u_d^{a^*}]$
	such that~${\mathcal{Z}}(u)<0$ for all~$[\underline{u},u^\sharp)$
	and~${\mathcal{Z}}(u^\sharp)=0$. This gives that~$
	{\mathcal{Z}}'(u^\sharp)\ge0$.
	But this inequality is in contradiction with~\eqref{ch24jncdkc6knsbd vc8}
	and therefore the proof of~\eqref{ch2TBP-SP-EL-34}
	is complete.
	
	The desired claim in~\eqref{ch21304b} follows easily
	from~\eqref{ch2TBP-SP-EL-34}, hence we focus now
	on the proof of~\eqref{ch21819b}.
	\medskip 
	
	To this end, we take $\rho>1$ and we 
	claim that, for every~$a^*>0$ and every~$a\in(0,a^*)$,
	\begin{equation}\label{ch2xcvbn881300-ALT56}\gamma_{a}(
	u_s^{a^*})\le\gamma_{a^*}( u_s^{a^*})\quad{\mbox{ with strict inequality when }}
	a^*\in\left(0,\frac{1}{c}\right).
	\end{equation}
	To prove this, we first notice that,
	if $a<a^*<\frac{1}{c}$, then
	\begin{equation*}
	u_s^{a^*}=\rho c \frac{1-a^*c}{1+\rho c} < \rho c \frac{1-ac}{1+\rho c} = u_s^a.
	\end{equation*}
	Hence by \eqref{ch2gamma1<r} in
	Lemma \ref{ch2lemma:vett_tg} we have
	\begin{equation*}
	\gamma_a(u_s^{a^*}) < \frac{u_s^{a^*}}{\rho c} = \gamma_{a^*}(u_s^{a^*}) \quad \text{for} \ a<a^*<\frac{1}{c}, 
	\end{equation*}
	and this establishes~\eqref{ch2xcvbn881300-ALT56}
	when~$a^*\in\left(0,\frac{1}{c}\right)$.
	Thus, we now focus on the case~$a^*\geq \frac{1}{c}$.
	In this situation, we have that~$u_s^{a^*}=0$
	and accordingly~$\gamma_a(u_s^{a^*})=
	\gamma_a(0)=\gamma_{a^*}(0) =
	\gamma_{a^*}(u_s^{a^*})$, that completes the proof
	of~\eqref{ch2xcvbn881300-ALT56}.
	
	In addition, by Propositions \ref{ch2lemma:M} and \ref{ch2M:p045} we have that
	\begin{equation}\label{ch21821}
	\gamma_a'(0)= \frac{a}{\rho-1+ac} \leq \frac{{a^*}}{\rho-1+{a^*}c}=\gamma_{a^*}'(0) \quad \text{for} \ a\in\left[\frac{1}{c}, {a^*}\right].
	\end{equation}
	Moreover, for $u>u_s^a$, if~$v=\gamma_a(u)$ we have that~$v>\gamma_a(u_s^a)=v_s^a$, thanks to the monotonicity of~$\gamma_a$,
	and, as a result,
	\begin{equation}\label{ch21804}
	u(1-u-v-ac)<u(1-u_s^a-v_s^a-ac)=0.
	\end{equation}
	
	Now we claim that, for all~$ a< {a^*}$,
	$u\in(u_s^{a^*}, u_d^{{a^*}})$ and~$  v=\gamma_{{a^*}}(u)$,
	we have
	\begin{equation}\label{ch2473-bniu-1}
	u(1-u-v)({a^*}-a)(u-\rho c v)<0
	.\end{equation}
	Indeed, 
	by the monotonicity of~$\gamma_{{a^*}}$,
	in this situation we have that~$v\le\gamma_{{a^*}}(u^{a^*}_d)$,
	and therefore,
	by \eqref{ch2ki87yh556g}, \begin{equation}\label{ch2473-bniu-2}
	1-u-v >1-u_d^{a^*}-\gamma_{{a^*}}(u_d^{a^*})=1-u_d^{a^*}-1+u_d^{a^*}=0. \end{equation}
	Moreover, by~\eqref{ch2gamma>r}
	in Lemma \eqref{ch2lemma:vett_tg},
	we have that~$ \gamma_{a^*}(u) > \frac{u}{\rho c}$, and
	hence $u-\rho c v> 0$. Combining this inequality with~\eqref{ch2473-bniu-2},
	we obtain~\eqref{ch2473-bniu-1}, as desired.
	
	Now, by~\eqref{ch21804},
	for all~$ a < {a^*}$,
	$u\in(u_s^{a}, u_d^{{a^*}})$
	and~$v=\gamma_{{a^*}}(u)$,
	$$ 
	0<
	-u(1-u-v-ac)=acu-u(1-u-v) <
	{a^*} cu-u(1-u-v)$$
	and then, by~\eqref{ch2473-bniu-1},
	\begin{equation}\label{ch21800}
	\begin{split}&
	\frac{au- \rho v(1-u-v)}{acu-u(1-u-v)} - \frac{{a^*} u- \rho v(1-u-v)}{{a^*} cu-u(1-u-v)} \\=\,&
	\frac{u(1-u-v)({a^*}-a)(u-\rho c v)}{\big(acu-u(1-u-v)\big)
		\big({a^*} cu-u(1-u-v)\big)}
	\\ <\,&0.
	\end{split}
	\end{equation}
	Now we recall the definition of~${\mathcal{Z}}$
	in~\eqref{ch23456784jncdkc6knsbd vc83456789}
	and we claim that
	\begin{equation}\label{ch2567890-4jncdkc6knsbd vc8}
	{\mbox{if~$u_o\in(u_s^{a^*}, u_d^{a^*})$ is such that~${\mathcal{Z}}(u_o)=0$,
			then~${\mathcal{Z}}'(u_o)<0$.}}
	\end{equation}
	To prove this, we let~$v_o:=\gamma_a(u_o)$,
	we notice that~$v_o=\gamma_{a^*}(u_o)$, we
	recall~\eqref{ch2989u:SMNDnb csn44}
	and apply it to a trajectory starting at~$(u_o,v_o)$,
	thus finding that
	\begin{eqnarray*}
		&&\rho v_o(1-u_o-v_a(t)) -au_o=
		\gamma_a'(u_o)\big( u_o(1-u_o-v_o) - acu_o\big)
		.\end{eqnarray*}
	This and~\eqref{ch21800} yield that
	\begin{eqnarray*}
		0>\frac{au- \rho v(1-u-v)}{acu-u(1-u-v)} - \frac{{a^*} u- \rho v(1-u-v)}{{a^*} cu-u(1-u-v)} =\gamma_a'(u_o)-\gamma_{a^*}'(u_o)={\mathcal{Z}}'(u_o),
	\end{eqnarray*}
	which proves the desired claim in~\eqref{ch2567890-4jncdkc6knsbd vc8}.
	
	We now point out that
	\begin{equation}\label{ch26879977xx124ff469}\begin{split}&
	{\mbox{there exists~$\underline{u}\in[u_s^{a^*}, u_d^{a^*}]$
			such that~${\mathcal{Z}}(\underline{u})<0$}}\\&{\mbox{and~${\mathcal{Z}}(u)\le0$ for every~$u\in[u_s^{a^*},\underline{u}]$.}}
	\end{split}\end{equation}
	Indeed, if~$a^*\in\left(0,\frac{1}{c}\right)$,
	this claim follows directly from~\eqref{ch2xcvbn881300}
	by choosing~$\underline{u}:=
	u_s^{a^*}$, while if~$a^*\ge\frac{1}{c}$,
	the claim follows from~\eqref{ch2xcvbn881300}
	and~\eqref{ch24jncdkc6knsbd vc8}
	by choosing~$\underline{u}:=u_s^{a^*}+\epsilon$
	with~$\epsilon>0$ sufficiently small.
	
	Now we claim that
	\begin{equation}\label{ch2jjjjdnfnfTBP-SP-EL-34}
	{\mathcal{Z}}(u)\le0\qquad{\mbox{for every }}u\in[u_s^{a^*}, u_d^{a^*}].
	\end{equation}
	Indeed, by~\eqref{ch26879977xx124ff469},
	we know that the claim is true for all~$u\in[u_s^{a^*},\underline{u}]$.
	Then, the claim for~$u\in(\underline{u}, u_d^{a^*}]$
	can be proved by contradiction,
	supposing that there exists~$u^\sharp\in(\underline{u}, u_d^{a^*}]$
	such that~${\mathcal{Z}}(u)<0$ for all~$[\underline{u},u^\sharp)$
	and~${\mathcal{Z}}(u^\sharp)=0$. This gives that~$
	{\mathcal{Z}}'(u^\sharp)\ge0$, which is in contradiction with~\eqref{ch24jncdkc6knsbd vc8}.
	
	Having completed the proof of~\eqref{ch2jjjjdnfnfTBP-SP-EL-34},
	one can use it to obtain the desired claim in~\eqref{ch21819b}.
\end{proof}

Now we perform the proof 
of Theorem~\ref{ch2thm:W}, analyzing separately the cases~$\rho=1$,~$\rho<1$
and~$\rho>1$.

\begin{proof}[Proof of Theorem~\ref{ch2thm:W}, case $\rho=1$]
	We notice that
	\begin{equation}\label{ch21959}
	\mathcal{V_{\mathcal{K}}}\subseteq 
	\mathcal{V_{\mathcal{A}}},
	\end{equation}
	since~$\mathcal{K}\subset \mathcal{A}$.
	
	Also, from Theorem~\ref{ch2thm:Vbound}, part (i), we get that~$\mathcal{V}_{\mathcal{A}}=\mathcal{S}_c$, where~$\mathcal{S}_c$ was defined in~\eqref{ch2def:S_c}. 
	On the other hand, by Lemma~\ref{ch2lemma:rho=1}, we know that for~$\rho=1$ and for all~$a>0$ we have~${\mathcal{E}(a)}=\mathcal{S}_c$.
	But since every constant~$a$ belongs to the set~$\mathcal{K}$, we have~$\mathcal{E}(a)\subseteq \mathcal{V}_{\mathcal{K}}$.
	This shows that~$\mathcal{V}_{\mathcal{A}} = {\mathcal{E}(a)}\subseteq \mathcal{V}_{\mathcal{K}}$, and together with~\eqref{ch21959} concludes the proof.
\end{proof}

\begin{proof}[Proof of Theorem~\ref{ch2thm:W}, case~$\rho<1$]
	We notice that
	\begin{equation}\label{ch200--fg5996}
	\mathcal{V_{\mathcal{K}}}\subseteq 
	\mathcal{V_{\mathcal{A}}},\end{equation} since~$\mathcal{K}\subset \mathcal{A}$.
	To prove that the inclusion is strict, we aim to find a point~$
	(\bar{u}, \bar{v})\in \mathcal{V}_{\mathcal{A}} \setminus  \mathcal{V}_{\mathcal{K}}$.
	Namely, we have to prove that there exists~$
	(\bar{u}, \bar{v})\in \mathcal{V}_{\mathcal{A}}$ such that,
	for all constant strategies~$a>0$, we have that~$(\bar{u}, \bar{v})\notin \mathcal{E}(a)$,
	that is, by the characterization in Proposition~\ref{ch2prop:char},
	it must hold true that~$\bar{v} \geq \gamma_a(\bar{u})$ and $\bar{u}\leq u_{\mathcal{M}}^a$.
	
	To do this, we define
	\begin{equation}\label{ch2def:f}
	f(u):= \frac{u}{c} +  \frac{1-\rho}{1+\rho c}\quad {\mbox{ and }}\quad
	m:= \min \left\{\frac{\rho c (c+1)}{1+\rho c}, 1 \right\}.
	\end{equation}
	By inspection, one can see that~$(u, f(u))\in[0,1]\times[0,1]$ if and only
	if~$u\in [0,m]$.
	We point out that, by~(ii) of Theorem~\ref{ch2thm:Vbound},
	for~$\rho <1$ and~$u\in [u_s^0, m]$,
	a point~$({u}, {v})$ belongs to~$\mathcal{V}_{\mathcal{A}}$
	if and only if~${v} < f({u})$. Here~$u_s^0$ is defined in~\eqref{ch2u0v0}. We underline that the interval~$[u_s^0, m]$ is non empty since
	\begin{equation}\label{ch22101}
	u_s^0=\frac{\rho c}{1+\rho c}<\min \left\{\frac{\rho c (c+1)}{1+\rho c}, 1 \right\}= m.
	\end{equation}	
	Now we point out that
	\begin{equation}\label{ch21633}
	m \leq u_{\mathcal{M}}^a .
	\end{equation}
	Indeed, by \eqref{ch2def:f} we already know that $m\leq 1$, thus if $u_{\mathcal{M}}^a=1$ the inequality in \eqref{ch21633} is true. On the other hand, when $u_{\mathcal{M}}^a<1$ we have that $(u_{\mathcal{M}}^a,1)\times(0,1)\subseteq{\mathcal{E}}(a)$.
	This and~\eqref{ch200--fg5996} give that~$
	(u_{\mathcal{M}}^a,1)\times(0,1)\subseteq\mathcal{V_{\mathcal{K}}}\subseteq 
	\mathcal{V_{\mathcal{A}}}$. Hence, in view of~\eqref{ch2bound:rho<1}, we deduce that~$\frac{\rho c(c+1)}{1+\rho c}\le u_{\mathcal{M}}^a$. In particular, we find that~$m\le u_{\mathcal{M}}^a$,
	and therefore~\eqref{ch21633} is true also in this case.
	
	With this notation, we claim the existence of a
	value~$\bar{v}\in(0,1]$
	such that for all~$a>0$ we have 
	$\gamma_a(m)\leq  \bar{v} < f(m)$.
	That is, we prove now that  there
	exists~$\theta>0$ such that
	\begin{equation}\label{ch20000}
	\gamma_a(m)+ \theta < f(m) \quad {\mbox{ for all }} a>0.
	\end{equation}
	The strategy is to study two cases separately, namely we prove~\eqref{ch20000}
	for sufficiently small values of~$a$ and then for the other 
	values of~$a$.
	
	To prove~\eqref{ch20000} for small values of~$a$, we start by
	looking at
	the limit function~$\gamma_0$ defined in~\eqref{ch2def:gamma0}.
	One observes that
	\begin{equation}\label{ch2wqffwoe3u8ry4}
	\gamma_0(u_s^0) = v_s^0= \frac{1}{1+\rho c} = \frac{\rho c}{c(1+\rho c)}+\frac{1-\rho}{1+\rho c}=
	f(u_s^0).\end{equation}
	Moreover, for all~$u\in(u_s^0, m]$,
	we have that 
	$$\gamma_0'(u) =\frac{v_s^0}{(u_s^0)^\rho} \,\rho u^{\rho-1}<
	\frac{v_s^0}{(u_s^0)^\rho} \,\rho (u_s^0)^{\rho-1}=\frac{\rho v_s^0}{u_s^0}=
	\frac{1}{c}= f'(u).$$
	Hence, using the fundamental theorem of calculus on the continuous functions $\gamma_{0}(u)$ and $f(u)$, we get 
	\begin{equation*}
	\gamma_{0}(m)= \gamma_{0}(u_s^0) +\int_{u_s^0}^{m} \gamma_{0}'(u)\,du <  f(u_s^0) +\int_{u_s^0}^{m} f'(u)\,du = f(m).
	\end{equation*}
	Then, the quantity  $$\theta_1:= \frac{f(m)-\gamma_0(m)}{4}$$
	is positive and
	we have
	\begin{equation}\label{ch21608}
	\gamma_0(m)+ 2\theta_1 < f(m).
	\end{equation}

	Now, by the uniform convergence of~$\gamma_a$ to~$\gamma_0$
	given by Lemma~\ref{ch2lemma:conv_gamma},
	we know that there exists~$\varepsilon\in \left(0,\frac1c\right)$ such that, if~$a\in(0,\varepsilon]$, 
	\begin{equation}\label{ch2KJ444S}\underset{u\in [u_s^0,m]}{\sup } |\gamma_a(u)-\gamma_0(u)| < {\theta_1}.\end{equation}
	By this and~\eqref{ch21608}, we obtain that
	\begin{equation}\label{ch20000BIS}
	\gamma_a(m) + {\theta_1} < f(m) \quad {\mbox{ for all }} a\in(0,\varepsilon] .
	\end{equation}
	We remark that formula~\eqref{ch20000BIS} will give the desired claim
	in~\eqref{ch20000} for conveniently small values of~$a$.
	
	We are now left with considering the case~$a> \varepsilon$.
	To this end, recalling~\eqref{ch2usvs}, \eqref{ch2ki87yh556g}, by the first statement in
	Lemma~\ref{ch2lemma:ord}, used here with~$a^*:=\varepsilon$,
	we get 
	\begin{equation}\label{ch21304}
	\gamma_a(u) \leq \gamma_{\varepsilon}(u) \quad \text{for all} \ a > \varepsilon \ \text{and for all} \ u\in[u_s^{\varepsilon}, u_d^{\varepsilon}].
	\end{equation}
	Now we observe that
	\begin{equation}\label{ch28j8j8fb8i903-1}
	u^a_d\ge u_s^\varepsilon.
	\end{equation}
	Indeed, suppose not, namely
	\begin{equation}\label{ch28j8j8fb8i903-2}
	u^a_d< u_s^\varepsilon.
	\end{equation}
	Then, by the monotonicity of~$\gamma_a$, we have that~$\gamma_a(u^a_d)\le\gamma_a( u_s^\varepsilon)$.
	This and~\eqref{ch21304} yield that~$\gamma_a(u^a_d)\le\gamma_\varepsilon( u_s^\varepsilon)$. Hence,
	the monotonicity of~$\gamma_\varepsilon$ gives that~$
	\gamma_a(u^a_d)\le\gamma_\varepsilon( u_d^\varepsilon)$.
	This and~\eqref{ch2ki87yh556g} lead to~$1-u^a_d\le1-u_d^\varepsilon$,
	that is~$u_d^\varepsilon\le u^a_d$. {F}rom this inequality,
	using again~\eqref{ch28j8j8fb8i903-2}, we deduce that~$u_d^\varepsilon<
	u_s^\varepsilon$. This is in contradiction with~\eqref{ch21306}
	and thus the proof of~\eqref{ch28j8j8fb8i903-1}
	is complete.
	
	We also notice that
	\begin{equation}\label{ch28j8j8fb8i903-11}
	u^a_d\ge u_d^\varepsilon.
	\end{equation}
	Indeed, suppose not, say
	\begin{equation}\label{ch28j8j8fb8i903-12}
	u^a_d< u_d^\varepsilon.\end{equation}
	Then, by~\eqref{ch28j8j8fb8i903-1}, we have that~$u^a_d\in[u_s^\varepsilon, u_d^\varepsilon]$ and therefore we can apply~\eqref{ch21304}
	to say that~$
	\gamma_a(u^a_d) \leq \gamma_{\varepsilon}(u^a_d)$.
	Also, by the monotonicity of~$\gamma_{\varepsilon}$,
	we have that~$\gamma_{\varepsilon}(u^a_d)\le \gamma_{\varepsilon}(u^\varepsilon_d)$.
	
	With these items of information and~\eqref{ch2ki87yh556g}, we find that
	$$ 1-u^a_d=\gamma_a(u^a_d) \leq
	\gamma_{\varepsilon}(u^\varepsilon_d)=1-u^\varepsilon_d,$$
	and accordingly~$u^a_d\ge u^\varepsilon_d$.
	This is in contradiction with~\eqref{ch28j8j8fb8i903-12}
	and establishes~\eqref{ch28j8j8fb8i903-11}.
	
	Moreover, by~\eqref{ch2usvs} and~\eqref{ch2u0v0},
	we know that~$u_s^0>u_s^{a^*}$, for every~$a^*>0$.
	Therefore, setting~$\tilde u_d^{a^*}:=\min
	\{u_d^{a^*},u_s^0\}$, we have that~$\tilde u_d^{a^*}\in
	[u_s^{a^*},u_d^{a^*}]$. Thus, we are in the position of
	using the first statement 
	in Lemma~\ref{ch2lemma:ord} with~$a:=\varepsilon$
	and deduce that
	\begin{equation}\label{ch290i3883jj889203}
	\gamma_\varepsilon (\tilde u_d^{a^*})\le\gamma_{a^*}(
	\tilde u_d^{a^*})\qquad{\mbox{for all}}\quad a^*<\varepsilon.
	\end{equation}
	We also remark that
	\begin{equation} \label{ch279ihkf843767676}u_d^{a^*}\to u_s^0\qquad{\mbox{as }}\;a^*\to0.\end{equation} Indeed, 
	up to a subsequence we can assume that~$u_d^{a^*}\to \tilde u$ as~$a^*\to0$, for some~$\tilde u\in[0,1]$. Also, by~\eqref{ch2ki87yh556g},
	$$ \gamma_{a^*} (u_d^{a^*})=1-u_d^{a^*},$$
	and then the uniform convergence of~$ \gamma_{a^*}$
	in Lemma~\ref{ch2lemma:conv_gamma} yields that
	$$ \gamma_{0} (\tilde u)=1-\tilde u.$$
	This and~\eqref{ch2ki87yh556g} lead to~$\tilde u=u_d^0$.
	Since
	\begin{equation}\label{ch2SQund0-dis0}
	u_d^0=u_s^0\end{equation} in virtue of~\eqref{ch2u0v0},
	we thus conclude that~$\tilde u=u_s^0$
	and the proof of~\eqref{ch279ihkf843767676}
	is thereby complete.
	
	As a consequence of~\eqref{ch279ihkf843767676},
	we have that~$\tilde u_d^{a^*}\to
	u_s^0$ as~$a^*\to0$. Hence,
	using again the uniform convergence of~$ \gamma_{a^*}$
	in Lemma~\ref{ch2lemma:conv_gamma},
	we obtain that~$\gamma_{a^*}(
	\tilde u_d^{a^*})\to\gamma_{0}(u^0_s)$.
	{F}rom this and~\eqref{ch290i3883jj889203}, we conclude that
	\begin{equation}\label{ch2SQund0-dis1}
	\gamma_\varepsilon (u^0_s)\le\gamma_{0}(
	u^0_s).\end{equation}
	Now we claim that
	\begin{equation}\label{ch29i9i9i78i9u-3934}
	u_d^{\varepsilon} > u_s^0 .\end{equation}
	Indeed, suppose, by contradiction,
	that
	\begin{equation}\label{ch29i9i9i78i9u-3934-0}u_d^{\varepsilon} \le u_s^0.\end{equation} Then, the monotonicity
	of~$\gamma_{\varepsilon} $, together with~\eqref{ch2SQund0-dis0}
	and~\eqref{ch2SQund0-dis1}, gives that
	$$ 1-u_d^{\varepsilon} =
	\gamma_{\varepsilon} (u_d^{\varepsilon}) \le
	\gamma_{\varepsilon} ( u_s^0)=1-u_s^0.$$
	{F}rom this and~\eqref{ch29i9i9i78i9u-3934-0} we deduce that~$u_d^{\varepsilon}=u_s^0$. In particular, we
	have that~$u^0_s\in(u_s^\varepsilon,u^\varepsilon_{\mathcal{M}})$.
	Accordingly, by~\eqref{ch2gamma<r},
	$$ 1- u^0_s=
	1-u_d^{\varepsilon}=
	\gamma_{\varepsilon}(u_d^{\varepsilon})=
	\gamma_{\varepsilon}(u^0_s)< \frac{u^0_s}{{\rho}c} 
	.$$
	As a consequence,
	$$ u^0_s>\frac{\rho c}{1+\rho c},$$
	and this is in contradiction with~\eqref{ch2u0v0}.
	The proof of~\eqref{ch29i9i9i78i9u-3934} is thereby complete.
	
	As a byproduct of~\eqref{ch2SQund0-dis0}
	and~\eqref{ch29i9i9i78i9u-3934}, we have that
	\begin{equation}\label{ch289ujfvdjjgjh599fghjkl6}
	v^{\varepsilon}_d=
	\gamma_{\varepsilon}(u^{\varepsilon}_d)=
	1-u_d^{\varepsilon} <1- u_s^0=1- u_d^0=\gamma_0(u^0_d)
	=\gamma_0(u^0_s)=v^0_s.\end{equation}
	Similarly, by means of~\eqref{ch28j8j8fb8i903-11},
	\begin{equation}\label{ch2Qcbvolr9fjevcanf9d4}
	v^a_d=
	\gamma_a( u^a_d)= 1-u^a_d\le1- u_d^\varepsilon=
	\gamma_\varepsilon(u_d^\varepsilon)=v_d^\varepsilon.
	\end{equation}
	In light of~\eqref{ch28j8j8fb8i903-11}, \eqref{ch29i9i9i78i9u-3934},
	\eqref{ch289ujfvdjjgjh599fghjkl6} and~\eqref{ch2Qcbvolr9fjevcanf9d4},
	we can write that
	\begin{equation}\label{ch21731}
	1>u_d^a \geq u_d^{\varepsilon} > u_s^0 >0 \quad \text{and} \quad 1> v_s^0 > v_d^{\varepsilon} \geq v_d^a >0.
	\end{equation}
	
	\begin{figure}[h] 
		\begin{subfigure}{.4\textwidth}
			\centering
			\includegraphics[scale=0.5]{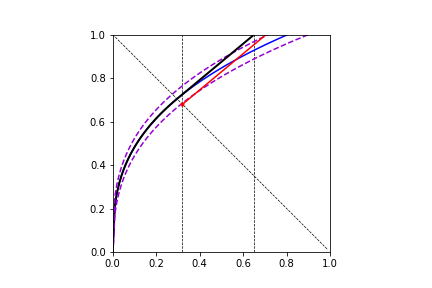}
		\end{subfigure}%
		\begin{subfigure}{.6\textwidth}
			\centering
			\includegraphics[scale=0.5]{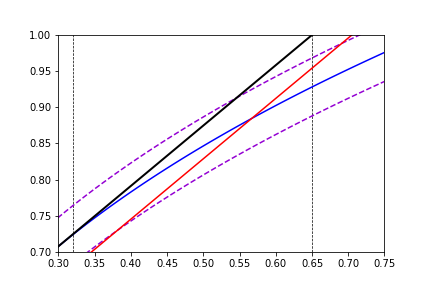}
		\end{subfigure}
		\caption{\it The figures illustrate the functions involved in the proof of Theorem \ref{ch2thm:W} for the case $\rho < 1$.
			The two vertical lines correspond to the values $u_d^{\varepsilon}$ and $m$. The thick black line represents the boundary of $\mathcal{V}_{\mathcal{A}}$; the blue line is the graph of $\gamma_0(u)$; the dark violet lines delimit the area where $\gamma_{a}(u)$ for $a\leq\varepsilon$ might be; the red line is the upper limit of $\gamma_a(u)$ for $a>\varepsilon$. The image was realized using a simulation in Python for the values $\rho=0.35$ and $c=1.2$. }
		\label{ch2fig:thm141}
	\end{figure}
	
	Now,
	to complete the proof of~\eqref{ch20000}
	when~$a>\varepsilon$,
	we consider two cases depending on the order of $m$ and $u_d^{\varepsilon}$. If $u_d^{\varepsilon}\geq m$, by \eqref{ch21731}  we have that~$m<1$ and $f(m)=1$. Then, 
	\begin{equation}\label{ch21802}
	\gamma_a(m) \leq \gamma_a(u_d^{\varepsilon}) \le
	\gamma_{\varepsilon}(u_d^{\varepsilon})= v_d^{\varepsilon} < 1 = f(m),
	\end{equation}
	thanks to the monotonicity of $\gamma_a$, \eqref{ch21304} and \eqref{ch21731}.
	We define $$\theta_2:=\frac{1-v_d^{\varepsilon
	}}{2},$$ which is positive  thanks to \eqref{ch21731}. {F}rom \eqref{ch21802}, we get that
	\begin{equation}\label{ch21815}
	\gamma_{a}(m) +\theta_2 \leq v_d^{\varepsilon} +\theta_2 <1= f(m).
	\end{equation}
	This formula proves the claim  in \eqref{ch20000} for $a>\varepsilon$ and $u_d^{\varepsilon}\geq m$.
	
	If instead~$u_d^{\varepsilon}< m$, then we proceed as follows.
	By \eqref{ch21731} we have
	\begin{equation}\label{ch21722}
	\gamma_a(u_d^{a}) =v_d^{a} \leq v_d^{\varepsilon} <  v_s^0 =  f(u_s^0).
	\end{equation}
	Now we set $$\theta_3:=\frac{f(u_d^{\varepsilon})-f(u_s^0)}{2}.$$
	Using the definition of $f$ in \eqref{ch2def:f}, we see that $$	\theta_3
	= \frac{u_d^{\varepsilon}-u_s^0}{2c} ,$$ 
	and accordingly~$\theta_3$ is positive, due to \eqref{ch21731}. 
	
	{F}rom \eqref{ch21722} we have
	\begin{equation}\label{ch21740}
	\gamma_a(u_d^{a}) + \theta_3 < f(u_s^0) +\theta_3 < f(u_d^{\varepsilon}).
	\end{equation}
	
	Now we show that, on any trajectory~$(u(t),v(t))$ lying
	on the graph of~$\gamma_{a}$, it holds that
	\begin{equation} \label{ch21640}
	\dot{v}(t) > \frac{\dot{u}(t)}{c} \quad \text{provided that} \ u(t)\in( u_d^a,u^a_{\mathcal{M}}) .
	\end{equation}
	To prove this, we first observe that  $u(t)>u_d^{a}> u_s^{a}$,
	thanks to~\eqref{ch21306}.
	Hence, we can exploit formula \eqref{ch2gamma<r} of Lemma~\ref{ch2lemma:vett_tg} and get that 
	\begin{equation}\label{ch28yh78749in2fd9}
	\gamma_a(u(t)) -
	\frac{u(t)}{\rho c}<0.\end{equation}
	Also, by the monotonicity of~$\gamma_a$
	and~\eqref{ch2ki87yh556g},
	$$\gamma_a(u(t))\ge \gamma_a(u_d^a) = 1-u_d^a > 1-u(t).$$
	{F}rom this and \eqref{ch28yh78749in2fd9} it follows that
	\begin{equation*}
	\left(\dot{v}(t) - \frac{\dot{u}(t)}{c} \right)= \rho \left(\gamma_a(u(t)) -
	\frac{u(t)}{\rho c} \right) (1-u(t)-\gamma_a(u(t))) > 0 
	\end{equation*}
	provided that~$ u(t)
	\in( u_d^a,u^a_{\mathcal{M}}) $, and this proves~\eqref{ch21640}.
	
	In addition, for such a trajectory~$(u(t),v(t))$ we have that
	\begin{equation*}\begin{split}&
	\dot{u}(t)=u(t)\,(1-u(t)-\gamma_a(u(t))- ac) \\&\qquad\qquad< u(t)\,(1-u(t)-\gamma_a(u_d^a))=u(t)\,(1-u(t)-1+u_d^a)<0,\end{split}
	\end{equation*}
	provided that~$ u(t)
	\in( u_d^a,u^a_{\mathcal{M}}) $.
	
	{F}rom this and \eqref{ch21640}, we get
	\begin{equation*}
	\gamma_a'(u(t))= \frac{\dot{v}(t)}{\dot{u}(t)} < \frac{1}{c} = f'(u(t))
	,\end{equation*}
	provided that~$ u(t)
	\in( u_d^a,u^a_{\mathcal{M}}) $.
	
	Consequently, taking as initial datum of the trajectory
	an arbitrary point~$(u,\gamma_a(u))$ with~$u\in
	( u_d^a,u^a_{\mathcal{M}}) $, we can write that, for all~$u\in( u_d^a,u^a_{\mathcal{M}})$,
	\begin{equation*}
	\gamma_a'(u)< f'(u).
	\end{equation*}
	As a result, integrating and using~\eqref{ch21304}, 
	for all~$u\in( u_d^a,u^a_{\mathcal{M}})$, we have
	\begin{equation*}
	\gamma_a(u)=
	\gamma_{a}(u_d^a)+ \int_{u_d^a}^{u}\gamma_a'(u)\,du
	< \gamma_{a}(u_d^a)+ \int_{u_d^a}^{u}f'(u)\,du=\gamma_{a}(u_d^a)+f(u)-f(u_d^a)
	.
	\end{equation*}
	Then, making use \eqref{ch21740}, for~$u\in( u_d^a,u^a_{\mathcal{M}})$,
	\begin{equation}\label{ch28781jh98172omOS} 	\gamma_a(u) + \theta_3 <  \gamma_{a}(u_d^a)+f(u)-f(u_d^a)  + \theta_3\le
	f(u)-f(u_d^a)+f(u_d^{\varepsilon}).
	\end{equation}
	Also, recalling~\eqref{ch21731}
	and the monotonicity of~$f$, we see that~$f(u_d^{\varepsilon})\le
	f(u_d^{a})$. Combining this and~\eqref{ch28781jh98172omOS},
	we deduce that
	\begin{equation}\label{ch28781jh98172omOS-0987654-PRE}
	\gamma_a(u) + \theta_3 <f(u)\qquad{\mbox{for all }}u\in( u_d^a,u^a_{\mathcal{M}})
	.\end{equation}
	We also observe that if~$u\in (u_d^\varepsilon, u_d^a]$,
	then the monotonicity of~$\gamma_a$ yields that~$\gamma_a(u)\le
	\gamma_a(u_d^a)$. It follows from this and~\eqref{ch21740}
	that~$\gamma_a(u)+\theta_3 < f(u_d^{\varepsilon})$.
	This and the monotonicity of~$f$ give that
	$$ \gamma_a(u)+\theta_3 < f(u)
	\qquad{\mbox{for all }}u\in(u_d^\varepsilon, u_d^a]
	.$$
	Comparing this with~\eqref{ch28781jh98172omOS-0987654-PRE},
	we obtain
	\begin{equation*}
	\gamma_a(u) + \theta_3 <f(u)\qquad{\mbox{for all }}u\in( u_d^\varepsilon,u^a_{\mathcal{M}})
	\end{equation*}
	and
	therefore
	\begin{equation}\label{ch28781jh98172omOS-0987654}
	\gamma_a(u) + \theta_3 \le f(u)\qquad{\mbox{for all }}u\in[u_d^\varepsilon,u^a_{\mathcal{M}}]
	.\end{equation}
	
	Now, in view of~\eqref{ch21633}, we have that~$m\in
	[u_d^\varepsilon,u^a_{\mathcal{M}}]$.
	Consequently, we can utilize~\eqref{ch28781jh98172omOS-0987654}
	with~$u:=m$ and find that
	\begin{equation}\label{ch2a}
	\gamma_a(m) + \theta_3 \le f(m)
	\end{equation}
	which gives \eqref{ch20000} in the case $a>\varepsilon$ and $u_d^{\varepsilon} \leq m$
	(say, in this case with~$\theta\le\theta_3/2$).
	
	That is, by~\eqref{ch20000BIS},~\eqref{ch21815}  and~\eqref{ch2a}
	we obtain that~\eqref{ch20000} holds true
	for
	\begin{equation*}
	\theta :=\frac12\, \min \left\{ \theta_1, \ \theta_2 , \ \theta_3   \right\}.
	\end{equation*}
	If we choose $\bar{v}:= f(m)-\frac{\theta}{2}$ we have that
	\begin{equation}\label{ch21643}
	0 < \gamma_{a}(m) \leq \bar{v} < f(m) \leq 1.
	\end{equation}
	This completes the proof of Theorem~\ref{ch2thm:W} when~$\rho<1$, in light of 
	the characterizations of $\mathcal{E}(a)$ and $\mathcal{V}_{\mathcal{A}}$ from Proposition \ref{ch2prop:char} and Theorem \ref{ch2thm:Vbound}, respectively.
\end{proof}

Now we focus on the case~$\rho>1$. 

\begin{proof}[Proof of Theorem~\ref{ch2thm:W}, case~$\rho>1$]
	As before, the inclusion~$\mathcal{V_{\mathcal{K}}}\subseteq \mathcal{V_{\mathcal{A}}}$ is trivial since~$\mathcal{K}\subset \mathcal{A}$.
	To prove that it is strict, we aim to find a point~$(\bar{u}, \bar{v})\in \mathcal{V}_{\mathcal{A}}$ such that~$(\bar{u}, \bar{v})\notin \mathcal{V}_{\mathcal{K}}$.
	Thus, we have to prove that there exists~$
	(\bar{u}, \bar{v})\in \mathcal{V}_{\mathcal{A}}$ such that,
	for all constant strategies~$a>0$, we have that~$(\bar{u}, \bar{v})\notin \mathcal{E}(a)$.
	
	To this end, using the characterizations given in Proposition~\ref{ch2prop:char} and Theorem~\ref{ch2thm:Vbound}, we claim that
	\begin{equation}\begin{split}\label{ch21219}
	&\mbox{there exists a point~$(\bar{u}, \bar{v})\in[0,1]\times[0,1]$ } \\ &\mbox{satisfying~$u_{\infty}\leq\bar{u}\leq u_{\mathcal{M}}^a$ and~$\gamma_a(\bar{u}) \leq \bar{v} < \zeta (\bar{u})$ for all~$a>0$.}   
	\end{split}\end{equation}
	For this, we let
	\begin{equation*}
	m:=  \min\left\{1, \frac{c}{(c+1)^{\frac{\rho-1}\rho}} \right\} .
	\end{equation*}
	By \eqref{ch2ZETADEF} one sees that
	\begin{equation}\label{ch21657}
	u_{\infty}<m.
	\end{equation}
	In addition, we point out that
	\begin{equation}\label{ch21633-0987654gfhyf}
	m \leq u_{\mathcal{M}}^a .
	\end{equation}
	Indeed, since $m\leq 1$, if $u_{\mathcal{M}}^a=1$ the desired inequality is obvious. If instead $u_{\mathcal{M}}^a<1$ we have that $(u_{\mathcal{M}}^a,1)\times(0,1)\subseteq{\mathcal{E}}(a)\subseteq\mathcal{V_{\mathcal{K}}}\subseteq 
	\mathcal{V_{\mathcal{A}}}$. Hence, by~\eqref{ch2bound:rho>1},
	it follows that~$\frac{c}{(c+1)^{\frac{\rho-1}\rho}}\le u_{\mathcal{M}}^a$,
	which leads to~\eqref{ch21633-0987654gfhyf}, as desired.
	
	Now we claim that there exists~$\theta >0$ such that
	\begin{equation}\label{ch20017}
	\gamma_a(m) + \theta < \zeta(m) \quad {\mbox{for all }}\; a>0.
	\end{equation}
	
	We first show some preliminary facts for~$\gamma_a(u)$.
	For all~$a>0$, we have that~$\mathcal{E}(a)\subseteq \mathcal{V}_{\mathcal{A}}$. Owing
	to the characterization of~$\mathcal{E}(a)$ from Proposition~\ref{ch2prop:char} and of~$\mathcal{V}_{\mathcal{A}}$ from Theorem~\ref{ch2thm:Vbound} (which can be used here, thanks
	to~\eqref{ch21657} and~\eqref{ch21633-0987654gfhyf}), we get that
	\begin{equation}\label{ch21826}
	\gamma_a(u)\le \frac{u}{c} \quad \text{for all } \ u\in(0, u_{\infty}]\ \text{ and }\ a>0.
	\end{equation}
	This is true in particular for~$u=u_{\infty}$.

	We choose 
	\begin{equation}\label{ch2pthm3:def:delta}
	\delta \in\left(0, \frac{\rho-1}{c}\right)\quad{\mbox{ and }}\quad M:= \max\left\{  \frac{1}{c}, \frac{\rho + \frac{1}{c}+\delta}{\delta c u_{\infty}}  \right\} ,
	\end{equation}
	and we prove~\eqref{ch20017} by treating separately the cases~$a>M$ and~$a\in(0, M]$.
	
	We first consider the case~$a>M$. 
	We let~$(u(t),v(t))$ be a trajectory for \eqref{ch2model} lying on~$\gamma_a$ and we show that
	\begin{equation}\label{ch21721}
	\dot{v}(t)- \left( \frac{1}{c}+\delta  \right)\dot{u} (t)> 0 \quad \text{provided that } \ u(t)>u_{\infty}\
	{\mbox{ and }}\ a>M.
	\end{equation}
	To check this, we observe that
	\begin{align*}
	\dot{v}(t)- \left( \frac{1}{c}+\delta  \right)\dot{u} (t)&=  \left[ \rho \gamma_{a}(u(t))
	-\left( \frac{1}{c}+\delta  \right) u(t)\right](1-u(t)-\gamma_{a}(u(t)))+ \delta acu(t) \\ 
	& \geq  - \left\vert \rho + \frac{1}{c} + \delta \right\vert
	+ \delta a c u_{\infty} >0,	
	\end{align*}
	where the last inequality is true thanks to the hypothesis~$a>M$ and the definition of~$M$
	in~\eqref{ch2pthm3:def:delta}. This proves~\eqref{ch21721}.
	
	Moreover, for~$a>M\geq \frac{1}{c}$ we have $\dot{u}<0$.
	{F}rom this, \eqref{ch21721} and the invariance of~$\gamma_a$ for the flow, we get
	\begin{equation}\label{ch21905}
	\gamma_a'(u(t))=\frac{\dot v(t)}{\dot u(t)}< 
	\frac{1}{c}+ \delta ,
	\end{equation}	
	provided that $u(t)>u_{\infty}$ and $a>M$.
	
	For this reason and~\eqref{ch21826}, we get
	\begin{equation}\label{ch2pthm3:1643}
	\gamma_a(u(t)) = \gamma_a(u_{\infty})  + \int_{u_{\infty}}^{u(t)} \gamma_a'(\tau)\,d\tau\le \frac{u_{\infty}}{c} + \left( \frac{1}{c}+\delta \right)(u(t)-u_{\infty}
	)
	\end{equation}
	provided that $u(t)>u_{\infty}$ and $a>M$.
	
	Furthermore, thanks to the choice of~$\delta$ in~\eqref{ch2pthm3:def:delta}, we have
	\begin{equation*}
	\zeta'(u)= \frac{\rho u^{\rho-1}}{c u_{\infty}^{\rho-1}}>\frac{\rho}{c} > \frac{1}{c}+\delta \quad \text{for all } \ u>u_{\infty}.
	\end{equation*} 
	Since also~$\zeta(u_{\infty})=\frac{u_{\infty}}{c}$, by \eqref{ch2pthm3:1643} we deduce that
	\begin{equation}\label{ch21621-PRE}
	\gamma_a(u(t))\leq    \frac{u_{\infty}}{c} + \left( \frac{1}{c}+\delta \right)(u(t)-u_{\infty}) < \zeta(u_{\infty}) + \int_{u_{\infty}}^{u(t)} \zeta'(\tau)\, d\tau = \zeta(u(t)),
	\end{equation}
	provided that $u(t)>u_{\infty}$ and $a>M$.
	
	In particular, given any~$u>u_\infty$, we can take a trajectory starting at~$(u,\gamma_a(u))$
	and deduce from~\eqref{ch21621-PRE} that
	\begin{equation*}
	\gamma_a(u)\leq    \frac{u_{\infty}}{c} + \left( \frac{1}{c}+\delta \right)(u-u_{\infty}) < \zeta(u_{\infty}) + \int_{u_{\infty}}^{u} \zeta'(\tau)\, d\tau = \zeta(u),
	\end{equation*}
	whenever $a>M$. We stress that, in light of~\eqref{ch21657}, we can take~$u:=m$
	in the above chain of inequalities, concluding that
	\begin{equation*}
	\gamma_a(m)\le	  \frac{u_{\infty}}{c} + \left( \frac{1}{c}+\delta \right)(m-u_{\infty}) < \zeta(m)
	.\end{equation*}
	We rewrite this in the form
	\begin{equation}\label{ch21621}
	\gamma_a(m)\le	\left( \frac{1}{c}+\delta  \right)m-{\delta}u_{\infty}< \zeta(m)
	.\end{equation}
	We define
	\begin{equation}\label{ch2pthm3:def:theta}
	\theta_1:=\frac{1}{2}\left[ \zeta(m)-\left( \frac{1}{c}+\delta  \right)m +{\delta}u_{\infty} \right],
	\end{equation}
	that is positive thanks to the last inequality in~\eqref{ch21621}. Then by the first inequality in~\eqref{ch21621} we have
	\begin{equation*}
	\gamma_a(m)+\theta_1 
	\le \left( \frac{1}{c}+\delta  \right)m-{\delta}u_{\infty}+\theta_1=
	\frac12\left[ \left( \frac{1}{c}+\delta  \right)m-{{\delta}u_{\infty}}\right]
	+\frac{ \zeta(m)}2.	\end{equation*}
	Hence, using again the last inequality in~\eqref{ch21621}, we obtain that
	\begin{equation}\label{ch22001}
	\gamma_a(m)+\theta_1<\zeta(m),\end{equation}
	which gives the claim in \eqref{ch20017} for the case~$a>M$.
	
	Now we treat the case~$a\in(0, M]$.
	We claim that
	\begin{equation}\label{ch21204}
	u_d^M > u_{\infty}.
	\end{equation}
	Here, we are using the notation~$u_d^M$ to denote the point~$u_d^a$ when~$a:=M$. To prove~\eqref{ch21204} we argue as follows.
	Since $M\geq \frac{1}{c}$, by Propositions \ref{ch2lemma:M} and~\ref{ch2M:p045} we have
	\begin{equation}\label{ch21719}
	\gamma_M'(0) = \frac{M}{\rho-1+Mc} < \frac{1}{c}.
	\end{equation}
	Moreover, since the graph of $\gamma_M(u)$ is a parametrization of a trajectory for \eqref{ch2model} with $a=M$, we have that $ \dot{v}(t)= \gamma_M'(u(t)) \dot{u}(t)$.
	Hence, at all points $(\bar{u}, \bar{v})$ with~$\bar{u}\in(0, u_{\infty})$ and $\bar{v}=\gamma_M(\bar{u})$ we have
	\begin{equation}\label{ch21734}
	\gamma_M ' (\bar{u}) = \frac{M \bar{u} - \rho \bar{v} (1-\bar{u}-\bar{v}) }{Mc \bar{u} - \bar{u}(1-\bar{u}-\bar{v})}.
	\end{equation}
	We stress that the denominator in the right hand side of~\eqref{ch21734}
	is strictly positive, since~$M\geq \frac{1}{c}$ and~$\bar u>0$.
	
	In addition, we have that  
	\begin{equation}\label{ch21757}
	\begin{split}
	\frac{1}{c} - \frac{M \bar{u} - \rho \bar{v} (1-\bar{u}-\bar{v}) }{Mc \bar{u} - \bar{u}(1-\bar{u}-\bar{v})} = \frac{(\rho c \bar{v}-\bar{u})(1-\bar{u}-\bar{v}) }{Mc^2 \bar{u} - c\bar{u}(1-\bar{u}-\bar{v})}.
	\end{split}
	\end{equation}
	Also,
	$$
	u_s^M=0<\bar{u}<u_\infty<m\le u^M_{\mathcal{M}},
	$$
	thanks to~\eqref{ch21657} and~\eqref{ch21633-0987654gfhyf}.
	Hence, we can exploit formula~\eqref{ch2gamma>r} in Lemma \ref{ch2lemma:vett_tg}
	with the strict inequality, thus obtaining that
	\begin{equation}\label{ch28ujINtdensnumeok3965}
	\rho c \bar{v}-\bar{u}=\rho c\gamma_M(\bar{u})-\bar{u} >0.\end{equation}
	Moreover, by \eqref{ch21826},
	$$ 1-\bar{u}-\bar{v} = 
	1-\bar{u}-\gamma_M(\bar{u})\ge
	1-\bar{u}-\frac{\bar{u}}{c} > 1-u_{\infty}- \frac{u_{\infty}}{c}=0. $$
	Therefore, using the latter estimate and~\eqref{ch28ujINtdensnumeok3965}
	into~\eqref{ch21757}, we get that
	\begin{equation*}
	\begin{split}
	\frac{1}{c} - \frac{M \bar{u} - \rho \bar{v} (1-\bar{u}-\bar{v}) }{Mc \bar{u} - \bar{u}(1-\bar{u}-\bar{v})} > 0.
	\end{split}
	\end{equation*}
	{F}rom this and \eqref{ch21734}, we have that
	\begin{equation*}
	\gamma_M'(u) < \frac{1}{c} \quad \text{for all} \ u\in(0,u_{\infty}).
	\end{equation*}
	This, together with \eqref{ch21719} and the fact that $\gamma_M(0)=0$, gives
	\begin{equation*}
	\gamma_M(u)
	=\gamma_M(u)-\gamma_M(0)=\int_0^u \gamma_M'(\tau)\,d\tau
	< \frac{u}{c} 
	\end{equation*}
	for all~$ u\in(0,u_{\infty}]$. 
	This inequality yields that
	\begin{equation}\label{ch21810}
	\gamma_M(u_{\infty}) <  \frac{u_{\infty}}{c}= 1-u_{\infty}.
	\end{equation}
	Now, to complete the proof of~\eqref{ch21204} we argue by contradiction
	and suppose that the claim in~\eqref{ch21204} is false, hence
	\begin{equation}\label{ch21814}
	u_d^M \leq u_{\infty}.
	\end{equation}
	Thus, by \eqref{ch21810}, the monotonicity of $\gamma_M(u)$ and the definition of $u_d^M$ given in \eqref{ch2ki87yh556g}, we get
	\begin{equation*}
	1-u_d^M  = \gamma_M(u_d^M) \le \gamma_M(u_{\infty}) < 1-u_{\infty} 
	\end{equation*}
	which is in contraddiction with \eqref{ch21814}. Hence, \eqref{ch21204} holds true, as desired.
	
	Also, by the second statement in Lemma \ref{ch2lemma:ord}, used here with~$a^*:=M$,
	\begin{equation}\label{ch21819}
	\gamma_a(u) \leq \gamma_{M}(u) \quad \text{for all } \ u\in[0, u_d^M].
	\end{equation}
	We claim that
	\begin{equation}\label{ch2181967890p67890-4567890456789}
	u_d^M\le u^a_d.
	\end{equation}
	Indeed, suppose, by contradiction, that
	\begin{equation}\label{ch2181967890p67890-4567890456789PRE}
	u_d^M>u^a_d.
	\end{equation}
	Then, by the monotonicity of~$\gamma_a$ and~\eqref{ch21819}, used here with~$u:=u^M_d$,
	we find that
	$$ 1-u^a_d=\gamma_a(u^a_d)\le\gamma_a(u^M_d) \leq \gamma_{M}(u^M_d)=1-u^M_d.$$
	This entails that~$u^a_d\ge u^M_d$, which is in contradiction with~\eqref{ch2181967890p67890-4567890456789PRE},
	and thus establishes~\eqref{ch2181967890p67890-4567890456789}.
	
	We note in addition that
	\begin{equation}\label{ch209876543988878-00-181967890p67890-4567890456789PRE}
	v_d^M =\gamma_M(u_d^M)=1-u_d^M
	<1-u_{\infty},
	\end{equation}
	thanks to the definition of $(u_d^M, v_d^M)$ and~\eqref{ch21204}.
	
	Similarly, by~\eqref{ch2181967890p67890-4567890456789},
	\begin{equation}\label{ch209876543988878-00-181967890p67890-4567890456789PRE098-08}
	v_d^a =\gamma_a(u_d^a)=1-u_d^a\le
	1-u_d^M=\gamma_M(u_d^M)
	= v_d^M.
	\end{equation}
	Collecting the pieces of information in~\eqref{ch21204}, \eqref{ch2181967890p67890-4567890456789},
	\eqref{ch209876543988878-00-181967890p67890-4567890456789PRE}
	and~\eqref{ch209876543988878-00-181967890p67890-4567890456789PRE098-08},
	we thereby conclude that, for all~$a\in(0,M]$,
	\begin{equation}\label{ch21834}
	0<	u_{\infty} < u_d^M \leq u_d^a<1 \qquad \text{and} \qquad 0< v_d^a \leq v_d^M < 1-u_{\infty} =:v_\infty<1.
	\end{equation}
	
	Now we consider two cases depending on the order of $m$ and $u_d^M$. If $u_d^M\geq m$,  by \eqref{ch21834} we have $m<1$ and $\zeta(m)=1$. Accordingly,
	for $a\in(0,M]$, by \eqref{ch21834} and \eqref{ch21819} we have 
	\begin{equation*}
	\gamma_a(m) \le \gamma_{a}(u_d^M) \leq  \gamma_{M}(u_d^M)=v_d^M < 1=\zeta (m).
	\end{equation*}
	Hence, we can define $$\theta_2:=\frac{1-v_d^M}{2},$$
	and observe that~$\theta_2$ is positive by \eqref{ch21834}, thus obtaining that
	\begin{equation}\label{ch21916}
	\gamma_{a}(m) + \theta_2 < \zeta(m).
	\end{equation}
	This is the desired claim in \eqref{ch20017} for $a\in(0,M]$ and $u^*\geq m$.
	
	\begin{figure} 
		\centering
		\includegraphics[scale=0.5]{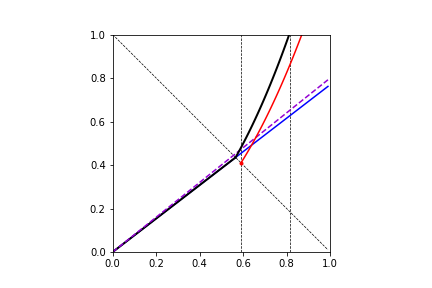}
		\caption{\it The figure illustrates the functions involved in the proof of Theorem \ref{ch2thm:W} for the case $\rho > 1$.
			The two vertical lines correspond to the values $u_d^{M}$ and $m$. The thick black line represents the boundary of $\mathcal{V}_{\mathcal{A}}$; the blue line is the graph of the line $v=\frac{u}{c}$; the dark violet line is the upper bound for $\gamma_{a}(u)$ for $a>M$; the red line is  $\phi(u)$. The image was realized using a simulation in Python for the values $\rho=2.3$ and $c=1.3$. }
		\label{ch2fig:thm142}
	\end{figure}
	
	If instead $u_d^M<m$, we consider the function
	\begin{equation*}
	\phi(u) := v_d^M\,\left(\frac{u}{u_d^M}\right)^{\rho} , \quad{\mbox{ for }} u\in[u_d^M, m]
	\end{equation*}and we claim that 
	\begin{equation}\label{ch21730}
	\gamma_a(u)\leq\phi(u) \quad \text{for all } \ a\in(0,M]\ {\mbox{ and }} \ u\in[u_d^M,m].
	\end{equation}
	To prove this, we recall \eqref{ch21834} and the fact that~$\gamma_a$ is an increasing function
	to see that
	\begin{equation}\label{ch21946}\gamma_a(u_d^M)\le 
	\gamma_a(u_d^a) =v_d^a \leq v_d^M = \phi(u_d^M) .
	\end{equation} 
	
	Now we remark that
	$$ \gamma_M(u_d^M)+u_d^M=1>1-Mc=\gamma_M(u_s^M)+u_s^M,$$
	and therefore~$u_d^M>u_s^M$. Notice also that~$u_d^M<m\le 
	u^M_{\mathcal{M}}$, thanks to~\eqref{ch21633-0987654gfhyf}.
	As a result, we find that~$\rho c \gamma_M(u_d^M) > u_d^M$ by inequality~\eqref{ch2gamma>r}
	in Lemma \ref{ch2lemma:vett_tg}. Therefore, if~$u\ge u_d^M$ and~$v=\phi(u)$, then
	\begin{equation*}\begin{split}&
	au \left(   1-\rho c \frac{v_d^M}{(u_d^M)^{\rho}} u^{\rho -1}   \right)
	= au \left(   1- \frac{\rho c\gamma_M(u_d^M) }{(u_d^M)^{\rho}} u^{\rho -1}   \right)\\&\qquad<
	au \left(   1- \left(\frac{u}{u_d^M} \right)^{\rho -1}   \right)
	\leq 0=\rho \left( v- \frac{v_d^M}{(u_d^M)^{\rho}} u^{\rho}  \right) (1-u-v).\end{split}
	\end{equation*}
	Using this and	\eqref{ch21804}, we deduce that, if~$a\in[0,M]$, $u\in[u_d^M, m]$ and~$v=\phi(u)$,
	\begin{equation}\label{ch21858}\begin{split}&
	\frac{au-\rho v (1-u-v)}{acu - u(1-u-v)}-\frac{v_d^M}{(u_d^M)^{\rho}} \rho u^{\rho-1}
	\\=\;&
	\frac{au-\rho v (1-u-v)
		-\big(acu - u(1-u-v) \big)\,\frac{v_d^M}{(u_d^M)^{\rho}} \rho u^{\rho-1}
	}{acu - u(1-u-v)}\\=\;&
	\frac{au\left(1-\rho c\frac{v_d^M}{(u_d^M)^{\rho}} u^{\rho-1}\right)-\rho (1-u-v)\left(v-
		\frac{v_d^M}{(u_d^M)^{\rho}}  u^{\rho}\right)
	}{acu - u(1-u-v)}\\ <\;&0.
	\end{split}
	\end{equation}
	Now we take~$ a\in(0,M]$, $u\in[u_d^M, m]$ and
	suppose that~$v=\phi(u)=\gamma_a(u)$,
	we consider an orbit~$(u(t),v(t))$ lying on~$\gamma_a$ with~$(u(0),v(0))=(u,v)$,
	and we notice that, by~\eqref{ch21804} and~\eqref{ch21858},
	\begin{equation}\label{ch21951}\begin{split}&
	\gamma_a'(u)=	\gamma_a'(u(0))=\frac{\dot v(0)}{\dot u(0)}
	=	\frac{au(0)-\rho v(0)\, (1-u(0)-v(0))}{acu (0)- u(0)(1-u(0)-v(0))}\\&\qquad
	=	\frac{au-\rho v\, (1-u-v)}{acu - u(1-u-v)}
	<
	\frac{v_d^M}{(u_d^M)^{\rho}} \rho u^{\rho-1}
	= \phi'(u).
	\end{split}\end{equation}
	
	To complete the proof of~\eqref{ch21730},
	we define
	$$ {\mathcal{H}}(u):=\gamma_a(u)-\phi(u)$$
	and we claim that for every~$a\in(0,M]$ there exists~$\underline{u}\in[u_d^M, m]$
	such that
	\begin{equation}\label{ch298989898kjkjkjkdfbv}
	{\mbox{${\mathcal{H}}(\underline u)<0$ and~${\mathcal{H}}(u)\le0$ for every~$u\in[u_d^M,\underline u]$.}}
	\end{equation}
	Indeed, by~\eqref{ch21946}, we know that~${\mathcal{H}}(u_d^M)\le0$.
	Thus, if~${\mathcal{H}}(u_d^M)<0$ then we can choose~$
	\underline u:=u_d^M$ and obtain~\eqref{ch298989898kjkjkjkdfbv}.
	If instead~${\mathcal{H}}(u_d^M)=0$, we have that~$
	\gamma_a(u_d^M)=\phi(u_d^M)$ and thus we can
	exploit~\eqref{ch21951} and find that~${\mathcal{H}}'(u_d^M)<0$, from which
	we obtain~\eqref{ch298989898kjkjkjkdfbv}.
	
	Now we claim that, for every~$ a\in(0,M]$ and $u\in[u_d^M, m]$,
	\begin{equation}\label{ch2KL:0ksf3566}
	{\mathcal{H}}(u)\le0.
	\end{equation}
	For this, given~$ a\in(0,M]$, we define
	$$ 
	{\mathcal{L}}:=\{ u_*\in [u_d^M, m]{\mbox{ s.t. }}{\mathcal{H}}(u)\le0 {\mbox{ for every }}u\in[u_d^M,u_*]\}
	\qquad{\mbox{and}}\qquad
	\overline u:=\sup {\mathcal{L}}.$$
	We remark that~$\underline u\in{\mathcal{L}}$, thanks to~\eqref{ch298989898kjkjkjkdfbv}
	and therefore~$\overline u$ is well defined.
	We have that
	\begin{equation}\label{ch212jnSikjm239gfvhb37}
	\overline u=m,
	\end{equation}
	otherwise we would have that~${\mathcal{H}}(\overline u)=0$
	and thus~${\mathcal{H}}'(\overline u)<0$, thanks to~\eqref{ch21951},
	which would contradict the maximality of~$\overline u$.
	Now, the claim in~\eqref{ch2KL:0ksf3566} plainly follows from~\eqref{ch212jnSikjm239gfvhb37}.
	
	We notice that by the inequalities in~\eqref{ch21834} we have
	\begin{equation}\label{ch22007}
	\zeta(u)= \frac{v_{\infty}}{(u_{\infty})^{\rho}} u^{\rho}>
	\frac{v_d^M}{(u_d^M)^{\rho}} u^{\rho}
	= \phi(u).
	\end{equation}
	Then, we define
	\begin{equation}\label{ch21958}
	\theta_3:= \frac{\zeta(m)-\phi(m)}{2},
	\end{equation}	
	that is positive thanks to \eqref{ch22007}.
	We get that
	\begin{equation}\label{ch21912}
	\phi(m)+\theta_3 < \zeta(m). 
	\end{equation}
	{F}rom this and~\eqref{ch21730}, we conclude that  
	\begin{equation}\label{ch20115}
	\gamma_a(m) + \theta_3 \leq \phi(m) + \theta_3 < \zeta(m) \quad \ \text{for} \ a\in(0,M].
	\end{equation}
	
	By~\eqref{ch22001},~\eqref{ch21916} and~\eqref{ch20115} we have that~\eqref{ch20017} is true for~$\theta = \min \{\theta_1, \ \theta_2, \ \theta_3  \}$. 
	
	This also establishes the claim in \eqref{ch21219}, and the proof is completed.  
\end{proof}

\subsection{Proof of Theorem~\ref{ch2thm:H}}

Now, we can complete the proof of Theorem~\ref{ch2thm:H} by building on the previous work.

\begin{proof}[Proof of Theorem~\ref{ch2thm:H}]
	Since the class of Heaviside functions~$\mathcal{H}$ is contained in the class of piecewise continuous functions~$\mathcal{A}$, we have that
	\begin{equation}
	\mathcal{V}_{\mathcal{H}}\subseteq \mathcal{V}_{\mathcal{A}},
	\end{equation}
	hence we are left with proving the converse inclusion. We treat separately the cases $\rho=1$, $\rho<1$ and $\rho>0$.
	
	If~$\rho=1$, the desired
	claim follows from Theorem \ref{ch2thm:W}, part (i).
	
	If~$\rho<1$, we deduce from~\eqref{ch2bound:rho<1} and~\eqref{ch28ujff994-p-1} that
	\begin{equation}\label{ch2iwfewuguew387627}
	\mathcal{V}_{\mathcal{A}}= \mathcal{F}_0 \cup \mathcal{P},
	\end{equation}
	where~$\mathcal{P}$ has been defined 
	in~\eqref{ch2PPDEFA}
	and~$\mathcal{F}_0$ in~\eqref{ch2qwertyuiolkjhgf}.
	
	Moreover, by~\eqref{ch28ujff994-p-3BIS}, we have that
	\begin{equation}\label{ch2iwfewuguew38762722} \mathcal{F}_0\subseteq \mathcal{V}_{\mathcal{K}}\subseteq \mathcal{V}_{\mathcal{H}}.	\end{equation}
	Also, in Proposition~\ref{ch2prop:construction} we construct a Heaviside winning strategy for every point in~$ \mathcal{P}$. Accordingly, it follows that~$ \mathcal{P} \subseteq \mathcal{V}_{\mathcal{H}}$. 
	This,~\eqref{ch2iwfewuguew387627} and~\eqref{ch2iwfewuguew38762722} entail that~$ \mathcal{V}_{\mathcal{A}} \subseteq \mathcal{V}_{\mathcal{H}}$, which completes the proof of
	Theorem~\ref{ch2thm:H} when~$\rho<1$.
	
	Hence, we now focus on the case~$\rho>1$.  By~\eqref{ch2bound:rho>1}
	and~\eqref{ch27hperpre923i5},
	\begin{equation}\label{ch28877SA}
	\mathcal{V}_{\mathcal{A}}= \mathcal{S}_{c} \cup \mathcal{Q},
	\end{equation}
	where~$\mathcal{S}_{c}$ was defined in~\eqref{ch2def:S_c} and~$\mathcal{Q}$ in~\eqref{ch2DEFQ}.
	
	For every point~$(u_0, v_0)\in\mathcal{S}_{c}$ there exists~$\bar{a}$ that is a constant winning strategy for~$(u_0, v_0)$, thanks to Proposition~\ref{ch2prop:bhva},
	therefore~$\mathcal{S}_{c}\subseteq
	\mathcal{V}_{\mathcal{H}}$.
	Moreover, in Proposition~\ref{ch2prop:construction} for every point~$(u_0, v_0)\in \mathcal{Q}$ we constructed a Heaviside winning strategy, whence~$ \mathcal{Q} \subseteq \mathcal{V}_{\mathcal{H}}$. In light of these observations
	and~\eqref{ch28877SA}, we see that also in this case~$ \mathcal{V}_{\mathcal{A}} \subseteq \mathcal{V}_{\mathcal{H}}$ and the proof is complete.
\end{proof}

\subsection{Bounds on winning initial positions under pointwise
	constraints for the possible strategies}

This subsection is dedicated to the analysis of~$\mathcal{V}_{\mathcal{A}}$ when we put some constraints on~$a(t)$. In particular, we consider $M\geq m \geq 0$ with $M>0$ and the set $\mathcal{A}_{m,M}$ of the functions $a(t)\in\mathcal{A}$ with $m\leq a(t)\leq M$ for all $t>0$.  We will prove Theorem \ref{ch2thm:limit} via
a technical proposition giving informative bounds on $\mathcal{V}_{{m,M}}$.

For this, we denote by~$(u_s^m,v_s^m)$ the point~$(u_s,v_s)$ introduced in~\eqref{ch2usvs}
when~$a(t)=m$ for all~$t>0$ (this when~$mc<1$, and we use the convention that~$(u_s^m,v_s^m)=(0,0)$
when~$mc\ge1$).
In this setting, we have the following result obtaining explicit
bounds on the favorable set~$\mathcal{V}_{{m,M}}$:

\begin{proposition} \label{ch2prop:limit}
	Let~$M\geq m\geq 0$ with $M>0$
	and
	\begin{equation}\label{ch2RANGEEP}
	\varepsilon\in\left(0,\,\min\left\{\frac{M(c+1)}{M+1},1\right\}\right).\end{equation}
	Then
	\begin{itemize}
		\item[(i)] If~$\rho<1$, we have
		\begin{equation}\label{ch28uj6tg574tygh}
		\begin{split}
		\mathcal{V}_{{m,M}} \subseteq \Big\{ (u,v)\in[0,1] \times [0,1] \;{\mbox{ s.t. }}\; v< f_{\varepsilon}(u)\Big\}
		\end{split}
		\end{equation}
		where~$f_{\varepsilon} : [0, u_{\mathcal{M}}]\to [0,1]$ is the continuous function given by
		\begin{equation*}
		f_{\varepsilon}(u)=\left\{ 
		\begin{array}{ll}
		\displaystyle		\frac{(u_s^m)^{1-\rho}u^{\rho}}{\rho c } & \text{if} \ u\in [0, u_s^m), \\		
		\displaystyle\frac{u}{\rho c}  & \text{if} \ u\in [u_s^m, u_s^0), \\
		\displaystyle\dfrac{u}{c}+\frac{1-\rho}{1+\rho c}  & \text{if} \ u\in [u_s^0, u_1), \\
		hu +p & \text{if} \ u\in [u_1, 1],
		\end{array}	\right.
		\end{equation*}
		with the convention that the first interval is empty if $m\geq \frac{1}{c}$, the second interval is empty if $m=0$,
		and $h$, $u_1$ and $p$ take the following values:
		\begin{align*}
		&h := \frac{1}{c}\left(1-\dfrac{\varepsilon^2(1-\rho)}{M (1+\rho c)(c+1-\varepsilon)^2 + \varepsilon (\rho c +\rho + \varepsilon-\varepsilon \rho)}\right), \\ &u_1:=\frac{c(\rho c+\rho+\varepsilon-\varepsilon \rho)}{(1+\rho c)(c+1-\varepsilon)}, \\
		&p :=\frac{c+1-hc(\rho c+\rho+\varepsilon-\varepsilon \rho)}{(1+\rho c)(c+1-\varepsilon)}.
		\end{align*}
		
		\item[(ii)]  If~$\rho>1$,  we have
		\begin{equation*}
		\begin{split}
		\mathcal{V}_{{m,M}} \subseteq \Big\{ (u,v)\in[0,1] \times [0,1] \;{\mbox{ s.t. }}\; v< g_{\varepsilon}(u)\Big\}
		\end{split}
		\end{equation*}
		where~$g_{\varepsilon} : [0, u_{\mathcal{M}}] \to [0,1]$ is the continuous function given by
		\begin{equation*}
		g_{\varepsilon}(u)=
		\begin{dcases}
		k\,u & \text{if} \ u\in [0, u_2), \\		
		\displaystyle\dfrac{u}{c} + q  & \text{if} \ u\in [u_2, u_3), \\
		\displaystyle\dfrac{(1-u_3)u^{\rho}}{(u_3)^{\rho}}  & \text{if} \ u\in [u_3, 1]
		\end{dcases}	
		\end{equation*}
		for the following values:
		\begin{align*}&
		k:= \frac{(c+1-\varepsilon)M}{(\rho -1)\varepsilon c + (c+1-\varepsilon) Mc}, \qquad q:=
		\frac{(kc-1)(1-\varepsilon)}{c(k-k\varepsilon+1)}, \\& u_2:=\frac{1-\varepsilon}{k-k\varepsilon+1}\qquad {\mbox{and}}\qquad
		u_3:=\frac{c+1-\varepsilon}{(c+1)(k-k\varepsilon +1)}.
		\end{align*}
	\end{itemize}
\end{proposition} 

We observe that it might be that for some $u\in[0,1]$ we have $f_{\varepsilon}(u)>1$ or $g_{\varepsilon}(u)>1$. In this case, the above proposition would produce the trivial result that $\mathcal{V}_{{m,M}} \cap (\{u\}\times[0,1]) \subseteq \{  u\}\times [0,1]$. On the other hand, a suitable choice of~$\varepsilon$ would lead to nontrivial
consequences entailing, in particular, the proof of
Theorem \ref{ch2thm:limit}.

\begin{proof}[Proof of Proposition~\ref{ch2prop:limit}]
	
	We start by proving the claim in~(i). For this, we will show that
	\begin{equation}\label{ch21642}
	\mathcal{D}:=\Big\{ (u,v)\in[0,1] \times [0,1] \;{\mbox{ s.t. }}\; v\geq  f_{\varepsilon}(u)\Big\} \subseteq \mathcal{V}_{m, M}^C.
	\end{equation}
	where~$\mathcal{V}_{m, M}^C$ is the complement of~$\mathcal{V}_{m, M}$ in the topology of~$[0,1]\times[0,1]$.
	We remark that once~\eqref{ch21642} is established, then the desired claim in~\eqref{ch28uj6tg574tygh}
	plainly follows by taking the complement sets.
	
	To prove~\eqref{ch21642} we first show that
	\begin{equation}\label{ch2rPjhnfvvcc}
	0 \leq u_s^m < u_s^0 < u_1 < 1.\end{equation}
	Notice, as a byproduct, that the above inequalities also give that~$f_{\varepsilon}$ is well defined.
	To prove~\eqref{ch2rPjhnfvvcc} we notice that, by \eqref{ch2usvs}, \eqref{ch2u0v0} and \eqref{ch2usvs2},
	\begin{equation*}
	0\leq 	u_s^m=\max \left\{ 0, \frac{1-mc}{1+\rho c}\,\rho c\right\} < \frac{\rho c}{1+\rho c} =u_s^0
	\end{equation*}
	(and actually the first inequality is strict if $m<\frac{1}{c}$). Next, one can check that, since~$\varepsilon>0$,
	$$  u_s^0-u_1=\frac{\rho c}{1+\rho c}
	-\frac{c(\rho c+\rho+\varepsilon-\varepsilon \rho)}{(1+\rho c)(c+1-\varepsilon)}=-\frac{c\varepsilon}{(1+\rho c)(c+1-\varepsilon)}
	<0.$$
	Furthermore, since~$\varepsilon<1$,
	$$u_1-1=
	\frac{c(\rho c+\rho+\varepsilon-\varepsilon \rho)}{(1+\rho c)(c+1-\varepsilon)}-1=\frac{(\varepsilon-1)(c+1)}{(1+\rho c)(c+1-\varepsilon)}<0.
	$$
	These observations prove~\eqref{ch2rPjhnfvvcc}, as desired.

	Now we point out that
	\begin{equation}\label{ch297896705689045-0}
	{\mbox{$f_{\varepsilon}$ is a continuous function. }}\end{equation}
	Indeed,
	\begin{equation}\label{ch297896705689045-1}
	\frac{(u_s^m)^{1-\rho}}{\rho c} (u_s^m)^\rho = \frac{u_s^m}{\rho c}\qquad{\mbox{ and }}\qquad
	\frac{u_s^0}{\rho c} = \frac{u_s^0}{ c}+ \frac{1-\rho}{1+\rho c}.\end{equation}
	Furthermore, by the definitions of~$p$ and~$u_1$ we see that
	\begin{equation}\label{ch2767thisbc0-i6yjh00}
	\begin{split} p\,&=
	\frac{c+1}{(1+\rho c)(c+1-\varepsilon)}
	-
	\frac{hc(\rho c+\rho+\varepsilon-\varepsilon \rho)}{(1+\rho c)(c+1-\varepsilon)}\\&
	=\frac{c+1}{(1+\rho c)(c+1-\varepsilon)}-hu_1.\end{split}\end{equation}
	Moreover, from the definition of~$u_1$,
	$$ \frac{u_1}{c}+\frac{1-\rho}{1+\rho c}  = \frac{c+1}{(1+\rho c)(c+1-\varepsilon)}.$$
	Combining this and~\eqref{ch2767thisbc0-i6yjh00}, we deduce that
	\begin{equation}\label{ch2indeh8idenf4596}
	\frac{u_1}{c}+\frac{1-\rho}{1+\rho c}  = h u_1+p.
	\end{equation}
	This observation and~\eqref{ch297896705689045-1}
	entail the desired claim in~\eqref{ch297896705689045-0}.

	Next, we show that
	\begin{equation}\label{ch21601}
	f_{\varepsilon}(u)>0 \quad \text{for} \ u>0.
	\end{equation}
	To prove this, we note that for $u\in(0,u_s^m)$ the function is an exponential times the positive constant $\frac{(u_s^m)^{1-\rho}}{\rho c}$, hence is positive. If $u\in[u_s^m, u_s^0)$ then $f_{\varepsilon}(u)$ is a linear function and it is positive since $\rho c >0$. On $[u_s^0, u_1)$, $f_{\varepsilon}(u)$ coincide with a linear function with positive angular coefficient, hence we have
	$$ f_{\varepsilon}(u) \geq \underset{u\in[u_s^0, u_1)}{\min} f_{\varepsilon}(u)= f_{\varepsilon}(u_s^0)= \frac{u_s^0}{\rho c} >0.  $$
	By inspection one can check that $h>0$. Hence, 
	in the interval $[u_1,1]$ we have
	$$ f_{\varepsilon}(u) \geq \underset{u\in[u_1, 1]}{\min} f_{\varepsilon}(u)= f_{\varepsilon}(u_1)\geq \frac{u_s^0}{\rho c} >0.  $$
	This completes the proof of \eqref{ch21601}.
	
	Let us notice that, as a consequence of~\eqref{ch21601}, 
	\begin{equation}\label{ch22344}
	\mathcal{D} \cap\big( (0,1]\times \{0\} \big)= \varnothing.
	\end{equation}
	
	Now we show that 
	\begin{equation}\label{ch22352}
	{ \mbox{for any strategy~$a\in\mathcal{A}_{m, M}$, no trajectory starting in~$\mathcal{D}$ leaves~$\mathcal{D}$.}  }
	\end{equation}
	To this end, we notice that, since~$\partial \mathcal{D} \cap \{v=0\}= \{(0,0)\}$, and the origin is an equilibrium,
	we already have that no trajectory can exit~$\mathcal{D}$ by passing through the points in~$\partial \mathcal{D}\cap \partial ([0,1]\times[0,1])$. Hence, we are left with
	considering the possibility of leaving $\mathcal{D}$ through~$\partial \mathcal{D}\cap ((0,1)\times(0,1))$.
	To exclude this possibility, we compute the velocity of a trajectory in the inward normal direction at~$\partial \mathcal{D}\cap ((0,1)\times(0,1))$.
	
	For every~$u\in[0, u_s^m)$ we have that this normal velocity is
	\begin{equation}\begin{split}\label{ch21614}&
	\dot{v}-	\frac{(u_s^m)^{1-\rho} \rho (u)^{\rho-1} \dot{u}}{\rho c }  \\
	&\qquad =\rho \left( v-  \frac{(u_s^m)^{1-\rho} \, u^{\rho} }{\rho c }  \right)  (1-u-v) -au\left(1- \frac{(u_s^m)^{1-\rho} }{u^{1-\rho}}  \right).
	\end{split}\end{equation}
	Notice that the term~$ 
	v-  \frac{(u_s^m)^{1-\rho} \, u^{\rho} }{\rho c } $
	vanishes on~$\partial \mathcal{D}\cap ((0,1)\times(0,1))$
	when~$u\in[0, u_s^m)$. Also, for all~$u\in[0, u_s^m)$ we have 
	\begin{equation*}
	1- \frac{(u_s^m)^{1-\rho} }{u^{1-\rho}}<0,
	\end{equation*} 
	thus the left hand side in~\eqref{ch21614} is positive. This observation
	rules out the possibility of leaving $\mathcal{D}$ through~$\partial \mathcal{D}\cap ((0,1)\times(0,1))$
	at points where~$u\in[0, u_s^m)$.
	
	It remains to exclude egresses at points of~$\partial \mathcal{D}\cap ((0,1)\times(0,1))$
	with~$u\in[ u_s^m,1)$. We first consider this type of points when~$(u_s^m,u_s^0)$.
	At these points, we have that the
	velocity in the
	inward normal direction on~$\{ v=\frac{u}{\rho c} \}$ is
	\begin{equation*}
	\dot{v}- \frac{\dot{u}}{\rho c}= \left( \rho v - \frac{u}{\rho c} \right)(1-u-v) + au\left( \frac{1}{\rho} -1\right)
	\end{equation*}  
	Expressing~$u$ with respect to~$v$ on~$\partial \mathcal{D}\cap ((0,1)\times(0,1))$
	with~$u\in( u_s^m,u_s^0)$, we have
	\begin{equation}\begin{split}\label{ch2Moiuyted645JN}
	\dot{v}- \frac{\dot{u}}{\rho c}&=v \left( \rho-1  \right)(1-\rho c v-v) + a \rho cv\,\frac{1-\rho}{\rho} \\
	&= v(1-\rho)(\rho c v + v-1 +ac).\end{split}
	\end{equation}  
	We also remark that, for these points,
	\begin{equation*}
	v>v_s^m= \frac{1-mc}{1+\rho c}\ge\frac{1-ac}{1+\rho c}	,\end{equation*} 
	thanks to~\eqref{ch2usvs}. This gives that the quantity in~\eqref{ch2Moiuyted645JN}
	is strictly positive and, as a consequence,
	we have excluded the possibility of
	exiting~$\mathcal{D}$
	at points of~$\partial \mathcal{D}\cap ((0,1)\times(0,1))$
	with~$u\in(u_s^m,u_s^0)$.
	
	It remains to consider the case~$u\in \{u_s^m\}\cup[u_s^0,1)$.
	We first focus on the range~$u\in (u_s^0,u_1)$.
	In this interval, the velocity of a trajectory starting at a point~$(u,v)\in\partial \mathcal{D}\cap ((0,1)\times(0,1))$
	lying on the line~$v=\frac{u}{c}+ \frac{1-\rho}{1+\rho c}$
	in the inward normal direction with respect to~$\partial \mathcal{D}$ is given by
	\begin{equation}\label{ch2tqwfe3857uvcjycer4cubt}
	\dot{v}- \frac{1}{c}\dot{u}= \left(\rho v - \frac{u}{c} \right)(1-u-v).
	\end{equation}
	We also observe that, in light of~\eqref{ch2u0v0},
	$$ u>u_s^0=\frac{\rho c}{1+\rho c}, $$
	and therefore, for any~$u\in( u_s^0,u_1)$ lying on the above line,
	\begin{equation*}
	1-u-v= 1-u-\frac{u}{c} - \frac{1-\rho}{1+\rho c} =(c+1)\left(\frac{\rho}{1+\rho c} -\frac{u}{c} \right) < 0
	\end{equation*}
	and
	\begin{equation*}
	\rho v - \frac{u}{c} = \frac{\rho u}{c}+ \frac{\rho (1-\rho)}{1+\rho c}- \frac{u}{c}  =(1-\rho)\left(  \frac{\rho}{1+\rho c} - \frac{u}{c} \right) < 0 .
	\end{equation*}
	Using these pieces of information in~\eqref{ch2tqwfe3857uvcjycer4cubt}, we conclude that the
	inward normal velocity of a trajectory starting at a point~$(u,v)\in\partial \mathcal{D}\cap ((0,1)\times(0,1))$
	with~$u\in (u_s^0,u_1)$ is strictly positive. This gives that no trajectory can exit~$\mathcal{D}$
	at this type of points, and we need to exclude the case~$u\in \{u_s^m, u_s^0\}\cup[u_1,1)$.
	
	We consider now the interval~$[u_1,1)$. 
	In this interval, the component of the velocity of a trajectory at a point on the straight
	line given by~$v=hu+p$ in the orthogonal inward pointing direction is 
	\begin{equation}\label{ch28gqwfOJHNsmeoout43906}\begin{split}&
	(\dot{u}, \dot{v}) \cdot \frac{(-h, 1)}{\sqrt{1+h^2}} = \frac{   (\rho v -h u)(1-u-v) -au(1-hc) }{\sqrt{1+h^2}}\\
	&\qquad =
	\frac{((1-\rho)hu-\rho p )(u+v-1) -au(1-hc)}{\sqrt{1+h^2}}
	\end{split}\end{equation}
	We observe that, if~$u\in[u_1,1)$,
	\begin{equation}\label{ch2jd723u9007432yhgvythgkliew}\begin{split}&(1-\rho)hu-\rho p \ge
	(1-\rho)hu_1-\rho p=hu_1-\rho (hu_1+p) \\ 
	&\qquad	= h u_1-\rho \left(\frac{u_1}{c}+\frac{1-\rho}{1+\rho c} \right) =
	h u_1-\rho \left(
	\frac{\rho c+\rho+\varepsilon-\varepsilon \rho}{(1+\rho c)(c+1-\varepsilon)}
	+\frac{1-\rho}{1+\rho c} \right) \\&\qquad=
	h u_1-\frac{\rho (c+1)}{(1+\rho c)(c+1-\varepsilon)},
	\end{split}\end{equation}
	thanks to~\eqref{ch2indeh8idenf4596}.
	
	We also remark that
	\begin{equation*}\begin{split}
	hu_1\,&=
	\left(1-\dfrac{\varepsilon^2(1-\rho)}{M (1+\rho c)(c+1-\varepsilon)^2 + \varepsilon (\rho c +\rho + \varepsilon-\varepsilon \rho)}\right)\,\frac{ \rho c+\rho+\varepsilon-\varepsilon \rho}{(1+\rho c)(c+1-\varepsilon)}, \\
	&=
	\frac{ \rho c+\rho+\varepsilon-\varepsilon \rho}{(1+\rho c)(c+1-\varepsilon)}\\&\qquad\qquad
	-
	\dfrac{\varepsilon^2(1-\rho)\big(\rho c+\rho+\varepsilon-\varepsilon \rho\big)}{\big(
		M (1+\rho c)(c+1-\varepsilon)^2 + \varepsilon (\rho c +\rho + \varepsilon-\varepsilon \rho)\big)(1+\rho c)(c+1-\varepsilon)}
	.
	\end{split}
	\end{equation*}	
	Accordingly,
	\begin{eqnarray*}&&
		h u_1-\frac{\rho (c+1)}{(1+\rho c)(c+1-\varepsilon)}=
		\frac{ \varepsilon(1- \rho)}{(1+\rho c)(c+1-\varepsilon)}\\&&\qquad\qquad
		-
		\dfrac{\varepsilon^2(1-\rho)\big(\rho c+\rho+\varepsilon-\varepsilon \rho\big)}{\big(
			M (1+\rho c)(c+1-\varepsilon)^2 + \varepsilon (\rho c +\rho + \varepsilon-\varepsilon \rho)\big)(1+\rho c)(c+1-\varepsilon)}\\
		&&\qquad=
		\frac{ \varepsilon(1- \rho)}{(1+\rho c)(c+1-\varepsilon)}\Bigg(1
		-
		\dfrac{\varepsilon \big(\rho c+\rho+\varepsilon-\varepsilon \rho\big)}{
			M (1+\rho c)(c+1-\varepsilon)^2 + \varepsilon (\rho c +\rho + \varepsilon-\varepsilon \rho)}\Bigg)\\&&\qquad=
		\frac{ \varepsilon(1- \rho)}{(1+\rho c)(c+1-\varepsilon)}\cdot
		\dfrac{M (1+\rho c)(c+1-\varepsilon)^2}{
			M (1+\rho c)(c+1-\varepsilon)^2 + \varepsilon (\rho c +\rho + \varepsilon-\varepsilon \rho)}\\&&\qquad=
		\dfrac{\varepsilon M(1- \rho)(c+1-\varepsilon)}{
			M (1+\rho c)(c+1-\varepsilon)^2 + \varepsilon (\rho c +\rho + \varepsilon-\varepsilon \rho)}
		.\end{eqnarray*}
	{F}rom this and~\eqref{ch2jd723u9007432yhgvythgkliew}, we gather that
	\begin{equation}\label{ch2ILpredmnow55}
	(1-\rho)hu-\rho p\ge
	\dfrac{\varepsilon M(1- \rho)(c+1-\varepsilon)}{
		M (1+\rho c)(c+1-\varepsilon)^2 + \varepsilon (\rho c +\rho + \varepsilon-\varepsilon \rho)}
	.\end{equation}
	
	Furthermore, we point out that, when~$[u_1, 1)$ and~$v=hu+p$,
	\begin{equation*}\begin{split}&
	u+v-1\ge u_1+hu_1+p-1 = 
	u_1+\frac{u_1}c+\frac{1-\rho}{1+\rho c}-1\\&\qquad
	=\frac{(c+1)(\rho c+\rho+\varepsilon-\varepsilon \rho)}{(1+\rho c)(c+1-\varepsilon)}
	-\frac{\rho(c+1)}{1+\rho c}
	=\frac{\varepsilon(c+1)}{(1+\rho c)(c+1-\varepsilon)} >\frac{\varepsilon}{c+1-\varepsilon},
	\end{split}\end{equation*}
	thanks to~\eqref{ch2indeh8idenf4596}.
	
	Combining this inequality and~\eqref{ch2ILpredmnow55}, we deduce that
	\begin{equation*}
	((1-\rho)hu-\rho p )(u+v-1) >
	\dfrac{\varepsilon^2 M(1- \rho)}{
		M (1+\rho c)(c+1-\varepsilon)^2 + \varepsilon (\rho c +\rho + \varepsilon-\varepsilon \rho)}.
	\end{equation*}	
	Therefore, noticing that~$h<\frac{1}{c}$,
	\begin{eqnarray*}&&
		((1-\rho)hu-\rho p )(u+v-1) -au(1-hc)\\&>&
		\dfrac{\varepsilon^2 M(1- \rho)}{
			M (1+\rho c)(c+1-\varepsilon)^2 + \varepsilon (\rho c +\rho + \varepsilon-\varepsilon \rho)}-Mu(1-hc)\\
		&=&
		\dfrac{\varepsilon^2 M(1- \rho)(1-u)}{
			M (1+\rho c)(c+1-\varepsilon)^2 + \varepsilon (\rho c +\rho + \varepsilon-\varepsilon \rho)},
	\end{eqnarray*}
	which is strictly positive.
	
	Using this information in~\eqref{ch28gqwfOJHNsmeoout43906},
	we can thereby exclude the possibility of leaving~${\mathcal{D}}$ through~$
	\partial \mathcal{D}\cap ((0,1)\times(0,1))$
	with~$u\in [u_1,1)$.
	As a result, it only remains to exclude the possibility
	of an egress from~${\mathcal{D}}$ through~$
	\partial \mathcal{D}\cap ((0,1)\times(0,1))$
	with~$u\in \{u^m_s,u^0_s\} $.
	
	For this, we perform a general argument of dynamics, as follows. We denote by~$P^m_s$
	and~$P^0_s$ the points on~$
	\partial \mathcal{D}\cap ((0,1)\times(0,1))$
	with~$u=u^m_s$ and~$u=u^0_s $, respectively (these points may also coincide,
	as it happens when~$m=0$). We stress that we already know by the previous arguments that
	\begin{equation}\label{ch26tGSHuj2fw7545}
	{\mbox{if a trajectory leaves~${\mathcal{D}}$ it must pass through~$\{P^m_s,P^0_s\}$.}}
	\end{equation}
	Our goal is to show that no trajectory leaves~${\mathcal{D}}$ and for this we argue by contradiction,
	supposing that there exist~$\bar P\in {\mathcal{D}}$ and~$T>0$ such that~$\phi^T(\bar P)$ lies in the complement
	of~${\mathcal{D}}$ in~$[0,1]\times[0,1]$. Here, we have denoted by~$\phi^T$ the flow associated to~\eqref{ch2model}.
	We let~$\bar Q:=\phi^T(\bar P)$ and, since the complement of~${\mathcal{D}}$ is open in~$[0,1]\times[0,1]$,
	we can find~$\rho>0$ such that~$B_\rho(\bar Q)\cap([0,1]\times[0,1])$ is contained in
	the complement of~${\mathcal{D}}$.
	
	Also, from~\eqref{ch26tGSHuj2fw7545}, there exists~$\bar t\in[0,T)$ such that~$\phi^{\bar t}(\bar P)\in\{P^m_s,P^0_s\}$.
	We suppose that~$\phi^{\bar t}(\bar P)=P^m_s$ (the case~$
	\phi^{\bar t}(\bar P)=P^0_s$ being completely analogous).
	We let~$\bar T:=T-\bar t$ and we notice that~$\phi^{\bar T}(P^m_s)= \phi^T(\bar P)=\bar Q$. Hence,
	by continuity with respect to the data, we can find~$r>0$ such that
	$$ \phi^{\bar T}\big( B_r(P^m_s)\cap ([0,1]\times[0,1]) \big)\subseteq B_\rho(\bar Q)\cap([0,1]\times[0,1]).$$
	We define~${\mathcal{U}}:=B_r(P^m_s)\cap{\mathcal{D}}$. We observe that
	\begin{equation}\label{ch21fe45-90jhwg3rg rewt57}
	{\mbox{${\mathcal{U}}$ has strictly positive Lebesgue measure,}} 
	\end{equation}
	since~$P^m_s\in\partial{\mathcal{D}}$ and~${\mathcal{D}}$ has boundary of H\"older class.
	In addition,
	$$ \phi^{\bar T}\big( {\mathcal{U}} \big)\subseteq B_\rho(\bar Q)\cap([0,1]\times[0,1])
	\subseteq\big( [0,1]\times[0,1]\big)
	\setminus{\mathcal{D}}.$$
	This and~\eqref{ch26tGSHuj2fw7545} give that for every~$P\in{\mathcal{U}}$ there exists~$t_P\in[0,\bar T]$ such that~$
	\phi^{t_P}(P)\in\{P^m_s,P^0_s\}$. In particular,
	$$ P\in \phi^{-t_P} \{P^m_s,P^0_s\}\subseteq \big\{ \phi^t (P^m_s), \;t\in [-\bar T,0]\big\}
	\cup\big\{ \phi^t (P^0_s), \;t\in [-\bar T,0]\big\}.$$
	Since this is valid for every~$P\in{\mathcal{U}}$, we conclude that
	\begin{equation}\label{ch29iKCIMWSsiofnooer}
	{\mathcal{U}}\subseteq\big\{ \phi^t (P^m_s), \;t\in [-\bar T,0]\big\}
	\cup\big\{ \phi^t (P^0_s), \;t\in [-\bar T,0]\big\}.
	\end{equation}
	Now we remark that~$\big\{ \phi^t (P^m_s), \;t\in [-\bar T,0]\big\}$ is an arc of a smooth curve,
	whence it has null Lebesgue measure, and a similar statement holds true for~$
	\big\{ \phi^t (P^0_s), \;t\in [-\bar T,0]\big\}$.
	Consequently, we deduce from~\eqref{ch29iKCIMWSsiofnooer} that~${\mathcal{U}}$
	has null Lebesgue measure, in contradiction with~\eqref{ch21fe45-90jhwg3rg rewt57}.
	
	In this way, we have shown that no trajectory can leave~${\mathcal{D}}$ and the proof
	of~\eqref{ch22352}
	is complete.
	
	By~\eqref{ch22344} and~\eqref{ch22352}, no trajectory starting in~$\mathcal{D}$ can arrive in~$(0,1]\times[0,1]$ when the bound~$m\leq a(t)\leq M$ holds, hence \eqref{ch21642} is true. 
	Therefore the statement~(i) in Proposition \ref{ch2prop:limit} is true.
	
	\medskip
	
	Now we establish the claim in~(ii). To this end,
	we point out that claim (ii) is equivalent to
	\begin{equation}\label{ch22359}
	\mathcal{G}:=\Big\{ (u,v)\in[0,1] \times [0,1] \;{\mbox{ s.t. }}\; v\geq  g_{\varepsilon}(u)\Big\} \subseteq \mathcal{V}_{m, M}^C,
	\end{equation}
	for all~$\varepsilon\in\left(0, 1\right)$,
	where~$\mathcal{V}_{m, M}^C$ is the complement of~$\mathcal{V}_{m, M}$ in the topology of~$[0,1]\times[0,1]$.
	
	First, we point out that
	\begin{equation}\label{ch2Ov54io0v9ik4gfvh}
	{\mbox{$g_{\varepsilon}$ is a well defined continuous function. }}\end{equation}
	Indeed, one can easily check for $\varepsilon\in(0,1)$ that 
	\begin{equation}\label{ch21409}\begin{split}&
	0 < u_2
	=\frac{1-\varepsilon}{k-k\varepsilon+1}-\frac{c+1-\varepsilon}{(c+1)(k-k\varepsilon +1)}+u_3
	=-\frac{c\varepsilon }{(c+1)(k-k\varepsilon +1)}+u_3\\&\qquad\qquad\qquad
	<u_3
	<\frac{c+1}{(c+1)(k-k\varepsilon +1)}<1.\end{split}
	\end{equation}
	Then, one checks that 
	\begin{align*}
	ku_2=\frac{u_2}{c}+q,
	\end{align*}
	hence $g_{\varepsilon}$ is continuous at the point $u_2$. In addition, one can check that $g_{\varepsilon}$ is continuous
	at the point~$u_3$ by observing that
	\begin{equation}\label{ch2isceocessvcpoo}\begin{split}&
	\frac{u_3}c+q-(1-u_3)=\frac{(c+1)u_3}{c}+q-1\\&\qquad=
	\frac{c+1-\varepsilon}{c(k-k\varepsilon +1)}+\frac{(kc-1)(1-\varepsilon)}{c(k-k\varepsilon+1)}-1\\&\qquad=
	\frac{c+1-\varepsilon+(kc-1)(1-\varepsilon)-c(k-k\varepsilon+1)}{c(k-k\varepsilon+1)}=0.
	\end{split}\end{equation}
	This completes the proof of \eqref{ch2Ov54io0v9ik4gfvh}.
	
	Now we show that
	\begin{equation}\label{ch21411}
	g_{\varepsilon}(u)>0 \quad \text{for every} \ u\in(0,1].
	\end{equation}
	We have that~$k>0$ for every~$\varepsilon<1$, and therefore~$g_{\varepsilon}(u)>0$ for all~$u\in(0, u_2)$.
	Also, since $g_{\varepsilon}(u_2)=ku_2>0$ and $g_{\varepsilon}$ is linear in $(u_2, u_3)$, we have that $g_{\varepsilon}(u)>0$ for all~$u\in(u_2, u_3)$. Moreover, in the interval~$\in[u_3,1]$ we have that~$g_{\varepsilon}$
	is an exponential function multiplied by a  positive constant, thanks to~\eqref{ch21409}, hence it is positive. These considerations prove~\eqref{ch21411}.
	
	As a consequence of~\eqref{ch21411}, we have that 
	\begin{equation}\label{ch22360}
	\mathcal{G} \cap \big((0,1]\times \{0\}\big) = \varnothing.
	\end{equation}
	Now we claim that 
	\begin{equation}\label{ch20002}
	{ \mbox{for any strategy~$a\in\mathcal{A}_{m,M}$, no trajectory starting in~$\mathcal{G}$ leaves~$\mathcal{G}$.}  }
	\end{equation}
	For this, 
	we observe that, in light of~\eqref{ch22360}, all the points on
	$$\partial \mathcal{G} \setminus \{ (u,g_{\varepsilon}(u)) \
	{\mbox{ with }}\ u\in[0,1]  \}$$ belong to~$\partial ([0,1]\times [0,1]) \setminus \{v=0 \}$,
	and these three sides of the square do not allow the flow to exit.
	Hence, to prove~\eqref{ch20002}
	it suffices to check that the trajectories starting on~$\partial \mathcal{G}\cap\big( (0,1)\times(0,1)\big)$ 
	enter~${\mathcal{G}}$. We do this
	by showing that the inner pointing derivative of the trajectory is nonnegative, according to the computation below.
	
	At a point on the line~$v=k u$, the velocity of a trajectory in the direction that is orthogonal to~$\partial \mathcal{G}$ for~$u\in[0,u_2)$ and pointing inward is:
	\begin{equation}\label{ch21741}
	(\dot{u}, \dot{v})\cdot \frac{(-k, 1)}{\sqrt{1+k^2}} =\frac{(\rho v- ku)(1-u-v)-au(1-kc) }{\sqrt{1+k^2}}    .   
	\end{equation}
	We also note that
	\begin{equation}\label{ch22018}
	kc
	= \frac{(c+1-\varepsilon)M}{(\rho -1)\varepsilon + (c+1-\varepsilon)M }
	<1,\end{equation} and therefore, at a point on~$v=k u$ with $u\in[0, u_2)$,
	\begin{equation*}\begin{split}&
	1-u-v \geq  1-u_2-k u_2= 
	1-	\frac{(1+k)(1-\varepsilon)}{k-k\varepsilon+1}
	=	\frac{\varepsilon}{k(1-\varepsilon)+1}\\&\qquad\qquad\qquad
	=\frac{\varepsilon c}{kc(1-\varepsilon)+c}
	> \frac{\varepsilon c}{1+c-\varepsilon}.\end{split}
	\end{equation*}
	This inequality entails that
	\begin{equation*}
	k= \frac{(1+c-\varepsilon)M}{(\rho-1)\varepsilon c+(1+c-\varepsilon)Mc } 
	=\frac{M}{\frac{(\rho-1)\varepsilon c}{1+c-\varepsilon}+Mc } 
	>  \frac{M}{(\rho-1)(1-u-v)+Mc}.
	\end{equation*}
	Consequently,
	\begin{equation*}
	(\rho-1)(1-u-v)k > M (1-kc).
	\end{equation*} 
	{F}rom this and~\eqref{ch21741}, one deduces that, for all~$u\in(0, u_2)$, $a\leq M$, and $v=k u$,
	\begin{equation*}\begin{split}&
	(\dot{u}, \dot{v})\cdot \frac{(-k, 1)}{\sqrt{1+k^2}} =\frac{ku(\rho - 1)(1-u-v)-au(1-kc) }{\sqrt{1+k^2}} \\&\qquad\qquad  >
	\frac{Mu (1-kc)-au(1-kc) }{\sqrt{1+k^2}}\ge0.\end{split}
	\end{equation*}
	This (and the fact that the origin is an equilibrium)
	rules out the possibility of exiting~${\mathcal{G}}$ from~$\{ u\in[0, u_2){\mbox{ and }} v=k u\}$.
	
	It remains to consider the portions of~$\partial\mathcal{G}\cap((0,1)\times(0,1))$ given by
	\begin{equation}\label{ch29u:9idkj:0oekdjfjfj81763yhrf}
	\left\{ u\in[ u_2,u_3){\mbox{ and }} v=\frac{ u}c+q\right\}\end{equation}
	and by
	\begin{equation}\label{ch29u:9idkj:0oekdjfjfj81763yhrf2}\left\{ u\in[u_3,1]{\mbox{ and }} v=\frac{(1-u_3)u^\rho}{(u_3)^\rho}\right\}.\end{equation}
	
	Let us deal with the case in~\eqref{ch29u:9idkj:0oekdjfjfj81763yhrf}.
	In this case,
	the velocity of a trajectory in the direction orthogonal to~$\partial \mathcal{G}$ for~$u\in[u_2,u_3)$ and pointing inward is
	\begin{equation}\label{ch22027}
	(\dot{u}, \dot{v})\cdot \frac{(-1, c)}{\sqrt{1+c^2}}=\frac{(\rho c v -u)(1-u-v)}{\sqrt{1+c^2}}.
	\end{equation}
	Recalling~\eqref{ch2RANGEEP}, we also observe that
	\begin{equation}\label{ch21536}\begin{split}&k-
	\frac{1}{\rho c}
	=\frac1c\left(
	\frac{(c+1-\varepsilon)M}{(\rho -1)\varepsilon + (c+1-\varepsilon) M}-\frac1\rho\right)\\&\qquad=
	\frac{(\rho-1)\big((c+1-\varepsilon)M
		-\varepsilon\big)}{\rho c\big( (\rho -1)\varepsilon + (c+1-\varepsilon) M\big)}
	>0.\end{split}
	\end{equation}
	Thus, on the line given by~$v=\frac{u}{c}+q$ we have that
	\begin{equation}\label{ch2Cnodizeos80p4}\begin{split}&
	\rho c v -u= (\rho-1)u +\rho c q
	\ge (\rho-1)u_2 +\rho c q\\&\qquad=
	\frac{(\rho-1)(1-\varepsilon)}{k-k\varepsilon+1}
	+\frac{\rho(kc-1)(1-\varepsilon)}{k-k\varepsilon+1}\\&\qquad
	= (1-\varepsilon)\frac{(\rho-1)+\rho(kc-1)}{k-k\varepsilon+1}=\frac{(1-\varepsilon)(\rho k c -1)}{k-k\varepsilon+1}>0,\end{split}
	\end{equation}
	where~\eqref{ch21536} has been used in the latter inequality.
	
	In addition, recalling~\eqref{ch2isceocessvcpoo},
	\begin{equation*}
	1-u-v > 1-u_3 - \frac{u_3}{c} -q = 1-u_3-1+u_3=0.
	\end{equation*}
	{F}rom this and~\eqref{ch2Cnodizeos80p4}, we gather that the velocity calculated in~\eqref{ch22027} is positive in~$[u_2, u_3)$ and this
	excludes the possibility of exiting~$\mathcal{G}$ from the boundary given in~\eqref{ch29u:9idkj:0oekdjfjfj81763yhrf}.

	Next, we focus
	on the portion of the boundary described in~\eqref{ch29u:9idkj:0oekdjfjfj81763yhrf2}
	by considering~$u\in[u_3, 1]$.
	That is, we now compute the component of the velocity at a point on~$\partial \mathcal{G}$ for ~$u\in[u_3, 1]$ in the direction that is orthogonal to~$\partial \mathcal{G}$ and pointing inward, that is
	\begin{equation}\label{ch21803}
	\begin{split}&
	(\dot{u}, \dot{v})\cdot \frac{(-\rho \frac{1-u_3}{(u_3)^{\rho}}u^{\rho-1}, 1)}{\sqrt{1+\rho^2\frac{(1-u_3)^2}{(u_3)^{2\rho}}u^{2\rho-2}}} \\=\,& \frac{\rho(1-u-v)\left(v-   \frac{1-u_3}{(u_3)^{\rho}} u^{\rho} \right) - au\left( 1-\rho c \frac{1-u_3}{(u_3)^{\rho}} u^{\rho-1}   \right) }{\sqrt{1+\rho^2\frac{(1-u_3)^2}{(u_3)^{2\rho}}u^{2\rho-2}}}
	\\=\,& \frac{ au\left( \rho c \frac{1-u_3}{(u_3)^{\rho}} u^{\rho-1} -1   \right) }{\sqrt{1+\rho^2\frac{(1-u_3)^2}{(u_3)^{2\rho}}u^{2\rho-2}}}\\ \ge\,&
	\frac{ au\left( \rho c \frac{1-u_3}{u_3} -1   \right) }{\sqrt{1+\rho^2\frac{(1-u_3)^2}{(u_3)^{2\rho}}u^{2\rho-2}}}.
	\end{split}
	\end{equation}
	Now we notice that
	\begin{eqnarray*}
		&&	\rho c (1-u_3) = \rho c \left(\frac{u_3}{c}+q \right) =
		\rho u_3+ \rho c q=\rho u_3+\frac{\rho (kc-1)(1-\varepsilon)(c+1) u_3}{c+1-\varepsilon},
	\end{eqnarray*} 
	thanks to~\eqref{ch2isceocessvcpoo}.
	
	As a result, using~\eqref{ch21536}, 
	\begin{eqnarray*}
		&&	\rho c (1-u_3) >\rho u_3+\frac{ (1-\rho)(1-\varepsilon)(c+1) u_3}{c+1-\varepsilon}
		\\&&\qquad
		=  \frac{ u_3}{c+1-\varepsilon}
		\Big(\rho (c+1-\varepsilon)+(1-\rho)(1-\varepsilon)(c+1) \Big)\\&&\qquad=
		\frac{ u_3\big( (1-\varepsilon)(c+1)+\varepsilon \rho c\big)}{c+1-\varepsilon}	\\&&
		\qquad= u_3+
		\frac{ \varepsilon c u_3( \rho-1)}{c+1-\varepsilon}>u_3.
	\end{eqnarray*}
	This gives that the quantity in  \eqref{ch21803} is positive.
	Hence, we have ruled
	out also the possibility of exiting~$\mathcal{G}$ from the boundary given in~\eqref{ch29u:9idkj:0oekdjfjfj81763yhrf2},
	and this ends the proof of~\eqref{ch20002}.
	
	Since no trajectory can exit~$\mathcal{G}$ for any~$a$ with~$m\leq a \leq M$, we get that no point~$(u,v)\in \mathcal{G}$ is mapped into~$(0,1]\times\{0\}$ because of~\eqref{ch22360}, thus~\eqref{ch22359} is true and the proof is complete.
\end{proof}

We end this section with the proof of Theorem \ref{ch2thm:limit}.

\begin{proof}[Proof of Theorem \ref{ch2thm:limit}]
	Since by definition $\mathcal{A}_{m,M}\subseteq \mathcal{A}$, we have that $\mathcal{V}_{{m,M}}\subseteq \mathcal{V}_{\mathcal{A}}$. Hence, we are left with proving that the latter inclusion is strict.
	
	We start with the case $\rho<1$. We choose
	\begin{equation}\label{ch21934567890-dfghjk4567890-fd11}
	\varepsilon\in\left(0,\,\min\left\{ \frac{ \rho c(c+1)}{1+\rho c}, \frac{M(c+1)}{M+1},1  \right\} \right). \end{equation}
	We observe that this choice is compatible with
	the assumption on~$\varepsilon$ in~\eqref{ch2RANGEEP}. We note that
	\begin{equation}\label{ch21911}
	u_1 < \min\left\{ \frac{ \rho c(c+1)}{1+\rho c}, 1  \right\},
	\end{equation}
	thanks to~\eqref{ch21934567890-dfghjk4567890-fd11}.
	Moreover, by \eqref{ch2indeh8idenf4596} and the fact that~$h<\frac1c$, it holds that
	\begin{equation}\label{ch21933}
	h u + p 
	=h (u-u_1)+hu_1 + p
	=h (u-u_1)+\frac{u_1}{c}+\frac{1-\rho}{1+\rho c}<\frac{u}{c}+ \frac{1-\rho }{1+\rho c} 
	\end{equation}
	for all~$u>u_1$.

	Now we choose $$\bar{u}\in \left( u_1,  \min\left\{ \frac{ \rho c(c+1)}{1+\rho c}, 1  \right\} \right),$$ which is possible thanks to \eqref{ch21911}, and 
	\begin{equation}\label{ch21925}
	\bar{v}: = \frac{1}{2}\left(  h \bar{u} + p   \right) + \frac{1}{2}\left( \frac{\bar{u} }{c}+ \frac{1-\rho }{1+\rho c}   \right).
	\end{equation}
	By \eqref{ch21933} we get that
	\begin{equation}\label{ch21937}
	h \bar{u} + p <\frac12\left(h \bar{u} + p\right)
	+ \frac{1}{2}\left( \frac{\bar{u} }{c}+ \frac{1-\rho }{1+\rho c}   \right)=
	\bar{v} < \frac{\bar{u}}{c}+ \frac{1-\rho }{1+\rho c}.
	\end{equation}
	Using Proposition \ref{ch2prop:limit} and \eqref{ch21937}, we deduce that  $(\bar{u}, \bar{v})\not\in \mathcal{V}_{{m,M}}$. By Theorem \ref{ch2thm:Vbound} and \eqref{ch21937} we obtain instead that $(\bar{u}, \bar{v})\in \mathcal{V}_{\mathcal{A}}$. Hence, the set  
	$\mathcal{V}_{{m,M}}$ is strictly included in~$\mathcal{V}_{\mathcal{A}}$
	when~$\rho<1$.

	Now we consider the case~$\rho>1$, using again the notation of Proposition \ref{ch2prop:limit}. We recall that~$u_2>0$
	and~$ u_{\infty}>0$, due to~\eqref{ch2ZETADEF} and~\eqref{ch21409}, hence we can choose
	$$ \bar{u} \in \left( 0, \min\{u_2, u_{\infty}\}   \right) .$$
	We also define 
	\begin{equation*}
	\bar{v} := \frac12\left( \frac{1}{c} +k   \right) \bar{u}. 
	\end{equation*}
	By \eqref{ch22018}, we get that 
	\begin{equation}\label{ch22031}
	k \bar{u} < \frac{k\bar{u}}2+\frac{\bar{u}}{2c}
	=	\bar{v} < \frac{ \bar{u}}{c}.
	\end{equation}
	Exploiting this and the characterization in Proposition \ref{ch2prop:limit}, it holds that  $(\bar{u}, \bar{v})\not\in \mathcal{V}_{{m,M}}$. On the other hand, by Theorem \eqref{ch2thm:Vbound} and \eqref{ch22031} we have instead that $(\bar{u}, \bar{v})\in \mathcal{V}_{\mathcal{A}}$. As a consequence, the set~$\mathcal{V}_{{m,M}}$
	is strictly contained in~$ \mathcal{V}_{\mathcal{A}}$ for $\rho>1$.
	This concludes the proof of Theorem~\ref{ch2thm:limit}.
\end{proof}

\subsection{Minimization of war duration: proof of Theorem~\ref{ch2thm:min}}

We now deal with the strategies leading to the quickest possible victory of the first population.

\begin{proof}[Proof of Theorem~\ref{ch2thm:min}]
	Our aim is to establish the existence 
	of the strategy leading to the quickest possible victory
	and to determine its range.
	For this, we consider
	the following minimization problem under constraints for~$x(t):=(u(t), v(t))$: 
	\begin{equation}\label{ch2sys:min}
	\left\{ 
	\begin{array}{ll}
	\dot{x}(t)=f(x(t), a(t) ), \\ x(0)=(u_0, v_0), \\ x(T_s)\in (0,1]\times\{0\}, \\
	\displaystyle\min_{a(t)\in [m, M]} \displaystyle\int_{0}^{T_s} 1 \,dt, 
	\end{array}
	\right.
	\end{equation}
	where 
	\begin{equation*}
	f(x, a) := \Big( u(1-u-v-ac), \ \rho v(1-u-v) -au    \Big).
	\end{equation*}
	Here~$T_s$ corresponds to the exit time introduced in~\eqref{ch2def:T_s}, in dependence of the strategy~$a(\cdot)$.
	
	Theorem 6.15 in~\cite{trelat2005controle} assures the existence of a minimizing solution~$(\tilde{a}, \tilde{x})$ with~$\tilde{a}(t)\in[m, M]$ for all~$t\in[0,T]$, and~$\tilde{x}(t)\in[0,1]\times[0,1]$ absolutely continuous, such that~$\tilde{x}(T)=(\tilde u(T), 0)$ with~$\tilde u(T)\in [0,1]$,
	where~$T$ is the exit time for~$\tilde{a}$.
	
	We now prove that
	\begin{equation}\label{ch290o-045}
	\tilde{u}(T)>0.\end{equation}
	Indeed, if this were false, then~$(\tilde{u}(T), \tilde{v}(T))=(0,0)$. Let us call~$d(t):
	= \tilde{u}^2(t)+ \tilde{v}^2(t)$. Then, we observe that
	the function~$d(t)$ satisfies the following differential inequality:
	\begin{equation}\label{ch21955}
	- \dot{d}(t) \le C d , \qquad \text{for}  \quad C:=4+4\rho+2Mc+M.
	\end{equation}
	To check this, we compute that
	\begin{align*}
	- \dot{d}  &= 2\left(  -\tilde{u}^2(1-\tilde{u}-\tilde{v}-\tilde ac) - \tilde{v}^2 \rho(1-\tilde{u}-\tilde{v}) + \tilde{u}\tilde{v}\tilde a     \right) \\
	& \le2\tilde{u}^2(2+Mc) +  4\rho\tilde{v}^2 + (\tilde{u}^2+\tilde{v}^2)M \\
	& \le C (\tilde{u}^2+\tilde{v}^2)\\&= C d,
	\end{align*}
	which proves~\eqref{ch21955}.
	
	{F}rom~\eqref{ch21955}, one has that
	\begin{equation*}
	0<(u_0^2+v_0^2 ) e^{-CT} \le d(T)=\tilde{u}^2(T)+ \tilde{v}^2(T)=\tilde{u}^2(T),
	\end{equation*}
	and this leads to~\eqref{ch290o-045}, as desired.	
	We remark that, in this way, we have found a trajectory~$\tilde{a}$ which
	leads to the victory of the first population in the shortest possible time.
	
	Theorem 6.15 in~\cite{trelat2005controle} assures that~$\tilde{a}(t)\in L^{1}[0,T]$, so~$\tilde{a}(t)$ is measurable.
	We have that the two vectorial functions~$F$ and~$G$, defined by
	\begin{equation*}
	F(u,v):= \left( 
	\begin{array}{c}
	u(1-u-v)\\
	\rho v (1-u-v)
	\end{array}
	\right)\qquad{\mbox{and}}\qquad G(u,v):= \left( 
	\begin{array}{c}
	-cu\\
	-u
	\end{array}
	\right),
	\end{equation*}
	and satisfying~$f(x(t), a(t))= F(x(t))+a(t)G(x(t))$, 
	are analytic. Moreover the set $\overline{\mathcal{V}}_{\mathcal{A}_{m,M}}$ is a subset of $\R^2$, therefore it can be seen as an analytic manifold with border which is also a compact set. For all $x_0\in{\mathcal{V}}_{\mathcal{A}_{m,M}}$ and $t>0$ we have that the trajectory starting from $x_0$ satisfies $x(\tau)\in\overline{\mathcal{V}}_{\mathcal{A}_{m,M}}$ for all $\tau\in[0,t]$. Then, by Theorem 3.1 in \cite{sussmann1987C}, there exists a couple $(\tilde{a}, \tilde{x})$ analytic a part from a finite number of points, such that $(\tilde{a}, \tilde{x})$ solves \eqref{ch2sys:min}.

	\medskip
	
	Now, to study the range of~$\tilde{a}$, we apply the Pontryagin Maximum Principle (see for example~\cite{trelat2005controle}). The Hamiltonian associated with system~\eqref{ch2sys:min} is
	\begin{equation*}
	H(x,p, p_0, a ): =  p\cdot f(x,a)  + p_0
	\end{equation*}
	where~$p=(p_u, p_v)$ is the adjoint to~$x=(u,v)$ and~$p_0$ is the adjoint to the cost function identically
	equal to~$1$.
	The Pontryagin Maximum Principle tells us that, since~$\tilde{a}(t)$ and~$\tilde{x}(t)=(\tilde{u}(t), \tilde{v}(t))$ give the optimal solution, there exist a vectorial function~$\tilde p : [0, T] \to \R^2$ and a scalar~$\tilde p_0\in(-\infty, 0]$ such that 
	\begin{equation}\label{ch2HJA}
	\left\{
	\begin{array}{ll}
	\dfrac{d\tilde{x}}{dt} (t)= \dfrac{\partial H}{\partial p} (\tilde{x}(t), \tilde p(t), \tilde p_0, \tilde{a}(t) ), & \text{for a.a.} \ t\in[0, T], \\
	\\
	\dfrac{d	\tilde{p}}{dt} (t)=- \dfrac{\partial H}{\partial x} (\tilde{x}(t), \tilde p(t), \tilde p_0, \tilde{a}(t) ), & \text{for a.a.} \ t\in[0, T],
	\end{array}
	\right.
	\end{equation} 
	and 
	\begin{equation}\label{ch22349}
	H(\tilde{x}(t), \tilde p(t), \tilde p_0, \tilde{a}(t) ) = \underset{a(\cdot)\in[m,M]}{\max} H(\tilde{x}(t), \tilde p(t), \tilde p_0, a ) \quad \text{for a.a.} \ t\in[0, T].
	\end{equation}
	Moreover, since the final time is free, we have
	\begin{equation}\label{ch21244}
	H(\tilde{x}(T), \tilde p(T),\tilde p_0, \tilde{a}(T) ) =0. 
	\end{equation}
	Also, since~$H(x,p,p_0,a)$ does not depend on~$t$, we get
	\begin{equation}\label{ch22343}
	H(\tilde{x}(t), \tilde p(t), \tilde p_0, \tilde{a}(t) ) ={\mbox{constant}}=0, \quad \text{for a.a.} \ t\in[0, T], 
	\end{equation}
	where the value of the constant is given by~\eqref{ch21244}.
	By substituting the values of~$f(x,a)$ in~$H(x,p,p_0,a)$ and using~\eqref{ch22343}, we get, for a.a.~$
	t\in[0, T]$,
	\begin{equation*}
	\tilde p_u \tilde{u}(1-\tilde{u}- \tilde{v}-\tilde{a}c)+ \tilde p_v\rho  \tilde{v}(1-\tilde{u}- \tilde{v}) -\tilde p_v \tilde{a} \tilde{u} + \tilde p_0 =0,
	\end{equation*}
	where~$\tilde p=(\tilde p_u,\tilde p_v)$.
	
	Also, by~\eqref{ch22349} we get that
	\begin{equation}\label{ch20oskdfee}
	\underset{a\in[m,M]}{\max} H(\tilde{x}(t), \tilde p(t), \tilde p_0, a )= \underset{a\in[m,M]}{\max} \Big[-a\tilde{u}(c\tilde p_u + \tilde p_v ) + \tilde p_u \tilde{u}(1-\tilde{u}- \tilde{v})+ \tilde p_v\rho  \tilde{v}(1-\tilde{u}- \tilde{v})+\tilde p_0\Big].
	\end{equation}Thus, to maximize the term in the square brackets we must choose appropriately the value of~$\tilde{a}$ depending on the sign of~$\varphi(t):=c\tilde p_u(t)+\tilde p_v(t)$, that is
	we choose
	\begin{equation}\label{ch21631}
	\tilde{a}(t):=
	\left\{
	\begin{array}{ll}
	m &{\mbox{ if }} \varphi(t)>0, \\
	M &{\mbox{ if }} \varphi(t)<0.
	\end{array}
	\right.
	\end{equation}
	When~$\varphi(t)=0$, we are for the moment free
	to choose~$	\tilde{a}(t):=a_s(t)$ for every~$a_s(\cdot)$
	with range in~$[m,M]$,
	without affecting the maximization problem in~\eqref{ch20oskdfee}.
	
	Our next goal is to determine that~$a_s(t)$ has the expression stated in~\eqref{ch2KSM94rt3rjjjdfe} for
	a.a.~$t\in[0,T]\cap \{\varphi=0\}$.
	
	To this end, we claim that
	\begin{equation}\label{ch29id0-3rgjj}
	{\mbox{$\dot\varphi(t)=0$ a.e.~$t\in[0,T]\cap \{\varphi=0\}$.}}
	\end{equation}
	Indeed, by~\eqref{ch2HJA}, we know that~$\tilde p$ is Lipschitz continuous in~$[0,T]$,
	hence almost everywhere differentiable, and thus the same holds for~$\varphi$.
	Hence, up to a set of null measure, given~$t\in[0,T]\cap \{\varphi=0\}$,
	we can suppose that~$t$ is not an isolated point in such a set,
	and that~$\varphi$ is differentiable at~$t$. That is, there exists an infinitesimal sequence~$h_j$
	for which~$\varphi(t+h_j)=0$ and
	$$ \dot\varphi(t)=\lim_{j\to+\infty}\frac{\varphi(t+h_j)-\varphi(t)}{h_j}
	=\lim_{j\to+\infty}\frac{0-0}{h_j}=0,$$
	and this establishes~\eqref{ch29id0-3rgjj}.
	
	Consequently, in light of~\eqref{ch29id0-3rgjj}, a.a.~$t\in[0,T]\cap \{\varphi=0\}$
	satisfies	\begin{equation*}\begin{split}
	&	0 =\dot\varphi(t)= c\frac{d\tilde p_u}{dt}(t)+ \frac{d\tilde p_v}{dt}(t) \\&\qquad = c\big[ -\tilde p_u(t)(1-2\tilde{u}(t)-\tilde{v}(t)-ca_s(t))+\tilde p_v(t) (\rho \tilde{v}(t)+a_s(t))  \big]\\&\qquad\qquad
	+ \tilde p_u(t)\tilde u(t)-\tilde p_v(t) \rho(1-\tilde{u}(t)-2\tilde{v}(t)).\end{split}
	\end{equation*}
	Now, since~$\varphi(t)=0$, we have that~$ \tilde p_v(t)=- c\tilde p_u(t)$; inserting this information in the last equation, we get
	\begin{equation}\label{ch20004}
	0= -\tilde p_u c (1-2\tilde u-\tilde v-a_s c) -\tilde p_u \rho c^2 \tilde v - \tilde p_u a_s c^2 + \tilde p_u \tilde u+ \tilde p_u \rho c (1-\tilde u-2\tilde v).
	\end{equation}
	Notice that if~$\tilde p_u=0$, then~$\tilde p_v=-c \tilde p_u=0$; moreover, by~\eqref{ch22343}, one gets~$\tilde p_0=0$. But by the Pontryagin Maximum Principle one cannot have~$(\tilde p_u, \tilde p_v, \tilde p_0)=(0,0,0)$, therefore one obtains~$\tilde p_u\neq 0$ in~$\{ \varphi=0 \}$.
	Hence, dividing~\eqref{ch20004} by~$\tilde p_u$ and rearranging the terms, one gets
	\begin{equation}\label{ch20007}
	\tilde{u}(2c+1-\rho c) + c\tilde{v}(1-\rho c-2\rho)+c(\rho-1)=0.
	\end{equation}
	Differentiating the expression in~\eqref{ch20007} with respect to time, we get
	\begin{equation*}
	\tilde{u} (2c+1-\rho c) (1-\tilde{u}-\tilde{v}-ac) + c(1-\rho c-2\rho) [  \rho \tilde{v} (1-\tilde{u}-\tilde{v}) -a\tilde{u} ]=0,
	\end{equation*}
	that yields
	\begin{equation}
	a_s = \frac{(1-\tilde{u}-\tilde{v}) (	\tilde{u} (2c+1-\rho c)+\rho c) }{2c\tilde{u}(c+1)},
	\end{equation}
	which is the desired expression. 
	By a slight abuse of notation, we define the function~$a_s(t)= a_s(\tilde{u}(t), \tilde{v}(t))$ for~$t\in[0,T]$. 
	Notice that since~$\tilde{u}(t)>0$ for~$t\in[0,T]$,~$a_s(t)$ is continuous for~$t\in[0,T]$.
\end{proof}






\chapter{\label{ch3}Decay estimates for evolution equations 	with classical and fractional time-derivatives} 

		Using energy methods, we prove some power-law and exponential decay estimates for classical and nonlocal evolutionary equations. 
		The results obtained are framed into a general setting, which
		comprise, among the others,
		equations involving both standard and Caputo time-derivative, complex valued magnetic operators,
		fractional porous media equations and
		nonlocal Kirchhoff operators.
		
		Both local and fractional space diffusion are taken into account, possibly in a nonlinear setting.
		The different quantitative behaviors,
		which distinguish polynomial decays from exponential ones,
		depend heavily on the structure of the time-derivative involved in the equation.

	The content of this chapter comes from the paper \cite{decay} in collaboration with Enrico Valdinoci and the paper \cite{matrix} in collaboration with Serena Dipierro and Enrico Valdinoci.
	
\section{Introduction and main results}
	
	\subsection{Setting of the problem}
	Fractional calculus is becoming popular thanks to both the deep mathematics that it involves and its adaptability to the modelization of several real-world phenomena. As a matter of fact, integro-differential operators can describe nonlocal interactions of various type and
	anomalous
	diffusion by using suitable kernels or fractional time-derivatives, see e.g.  \cite{ultraslow}. 
	Integro-differential equations and fractional derivatives have been involved in designing, for example, wave equations, magneto-thermoelastic heat conduction, hydrodynamics, quantum physics, porous medium equations.
	
	A wide literature is devoted to the study of existence, uniqueness, regularity and asymptotic theorems.  Here we study the behaviour of the Lebesgue norm of solutions of integro-differential equations on bounded domains, extending the method of \cite{SD.EV.VV}
	to a very broad class of nonlocal equations
	and obtaining a power-law decay in time
	of the $L^s$ norm with $s\geq 1$. 
	Also, for the case of classical time-derivatives,
	we obtain exponential decays in time. The difference between
	polynomial and exponential decays in time is thus related to
	the possible presence of a fractional derivative in the operator involving the time variable.
	
	The setting in which we work
	takes into account a
	parabolic evolution of a function under the action of a spatial diffusive operator,
	which possesses suitable ``ellipticity'' properties, can be either
	classical or fractional, and can also be of nonlinear type.
	We work in a very general framework that adapts to both local
	and nonlocal operators. We comprise in this analysis also the case of
	complex valued operators and of a combination
	of fractional and classical time-derivatives.
	\medskip
	
	The main assumptions that we take is an ``abstract'' hypothesis
	which extends a construction made in \cite{SD.EV.VV}, 
	and which, roughly speaking, can be seen as a quantitative
	counterpart of the uniform ellipticity of the spatial diffusive operators.
	In~\cite{SD.EV.VV}, several time-decay estimates have been given
	covering the cases in which the time-derivative is of fractional
	type and the spatial operator is either the Laplacian,
	the fractional Laplacian, the $p-$Laplacian and
	the mean curvature equation. In this chapter,
	we deal with the cases in which {\em the time-derivative can be
		either classical or fractional, or a convex combination of the two},
	and we deal with new examples of spatial diffusive operators,
	which include the case of a
	complex valued operators. In particular,
	we present applications to the {\em fractional porous medium equation}, 
	to the {\em classical} and {\em fractional Kirchhoff equations}, to
	the {\em classical} and {\em fractional magnetic operators}.
	Referring to \cite{SD.EV.VV} for the corresponding results, we also present in Table \ref{TA1} the decay results for the \emph{$p-$Laplacian}, the \emph{nonlinear diffusion operator}, the \emph{graphical mean curvature operator}, the \emph{fractional $p-$Laplacian}, the \emph{anisotropic fractional $p-$Laplacian}, a second version  of \emph{fractional porous medium} (unfortunately, two different operators are known under the same name), and the \emph{fractional graphical mean curvature}.
	\medskip
	
	We recall that the Caputo derivative of order $\alpha\in(0,1)$ is given by
	\begin{equation*}
	\partial_t^\alpha u(t)  := \dfrac{d}{dt} \int_{0}^{t} \dfrac{u(\tau)-u(0)}{(t-\tau)^\alpha} d\tau
	\end{equation*}
	up to a normalizing constant (that we omit here for the sake of simplicity). 
	
	Let also $\lambda_1, \lambda_2 \geq 0$ be fixed. We suppose, for concreteness,
	that $$\lambda_1 + \lambda_2=1,$$
	but up to a rescaling of the operator we can take $\lambda_1, \lambda_2$
	any nonnegative number with positive sum. Let $\Omega \subset \R^n$ be a
	bounded open set and let $u_0\in L^{\infty}(\R^n)$ such that $\text{supp} \,u_0 \subset \Omega$. Consider the Cauchy problem
	\begin{equation} \label{ch3sys:generalform}
	\left\{ \begin{array}{lr}
	(\lambda_1 \partial_t^{\alpha} + \lambda_2 \partial_t) [u] + \mathcal{N}[u]=0, & {\mbox{for all }}x\in \Omega, \ t>0, \\
	u(x,t)=0, & {\mbox{for all }}x\in \R^n \setminus \Omega , \ t>0, \\
	u(x,0)=u_0(x), & {\mbox{for all }}x\in \R^n ,
	\end{array} \right.
	\end{equation}
	where $\mathcal{N}$ is a possibly nonlocal operator.
	Given $s\in[1, +\infty)$ we want to find some estimates on the ${L}^s(\Omega)$
	norm of $u$. To this end,
	we exploit analytical techniques relying on
	energy methods, exploiting also some 
	tools 
	that have been recently developed in \cite{KSVZ, VZ15, SD.EV.VV}.
	Namely, as in \cite{SD.EV.VV}, we want to compare the $L^{s}$ norm of
	the solution $u$ with an explicit function that has a power law decay,
	and to do this we take advantadge of a suitable comparison result
	and of the study of auxiliary fractional parabolic equations as
	in \cite{KSVZ, VZ15}. 
	
	\subsection{Notation and structural assumptions}
	Let us recall that for a complex valued function $v:\Omega\to\C$ the Lebesgue norm is
	\begin{equation*}
	\Vert v \Vert_{L^s(\Omega)} = \left( \int_{\Omega} |v(x)|^s \; dx \right)^{\frac{1}{s}}
	\end{equation*}
	for any $s\in[1, +\infty)$. Also, we will call $\Re \{ z\}$ the real part of $z\in\C$.
	The main assumption we take is the following: there exist $\gamma \in (0,+\infty) $ and $C\in (0,+\infty)$ such that
	\begin{equation} \label{ch3cond:complexstr}
	\Vert u(\cdot,t) \Vert_{L^{s}(\Omega) }^{s-1+\gamma} \leq C \int_{\Omega} |u(x,t)|^{s-2} \Re \{ \bar{u}(x,t)\mathcal{N} [u](x,t)\} \; dx.
	\end{equation}
	The constants $\gamma$ and $C$ and their dependence from the parameters of the problem may vary from case to case. This structural assumption says, essentially, that $\mathcal{N}$ has
	an elliptic structure and it is also related (via an integration by parts)
	to a general form of the Sobolev inequality
	(as it is apparent in the basic
	case in which~$u$ is real valued, $s:=2$ and $\mathcal{N}u:=-\Delta u$). 
	
	In our setting, the structural inequality in~\eqref{ch3cond:complexstr}
	will be the cornerstone to obtain general energy estimates,
	which, combined with appropriate barriers, in turn
	produce time-decay estimates. The results
	obtained in this way are set in a general framework, and then we
	make concrete examples of operators that satisfy the structural
	assumptions, which is sufficient to establish asymptotic bounds that fit
	to the different cases of interest and take into account the peculiarities
	of each example in a quantitative way.
	
	Our general result also includes Theorem 1 of \cite{SD.EV.VV}
	as a particular case, since, if $\mathcal{N}$ and $u$ are real valued, the
	\eqref{ch3cond:complexstr} boils down to hypothesis (1.3) of \cite{SD.EV.VV} (in any case, the applications and examples covered here
	go beyond the ones presented in \cite{SD.EV.VV} both for
	complex and for real valued operators).
	
	\subsection{Main results}
	
	The ``abstract'' result that we establish here is the following:
	
	\begin{theorem} \label{ch3thm:complex}
		Let $u$ be a solution of the Cauchy problem \eqref{ch3sys:generalform}, with $\mathcal{N}$ possibly complex
		valued. Suppose that there exist $s\in[1, +\infty)$, $\gamma\in(0,+\infty)$ and $C\in(0,+\infty)$ such that $u$ satisfies \eqref{ch3cond:complexstr}.
		Then 
		\begin{equation} \label{ch3claim1gen}
		(\lambda_1\partial_t^{\alpha} + \lambda_2\partial_t) \Vert u(\cdot,t) \Vert_{L^{s}(\Omega) } \leq -\dfrac{\Vert u(\cdot,t) \Vert_{L^{s}(\Omega) }^{\gamma}}{C},
		\qquad{\mbox{ for all }}t>0,\end{equation}
		where $C$ and $\gamma$ are the constants appearing in~\eqref{ch3cond:complexstr}. 
		Furthermore,
		\begin{equation} \label{ch3claim2gen}
		\Vert u(\cdot,t) \Vert_{L^{s}(\Omega) } \leq
		\dfrac{C_*}{1+t^{\frac{\alpha}{\gamma}}},\qquad{\mbox{ for all }}t>0,	
		\end{equation}
		for some~$C_*>0$, depending only on $C$, $\gamma$, $\alpha$
		and $\Vert u_0(\cdot) \Vert_{L^{s}(\R^n)}$.
	\end{theorem}
	
	Theorem~\ref{ch3thm:complex} here
	comprises previous results in \cite{SD.EV.VV},
	extending their applicability
	to a wider class of equations, which include the
	cases of both standard and fractional 
	time-derivatives and complex valued operators.
	
	We also recall that the power-law decay in~\eqref{ch3claim2gen}
	is due to
	the behaviour of the solution of the equation 
	\begin{equation} \label{ch3mittagleffler}
	\partial_t^{\alpha} e(t)=-e(t),
	\end{equation}
	for $t\in(0, +\infty)$. Indeed, the solution
	of~\eqref{ch3mittagleffler} is explicit in terms of the Mittag-Leffler function and it is asymptotic to $\frac{1}{t^{\alpha}}$ as $t\rightarrow +\infty$ (see \cite{Mittag-Leffler}, \cite{Mittag-Leffler_asymt}); notice that
	the latter decay corresponds to the one
	in~\eqref{ch3mittagleffler} when $\gamma=1$. 
	
	As pointed out in \cite{KSVZ}, 
	the power law decay for solutions of
	time-fractional equations is, in general, unavoidable.
	On the other hand, solutions of equations
	driven by the standard time-derivative
	of the type $$ \partial_t v(t) + \mathcal{N}[v](t)=0$$ often have a faster decay in many concrete examples, for instance for $\mathcal{N}=-\Delta$ where exponential decay is attained. This particular feature of
	the classical heat equation is in fact a special
	case of a general phenomenon, described
	in details in the following result:
	
	\begin{theorem}\label{ch3thm:classic}
		Let $u$ be a solution of the Cauchy problem \eqref{ch3sys:generalform} with only classical derivative ($\lambda_1=0$) and $\mathcal{N}$ possibly complex
		valued. Suppose that there exist $s\in[1, +\infty)$, $\gamma\in(0,+\infty)$ and $C\in(0,+\infty)$ such that $u$ satisfies \eqref{ch3cond:complexstr}.
		Then, for some~$C_*>0$, depending only on the constants~$C$ and~$\gamma$
		in~\eqref{ch3cond:complexstr}, 
		and on~$\Vert u_0(\cdot ) \Vert_{L^{s}(\R^n)}$, we have that:
		\begin{itemize}
			\item	if $0<\gamma \leq 1$ the solution $u$ satisfies
			\begin{equation} \label{ch3claim3}
			\Vert u(\cdot ,t) \Vert_{L^{s}(\Omega) } \leq
			C_* \, e^{-\frac{t}{C}},\qquad{\mbox{for all }}t>0;	
			\end{equation}
			\item if $ \gamma>1$, the solution $u$ satisfies
			\begin{equation} \label{ch3claim4}
			\Vert u(\cdot ,t) \Vert_{L^{s}(\Omega) } \leq
			\dfrac{C_*}{1+t^{\frac{1}{\gamma-1}}},\qquad{\mbox{for all }}t>0.	
			\end{equation}
		\end{itemize} 
	\end{theorem} 
	
	We stress that Theorem~\ref{ch3thm:classic}
	is valid for a very general class of diffusive
	operators~${\mathcal{N}}$, including also
	the ones which take into account
	fractional derivatives in the space-variables.
	In this sense, the phenomenon described in
	Theorem~\ref{ch3thm:classic} is that:
	\begin{itemize} \item on the one hand,
		the fractional behaviour induces power-law
		decay,
		\item on the other hand, for long times,
		the interactions between different derivatives
		``decouple'': for instance, a space-fractional
		derivative, which would naturally induce
		a polynomial decay, does not asymptotically ``interfere''
		with a classical time-derivative
		in the setting of Theorem~\ref{ch3thm:classic},
		and the final result is that the decay in time is
		of exponential, rather than polynomial, type.
	\end{itemize}
	
	The fact that long-time asymptotics
	of mixed type (i.e. classical time-derivatives
	versus fractional-space diffusion) reflect the
	exponential decay of linear ordinary differential
	equations was also observed in~\cite{MR3703556}
	for equations inspired by the Peierls-Nabarro model
	for atom dislocations in crystal.\medskip
	
	As we will see in the proof
	of Theorem~\ref{ch3thm:classic}, the idea is to find a supersolution of \eqref{ch3claim1gen} and use a comparison principle in order to estimate the decay of the solution $u$. For the case of mixed derivatives, Vergara and Zacher \cite{VZ15} find both a supersolution and a subsolution decaying as $t^{-\frac{\alpha}{\gamma}}$. When $\alpha\rightarrow 1$, thus when the mixed derivative is approaching the classical one, the subsolution tends to 0. This allows possibly better decays, which are in fact proven. On the other side, the supersolution gains some extra decay, possibly reaching an exponential decay.
	\medskip
	
	The optimality of the decay estimates
	obtained in our results
	and some further comparisons with the existing literature are discussed in Subsection \ref{ch3sss:comparison}. 
	
	\subsection{Applications} \label{ch3applications}
	
	We now present several applications of Theorem \ref{ch3thm:complex} to some concrete examples. 
	
	\paragraph{The case of the fractional porous medium equation.}
	Let $0<\sigma<1$
	and 
	\begin{equation}\label{ch3kappa}
	K:\R^n \rightarrow \R^n
	\end{equation}
	be the positive function
	\begin{equation*}
	K(x):= c(n,\sigma) |x|^{-(n-2\sigma)},
	\end{equation*}
	being $c(n,\sigma)$ a constant. 
	The fractional\footnote{As a matter of fact,
		as clearly explained in 
		\url{https://www.ma.utexas.edu/mediawiki/index.php/Nonlocal_porous_medium_equation},
		the fractional porous medium equation
		is ``the name currently given to two very different equations''.
		The one introduced in~\cite{MR2737788}
		has been studied in details in \cite{SD.EV.VV}
		in terms of decay estimates. We focus here
		on the equation introduced in~\cite{porous}.
		As discussed in the above mentioned mediawiki page,
		the two equations have very different structures and typically
		exhibit different behaviors, so we think that it is a nice feature that,
		combining the results here with those in~\cite{SD.EV.VV},
		it follows that a complete set of decay estimates is valid for {\em both} the
		fractional porous medium equations at the same time.}
	porous medium operator (as defined in \cite{porous}) is 
	\begin{equation} \label{ch3op:porous}
	\mathcal{N}[u]:=-\nabla \cdot (u \nabla \mathcal{K}(u)), \qquad {\mbox{where}}\qquad\mathcal{K}(u):=u \star K
	\end{equation}
	where $\star$ denotes the convolution. This operator is used to describe the diffusion of a liquid under pressure in a porous environment in presence of memory effects and long-range
	interactions, and also has some application in biological models, see~\cite{porous}.\medskip
	
	In this framework, the following result holds:
	
	\begin{theorem} \label{ch3thm:porous}
		Take $u_0(x) \in L^{\infty}(\R^n)$ and let $u$ be a solution in $\Omega \times (0, + \infty)$ to \eqref{ch3sys:generalform} with $\mathcal{N}$ the fractional porous medium operator as in \eqref{ch3op:porous}. Then for all $s\in (1, +\infty)$ there exists $C_*>0$ depending on $n,\ s,\ \sigma,\ \Omega$ such that 
		\begin{equation*}
		\Vert u(\cdot, t) \Vert_{L^{s}(\Omega) } \leq \dfrac{C_*}{1+t^{\alpha /2}}.
		\end{equation*}
		Also, in the case of only classical derivative ($\lambda_1=0$), we have 
		\begin{equation*}
		\Vert u(\cdot, t) \Vert_{L^{s}(\Omega) } \leq \dfrac{C_*}{1+t}
		\end{equation*}
		where $C_*>0$, possibly different than before,  depends on $n,\ s,\ \sigma,\ \Omega$.	
	\end{theorem}
	
	\paragraph{The case of the Kirchhoff operator and the fractional Kirchhoff operator.}
	The Kirchhoff equation describes
	the movement of an elastic string that is constrained at the extrema,
	taking into account a possible growth of 
	the tension of the vibrating string in view of its extension. It was first introduced by
	Gustav Robert Kirchhoff in~1876, see
	\url{https://archive.org/details/vorlesungenberm02kircgoog},
	and fully addressed from the mathematical
	point of view in the 20th century, see~\cite{MR0002699}.
	
	Parabolic equations of Kirchhoff type have been widely studied during the '90s (see for example \cite{gobbino} and the reference therein). Recently a fractional counterpart to the Kirchhoff operator has been introduced by Fiscella and Valdinoci \cite{kirchhoff}.\medskip
	
	The setting that we consider here is the following.
	Let $m:[0,+\infty)\to[0,+\infty)$ be an nondecreasing function. A typical example is 
	\begin{equation}\label{ch3def:m}
	m(\xi)=m_0 +b\xi
	\end{equation}
	where $b> 0$ and $m_0 \geq 0$. We consider here both the cases\footnote{The case $m_0=0$ for \eqref{ch3def:m} is usually called the degenerate case and it presents several additional
		difficulties with respect with the non-degenerate case. } in which $m(0)>0$ and in which $m$ takes the form in~\eqref{ch3def:m} with $m_0=0$. In this setting, the Kirchhoff operator that we take into account is 
	\begin{equation}\label{ch3KKOP} \mathcal{N}[u]:= m \left(\Vert \nabla u \Vert_{L^2(\Omega)}^2\right)  (-\Delta)u =0.  \end{equation}
	Then, we obtain the following decay estimates:
	
	\begin{theorem}\label{ch3thm:cl_kirchhoff}
		Let $u$ be the solution of problem \eqref{ch3sys:generalform} with $\mathcal{N}$ the Kirchhoff operator in \eqref{ch3KKOP}. Then  
		there exist $\gamma>0$ and $C>0$ depending on $n,\ s,\ \Omega,\ \inf m(t)$ such that
		\begin{equation*}
		\Vert u(\cdot ,t) \Vert_{L^{s}(\Omega) } \leq \dfrac{C}{1+t^{\frac{\alpha}{\gamma}}},\qquad{\mbox{for all }}t>0,
		\end{equation*}
		in the following cases:
		\begin{itemize}
			\item[(i)\;\;] for all $s\in[1, +\infty)$ when $m$ is non-degenerate; in particular, in this case~$\gamma=1$.
			\item[(ii)\;\;] for all $s\in[1,+ \infty)$ when $m$ is degenerate and $n\leq 4$; in particular, in this case~$\gamma=3$.
			\item[(iii)\;\;] for $s\leq\frac{2n}{n-4}$ when $m$ is degenerate and $n>4$; in particular, in this case~$\gamma=3$.
		\end{itemize}
		Moreover, if we take $\lambda_1=0$,
		then there exists~$C_*>0$, $C'>0$
		depending on $n,\ s,\ \Omega,\ \inf m(t)$, 
		for which the following statements hold true:
		\begin{itemize}
			\item in case (i) we have
			\begin{equation*}
			\Vert u(\cdot,t) \Vert_{L^{s}(\Omega) } \leq
			C_* \, e^{-\frac{t}{C'}},\qquad{\mbox{for all }}t>0,
			\end{equation*}
			\item in cases (ii) and (iii) we have
			\begin{equation*}
			\Vert u(\cdot,t) \Vert_{L^{s}(\Omega) } \leq \dfrac{C_*}{1+t^\frac{1}{2}},\qquad{\mbox{for all }}t>0.
			\end{equation*}\end{itemize}
	\end{theorem}
	
	Next, we consider the case of the fractional Kirchhoff operator.
	We take a nondecreasing positive function
	$M:[0,+\infty)\rightarrow[0,+\infty)$. As for the classic Kirchhoff
	operator, we consider either the case when $M(0)>0$
	or the case $M(\xi)=b\xi$ with $b>0$.
	We fix $\sigma\in(0,1)$. We define the norm 
	\begin{equation}\label{ch3FKPO-1}
	\Vert u(\cdot , t) \Vert_{Z} = \left( \int_{\R^{2n}} \frac{|u(x,t)-u(y,t)|^2 }{|x-y|^{n+2\sigma}} \; dxdy \right)^{\frac{1}{2}} .
	\end{equation}
	Finally, the fractional Kirchhoff operator reads
	\begin{equation}\label{ch3FKPO}
	\mathcal{N}[u](x,t):= -M\left( \Vert u(\cdot , t)\Vert_{Z}^2 \right) \int_{\R^n}  \frac{ u(x+y,t) + u(x-y,t) -2u(x,t)}{|x-y|^{n+2\sigma}} \; dy.
	\end{equation}
	In this setting, our result is the following:
	
	\begin{theorem} \label{ch3thm:Kirchhoff}
		Let $u$ be the solution of problem \eqref{ch3sys:generalform} with $\mathcal{N}$ the fractional Kirchhoff operator in~\eqref{ch3FKPO}. Then  
		there exist $\gamma>0$ and $C>0$, depending on $K$, $n$, $s$, $\Omega$
		and~$\inf M(\xi)$, such that
		\begin{equation*}
		\Vert u(\cdot ,t) \Vert_{L^{s}(\Omega) } \leq
		\dfrac{C}{1+t^{\frac{\alpha}{\gamma}}}
		,\qquad{\mbox{for all }}t>0,
		\end{equation*}
		in the following cases:
		\begin{itemize}
			\item[(i)\;\;] for all $s\in[1, +\infty)$ when $M$ is non-degenerate; in particular, in this case~$\gamma=1$.
			\item[(ii)\;\;] for all $s\in[1,+ \infty)$ when $M$ is degenerate and $n\leq 4\sigma$; in particular, in this case~$\gamma=3$.
			\item[(iii)\;\;] for $s\leq\frac{2n}{n-4\sigma}$ when $M$ is degenerate and $n>4\sigma$; in particular, in this case~$\gamma=3$.
		\end{itemize}
		Moreover, if we take $\lambda_1=0$,
		then there exists $C_*>0$, depending on $n,\ s,\ \Omega,\ \inf M(t)$, such that:
		\begin{itemize} \item in case (i) we have
			\begin{equation*}
			\Vert u(\cdot,t) \Vert_{L^{s}(\Omega) } \leq
			C_* \, e^{-\frac{t}{C'}}	,\qquad{\mbox{for all }}t>0,
			\end{equation*}
			for some~$C'>0$, \item in cases (ii) and (iii) we have
			\begin{equation*}
			\Vert u(\cdot,t) \Vert_{L^{s}(\Omega) } \leq \dfrac{C_*}{1+t^\frac{1}{2}},\qquad{\mbox{for all }}t>0.
			\end{equation*}\end{itemize}
	\end{theorem}
	
	It is interesting to remark that
	the cases~(i), (ii) and~(iii) in Theorem~\ref{ch3thm:Kirchhoff} formally
	reduce to those in Theorem~\ref{ch3thm:cl_kirchhoff}
	when~$\sigma\to1$.
	
	\paragraph{The case of the magnetic operator and the fractional magnetic operator.}
	
	We consider here an operator similar to Schr\"{o}dinger equation
	with a magnetic potential (see e.g.~\cite{MR0142894} and the references therein), that is
	\begin{equation}\label{ch3NuMAG}
	\mathcal{N}[u]:= -(\nabla -iA)^2 u(x,t)= -\Delta u + |A|^2u -iA\cdot\nabla u -\nabla \cdot (iAu)
	\end{equation} 
	where $A: \R^n \rightarrow \R^n$ has the physical meaning of a magnetic field
	(in this case, one usually studies the three-dimensional case~$n=3$, but our
	approach is general).
	The goal of these pages is to apply Theorem~\ref{ch3thm:complex}
	to the magnetic operator in~\eqref{ch3NuMAG}, thus obtaining decay estimates in time in this framework.
	
	It is interesting to remark that the operator in~\eqref{ch3NuMAG}
	is structurally very different from the linear Schr\"{o}dinger operator,
	which corresponds to the choice
	\begin{equation}\label{ch3NuMAG:SC}
	\mathcal{N}[u]= -i(\Delta +V)u.
	\end{equation}
	Indeed, for the operator in~\eqref{ch3NuMAG:SC} decay estimates 
	in time do not\footnote{Indeed, if~$V\in\R$ and~$u$
		is a solution of the Schr\"{o}dinger parabolic equation~$
		\partial_t u+i(\Delta +V)u=0$ in~$\Omega$
		with homogeneous data along~$\partial\Omega$, the conjugated equation reads~$
		\partial_t \bar u-i(\Delta +V)\bar u=0$, and therefore
		\begin{eqnarray*}&& \partial_t\int_\Omega |u(x,t)|^2\,dx=
			\int_\Omega u(x,t)\,\partial_t\bar u(x,t)+\bar u(x,t)\,\partial_t u(x,t)\,dx\\&&\qquad
			=i\int_\Omega u(x,t)\,\Delta\bar u(x,t)-\bar u(x,t)\,\Delta u(x,t)\,dx
			\\&&\qquad=\int_\Omega \nabla\cdot\big( u(x,t)\,\nabla\bar u(x,t)-\bar u(x,t)\,\nabla u(x,t)\big)\,dx
			=0,\end{eqnarray*}
		where the last identity follows from the Divergence Theorem and the boundary conditions.
		This shows that decay estimates in time are in general
		not possible in this setting, thus highlighting an interesting
		difference between the Schr\"{o}dinger operator in~\eqref{ch3NuMAG:SC}
		and the magnetic operator in~\eqref{ch3NuMAG}.
		
		This difference, as well as the computation above,
		has a natural physical meaning,
		since in the Schr\"odinger equation the squared modulus of
		the solution represents
		the probability density of a wave function,
		whose total amount remains constant if no dissipative forces appear in the equation.} hold in general, not even in the case of classical 
	time-derivatives.   \medskip
	
	The decay estimate for the classical magnetic operator is the following:
	
	\begin{theorem} \label{ch3thm:cl_magnetic}
		Let $u$ be the solution of problem \eqref{ch3sys:generalform} with $\mathcal{N}$ the magnetic operator in~\eqref{ch3NuMAG}. 
		Then for all $s \in [1, +\infty)$ there exist $C_1>0$ depending on $A$, $n$, $s$ and $\sigma$ such that
		\begin{equation*}
		\Vert u(\cdot ,t) \Vert_{L^{s}(\Omega) } \leq \dfrac{C_1}{1+t^{{\alpha}}}\qquad{\mbox{for all }}t>0.
		\end{equation*}	
		Moreover, in the case of classical derivatives ($\lambda_1=0$), we have
		\begin{equation*}
		\Vert u(\cdot ,t) \Vert_{L^{s}(\Omega) } \leq C_2\,e^{-\frac{t}{C_3}}\qquad{\mbox{for all }}t>0
		\end{equation*}
		for some $C_2$, $C_3>0$, depending on $A$, $n$, $ s$ and $\sigma$.
	\end{theorem}
	
	In \cite{groundstates} D'Avenia and Squassina introduced a fractional operator where a magnetic field $A: \R^n \rightarrow \R^n $ appears. Their aim was to study the behaviour of free particles interacting with a magnetic field. For a fixed $\sigma \in (0,1)$, such an operator in dimension $n$ reads
	\begin{equation}\label{ch3SQ}
	\mathcal{N}[u](x,t):= \int_{\R^n} \frac{u(x,t)-e^{i(x-y)A(\frac{x+y}{2})} u(y,t)}{|x-y|^{n+2\sigma}} \;dy.
	\end{equation}
	In the appropriate framework, the fractional magnetic operator in~\eqref{ch3SQ}
	recovers the classical magnetic operator in~\eqref{ch3NuMAG}
	as~$\sigma\to1$, see~\cite{MR3528339} (see also~\cite{MR3794886}
	for a general approach involving also nonlinear operators).
	
	In the setting of the fractional magnetic operator, we present the following result:
	
	\begin{theorem} \label{ch3thm:magnetic}
		Let $u$ be the solution of problem \eqref{ch3sys:generalform} with $\mathcal{N}$ the
		fractional magnetic operator in~\eqref{ch3SQ}. Then for all $s \in [1, +\infty)$ there exist $C_1>0$ depending on $n$, $s$ and $\sigma$ such that
		\begin{equation*}
		\Vert u(\cdot ,t) \Vert_{L^{s}(\Omega) } \leq \dfrac{C_1}{1+t^{{\alpha}}}\qquad{\mbox{for all }}t>0.
		\end{equation*}	
		Moreover, in the case of classical derivatives ($\lambda_1=0$), we have
		\begin{equation*}
		\Vert u(\cdot ,t) \Vert_{L^{s}(\Omega) } \leq C_2\,e^{-\frac{t}{C_3}}\qquad{\mbox{for all }}t>0,
		\end{equation*}
		for some $C_2$, $C_3>0$ depending on $n$, $s$ and $\sigma$.
	\end{theorem}
	
	The magnetic operators present a crucial difference with respect to
	the other operators considered in the previous applications, since they are complex valued operators.\medskip

	\paragraph{Other operators.}
	
	We point out that condition \eqref{ch3cond:complexstr} has already been checked in many cases in \cite{SD.EV.VV}. We present here very briefly the operators treated there that may need an introduction. The list includes the cases of the classical~$p$-Laplacian
	and porous media diffusion (see~\cite{dibenedetto2012degenerate,vazquez2007porous})
	$$ \Delta_p u^m := {\rm div} (|\nabla u^m|^{p-2}\nabla u^m),
	\qquad{\mbox{with $p\in(1,+\infty)$ and $m\in(0,+\infty)$,}}$$
	the case of graphical mean curvature, given in formula (13.1) of \cite{giusti}, 
	$$ {\rm div}\left( \frac{\nabla u}{\sqrt{1+|\nabla u|^2}}\right), $$
	the case of the fractional $p$-Laplacian (see e.g.~\cite{brasco2018higher})
	\begin{eqnarray*}&&(-\Delta)^s_pu(x):=
		\int_{\R^n}\frac{|u(x)-u(y)|^{p-2}(u(x)-u(y))}{|x-y|^{n+sp}}\,dy,\\&&\qquad{\mbox{with $
				p\in(1,+\infty)$ and $s\in(0,1)$,}}\end{eqnarray*}
	and possibly even the sum of different nonlinear operator of this type, with coefficients $\beta_j>0$, 
	$$  \sum_{j=1}^N \beta_j (-\Delta)^{s_j}_{p_j} u, \qquad \text{with} \
	p_j\in(1,+\infty) \ \text{and}  \ s_j\in(0,1),  $$
	the case of the anisotropic fractional Laplacian, that is the sum of fractional directional derivatives in the directions of the space $e_j$, given by
	$$(-\Delta_{\beta})^{\sigma} u(x)= \sum_{j=1}^{n} \beta_j (-\partial_{x_j}^2)^{\sigma_j} u(x) $$
	for $\beta_j>0$, $\beta=(\beta_1, \dots, \beta_n)$ and  $\sigma=(\sigma_1, \dots, \sigma_n)$, where 
	$$ (-\partial_{x_j}^2)^{\sigma_j} u(x) = \int_{\R} \frac{u(x)- u(x+\rho e_j)}{\rho^{1+2\sigma_j}} d\rho, $$
	considered for example in \cite{farina2017regularity}. 
	The list of possible diffusion operators continues with
	a fractional porous media operators (see~\cite{MR2737788})
	\begin{equation*}
		{\mathcal{P}}_{1,s}(u):=(-\Delta)^s u^m
		\qquad{\mbox{with $s\in(0,1)$ and $m\in(0,+\infty)$,}}
	\end{equation*}
	and the graphical fractional mean curvature operator (see~\cite{barrera2014bootstrap})
	\begin{eqnarray*}&& {\mathcal{H}}^s(u)(x):=\int_{\R^n} F\left(\frac{u(x)-u(x+y)}{|y|}\right)\frac{dy}{|y|^{n+s}},\\&&\qquad\qquad{\mbox{with
				$s\in(0,1)$ and }}F(r):=\int_0^r \frac{d\tau}{(1+\tau^2)^{\frac{n+1+s}{2}}},\end{eqnarray*}
	
	For the sake of brevity, we recall the corresponding results in Table \ref{TA1}.

	{
		\begin{center}\begin{table}
				\fontsize{8}{12}\selectfont
				\begin{tabular}{l | c | c | c | c | }
					$\,$ & {\bf Operator ${\mathcal{N}}$} & {\bf Values of~$\lambda_1$, $\lambda_2$} & {\bf Range of~$\ell$} & {\bf Decay rate $\Theta$} 
					\\ \hline \\
					\begin{minipage}[l]{0.15\columnwidth}
						{\bf Nonlinear classical diffusion 
					}\end{minipage}
					& $\Delta_p u^m$ &
					$\lambda_1\in(0,1]$, $\lambda_2\in[0,1)$ & $\ell\in[1,+\infty)$ &
					$\Theta(t)=\frac1{t^{\frac\alpha{m(p-1)}}}$ 
					\\ \hline \\
					\begin{minipage}[l]{0.15\columnwidth}
						{\bf Nonlinear classical diffusion 
					}\end{minipage}
					& $\Delta_p u^m$ {with}  $(m,p)\neq(1,2)$ &
					$\lambda_1=0$, $\lambda_2=1$ & $\ell\in[1,+\infty)$ &
					$\Theta(t)=\frac1{t^{\frac1{m(p-1)-1}}}$ 
					\\ \hline \\
					\begin{minipage}[l]{0.15\columnwidth}
						{\bf Bi-Laplacian }\end{minipage}
					& $\Delta_2 u$  &
					$\lambda_1=0$, $\lambda_2=1$ & $\ell\in[1,+\infty)$ &
					$\Theta(t)=e^{-\frac{t}{C}}$ 
					\\ \hline \\
					\begin{minipage}[l]{0.15\columnwidth}
						{\bf Graphical mean curvature}\end{minipage}
					& $ {\rm div}\left( \frac{\nabla u}{\sqrt{1+|\nabla u|^2}}\right)$ &
					$\lambda_1\in(0,1]$, $\lambda_2\in[0,1)$ & $\ell\in[1,+\infty)$ &
					$\Theta(t)=\frac1{t^\alpha}$ 
					\\ \hline \\
					\begin{minipage}[l]{0.15\columnwidth}
						{\bf Graphical mean curvature}\end{minipage}
					& $ {\rm div}\left( \frac{\nabla u}{\sqrt{1+|\nabla u|^2}}\right)$ &
					$\lambda_1=0$, $\lambda_2=1$ & $\ell\in[1,+\infty)$ &
					$\Theta(t)=e^{-\frac{t}{C}}$ 
					\\ \hline \\
					\begin{minipage}[l]{0.15\columnwidth}
						{\bf Fractional $p{\mbox{-}}$Laplacian}\end{minipage}
					& $(-\Delta)^s_pu$ &
					$\lambda_1\in(0,1]$, $\lambda_2\in[0,1)$ & $\ell\in[1,+\infty)$ &
					$\Theta(t)=\frac1{t^{\frac\alpha{p-1}}}$ 
					\\ \hline \\
					\begin{minipage}[l]{0.15\columnwidth}
						{\bf Fractional $p{\mbox{-}}$Laplacian}\end{minipage}
					& $(-\Delta)^s_pu$, $p> 2$ &
					$\lambda_1=0$, $\lambda_2=1$ & $\ell\in[1,+\infty)$ &
					$\Theta(t)=\frac1{t^{\frac1{p-2}}}$ 
					\\ \hline \\
					\begin{minipage}[l]{0.15\columnwidth}
						{\bf Fractional $p{\mbox{-}}$Laplacian}\end{minipage}
					& $(-\Delta)^s_pu$, $p\leq 2$ &
					$\lambda_1=0$, $\lambda_2=1$ & $\ell\in[1,+\infty)$ &
					$\Theta(t)=e^{-\frac{t}{C}}$ 
					\\ \hline \\
					\begin{minipage}[l]{0.15\columnwidth} 
						{\bf Superposition of fractional $p{\mbox{-}}$Laplacians}\end{minipage}
					&{ $ \sum_{j=1}^N \beta_j (-\Delta)^{s_j}_{p_j} u$, $\beta_j>0$}&
					$\lambda_1\in(0,1]$, $\lambda_2\in[0,1)$ & $\ell\in[1,+\infty)$ &
					$\Theta(t)=\frac1{t^{\frac\alpha{p_{\max}-1}}}$ 
					\\ \hline \\
					\begin{minipage}[l]{0.15\columnwidth}
						{\bf Superposition of fractional $p{\mbox{-}}$Laplacians}\end{minipage}
					& \makecell{$ \sum_{j=1}^N \beta_j (-\Delta)^{s_j}_{p_j} u$, \\ with $\beta_j>0$ and $p_{\max}>2$} &
					$\lambda_1=0$, $\lambda_2=1$ & $\ell\in[1,+\infty)$ &
					$\Theta(t)=\frac1{t^{\frac1{p_{\max}-2}}}$ 
					\\ \hline \\
					\begin{minipage}[l]{0.15\columnwidth}
						{\bf Superposition of fractional $p{\mbox{-}}$Laplacians}\end{minipage}
					& \makecell{$ \sum_{j=1}^N \beta_j (-\Delta)^{s_j}_{p_j} u$, \\ with $\beta_j>0$ and $p_{\max}\leq2$}&
					$\lambda_1=0$, $\lambda_2=1$ & $\ell\in[1,+\infty)$ &
					$\Theta(t)=e^{-\frac{t}{C}}$ 
					\\ \hline \\
					\begin{minipage}[l]{0.15\columnwidth}
						{\bf Superposition of anisotropic fractional Laplacians}\end{minipage}
					& $ \sum_{j=1}^N \beta_j (-\partial_{x_j}^2)^{s_j} u$, $\beta_j>0$&
					$\lambda_1\in(0,1]$, $\lambda_2\in[0,1)$ & $\ell\in[1,+\infty)$ &
					$\Theta(t)=\frac1{t^{\alpha}}$ 
					\\ \hline \\
					\begin{minipage}[l]{0.15\columnwidth}
						{\bf Superposition of anisotropic fractional Laplacians}\end{minipage}
					& $ \sum_{j=1}^N \beta_j (-\partial_{x_j}^2)^{s_j} u$, $\beta_j>0$&
					$\lambda_1=0$, $\lambda_2=1$ & $\ell\in[1,+\infty)$ &
					$\Theta(t)=e^{-\frac{t}{C}}$ 
					\\ \hline \\
					\begin{minipage}[l]{0.15\columnwidth} 
						{\bf Fractional porous media~I}\end{minipage}
					& $ {\mathcal{P}}_{1,s}(u)$&
					$\lambda_1\in(0,1]$, $\lambda_2\in[0,1)$ & $\ell\in[1,+\infty)$ &
					$\Theta(t)=\frac1{t^{\frac\alpha{m}}}$ 
					\\ \hline \\
					\begin{minipage}[l]{0.15\columnwidth}
						{\bf Fractional porous media~I}\end{minipage}
					& $ {\mathcal{P}}_{1,s}(u)$, $m>1$ &
					$\lambda_1=0$, $\lambda_2=1$ & $\ell\in[1,+\infty)$ &
					$\Theta(t)=\frac1{t^{\frac1{m-1}}}$ 
					\\ \hline \\
					\begin{minipage}[l]{0.15\columnwidth}
						{\bf Fractional porous media~I}\end{minipage}
					& $ {\mathcal{P}}_{1,s}(u)$, $m\leq1$ &
					$\lambda_1=0$, $\lambda_2=1$ & $\ell\in[1,+\infty)$ &
					$\Theta(t)=e^{-\frac{t}{C}}$ 
					\\ \hline \\
					\begin{minipage}[l]{0.15\columnwidth}
						{\bf Fractional graphical mean curvature}\end{minipage}
					& $ {\mathcal{H}}^s(u)$&
					$\lambda_1\in(0,1]$, $\lambda_2\in[0,1)$ & $\ell\in[1,+\infty)$ &
					$\Theta(t)=\frac1{t^{\alpha}}$ 
					\\ \hline \\
					\begin{minipage}[l]{0.15\columnwidth}
						{\bf Fractional graphical mean curvature}\end{minipage}
					& $ {\mathcal{H}}^s(u)$&
					$\lambda_1=0$, $\lambda_2=1$ & $\ell\in[1,+\infty)$ &
					$\Theta(t)=e^{-\frac{t}{C}}$ 
					\\ \hline
				\end{tabular}\medskip
				\caption{Results from~\cite{SD.EV.VV}.}\label{TA1}
		\end{table}\end{center}
		
	}

	The examples provided here show that the ``abstract''
	structural hypothesis \eqref{ch3cond:complexstr} is reasonable and can be explicitly
	checked in several cases of interest.
	We are also confident that other interesting
	examples fulfilling such an assumption
	can be found, therefore Theorem \ref{ch3thm:complex} turns out to play a pivotal
	role in the asymptotics of real and complex valued,
	possibly nonlinear, and possibly fractional, operators.

	\subsubsection{Comparison with the existing literature} \label{ch3sss:comparison}
	
	In general, in problems of the type \eqref{ch3sys:generalform} it is very difficult to provide
	explicit solutions and often the system has no unique solution, see e.g.~\cite{biler3al}. Therefore, even partial information on the solutions is important. 
	
	In the case of a Kirchhoff parabolic equation with purely classical time-derivative in the degenerate case $m(0)=0$, Ghisi and Gobbino \cite{gobbino} found the time-decay estimate
	\begin{equation}\label{ch3GOBLA1}  c (1+t)^{-1} \leq \Vert \nabla u(\cdot, t)
	\Vert_{L^2(\Omega)}^2 \leq C (1+t)^{-1} \qquad{\mbox{for all }}t>0.\end{equation}
	for some costants $C, c>0$ depending on initial data. {F}rom this,
	performing an integration of the gradient along paths\footnote{More precisely,
		the fact that~\eqref{ch3GOBLA1} implies~\eqref{ch3GOBLA2}
		can be seen as a consequence of the following observation:
		for every~$u\in C^\infty_0(\Omega)$,
		\begin{equation}\label{ch3L2GRAs}
		\|u\|_{L^2(\Omega)}\le
		C\,\|\nabla u\|_{L^2(\Omega)},
		\end{equation}
		where~$C>0$ depends on~$n$ and~$\Omega$.
		Indeed, 
		fix~$x_0\in\R^n$ such that~$B_1(x_0)\subset\R^n\setminus\Omega$
		and~$\Omega\subset B_R(x_0)$, for some~$R>1$.
		Then, for every~$x\in\Omega$ we have that~$|x-x_0|\in[1,R]$ and thus
		\begin{eqnarray*}
			&& |u(x)|^2=|u(x)-u(x_0)|^2=\left| 
			\int_0^1 \nabla u(x_0+t(x-x_0))\cdot(x-x_0)\,dt
			\right|^2\\
			&&\qquad\le |x-x_0|^2\,\int_0^1 |\nabla u(x_0+t(x-x_0))|^2\,dt\le
			R^2\,\int_0^1 |\nabla u(x_0+t(x-x_0))|^2\,dt.
		\end{eqnarray*}
		On the other hand, if~$t\in[0,1/R)$ we have that
		$$ \big|t(x-x_0)\big|< \frac{|x-x_0|}{R}\le1$$
		and so~$x_0+t(x-x_0)\in B_1(x_0)\subset\R^n\setminus\Omega$,
		which in turn implies that~$\nabla u(x_0+t(x-x_0))=0$. This gives that
		$$ |u(x)|^2\le R^2\,\int_{1/R}^1 |\nabla u(x_0+t(x-x_0))|^2\,dt.$$
		Hence, using the substitution~$x\mapsto y:=x_0+t(x-x_0)$, we conclude that
		\begin{eqnarray*}
			&&\int_\Omega |u(x)|^2\,dx\le
			R^2\,\int_{1/R}^1 \left[\int_{\R^n}|\nabla u(x_0+t(x-x_0))|^2\,dx\right]\,dt
			=
			R^2\,\int_{1/R}^1 \left[\int_{\R^n}|\nabla u(y)|^2\,\frac{dy}{t^n}\right]\,dt\\&&\qquad
			\leq R^{n+2}\,\int_{1/R}^1 \left[\int_{\R^n}|\nabla u(y)|^2\,dy\right]\,dt
			\leq
			R^{n+2}\,\|\nabla u\|^2_{L^2(\R^n)}=R^{n+2}\,\|\nabla u\|^2_{L^2(\Omega)},\end{eqnarray*}
		which proves~\eqref{ch3L2GRAs}.}, one can find the estimate
	\begin{equation}\label{ch3GOBLA2}
	\Vert u(\cdot, t) \Vert_{L^2(\Omega)} \leq C (1+t)^{-\frac{1}{2}}
	\qquad{\mbox{for all }}t>0. \end{equation}
	The latter is exactly the estimate we found
	in Theorem~\ref{ch3thm:cl_kirchhoff} as a particular case of our analysis.
	
	The fractional porous medium equation with classical derivative has been studied 
	by Biler, Karch and Imbert
	in \cite{biler3al}, establishing some decay estimates of the $ L^s$ norm, such as
	\begin{equation}\label{ch3bilerdecay}
	\Vert u(\cdot, t) \Vert_{L^{s}(\Omega) } \leq t^{-\frac{n}{n+2-2\sigma} \left(1-\frac{1}{s} \right)}.
	\end{equation}
	As a matter of fact, this decay
	is slower than what we find in Theorem \ref{ch3thm:porous}, which is asymptotic to $t^{-1}$
	(in this sense, Theorem \ref{ch3thm:porous} here
	can be seen as an improvement of the estimates
	in~\cite{biler3al}).
	
	On the other hand, in~\cite{biler3al} the Authors also provide a {\em weak} solution that has exactly the decay in~\eqref{ch3bilerdecay}, 
	thus showing the optimality of~\eqref{ch3bilerdecay}
	in this generality,
	while our result holds for {\em strong} solutions. Then, comparing Theorem \ref{ch3thm:porous} here with the results in~\eqref{ch3bilerdecay}
	we obtain that a better decay is valid for regular solutions with respect to the one which is valid
	also for irregular ones.
	
	\section{Proofs}
	
	This section contains the proofs of our main results. We start with the proof of Theorem \ref{ch3thm:complex}.
	
	In order to prove Theorem \ref{ch3thm:complex}, we need a comparison result for the equation involving the mixed time-derivative. As a matter of fact,
	comparison results
	for the case of the Caputo derivative are available
	in the literature, see e.g. Lemma 2.6 of \cite{VZ15}. In our arguments
	we employ the differentiability of $u$ and the fact that $u$ is a strong solution, and we obtain:
	
	\begin{lemma} \label{ch3lemma:comparison}
		Let $T\in(0,+\infty)\cup\{+\infty\}$
		and~$w, \ v: [0,T) \rightarrow [0,+\infty) $ be two Lipschitz continuous
		functions.
		Assume that~$w$ is a supersolution and $v$ is a subsolution at each differentiability point for the equation
		\begin{equation}\label{ch30919}
		\lambda_1 \partial_t^{\alpha} u(t) + \lambda_2 \partial_t u(t) =-ku^{\gamma}(t)
		\end{equation}
		with $\lambda_1$, $\lambda_2$, $\gamma$, $k >0$.
		
		Then:
		if \begin{equation}\label{ch3X00}
		w(0)> v(0), \end{equation}
		we have that \begin{equation}\label{ch3X01}
		w(t)>v(t)\qquad{\mbox{ for all }}t\in(0,T).\end{equation}
	\end{lemma}
	
	\begin{proof}
		By contradiction, let us suppose that for some time $t \in(0,T)$ we have $w(t)=v(t)$, and let us call $\tau$ the first time for which the equality is reached. Then,
		since~$w$ is a supersolution and $v$ is a subsolution of~\eqref{ch30919},
		we obtain that
		\begin{equation}\label{ch3QUA1}
		\lambda_1 \partial_t^{\alpha} (w-v)(\tau) + \lambda_2 \partial_t (w-v)(\tau) \geq -k [w^{\gamma}(\tau) - v^{\gamma}(\tau)]=0.
		\end{equation}
		Now we distinguish two cases, depending on whether or not~$w-v$ is differentiable at~$\tau$.
		To start with,	suppose that $w-v$ is differentiable at~$\tau$.
		Since~$w\ge v$ in~$(0,\tau)$, we have that
		\begin{equation*}
		\partial_t (w-v)(\tau) \le 0.\end{equation*}
		{F}rom this and~\eqref{ch3QUA1}, we obtain that
		\begin{eqnarray*}
			0&\le&\partial_t^{\alpha} (w-v)(\tau) \\&=&
			\frac{(w-v)(\tau) - (w-v)(0)}{\tau^{\alpha}} +\alpha \int_0^{\tau} \frac{(w-v)(\tau) - (w-v)(\rho)}{(\tau-\rho)^{1+\alpha}} d\rho\\&=&
			-\frac{ (w-v)(0)}{\tau^{\alpha}} -\alpha \int_0^{\tau} \frac{ (w-v)(\rho)}{(\tau-\rho)^{1+\alpha}} d\rho\\&\le&
			-\frac{ (w-v)(0)}{\tau^{\alpha}} 
			.\end{eqnarray*}
		This is in contradiction with~\eqref{ch3X00} and so it proves~\eqref{ch3X01} in this case.
		
		Now we focus on the case in which~$w-v$ is not differentiable at~$\tau$.
		Then, there exists a sequence~$t_j\in(0,\tau)$
		such that $w-v$ is differentiable at~$t_j$, with~$\partial_t(w-v)(t_j)\le 0$
		and~$t_j\to\tau$ as~$j\to+\infty$.
		Consequently, 
		since~$w$ is a supersolution and $v$ is a subsolution of~\eqref{ch30919},
		we obtain that
		\begin{equation}\label{ch3QUA2}
		\begin{split}
		&\frac{(w-v)(t_j) - (w-v)(0)}{t_j^{\alpha}} +\alpha \int_0^{t_j} \frac{(w-v)(t_j)
			- (w-v)(\rho)}{(t_j-\rho)^{1+\alpha}} d\rho
		\\ =\;&
		\partial_t^{\alpha} (w-v)(t_j)
		\\ \ge\;&
		\partial_{t}^{\alpha} (w-v)(t_j) + \frac{\lambda_2}{\lambda_1} \partial_t (w-v)(t_j) 
		\\ \geq\;& -\frac{k}{\lambda_1}\, [w^{\gamma}(t_j) - v^{\gamma}(t_j)]
		.\end{split}\end{equation}
		Now we observe that if~$f$ is a Lipschitz function and~$t_j\to\tau>0$
		as~$j\to+\infty$, then
		\begin{equation}\label{ch3VITALI}
		\lim_{j\to+\infty}
		\int_0^{t_j} \frac{f(t_j)
			- f(\rho)}{(t_j-\rho)^{1+\alpha}} d\rho
		=\int_0^{\tau} \frac{f(\tau)
			- f(\rho)}{(\tau-\rho)^{1+\alpha}} d\rho.
		\end{equation}
		To check this, let
		$$ F_j(\rho):=\chi_{(0,t_j)}(\rho)\,
		\frac{f(t_j)
			- f(\rho)}{(t_j-\rho)^{1+\alpha}},$$
		and let~$E\subset(0,+\infty)$ be a measurable set, with measure~$|E|$ less than
		a given~$\delta>0$. Let also~$q:=\frac{1+\alpha}{2\alpha}>1$
		and denote by~$p$ its conjugated exponent. Then, by H\"older inequality, for large~$j$
		we have that
		\begin{eqnarray*}
			\int_E| F_j(\rho)|\,d\rho&\le& |E|^{1/p}\,\left( \int_0^{+\infty} |F_j(\rho)|^q\,d\rho\right)^{1/q}\\
			&\le& \delta^{1/p}\,\left( \int_0^{t_j} 
			\frac{|f(t_j) - f(\rho)|^q}{(t_j-\rho)^{(1+\alpha)q}}
			\,d\rho\right)^{1/q}\\
			&\le& L\,\delta^{1/p}\,\left( \int_0^{t_j} 
			\frac{d\rho}{(t_j-\rho)^{\alpha q}}\right)^{1/q}\\
			&=& L\,\delta^{1/p}\,\left( \int_0^{t_j} 
			\frac{d\rho}{(t_j-\rho)^{(1+\alpha)/2}}\right)^{1/q}\\
			&=& L\,\delta^{1/p}\,\left( \frac{2 t_j^{(1-\alpha)/2}}{1-\alpha}\right)^{1/q}
			\\&\le&
			L\,\left( \frac{2 (\tau+1)^{(1-\alpha)/2}}{1-\alpha}\right)^{1/q}
			\,\delta^{1/p}
			,\end{eqnarray*}
		where~$L$ is the Lipschitz constant of~$f$. Consequently, by the
		Vitali Convergence Theorem,
		we obtain that
		$$ \lim_{j\to+\infty}\int_0^{+\infty} F_j(\rho)\,d\rho=
		\int_0^{+\infty} \lim_{j\to+\infty}F_j(\rho)\,d\rho,$$
		which gives~\eqref{ch3VITALI}, as desired.
		
		Now, we take the limit as~$j\to+\infty$ in~\eqref{ch3QUA2},
		exploiting~\eqref{ch3VITALI} and the fact that~$w(\tau)=v(\tau)$. In this
		way, we have that
		$$
		-\frac{ (w-v)(0)}{\tau^{\alpha}} -\alpha \int_0^{\tau} \frac{
			(w-v)(\rho)}{(\tau-\rho)^{1+\alpha}} d\rho
		\ge0.$$
		Since~$w\ge v$ in~$(0,\tau)$, the latter inequality implies that
		$$
		-\frac{ (w-v)(0)}{\tau^{\alpha}} 
		\ge0.$$
		This is in contradiction with~\eqref{ch3X00} and so it completes the proof of~\eqref{ch3X01}.
	\end{proof}
	
	It is also useful to observe that Lemma~\ref{ch3lemma:comparison}
	holds true also for the classical derivative (i.e. when~$\lambda_1=0$).
	We give its statement and proof for the sake of completeness:
	
	\begin{lemma} \label{ch3lemma:comparison2}
		Let $T\in(0,+\infty)\cup\{+\infty\}$,
		$w, \ v: [0,T) \rightarrow [0,+\infty)$ be two Lipschitz continuous
		functions.
		Assume that~$w$ is a supersolution and $v$ is a subsolution at each differentiability point for the equation
		\begin{equation}\label{ch309192}
		\partial_t u(t) =-ku^{\gamma}(t)
		\end{equation}
		with $\gamma$, $k >0$.
		
		Then:
		if \begin{equation}\label{ch3X002}
		w(0)> v(0), \end{equation}
		we have that \begin{equation}\label{ch3X012}
		w(t)>v(t)\qquad{\mbox{ for all }}t\in(0,T).\end{equation}
	\end{lemma}
	
	\begin{proof} Suppose that~\eqref{ch3X012} is false. Then there exists~$\tau\in(0,T)$
		such that~$w>v$ in~$(0,\tau)$ and
		\begin{equation}\label{ch3TAU0}w(\tau)=v(\tau).\end{equation}
		We fix~$\e>0$, to be taken as small as we wish in the sequel,
		and define
		\begin{equation}\label{ch3TAU3}
		f(t):=w(t)-v(t)+\e\,(t-\tau).
		\end{equation}
		We observe that
		$$f(0)=w(0)-v(0)-\e\tau\ge\frac{w(0)-v(0)}2>0,$$
		as long as~$\e$ is sufficiently small,
		and~$f(\tau)=w(\tau)-v(\tau)=0$.
		Therefore there exists~$\tau_\e\in(0,\tau]$ such that
		\begin{equation}\label{ch3TAU2}
		{\mbox{$f>0$ in~$(0,\tau_\e)$
				and~$f(\tau_\e)=0$.}}\end{equation}
		We claim that
		\begin{equation}\label{ch3TAU}
		\lim_{\e\to0^+}\tau_\e=\tau.
		\end{equation}
		Indeed, suppose, by contradiction, that, up to a subsequence, $\tau_\e$
		converges to some~$\tau_0\in[0,\tau)$ as~$\e\to0^+$. Then we have that
		$$ 0=\lim_{\e\to0^+} f(\tau_\e)=\lim_{\e\to0^+}
		w(\tau_\e)-v(\tau_\e)+\e\,(\tau_\e-\tau)=w(\tau_0)-v(\tau_0).$$
		This is in contradiction with the definition of $\tau$
		and so~\eqref{ch3TAU} is proved.
		
		Now, from~\eqref{ch3TAU2}, we know that
		there exists a sequence~$t_j\in(0,\tau_\e]$ such that~$f$
		is differentiable at~$t_j$, $\partial_t f(t_j)\le0$
		and~$t_j\to\tau_\e$ as~$j\to+\infty$.
		
		Accordingly, we deduce from~\eqref{ch309192} and~\eqref{ch3TAU3} that
		$$ 0 \ge \partial_t f(t_j)=
		\partial_t (w-v)(t_j)+\e
		\ge-k\big( w^{\gamma}(t_j)-v^\gamma(t_j)\big)+\e.$$
		Hence, taking the limit as~$j\to+\infty$,
		\begin{equation}\label{ch3TAY1}
		\frac{\e}{k}\le w^{\gamma}(\tau_\e)-v^\gamma(\tau_\e)=
		\big( v(\tau_\e)+\e\,(\tau-\tau_\e)\big)^\gamma-v^\gamma(\tau_\e).
		\end{equation}
		We claim that 
		\begin{equation}\label{ch3TAY2}\liminf_{\e\to0^+}
		v(\tau_\e)>0.
		\end{equation}
		Indeed, if not, by~\eqref{ch3TAU0} and~\eqref{ch3TAU},
		\begin{equation}\label{ch3TAY3}
		0=\liminf_{\e\to0^+} v(\tau_\e)=v(\tau)=w(\tau).
		\end{equation}
		We observe that this implies that
		\begin{equation}\label{ch31GAMMA}
		\gamma\in(0,1).
		\end{equation}
		Indeed, since~$w$ is a supersolution of~\eqref{ch309192},
		we have that
		\begin{eqnarray*}&&
			w(t)\ge w(0)\,e^{-kt},\qquad{\mbox{when }}\gamma=1\\
			{\mbox{and }}&&w(t)\ge\frac{1}{\left(
				\frac1{w^{\gamma-1}(0)}+k(\gamma-1)t
				\right)^{\frac1{\gamma-1}}},
			\qquad{\mbox{when }}\gamma>1,\end{eqnarray*}
		as long as~$w(t)>0$, and so for all~$t>0$.
		In particular, we have that~$w(\tau)>0$, in contradiction with~\eqref{ch3TAY3},
		and this proves~\eqref{ch31GAMMA}.
		
		Then, we use that~$v$ is a subsolution of~\eqref{ch309192} and~\eqref{ch3TAY3}
		to write that, for any~$t\in(0,\tau)$,
		$$ 
		-\frac{v^{1-\gamma}(t)}{1-\gamma}=
		\frac{v^{1-\gamma}(\tau)-v^{1-\gamma}(t)}{1-\gamma}=\frac1{1-\gamma}
		\int_t^\tau \partial_\rho (v^{1-\gamma}(\rho))\,d\rho
		=\int_t^\tau \frac{\partial_t v(\rho)}{v^\gamma(\rho)}\,d\rho\le-k(\tau-t).$$
		Therefore, recalling~\eqref{ch31GAMMA},
		$$ v^{1-\gamma}(t)\ge k(1-\gamma)(\tau-t),$$
		and thus
		\begin{equation}\label{ch37ygfugv}
		v(t)=v(t)\ge \big( k(1-\gamma)(\tau-t)\big)^{1/(1-\gamma)}.\end{equation}
		Similarly, using that~$w$ is a supersolution of~\eqref{ch309192} and~\eqref{ch3TAY3}
		we obtain that, for any~$t\in(0,\tau)$,
		$$ w(t)\le \big( k(1-\gamma)(\tau-t)\big)^{1/(1-\gamma)}.$$
		Comparing this and~\eqref{ch37ygfugv}, we conclude that
		$$ w(0)\le\big( k(1-\gamma)\tau\big)^{1/(1-\gamma)}\le v(0),$$
		which is in contradiction with~\eqref{ch3X002}, and so the proof of~\eqref{ch3TAY2}
		is complete.
		
		Then, using~\eqref{ch3TAY1}
		and~\eqref{ch3TAY2}, a Taylor expansion gives that
		\begin{eqnarray*}
			\frac{1}{k}&\le& \frac{v^\gamma(\tau_\e)}{\e}\,\left[
			\left( 1+\frac{\e\,(\tau-\tau_\e)}{v(\tau_\e)}\right)^\gamma-1\right]
			\\&=&\frac{
				v^\gamma(\tau_\e)}{\e}\left[
			\frac{\gamma\e\,(\tau-\tau_\e)}{v(\tau_\e)}+
			O\left( \frac{\e^2\,(\tau-\tau_\e)^2}{v^2(\tau_\e)}\right)
			\right]\\&=&
			\frac{\gamma\,(\tau-\tau_\e)}{v^{1-\gamma}(\tau_\e)}+
			O\left( \frac{\e\,(\tau-\tau_\e)^2}{v^{2-\gamma}(\tau_\e)}\right)
			.\end{eqnarray*}
		Then, sending~$\e\to0^+$ and recalling~\eqref{ch3TAU}
		and~\eqref{ch3TAY2}, we conclude that~$\frac1k\le0$.
		This is a contradiction and the proof of~\eqref{ch3X012} is thereby complete.
	\end{proof}
	
	With this preliminary work, we are in the position of proving the general
	claim stated in Theorem \ref{ch3thm:complex}.
	
	\begin{proof}[Proof of Theorem \ref{ch3thm:complex}] 	
		First, notice that
		\begin{equation} \label{ch30721}
		{\partial_t |u|^{s}}= s |u|^{s-1} \left(\frac{\Re(u) \partial_{t} \Re(u)+ \Im(u) \partial_{t} \Im(u)}{|u|}  \right) = s|u|^{s-2} \Re\{\bar{u} \, \partial_t u\}.
		\end{equation}
		Using \eqref{ch30721} and exchanging the order of the integral and the derivative, we have
		\begin{equation}\label{ch3ineq:complex0}
		\begin{split}
		\int_{\Omega} |u|^{s-2} \Re\{\bar{u} \, \partial_t u\} \; dx &= \int_{\Omega}  \frac{\partial_t |u|^{s}}{s} \; dx 
		=\frac{1}{s} \partial_t \int_{\Omega} |u|^s \; dx =  \frac{1}{s} \partial_t \Vert u(\cdot, t) \Vert_{L^{s} (\Omega) }^{s} \\
		& =\Vert u(\cdot, t) \Vert_{L^{s} (\Omega) }^{s-1} \partial_t \Vert u(\cdot, t) \Vert_{L^{s} (\Omega)}.
		\end{split}
		\end{equation}
		Now we claim that
		\begin{equation}\label{ch3ineq:complex}
		\Vert u(\cdot ,t) \Vert_{L^{s}(\Omega) }^{s-1} \partial_t^{\alpha} (\Vert u(\cdot ,t) \Vert_{L^{s}(\Omega) }) \leq \int_{\Omega} |u(x,t)|^{s-2} \Re \{ \bar{u}(x,t) \partial_t^{\alpha} (u(x,t)) \} \, dx.
		\end{equation}
		This formula is similar to one given in Corollary 3.1 of \cite{VZ15}
		for general kernels. In our setting,
		we provide
		an easier proof for the case of the Caputo derivative, comprising also the case
		of complex valued operators. 
		To prove~\eqref{ch3ineq:complex},
		using the definition of Caputo derivative we see that
		\begin{equation*}
		\begin{split}
		\int_{\R^n} & |u(x,t)|^{s-2} \Re \{ \bar{u}(x,t)\partial_t^{\alpha} u(x,t) \} \; dx \\
		&=\int_{\Omega} |u(x,t)|^{s-2} \Re \left\{ \bar{u}(x,t) \left[ 
		\dfrac{u(x,t)-u(x,0)}{t^{\alpha}} + \alpha \int_{0}^{t} \dfrac{u(x,t)-u(x,\tau)}{(t-\tau)^{1+\alpha}} \;d\tau
		\right]
		\right\} \; dx \\
		&= \int_{\Omega} |u(x,t)|^{s-2} \bigg( \frac{|u(x,t)|^2 - \Re \{ \bar{u}(x,t) u(x,0) \} }{t^{\alpha}}  \\
		& \hspace{1em} + \alpha \int_{0}^{t} \frac{|u(x,t)|^2 - \Re \{ \bar{u}(x,t) u(x,\tau) \} }{(t-\tau)^{1+\alpha}} \; d\tau \bigg) dx.
		\end{split}
		\end{equation*}	
		Hence, by using the H\"older inequality, we get
		\begin{equation*}
		\begin{split}	
		\int_{\R^n} & |u(x,t)|^{s-2} \Re \{ \bar{u}(x,t)\partial_t^{\alpha} u(x,t) \} \; dx \\
		& \geq \frac{ \Vert u(\cdot, t) \Vert_{L^s(\Omega)}^{s} - \Vert u(\cdot, t) \Vert_{L^s(\Omega)}^{s-1} \Vert u(\cdot, 0) \Vert_{L^s(\Omega)}  }{t^{\alpha}} +\alpha \int_{0}^{t} \frac{\Vert u(\cdot, t) \Vert_{L^s(\Omega)}^{s}}{(t-\tau)^{1+\alpha}} \; d\tau \\
		& \hspace{1em} - \alpha \int_{0}^{t} \dfrac{\Vert u(\cdot, t) \Vert_{L^s(\Omega)}^{s-1} \Vert u(\cdot, \tau) \Vert_{L^s(\Omega)}}{(t-\tau)^{1+\alpha}} \; d\tau \\
		& = \Vert u(\cdot, t) \Vert_{L^s(\Omega)}^{s-1} \bigg[ \dfrac{\Vert u(\cdot, t) \Vert_{L^s(\Omega)}-\Vert u(\cdot, 0) \Vert_{L^s(\Omega)}}{t^{\alpha}} \\
		& \hspace{1em} + \alpha \int_{0}^{t} \dfrac{\Vert u(\cdot, t) \Vert_{L^s(\Omega)} - \Vert u(\cdot, \tau) \Vert_{L^s(\Omega)}}{(t-\tau)^{1+\alpha}}\; d\tau  \bigg] \\
		& = \Vert u(\cdot, t) \Vert_{L^s(\Omega)}^{s-1} \partial_t^{\alpha} \Vert u(\cdot, t) \Vert_{L^s(\Omega)}.
		\end{split}
		\end{equation*}
		This completes the proof of \eqref{ch3ineq:complex}.
		
		Now,	to make the notation simpler, we set $v(t):= \Vert u(\cdot, t) \Vert_{L^{s} (\Omega) }$. By combining \eqref{ch3ineq:complex0} and~\eqref{ch3ineq:complex}, we find that 
		\begin{equation*}
		v^{s-1}(t) \left( \lambda_1 \partial_t^{\alpha} v(t) + \lambda_2 \partial_t v(t) \right) \leq \int_{\Omega} |u|^{s-2}(x,t) \Re \left\{\bar{u}(x,t) \left( \lambda_1 \partial_t^{\alpha} u(x,t) +\lambda_2 \partial_t u(x,t) \right) \right\} dx
		\end{equation*}
		and so, using the fact that $u$ is a solution of~\eqref{ch3sys:generalform}, we conclude that
		\begin{equation*}
		v^{s-1}(t) \left( \lambda_1 \partial_t^{\alpha} v(t) + \lambda_2 \partial_t v(t) \right) \leq - \int_{\Omega} |u|^{s-2}(x,t) \Re \{\bar{u}(x,t)  \mathcal{N}[u](x,t)\}  dx.
		\end{equation*}
		{F}rom this, we use the structural hypothesis \eqref{ch3cond:complexstr} and we obtain that
		\begin{equation*}
		v^{s-1}(t) \left( \lambda_1 \partial_t^{\alpha} v(t) + \lambda_2 \partial_t v(t) \right) \leq - \frac{v^{s-1+\gamma}(t)}{C}.
		\end{equation*}
		Hence, we have established
		the claim in~\eqref{ch3claim1gen} 
		for all $t>0$ such that $v(t)>0$.
		Then, suppose that for some $\bar{t}>0$ we have $v(\bar{t})=0$. Since~$v$ is nonnegative,
		it follows that
		\begin{equation}\label{ch3210718}
		\partial_t v(\bar{t})=0.
		\end{equation}
		On the other hand, if $v(t)=0$, then
		\begin{equation}\label{ch31321}
		{ \partial_t^{\alpha} v(t)
			\le 0},\end{equation}
		because
		\begin{equation*}
		\partial_t^{\alpha} v(t)= \frac{v(t)-v(0)}{t^{\alpha}} + \int_{0}^{t} \frac{v(t)-v(\tau)}{(t-\tau)^{1+\alpha}} d\tau \leq -\frac{v(0)}{t^{\alpha}} - \int_{0}^{t} \frac{v(\tau)}{(t-\tau)^{1+\alpha}} d\tau \leq 0.
		\end{equation*}  
		So, by~\eqref{ch3210718} and~\eqref{ch31321},
		$\left( \lambda_1 \partial_t^{\alpha} v(\bar t) + \lambda_2 \partial_t v(\bar t) \right) \leq 0$, which gives~\eqref{ch3claim1gen} also in this case, as desired.
		
		Now we exhibit a supersolution $w(t)$ of the equation $(\lambda_1 \partial_t^{\alpha} + \lambda_2 \partial_t) v(t) = -\nu v^{\gamma}(t)$, where $\nu:=\frac{1}{C}$. 
		For this, we recall
		Section 7 of \cite{VZ15}, and we have that the function 
		\begin{equation*}
		w(t):= \left\{
		\begin{array}{ll}
		u_0 & {\mbox{if }} t\in [0,t_0], \\
		Kt^{-\frac{\alpha}{\gamma}} & {\mbox{if }}t\geq t_0, 
		\end{array}
		\right.
		\end{equation*}
		with $K:=u_0t_0^{\frac{\alpha}{\gamma}}$ is a supersolution of~$
		\partial_t^{\alpha} w(t) = -\nu w^{\gamma}(t)$
		as long as
		\begin{equation*}
		t_0 \geq \dfrac{u_0^{1-\gamma}}{\nu} \left(\frac{2^{\alpha}}{\Gamma(1-\alpha)} + \frac{\alpha}{\gamma} \frac{2^{\alpha + \frac{\alpha}{\gamma}}}{\Gamma(2-\alpha)}  \right).
		\end{equation*}
		We claim that $\partial_t w(t) \geq -\nu w^{\gamma} (t)$.
		To prove this, it is equivalent to check that
		\begin{equation*}
		\frac{\alpha}{\gamma} u_0 \, t_0^{\frac{\alpha}{\gamma}} \, t^{-\frac{\alpha}{\gamma}-1} \leq \nu \, u_0^{\gamma} \, t_0^{{\alpha}} \, t^{-\alpha}
		,\end{equation*}
		which is in turn equivalent to
		\begin{equation*}
		\frac{\alpha}{\gamma \, \nu} u_0^{1-\gamma} \, t_0^{ \frac{\alpha}{\gamma} -\alpha} \leq t^{1+\frac{\alpha}{\gamma} -\alpha},
		\end{equation*}
		and the latter equation holds if $$t_0 \geq \max \left\{ 1, \frac{\alpha}{\gamma  \nu} u_0^{1-\gamma} \right\}. $$ Therefore for $t_0$ big enough we have that $w(t)$ is a supersolution of the equation $(\lambda_1 \partial_t^{\alpha} + \lambda_2 \partial_t) v(t) = -\nu v^{\gamma}(t)$. Also, $w(t)$ satisfies
		\begin{equation*}
		w(t)\leq \frac{c}{1+t^{\frac{\alpha}{\gamma}}}
		\end{equation*}
		for some $c>0$ depending only on $\nu,\ \gamma, \ \alpha$ and $w(0)$. Hence by the comparison principle in Lemma \ref{ch3lemma:comparison}, we infer that $v(t) \leq w(t)$, which
		completes the proof of the desired result in~\eqref{ch3claim2gen}.	
	\end{proof}
	
	\begin{proof}[Proof of Theorem \ref{ch3thm:classic}]
		The proof is identical to the one of Theorem \ref{ch3thm:complex} a
		part from the construction of the supersolution
		(and from the use of the comparison principle in Lemma~\ref{ch3lemma:comparison2}
		rather than in Lemma~\ref{ch3lemma:comparison}). Our aim is now to find a supersolution to the equation \eqref{ch3claim1gen} in the case $\lambda_1=0$, that we can write as
		\begin{equation}\label{ch3CAS1}
		v'(t) = -\frac{1}{C}v^{\gamma}(t)
		\end{equation}
		where $C$ is the constant given in the hypothesis. To construct this supersolution,
		we distinguish the cases~$0<\gamma \leq 1$
		and~$\gamma>1$.
		
		We define
		\begin{equation}\label{ch379}
		w_0:=\Vert u_0(\cdot) \Vert_{L^{s} (\Omega) },
		\end{equation}
		
		\begin{equation}
		t_0:=\left\{\begin{matrix}
		0 & {\mbox{if }}\gamma=1,\\ \max\left\{0, \ \frac{C}{1-\gamma}(w_0^{1-\gamma}-1) \right\}& {\mbox{if }}0<\gamma<1,
		\end{matrix}
		\right.\end{equation}
		and 
		\begin{equation}\label{ch3theta1}
		\theta_0= \left(w_0-\dfrac{(1-\gamma)}{C}t_0\right).
		\end{equation}
		Notice that, for $0<\gamma<1$ 
		\begin{equation}\label{ch3theta}
		\theta_0 \leq 1.
		\end{equation}
		In fact, 
		\begin{equation*}
		\frac{C}{1-\gamma} (w_0^{1-\gamma}-1) \leq t_0
		\end{equation*}
		implies
		\begin{equation*}
		\left( w_0^{1-\gamma} - \frac{(1-\gamma)}{C}t_0  \right) \leq 1
		\end{equation*}
		and that proves \eqref{ch3theta}.
		Then, we see that the function
		\begin{equation}\label{ch397}
		w(t):= \left\{ 
		\begin{array}{lr}
		\left(w_0^{1-\gamma}-\dfrac{(1-\gamma)t}{C} \right)^{\frac{1}{1-\gamma}}, & {\mbox{if }}t\in[0,t_0] \\
		\theta_0  \,e^{\frac{t_0-t}{C}}, & {\mbox{if }}t\in(t_0, +\infty)
		\end{array}
		\right.
		\end{equation}
		is a continuous and Lipschitz function, moreover it is a solution of~\eqref{ch3CAS1}
		in the case $\gamma=1$ and a supersolution of~\eqref{ch3CAS1} in the case $0 <\gamma<1$. Indeed, to check this, we observe that, for $t\in[0, t_0]$,
		\begin{eqnarray*}
			&& \hspace{-1em}  w'(t)+\frac1{C}w^{\gamma}(t)) \\
			&& \hspace{3em} = -\dfrac{1}{C}\left( w_0^{1-\gamma} -\dfrac{(1-\gamma)t}{C} \right)^{\frac{\gamma}{1-\gamma}}  + \dfrac1{C}\left( w_0^{1-\gamma} -\dfrac{(1-\gamma)t}{C} \right)^{\frac{\gamma}{1-\gamma}} \\
			&& \hspace{3em}  =0,
		\end{eqnarray*}
		while
		for all~$t>t_0$,
		\begin{eqnarray*}
			&& C\left( w'(t)+\frac1{C}w^\gamma(t)\right)=
			-\theta_0 e^{\frac{(t_0-t)}{C}}
			+\theta_0^\gamma  e^{\frac{\gamma(t_0-t)}{C}}=
			\theta_0^\gamma  e^{\frac{\gamma(t_0-t)}{C}}\left(1
			-\theta_0^{1-\gamma} e^{\frac{(1-\gamma)(t_0-t)}{C}}\right)\\&&\qquad\ge
			\theta_0^\gamma  e^{\frac{\gamma(t_0-t)}{C}}\left(1
			-\theta_0^{1-\gamma} \right)\ge0,
		\end{eqnarray*}
		where the inequality holds thanks to \eqref{ch3theta}. 
		Notice also that the function $w$ is Lipschitz since it is piecewise continuous and derivable and it is continuous in the point $t=t_0$ because of the definition of $\theta$ given in \eqref{ch3theta1}.
		These observations establish
		the desired supersolution properties for the
		function in~\eqref{ch397} for~$0<\gamma\le1$.
		{F}rom this and the comparison result
		in Lemma \ref{ch3lemma:comparison2}, used here with $w(t)$ and $v(t):= \Vert u(\cdot, t) \Vert_{L^{s} (\Omega) } $, we obtain that
		$v(t)\le w(t)$ for any~$t\ge0$, and in particular,
		\begin{equation}\label{ch399}
		\Vert u(\cdot, t) \Vert_{L^{s} (\Omega) }\le
		K e^{-\frac{t}{C}}
		\qquad{\mbox{for any~$t>t_0$}}
		\end{equation}
		for $K:= \theta_0 e^{\frac{t_0}{C}}$.
		This proves~\eqref{ch3claim3}. 
		\medskip
		
		Now we deal with the case $\gamma>1$. In this case, we set $$ w_0:=\max \left\{\Vert u_0(\cdot) \Vert_{L^{s} (\Omega) }, \Big( \frac{C}{\gamma -1} \Big)^{\frac{1}{\gamma-1}}  \right\}. $$
		Then the function
		\begin{equation}\label{ch3992}
		w(t):= \left\{ 
		\begin{array}{lr}
		w_0  , & {\mbox{if }}t\in[0,1] \\
		w_0  t^{-\frac{1}{\gamma-1}}, & {\mbox{if }}t>1
		\end{array}
		\right.
		\end{equation}
		is a supersolution of~\eqref{ch3CAS1}. Indeed, if~$t>1$,
		\begin{eqnarray*}
			C\left( w'(t)+\frac1{C}w^\gamma(t)\right)=
			-\frac{C}{\gamma-1} w_0 t^{-\frac{\gamma}{\gamma-1}}
			+w_0^\gamma  t^{-\frac{\gamma}{\gamma-1}}=
			w_0 t^{-\frac{\gamma}{\gamma-1}}\left(
			w_0^{1-\gamma}-\frac{C}{\gamma-1}
			\right)\ge0,
		\end{eqnarray*}
		while, if~$t\in(0,1)$,
		$$ w'(t)+\frac1{C}w^\gamma(t)=\frac1{C}w^\gamma(t)\ge0.$$
		This gives that the function in~\eqref{ch3992} has the desired supersolution property and consequently we can apply the comparison result
		in Lemma \ref{ch3lemma:comparison2} with $w(t)$ and $v(t):= \Vert u(\cdot, t) \Vert_{L^{s} (\Omega) } $. In this way, we obtain that for all~$t\ge1$
		$$ \Vert u(\cdot, t) \Vert_{L^{s} (\Omega) }\le 
		w_0  t^{-\frac{1}{\gamma-1}},$$
		and so
		the proof of~\eqref{ch3claim4}
		is complete.
	\end{proof}
	
	Now, we present the applications of the abstract results to the operators introduced in Section \ref{ch3applications}.
	
	We start with the case of the fractional porous medium equation.
	
	\begin{proof}[Proof of Theorem~\ref{ch3thm:porous}]
		In order to prove Theorem \ref{ch3thm:porous}, our strategy is to verify the
		validity of inequality \eqref{ch3cond:complexstr} with~$\gamma:=2$
		for the porous medium operator,
		which would put us in the position of exploiting Theorems~\ref{ch3thm:complex}
		and~\ref{ch3thm:classic}.
		
		To this end, by elementary computations, up to changes of the positive constant $c$ depending on $n, \ s,$ and $ \sigma$, we see that
		\begin{equation}\label{ch3110}
		\begin{split}
		\int_{\Omega} u^{s-1}(x,t)\mathcal{N}[u](x,t) \; dx   &=\int_{\Omega} - u^{s-1} \nabla \cdot (u \nabla \mathcal{K} u)(x,t) \; dx \hspace{10em} \\
		&= \int_{\Omega} (s-1) u^{s-1}(x,t) \nabla u(x,t) \cdot \nabla \mathcal{K} u (x,t) dx \\
		&= \int_{\Omega} \nabla u^{s}(x,t) \cdot  \nabla \mathcal{K} u (x,t) \,dx
		\end{split}
		\end{equation}
		Now, define for $\e>0$, the regularized operator
		\begin{equation}
		\mathcal{K}_{\e}u= \int_{\Omega}c(n,\sigma) \frac{u(x-y,t)}{(|y|^2+\e^2)^{\frac{n-2\sigma}{2}}} dy.
		\end{equation}
		where $c(n,\sigma)$ is the same constant that appears in the definition of $\mathcal{K}$ in~\eqref{ch3kappa}.
		Notice that, since $u$ is regular, we have 
		\begin{multline}\label{ch3conv}
		\int_{\Omega} \nabla u^{s}(x,t) \cdot \nabla \mathcal{K}_{\e} u (x,t) \,dx \\
		\leq \iint_{\R^n \times \R^n}
		\frac{\chi_{\Omega}(x) \underset{x\in\Omega}{\sup} |\nabla u^s(x,t) | \, \chi_{\Omega}(x-y)\underset{(x-y)\in\Omega}{\sup} |\nabla u(x-y,t)|}{|y|^{n-2\sigma}} dxdy
		\end{multline}
		where $\chi$ is the characteristic function. Thus, thanks to~\eqref{ch3conv} we can apply the
		Dominated Convergence Theorem and obtain
		\begin{equation} \label{ch3limit}
		\underset{\e\rightarrow 0}{\lim} \int_{\Omega} \nabla u^{s}(x,t) \cdot  \nabla \mathcal{K}_{\e} u (x,t) \,dx  = \int_{\Omega} \nabla u^{s}(x,t) \cdot  \nabla \mathcal{K} u (x,t) \,dx.
		\end{equation}
		So, using~\eqref{ch3110} and~\eqref{ch3limit}, we have
		\begin{equation}\label{ch3lim}
		\begin{split}
		\int_{\Omega} u^{s-1}(x,t)\mathcal{N}[u](x,t) \; dx   &=\underset{\e\rightarrow 0}{\lim} \int_{\Omega} \nabla u^{s}(x,t) \cdot  \nabla \mathcal{K}_{\e} u (x,t) \,dx \\
		&=\underset{\e\rightarrow 0}{\lim} \int_{\Omega} \nabla u^{s}(x,t) \cdot \int_{\Omega}  \frac{(-n+2\sigma)c(n,\sigma)u(y)}{(|x-y|^2+\e^2)^{\frac{n-2\sigma+2}{2}}} (x-y) dy  \,dx \\
		&=\underset{\e\rightarrow 0}{\lim} \iint_{\Omega\times\Omega}  \dfrac{c(n,\sigma) u(y,t)\nabla u^{s}(x,t) \cdot (y-x) }{(|x-y|^2+\e^2)^{\frac{n-2\sigma+2}{2}}} \; dy \,dx,
		\end{split}	
		\end{equation}
		up to changes of the positive constant $c(n,\sigma)$.	
		Now we adapt a method that was introduced in \cite{porous} to obtain $L^p$ estimates. We exchange the order of integration and have that
		\begin{equation*}
		\begin{split}
		\iint_{\R^n}  c\, u(y,t) \dfrac{\nabla u^{s}(x,t) \cdot (y-x) }{(|x-y|^2+\e^2)^{\frac{n-2\sigma+2}{2}}} \; dx \,dy \\
		&  \hspace{-12em} = \iint_{\R^n}  c\, u(y,t) \dfrac{\nabla (u^{s}(x,t)-u^s(y,t)) \cdot (y-x) }{(|x-y|^2+\e^2)^{\frac{n-2\sigma+2}{2}}} \; dx \,dy \\
		&  \hspace{-12em} = \iint_{\R^n}  -c {(u^{s}(x,t)-u^s(y,t))u(y,t)} \Bigg[\dfrac{-n}{(|x-y|^2+\e^2)^{\frac{n-2\sigma+2}{2}}} \\ &
		\hspace{-9em }+\frac{(n-2\sigma+2)|x-y|^2}{(|x-y|^2+\e^2)^{\frac{n-2\sigma+4}{2}}}  \Bigg] dx \,dy \\
		&  \hspace{-12em} = \iint_{\R^n} c \frac{(u^{s}(x,t)-u^s(y,t))(u(x,t)-u(y,t))}2\Bigg[\dfrac{-n}{(|x-y|^2+\e^2)^{\frac{n-2\sigma+2}{2}}} \\ &
		\hspace{-9em }+\frac{(n-2\sigma+2)|x-y|^2}{(|x-y|^2+\e^2)^{\frac{n-2\sigma+4}{2}}}  \Bigg] dx \,dy.
		\end{split}
		\end{equation*}
		We observe now that, since $(u^{s}(x,t)-u^s(y,t))(u(x,t)-u(y,t))$ is always positive,  
		\begin{equation*}
		\begin{split}
		& \hspace{-3em}\iint_{\R^n} c \frac{(u^{s}(x,t)-u^s(y,t))(u(x,t)-u(y,t))}2\Bigg[\dfrac{-n}{(|x-y|^2+\e^2)^{\frac{n-2\sigma+2}{2}}}  \\
		&+\frac{(n-2\sigma+2)|x-y|^2}{(|x-y|^2+\e^2)^{\frac{n-2\sigma+4}{2}}}  \Bigg] dx \,dy \\
		&  \leq \iint_{\R^n} c \frac{(u^{s}(x,t)-u^s(y,t))(u(x,t)-u(y,t))(2-2\sigma)}{2|x-y|^{n+2(1-\sigma)}} dx \,dy.
		\end{split}
		\end{equation*}
		Thus, again by the Dominated Convergence Theorem, we can pass to the limit in \eqref{ch3lim} and obtain
		\begin{equation}\label{ch3111}
		\begin{split}
		&\hspace{-1.5em}\int_{\Omega} u^{s-1}(x,t)\mathcal{N}[u](x,t) \; dx \\
		& \hspace{1.5em}=\iint_{\R^n} c \frac{(u^{s}(x,t)-u^s(y,t))(u(x,t)-u(y,t))(2-2\sigma)}{2|x-y|^{n+2(1-\sigma)}} dx \,dy.
		\end{split}
		\end{equation}
		Now, we define $v(x,t)=u^{\frac{s+1}{2}}(x,t)$. Then, by inequality (2.15) of \cite{SD.EV.VV} we have, for some $C>0$,
		\begin{equation*}
		C (u^s(x,t)-u^s(y,t))(u(x,t)-u(y,t) ) \geq |v(x,t)-v(y,t)|^2.
		\end{equation*}
		{F}rom this, \eqref{ch3110} and~\eqref{ch3111} we obtain that
		\begin{equation}\label{ch3112}
		\begin{split}&
		C	\int_{\Omega} u^{s-1}(x,t)\mathcal{N}[u](x,t) \; dx  \\
		&\qquad=
		\iint_{\R^n} c\,C\, \dfrac{2-2s}{2} \dfrac{(u^{s}(x,t)-u^s(y,t))(u(x,t)-u(y,t))}{|x-y|^{n+2(1-\sigma)}} \;dx \,dy\\&\qquad
		\ge
		\iint_{\R^n} c\, \dfrac{2-2s}{2} \dfrac{|v(x,t)-v(y,t)|^2}{|x-y|^{n+2(1-\sigma)}} \;dx \,dy.\end{split}\end{equation}
		Now we set $z:=(1-s)$; then $z\in(0,1)$ and $n\geq 2z$. Let also
		$$p_z:= \dfrac{2n}{n-2z} \geq 2. $$
		Then for any $q\in [2, p_z]$ we can apply the Gagliardo-Sobolev-{S}lobodetski\u\i \
		fractionary inequality (compare \cite{guida}, Theorem 6.5) and obtain
		\begin{equation}\label{ch3113}
		\left( \int_{\Omega} u^{\frac{s+1}{2}q} \right)^{\frac{2}{q}}
		= \Vert v \Vert_{L^{q} (\Omega) }^2
		\leq C \iint \dfrac{|v(x,t)-v(y,t)|^2}{|x-y|^{n+2z}} \; dxdy
		\end{equation}
		with $C$ depending only on $\Omega,\ n,\ z$ and $q$. 
		In particular, choosing $q=2$, we deduce from~\eqref{ch3113} that
		\begin{equation}\label{ch3091a}
		\Vert u(\cdot, t) \Vert_{L^{s+1}(\Omega) }^{s+1}
		\leq C \iint \dfrac{|v(x,t)-v(y,t)|^2}{|x-y|^{n+2z}} \; dxdy
		\end{equation}
		On the other hand, using the H\"{o}lder inequality, one has that
		\begin{equation*}
		\Vert u(\cdot, t) \Vert_{L^{s}(\Omega) }^{s+1} \leq \Vert u(\cdot, t) \Vert_{L^{s+1}(\Omega)}^{s+1} |\Omega|^{1/s}.
		\end{equation*}
		Combining this and~\eqref{ch3091a}, we obtain
		\begin{equation*}
		\Vert u(\cdot, t) \Vert_{L^{s+1}(\Omega) }^{s}
		\leq C \iint \dfrac{|v(x,t)-v(y,t)|^2}{|x-y|^{n+2z}} \; dxdy,
		\end{equation*}
		up to renaming~$C>0$.
		This and~\eqref{ch3112} establish the validity
		of \eqref{ch3cond:complexstr} for $\gamma:=2$, as desired.
	\end{proof}
	
	Now we focus on the Kirchhoff equation, first dealing with the case
	of classical derivatives.
	
	\begin{proof}[Proof of Theorem \ref{ch3thm:cl_kirchhoff}]
		Our objective here is to verify the
		validity of inequality \eqref{ch3cond:complexstr} for suitable values
		of~$\gamma$, and then make use of Theorems~\ref{ch3thm:complex}
		and~\ref{ch3thm:classic}.
		
		First we present the proof for the non-degenerate case, that takes place when $m(\xi)$ has a positive minimum. Let us call $m_0:=\min m(\xi)$, then
		\begin{equation}\label{ch3331}
		m \left(\Vert \nabla u \Vert_{L^2(\Omega)}\right) \int_{\Omega} |u|^{s-2}u (-\Delta)u \; dx \geq m_0	\int_{\Omega} |u|^{s-2}u (-\Delta)u \; dx.
		\end{equation} 
		In Theorem 1.2 of~\cite{SD.EV.VV}, the case of the Laplacian was considered:
		there it was found that, for some $C>0$ depending on $s,\ n,\ \Omega$,
		\begin{equation*}
		\int_{\Omega} |u|^{s-2}u (-\Delta)u \; dx \geq C \Vert u \Vert_{L^{s} (\Omega) }^s.
		\end{equation*}
		Combining this with~\eqref{ch3331}
		we see that \eqref{ch3cond:complexstr} holds true for $\gamma=1$ and $C> 0$ depending on $s,\ n,\ \Omega, \ \min m(\xi)$. 
		
		Now we deal with
		the degenerate case, which requires the use of
		finer estimates. In this case, we have that
		\begin{equation}\label{ch30909}\begin{split}
		& \hspace{-3em} b \Vert \nabla u \Vert_{L^2(\Omega)}^2  \int_{\Omega} |u(x,t)|^{s-2} u(x,t)(-\Delta)u (x,t) \; dx\\
		& \hspace{3em} = b \Vert \nabla u \Vert_{L^2(\Omega)}^2 \int_{\Omega} |u(x,t)|^{s-2} |\nabla u (x,t)|^2 \; dx \\
		& \hspace{3em} \geq C \left(\int_{\Omega} |u(x,t)|^{\frac{s-2}{2}} |\nabla u(x,t)|^2 \; dx\right)^2
		,\end{split}\end{equation}
		where the first passage is an integration by parts and the last inequality holds in view of the Cauchy-Schwarz inequality. 
		
		Now define \begin{equation}\label{ch3992k}
		v(x,t):=|u|^{\frac{s+2}{4}}(x,t).\end{equation} We have that
		$$ |\nabla v|^2 = \left( \frac{s+2}{4} \right)^2 |u|^{\frac{s-2}{2}} |\nabla u|^2. $$
		This and~\eqref{ch30909} give that
		\begin{equation}\label{ch3781-b}
		\begin{split}
		&\left( \frac{s+2}{4} \right)^4
		b \Vert \nabla u \Vert_{L^2(\Omega)}^2  \int_{\Omega} |u(x,t)|^{s-2} u(x,t)(-\Delta)u (x,t) \; dx\\
		\ge\,&
		C \left(\int_{\Omega} \left( \frac{s+2}{4} \right)^2|u(x,t)|^{\frac{s-2}{2}} |\nabla u(x,t)|^2 \; dx\right)^2\\
		=\,&
		C \left(\int_{\Omega} |\nabla v(x,t)|^2 \; dx\right)^2.
		\end{split}\end{equation}
		We now use Sobolev injections (in the form
		given, for instance, in formula (2.9) of~\cite{SD.EV.VV}), remembering that $v$ is zero outside $\Omega$. The inequality
		\begin{equation}\label{ch3781-a}
		\Vert \nabla v \Vert_{L^2(\Omega)} \geq C \Vert v \Vert_{L^q(\Omega)} \end{equation}
		holds 
		\begin{equation}\label{ch3PER781-a}
		{\mbox{for all $q\geq 1$ if $n\in\{1, 2\}$,
				and for all~$q\in\left[1,\displaystyle\frac{2n}{n-2}\right]$ if $n>2$.}}\end{equation}
		Therefore, we set
		\begin{equation}\label{ch3PER781-b} q:=\frac{4s}{s+2}.\end{equation}
		Recalling the ranges of~$s$ in claim~(iii) of Theorem~\ref{ch3thm:cl_kirchhoff},
		when~$n>2$ we have that
		$$ (n-2) q-2n=\frac{4s(n-2)}{s+2}-2n=\frac{2}{s+2}\,\big(
		(n-4)s-2n
		\big)\le0,$$
		which shows that the definition in~\eqref{ch3PER781-b}
		fulfills the conditions in~\eqref{ch3PER781-a}, and so~\eqref{ch3781-a}
		is valid in this setting.
		
		Hence, making use of~\eqref{ch3992k}, \eqref{ch3781-b} and \eqref{ch3781-a}, up to renaming~$C$ line after line, we deduce that
		\begin{eqnarray*}
			&&b \Vert \nabla u \Vert_{L^2(\Omega)}^2  \int_{\Omega} |u(x,t)|^{s-2} u(x,t)(-\Delta)u (x,t) \; dx\\
			&&\qquad\ge
			C \|\nabla v(\cdot,t)\|^4_{L^2(\Omega)}
			\ge C\|v\|_{L^q(\Omega)}^4= C\|u\|_{L^s(\Omega)}^{{s+2}}.
		\end{eqnarray*}
		These observations imply that
		condition \eqref{ch3cond:complexstr} is satisfied here
		with $\gamma=3$ and $C$ depending on $s,\ m(\xi)$ and $\Omega$.
	\end{proof}
	
	Now we deal with the case of the fractional Kirchhoff equation.
	
	\begin{proof}[Proof of Theorem~\ref{ch3thm:Kirchhoff}]
		As in the case of classical space-derivatives
		dealt with in the proof of
		Theorem~\ref{ch3thm:cl_kirchhoff},
		a quick proof for the non-degenerate case is available. Indeed,
		\begin{equation*}
		\int_{\Omega} |u|^{s-2} u \mathcal{N}[u] \, dx =
		m \left(\Vert \nabla u\Vert_{L^{2} (\Omega) }^2 \right) \int_{\Omega} |u|^{s-2} u (-\Delta)^{\sigma}u \, dx  \geq \int_{\Omega} m_0 |u|^{s-2}u(-\Delta)^{\sigma}u \, dx
		\end{equation*}
		and in \cite{SD.EV.VV} it was shown that 
		\begin{equation*}
		\int_{\Omega} m_0 |u|^{s-2}u(-\Delta)^{\sigma}u \, dx \geq \Vert u \Vert_{L^{s} (\Omega) }^s.
		\end{equation*}
		Thus, the validity of inequality \eqref{ch3cond:complexstr} with~$\gamma=1$
		is established in this case.
		
		We now deal with the degenerate case. 
		We fix 
		\begin{equation}\label{ch3pge2si}
		p\in[2, +\infty)\end{equation}
		and we define
		\begin{equation}\label{ch318}
		r:= \frac{s+2}{2p}\qquad{\mbox{and}}\qquad
		v(x,t):=|u(x,t)|^r.\end{equation} We claim that
		\begin{equation} \label{ch3kirch:claim1}
		\begin{split}
		&|v(x,t)-v(y,t)|^p  \\
		& \hspace{2em}\leq c_0 |u(x,t)-u(y,t)|\sqrt{(u(x,t)-u(y,t)) (|u(x,t)|^{s-2}u(x,t) -|u(y,t)|^{s-2}u(y,t))}
		\end{split}
		\end{equation}
		for some $c_0> 0$, independent of $u$. To prove this, we first observe
		that the radicand in~\eqref{ch3kirch:claim1} is well defined, since, for every~$a$, $b\in\R$ we have that
		\begin{equation}\label{ch3DO1}
		(a-b) (|a|^{s-2}a -|b|^{s-2}b)\ge0.
		\end{equation}
		To check this, up to exchanging~$a$ and~$b$, we can suppose that~$a\ge b$.
		Then, we have three cases to take into account: either~$a\ge b\ge0$,
		or~$a\ge 0\ge b$, or~$0\ge a\ge b$.
		If~$a\ge b\ge0$, we have that
		$$ |a|^{s-2}a -|b|^{s-2}b= a^{s-1} -b^{s-1}\ge0,$$
		and so~\eqref{ch3DO1} holds true. If instead~$a\ge 0\ge b$, we have that
		$$ |a|^{s-2}a -|b|^{s-2}b=|a|^{s-1} +|b|^{s-1}\ge0,$$
		which gives~\eqref{ch3DO1}
		in this case. Finally, if~$0\ge a\ge b$,
		$$ |a|^{s-2}a -|b|^{s-2}b=-|a|^{s-1} +|b|^{s-1}\ge0,$$
		again since~$-|a|=a\ge b=-|b|$, thus completing the proof of~\eqref{ch3DO1}.
		
		Then, by~\eqref{ch3DO1}, we have that~\eqref{ch3kirch:claim1}
		is equivalent to 
		\begin{equation}\label{ch3kirch:claimE}
		|v(x,t)-v(y,t)|^{2p} \leq c_1 (u(x,t)-u(y,t))^3{ (|u(x,t)|^{s-2}u(x,t) -|u(y,t)|^{s-2}u(y,t))}.
		\end{equation}
		We also note that when $u(x,t)=u(y,t)$ the inequality in~\eqref{ch3kirch:claimE}
		is trivially satisfied.
		Hence, without loss of generality we can suppose that
		\begin{equation}\label{ch3WAG01la}
		{\mbox{$|u(x,t)|>|u(y,t)|$,\; for fixed $x,\ y \in \R^n$.}}\end{equation}
		We define
		the function
		\begin{equation}\label{ch3EA2}
		(-1,1)\ni\lambda\mapsto 	g (\lambda)=\frac{(1- |\lambda|^{\frac{s+2}{2p}})^{2p}}{(1-\lambda)^3(1-|\lambda|^{s-2}\lambda)}
		\end{equation}
		and we claim that
		\begin{equation}\label{ch3EA20}
		\sup_{(-1,1)}g(\lambda) < +\infty.
		\end{equation}
		To this end, we point out that
		$g$ is regular
		for all $\lambda\in (-1,1)$, so, to establish~\eqref{ch3EA20}, we only have to study the
		limits of~$g$ for $\lambda \rightarrow -1^+$ and $\lambda \rightarrow 1^-$.
		
		When~$\lambda \rightarrow -1^+$, this limit is immediate and $g(-1)=0$. On the other hand,
		when~$\lambda \rightarrow 1^-$, we see that
		\begin{equation*}
		\begin{split}
		\underset{\lambda\rightarrow 1^-}{\lim} g(\lambda) &= \underset{\varepsilon\rightarrow 0^+}{\lim} \frac{(1- (1-\e)^{\frac{s+2}{2p}})^{2p}}{(1-(1-\varepsilon))^3(1-(1-\varepsilon)^{s-1})} \\
		&= \underset{\varepsilon\rightarrow 0^+}{\lim} \frac{\left(
			\frac{s+2}{2p}\varepsilon+O(\varepsilon^2)\right)^{2p}}{\varepsilon^3((s-1)\varepsilon + O(\varepsilon^2))} \\
		&=\underset{\varepsilon\rightarrow 0^+}{\lim} 
		\frac{\varepsilon^{2p-4}\,\left(
			\frac{s+2}{2p}+O(\varepsilon)\right)^{2p}}{(s-1 + O(\varepsilon))},
		\end{split}
		\end{equation*}
		which is finite, thanks to~\eqref{ch3pge2si}.
		Then~\eqref{ch3EA20} holds true, as desired.
		
		Then, using~\eqref{ch3EA20} with~$\lambda:=\frac{b}{a}$, we have that 
		\begin{equation}\label{ch3WAG01la2}
		\begin{split}
		&{\mbox{for any~$a$, $b\in\R$ with $|a|>|b|$,}}\\
		&\frac{\left(
			\left|a\right|^{\frac{s+2}{2p}}
			- \left|b\right|^{\frac{s+2}{2p}}\right)^{2p}}{
			\left(a-b\right)^3\left(
			\left|a\right|^{s-2}a
			-\left|b\right|^{s-2}b\right)}
		\;=\;
		\frac{|a|^{s+2}}{|a|^{s-2}\;a^4}\cdot \frac{\left(1- \left|\frac{b}{a}\right|^{\frac{s+2}{2p}}\right)^{2p}}{
			\left(1-\frac{b}{a}\right)^3\left(1-\left|\frac{b}{a}\right|^{s-2}
			\frac{b}{a}\right)}\\&\qquad
		\; =\;
		\frac{(1- |\lambda|^{\frac{s+2}{2p}})^{2p}}{
			(1-\lambda)^3(1-|\lambda|^{s-2}\lambda)}
		\; =\;g(\lambda)\le C,
		\end{split}\end{equation}
		for some~$C>0$. Then, in view of~\eqref{ch3WAG01la},
		we can exploit~\eqref{ch3WAG01la2}
		with~$a:=u(x,t)$ and~$b:=u(y,t)$, from which we obtain that
		\begin{eqnarray*}&& 
			\left|
			\left|u(x,t)\right|^{\frac{s+2}{2p}}
			- \left|u(y,t)\right|^{\frac{s+2}{2p}}\right|^{2p}\;=\;
			\left(
			\left|u(x,t)\right|^{\frac{s+2}{2p}}
			- \left|u(y,t)\right|^{\frac{s+2}{2p}}\right)^{2p}\\ &&\qquad\;\leq \;C\,
			\left(u(x,t)-u(y,t)\right)^3\left(
			\left|u(x,t)\right|^{s-2}u(x,t)
			-\left|u(y,t)\right|^{s-2}u(y,t)\right).\end{eqnarray*}
		This and~\eqref{ch318} 
		imply~\eqref{ch3kirch:claim1}, as desired.
		
		Now, fixed
		$p$ as in~\eqref{ch3pge2si},
		we set
		\begin{equation}\label{ch38iswdjc8383}
		z:=\frac{ 2\sigma}{p}\in(0,\sigma]\subset(0,1).\end{equation}
		We apply the Gagliardo-Sobolev-{S}lobodetski\u\i \
		fractional immersion (for instance, in the version given
		in formula~(2.18) of~\cite{SD.EV.VV}) to $v$. 
		In this way,
		\begin{equation}\label{ch3202020}
		{\mbox{for all $q\in[1,+\infty)$ when $n\le zp$, and for all $q\in\left[1,
				\displaystyle\dfrac{np}{n-zp}\right]$ when~$n>zp$,}}
		\end{equation}
		we have that
		\begin{equation}\label{ch320}\begin{split}
		\Vert u(\cdot,t) \Vert_{L^{\frac{(s+2)q}{2p}}(\Omega)}^{\frac{s+2}2}=
		\Vert v(\cdot,t) \Vert_{L^q(\Omega)}^p &\,\leq C \iint_{\R^{2n}}
		\frac{|v(x,t)-v(y,t)|^p}{|x-y|^{n+zp}} dxdy\\&=
		C \iint_{\R^{2n}}
		\frac{|v(x,t)-v(y,t)|^p}{|x-y|^{n+2\sigma}} dxdy,\end{split}
		\end{equation}
		where the first equality comes from~\eqref{ch318} and the
		latter equality is a consequence of~\eqref{ch38iswdjc8383}.
		
		Now we choose
		\begin{equation}\label{ch3LAq}
		p:=\max\left\{ 2,\,\frac{s+2}{2}\right\}\qquad{\mbox{and}}\qquad
		q:=\frac{2sp}{s+2}.\end{equation}
		Notice that condition~\eqref{ch3pge2si} is fulfilled in this setting.
		Furthermore, recalling~\eqref{ch38iswdjc8383} and
		the assumptions in point~(iii)
		of Theorem~\ref{ch3thm:Kirchhoff}, we have that, when~$n>2\sigma=zp$,
		we have
		\begin{eqnarray*}
			&&(n-zp)q-np=
			\frac{2(n-2\sigma)sp}{s+2}-np=
			\frac{p}{s+2}\,\big(
			2(n-2\sigma)s-n(s+2)
			\big)\\
			&&\qquad=
			\frac{p}{s+2}\,\big(
			s(n-4\sigma)-2n
			\big)\le0.
		\end{eqnarray*}
		As a consequence, we have that condition~\eqref{ch3202020}
		is fulfilled the setting prescribed by~\eqref{ch3LAq},
		hence we can exploit~\eqref{ch320} in this framework.
		
		Then, from~\eqref{ch3LAq} we have that
		$$ \frac{(s+2)q}{2p}=s,$$
		and so~\eqref{ch320} gives that
		$$ \Vert u(\cdot,t) \Vert_{L^{s}(\Omega)}^{\frac{s+2}2}\le
		C \iint_{\R^{2n}}
		\frac{|v(x,t)-v(y,t)|^p}{|x-y|^{n+2\sigma}} dxdy.$$
		Hence, recalling~\eqref{ch3kirch:claim1}, up to renaming~$C>0$,
		we have that
		\begin{equation}\label{ch3COM234520Aiekwd}
		\begin{split} &\Vert u(\cdot,t) \Vert_{L^{s}(\Omega)}^{s+2}\\ \le\,&
		C\big( \iint_{\R^{2n}}
		\frac{|u(x,t)-u(y,t)|\sqrt{(u(x,t)-u(y,t))
				}}{|x-y|^{n+2\sigma}} \\
			& \hspace{2em}\times (|u(x,t)|^{s-2}u(x,t) -|u(y,t)|^{s-2}u(y,t))
			 dxdy\big)^2\\
		\le\,& C
		\iint_{\R^{2n}}
		\frac{|u(x,t)-u(y,t)|^2}{|x-y|^{n+2\sigma}} dxdy\\
		&\times \hspace{2em}
		\iint_{\R^{2n}}
		\frac{{(u(x,t)-u(y,t))
				(|u(x,t)|^{s-2}u(x,t) -|u(y,t)|^{s-2}u(y,t))}}{|x-y|^{n+2\sigma}} dxdy.
		\end{split}\end{equation}	
		Notice also that, in the degenerate case, we deduce from~\eqref{ch3FKPO-1}
		and~\eqref{ch3FKPO} that
		\begin{equation}\label{ch3ghUAJ:a9ok01}
		\begin{split}
		&\int_{\R^n} {\mathcal{N}}[u](x,t)\,|u(x,t)|^{s-2}\,u(x,t)\,dx\\
		=\;&
		-M_u\,
		\iint_{\R^{2n}} 
		\Big( u(x+y,t) + u(x-y,t) -2u(x,t) \Big) 
		\,|u(x,t)|^{s-2}\,u(x,t)\,\frac{dx\,dy}{|y|^{n+2\sigma}}\\
		=\;&
		-2M_u\,\iint_{\R^{2n}} 
		\big( u(y,t)-u(x,t) \big) 
		\,|u(x,t)|^{s-2}\,u(x,t)\,\frac{dx\,dy}{|x-y|^{n+2\sigma}}\\
		=\;&
		2M_u\,\iint_{\R^{2n}} 
		\big( u(x,t)-u(y,t) \big) 
		\,|u(x,t)|^{s-2}\,u(x,t)\,\frac{dx\,dy}{|x-y|^{n+2\sigma}}\\
		=\;&
		M_u\,\iint_{\R^{2n}} 
		\big( u(x,t)-u(y,t) \big) 
		\,\big(|u(x,t)|^{s-2}\,u(x,t)-|u(y,t)|^{s-2}\,u(y,t)\big)\,\frac{dx\,dy}{|x-y|^{n+2\sigma}}
		,\end{split}\end{equation}
		with
		\begin{equation}\label{ch3ghUAJ:a9ok02}\begin{split} M_u\,&:=
		M
		\left( \int_{\R^{2n}}\frac{ |u(x,t)-u(y,t)|^2 }{|x-y|^{n+2\sigma}}
		\,dx\,dy \right)\\&\ge
		b\,\int_{\R^{2n}}\frac{ |u(x,t)-u(y,t)|^2 }{|x-y|^{n+2\sigma}}\,dx\,dy,
		\end{split}\end{equation}
		with~$b>0$.
		
		Then, from~\eqref{ch3ghUAJ:a9ok01}
		and~\eqref{ch3ghUAJ:a9ok02},
		\begin{eqnarray*}
			&&\int_{\R^n} {\mathcal{N}}[u](x,t)\,|u(x,t)|^{s-2}\,u(x,t)\,dx
			\ge b\,\int_{\R^{2n}}\frac{ |u(x,t)-u(y,t)|^2 }{|x-y|^{n+2\sigma}}\,dx\,dy\\
			\qquad&&\times
			\iint_{\R^{2n}} 
			\big( u(x,t)-u(y,t) \big) 
			\,\big(|u(x,t)|^{s-2}\,u(x,t)-|u(y,t)|^{s-2}\,u(y,t)\big)\,\frac{dx\,dy}{|x-y|^{n+2\sigma}}.
		\end{eqnarray*}
		Comparing this with~\eqref{ch3COM234520Aiekwd}, we conclude that
		$$ \Vert u(\cdot,t) \Vert_{L^{s}(\Omega)}^{s+2}\le C\int_{\R^n} {\mathcal{N}}[u](x,t)\,|u(x,t)|^{s-2}\,u(x,t)\,dx,$$
		up to renaming~$C$.
		This gives that hypothesis \eqref{ch3cond:complexstr} is fulfilled
		in this case with $\gamma=3$.
	\end{proof} 
	
	
	Now we deal with the case of the magnetic operators.
	We start with the case of classical space-derivatives.
	For this, we exploit an elementary, but useful, inequality,
	stated in the following auxiliary result:
	
	\begin{lemma}
		Let~$a$, $b\in\R$, and~$\alpha$, $\beta$, $t\in\R^n$.
		Then
		\begin{equation}\label{ch3ST:00}
		(a^2+b^2)\Big(|a t-\beta|^2 + |bt+\alpha|^2 \Big)\ge 
		|a\alpha+b\beta|^2.\end{equation}
	\end{lemma}
	
	\begin{proof} For any~$t\in\R^n$, we define
		\begin{equation} \label{ch3872wj2xz}
		f(t):=
		(a^2+b^2)\Big(|a t-\beta|^2 + |bt+\alpha|^2 \Big)- 
		|a\alpha+b\beta|^2.\end{equation}
		We observe that
		\begin{equation}\label{ch3ST:01}
		\begin{split}
		f(0)\,&=
		(a^2+b^2)(\alpha^2+\beta^2) - |a\alpha+b\beta|^2\\ &=
		a^2\alpha^2+
		a^2\beta^2+
		b^2\alpha^2+
		b^2\beta^2
		- ( a^2\alpha^2+b^2\beta^2 +2ab\alpha\beta)
		\\ &=
		a^2\beta^2+
		b^2\alpha^2
		- 2ab\alpha\beta\\
		&= |a\beta-b\alpha|^2.
		\end{split}\end{equation}
		Moreover
		\begin{equation}\label{ch3ST:02}
		\lim_{|t|\to+\infty} f(t)=\left\{
		\begin{matrix}
		+\infty & {\mbox{ if }} a^2+b^2>0,\\
		0 & {\mbox{ otherwise.}}
		\end{matrix}
		\right.
		\end{equation}
		Now we claim that
		\begin{equation}\label{ch3ST:03}
		f(t)\ge0,\end{equation}
		for all~$t\in\R^n$. To prove~\eqref{ch3ST:03} we argue by contradiction
		and assume that
		$$ \inf_{\R^n} f<0.$$
		Then, in view of~\eqref{ch3ST:01} and~\eqref{ch3ST:02}, we have that
		\begin{equation}\label{ch3ST:04}
		f(\bar t)=\inf_{\R^n} f<0,\end{equation}
		for some~$\bar t\in \R^n$. As a consequence,
		$$ 0=\nabla f(\bar t)=
		2(a^2+b^2)\Big(a(a \bar t-\beta) + b(b\bar t+\alpha) \Big)=
		2(a^2+b^2)\Big((a^2+b^2) \bar t-a\beta+b\alpha \Big),$$
		which implies that
		$$ \bar t=
		\frac{a\beta-b\alpha}{a^2+b^2}.$$
		Thus, we substitute this information into~\eqref{ch3872wj2xz}
		and we obtain that
		\begin{eqnarray*}
			f(\bar t) &=&
			(a^2+b^2)\left(\left|
			\frac{a^2\beta-ab\alpha}{a^2+b^2}-\beta\right|^2 + 
			\left|\frac{ab\beta-b^2\alpha}{a^2+b^2}+\alpha\right|^2 \right)- 
			|a\alpha+b\beta|^2\\&=&
			(a^2+b^2)\left(\left|
			\frac{b^2\beta+ab\alpha}{a^2+b^2}\right|^2 + 
			\left|\frac{ab\beta+a^2\alpha}{a^2+b^2}\right|^2 \right)- 
			|a\alpha+b\beta|^2\\&=&
			(a^2+b^2)\left(b^2\left|
			\frac{b\beta+a\alpha}{a^2+b^2}\right|^2 + 
			a^2\left|\frac{b\beta+a\alpha}{a^2+b^2}\right|^2 \right)- 
			|a\alpha+b\beta|^2\\&=&
			(a^2+b^2)(a^2+b^2)\left|
			\frac{b\beta+a\alpha}{a^2+b^2}\right|^2 - 
			|a\alpha+b\beta|^2
			\\&=&0.\end{eqnarray*}
		This is in contradiction with~\eqref{ch3ST:04} and so it proves~\eqref{ch3ST:03},
		which in turn implies~\eqref{ch3ST:00}, as desired.
	\end{proof}
	
	With this, we are now in the position of completing the proof
	of Theorem \ref{ch3thm:cl_magnetic} and obtain the desired decay estimates
	for the classical magnetic operator.
	
	\begin{proof}[Proof of Theorem \ref{ch3thm:cl_magnetic}]
		We want to prove inequality \eqref{ch3cond:complexstr} for the classical magnetic operator in order to apply Theorem \ref{ch3thm:complex}.
		To this end,
		we aim at proving that 
		\begin{equation}\label{ch3FU:MAGN}
		\Re\big\{ \bar u{\mathcal{N}} u\big\}+|u|\Delta|u|\ge0.
		\end{equation}
		To check this, we observe\footnote{For an alternative proof based
			on fractional arguments, see the forthcoming footnote~\ref{ch3FU:MAGN3}.}
		that we can make the computations
		in the vicinity of a point~$x$ for which~$|u(x)|>0$.
		Indeed, if~\eqref{ch3FU:MAGN} holds true at~$\{|u|>0\}$,
		we can fix~$\epsilon>0$ and consider the function~$u_\epsilon:=u+\epsilon$.
		In this way, $u_\epsilon(x)=\epsilon>0$, hence we can apply~\eqref{ch3FU:MAGN}
		to~$u_\epsilon$ and conclude that
		\begin{equation}\label{ch390:91}
		\begin{split}
		0 \,&\le \Re\big\{ \bar u_\epsilon(x){\mathcal{N}} u_\epsilon(x)\big\}+|u_\epsilon(x)|\Delta|u_\epsilon(x)|\\
		&=
		\Re\big\{ (\bar u(x)+\epsilon){\mathcal{N}}u(x)\big\}+
		|u(x)+\epsilon|\Delta|u_\epsilon(x)|.
		\end{split}\end{equation}
		Notice that, for any test function~$\varphi\in C^\infty_0(\Omega)$,
		we have that
		$$ \lim_{\epsilon\to0}\int_\Omega
		\Delta|u_\epsilon(y)|\,\varphi(y)\,dy=
		\lim_{\epsilon\to0}\int_\Omega
		|u_\epsilon(y)|\,\Delta\varphi(y)\,dy=
		\int_\Omega
		|u(y)|\,\Delta\varphi(y)\,dy,$$
		and so (in the distributional sense)
		$$ \lim_{\epsilon\to0} \Delta|u_\epsilon|=
		\Delta|u|.$$
		Hence, we can pass to the limit in~\eqref{ch390:91}
		and obtain~\eqref{ch3FU:MAGN}.
		
		Accordingly, to prove~\eqref{ch3FU:MAGN}, from now on we will
		focus on the case in which~$|u|>0$. We write~$u=a+ib$ and
		we observe that 
		\begin{equation}\label{ch3Bvah}
		\begin{split}
		&\Re \{ -\bar{u}(\nabla-iA)^2 u   \} \\
		=\;& \Re \left\{ -\bar{u}(\Delta u - |A|^2 u -iA \cdot \nabla u -\nabla \cdot (iAu) )   \right\} \\
		=\;& \Re \left\{ -\bar{u} \Delta u + |A|^2 |u|^2 +2\bar{u}iA \cdot \nabla u +i(\nabla \cdot A)|u|^2 \right\} \\
		=\;& \Re \left\{ (-a+ib)( \Delta a+i\Delta b)
		+ |A|^2 (a^2+b^2) +2(b+ia)A \cdot( \nabla a+i\nabla b)
		+i(\nabla \cdot A)|u|^2 \right\}\\
		=\;&  -a\Delta a-b\Delta b
		+ |A|^2 (a^2+b^2) +2b\nabla a\cdot A -2a\nabla b\cdot A
		,
		\end{split}\end{equation}
		where we used the fact that $A$ is
		real valued.
		
		On the other hand, at points where~$|u|\ne0$,
		\begin{eqnarray*}&& \Delta |u|^2=2|u|\Delta |u|+2|\nabla |u||^2\\ 
			{\mbox{and }}&&\nabla |u|= \frac{a \nabla a + b \nabla b}{|u|},
		\end{eqnarray*}
		therefore
		\begin{eqnarray*}
			|u|\Delta |u|&=&\frac12\, \Delta |u|^2-|\nabla |u||^2\\
			&=& \frac12\,\Delta(a^2+b^2)-\frac{|a \nabla a + b \nabla b|^2}{|u|^2}\\
			&=& a\Delta a+b\Delta b+|\nabla a|^2+|\nabla b|^2-\frac{|a \nabla a + b \nabla b|^2}{a^2+b^2}.
		\end{eqnarray*}
		{F}rom this and~\eqref{ch3Bvah}, we conclude that
		\begin{equation}\label{ch39384-0348}
		\begin{split}&
		\Re\big\{ \bar u{\mathcal{N}} u\big\}+|u|\Delta|u|
		\\ =\;& |\nabla a|^2+|\nabla b|^2-\frac{|a \nabla a + b \nabla b|^2}{a^2+b^2}
		+ |A|^2 (a^2+b^2) +2b\nabla a\cdot A -2a\nabla b\cdot A
		\\ =\;& \big| aA-\nabla b\big|^2+\big| bA+\nabla a\big|^2
		-\frac{|a \nabla a + b \nabla b|^2}{a^2+b^2},
		\end{split}\end{equation}
		and the latter term is nonnegative, thanks to~\eqref{ch3ST:00}
		(applied here with~$t:=A$, $\alpha:=\nabla a$ and~$\beta:=\nabla b$).
		This completes the proof of~\eqref{ch3FU:MAGN}.
		
		Then, from~\eqref{ch3FU:MAGN} here and~\cite{SD.EV.VV}
		(see in particular the formula before~(2.12) in~\cite{SD.EV.VV},
		exploited here with~$p:=2$ and~$m:=2$),
		\begin{eqnarray*}
			&& \int_\Omega |u|^{s-2}
			\Re\big\{ \bar u{\mathcal{N}} u\big\}\,dx\ge
			-\int_\Omega |u|^{s-1}\Delta|u|\,dx\\
			&&\qquad=\int_\Omega\nabla |u|^{s-1}\cdot\nabla|u|\,dx
			\ge C\,\Vert u \Vert_{L^s{(\Omega)}}^{s},
		\end{eqnarray*}
		for some~$C>0$.
		This establishes inequality~\eqref{ch3cond:complexstr} in this case,
		with~$\gamma=1$. Hence, Theorem \ref{ch3thm:cl_magnetic} follows
		from Theorems~\ref{ch3thm:complex} and~\ref{ch3thm:classic}.
	\end{proof}
	
	Now we deal with the fractional magnetic operator.
	
	\begin{proof}[Proof of Theorem \ref{ch3thm:magnetic}]
		We have to verify the structural hypothesis \eqref{ch3cond:complexstr}. We already know
		that the desired inequality holds for the fractional Laplacian $(-\Delta)^{\sigma} v$ for $\sigma \in (0,1)$ and $v\geq 0$ (compare Theorem 1.2 of \cite{SD.EV.VV}). We notice that
		\begin{equation}\label{ch3FU:MAGN2}
		\begin{split}
		\Re& \left\{ \frac{\bar{u}(x,t) \left( u(x,t) - e^{i(x-y)A(\frac{x+y}{2})}u(y,t) \right)}{|x-y|^{n+2\sigma}}  \right\} \\
		& \hspace{3em} = \frac{|{u}(x,t)|^2  - \Re \left\{ e^{i(x-y)A(\frac{x+y}{2})}u(y,t) \bar{u}(x,t )\right\}  }{|x-y|^{n+2\sigma}} \\
		& \hspace{3em} \geq |u(x,t)| \frac{|{u}(x,t)|  -  |u(y,t)| }{|x-y|^{n+2\sigma}}, 
		\end{split}	
		\end{equation}
		and therefore\footnote{Interestingly, \label{ch3FU:MAGN3}
			integrating and taking the limit as~$\sigma\to1$ in~\eqref{ch3FU:MAGN2},
			one obtains an alternative (and conceptually simpler)
			proof of~\eqref{ch3FU:MAGN}.
			This is a nice example of analysis in a nonlocal setting
			which carries useful information to the classical case.}
		\begin{equation}\label{ch38i9ik9iok92}
		\int_{\Omega} |u(x,t)|^{s-2} \Re \{ \bar{u}(x,t)\mathcal{N} [u](x,t)\} \; dx 
		\geq \int_{\Omega} |u(x,t)|^{s-1} (-\Delta)^{\sigma}|u|(x,t) \; dx.\end{equation}
		Also, since $|u|$ is a real and positive function, we can exploit
		formula~(2.25) in~\cite{SD.EV.VV} (used here with~$p:=2$) and write that
		$$ \int_{\Omega} |u(x,t)|^{s-1} (-\Delta)^{\sigma}|u|(x,t) \; dx\ge
		{C} \Vert u \Vert_{L^{s}(\Omega) }^{s}.$$
		{F}rom this and~\eqref{ch38i9ik9iok92} we infer that
		condition~\eqref{ch3cond:complexstr} is satisfied in this case
		with~$\gamma=1$. 
		Then, the desired conclusion in Theorem \ref{ch3thm:magnetic}
		follows from Theorems \ref{ch3thm:complex}
		and~\ref{ch3thm:classic}.
	\end{proof}


	
	

	\bibliographystyle{abbrv}
	\bibliography{bibliog}

\end{document}